\documentclass[]{article}

\usepackage{amsmath}
\usepackage{amsfonts}
\usepackage{color}

\usepackage[margin=0.8in]{geometry}
\usepackage{graphicx}
\usepackage{caption}
\usepackage[labelformat=parens,labelfont={}]{subcaption}

\usepackage{url}

\newcommand{\bs}{\boldsymbol}

\newcommand{\dif}{\hspace{0.1cm}\mathrm{d}}

\makeatletter
\newcommand*{\rmnum}[1]{\expandafter\@slowromancap\romannumeral #1@}
\makeatother

\allowdisplaybreaks
\providecommand{\keywords}[1]{\textbf{\textit{keywords---}} #1}

\begin{document}
\title{A numerical study of branching and stability of solutions to three-dimensional martensitic phase transformations using gradient-regularized, non-convex, finite strain elasticity}
\author{K. Sagiyama\thanks{Mechanical Engineering, University of Michigan}, S. Rudraraju\thanks{Mechanical Engineering, University of Michigan} and K. Garikipati\thanks{Mechanical Engineering and Mathematics, University of Michigan, corresponding author \textsf{krishna@umich.edu}}}
\maketitle
\abstract{
In the setting of continuum elasticity, phase transformations involving martensitic variants are modeled by a free energy density function that is non-convex in strain space.  
Here, we adopt an existing mathematical model in which we regularize the non-convex free energy density function by higher-order gradient terms at finite strain and derive boundary value problems via the standard variational argument applied to the corresponding total free energy, inspired by Toupin's theory of gradient elasticity.  
These gradient terms are to preclude existence of arbitrarily fine microstructures, while still allowing for existence of multiple solution branches corresponding to local minima of the total free energy; these are classified as \emph{metastable} solution branches.  
The goal of this work is to solve the boundary value problem numerically in three dimensions, observe solution branches, and assess stability of each branch by numerically evaluating the second variation of the total free energy.  
We also study how these microstructures evolve as the length-scale parameter, the coefficient of the strain gradient terms in the free energy, approaches zero. 
}

$  $\\

\section{Introduction}\label{S:introduction}
Many multi-component solids, such as shape memory alloys (\textsf{NiTi}), involve phase transformations from cubic austenite to tetragonal martensite crystal structures.  
The tetragonal lattice is characterized by transformation strains relative to the undistorted reference cubic structure.  The strain splits the symmetry group of the cubic lattice into three equivalent sub-groups, each of which corresponds to a tetragonal lattice oriented along one of the cubic crystal axes.  
These tetragonal variants accommodate themselves in a body to achieve configurations that are local energy minimizers while maintaining kinematic compatibility.  
As a result, \emph{twin} microstructures form with a tiled appearance due to the near constancy of strain within each twinned sub-domain. 

The underlying phenomenology can be described by a free energy density function that is non-convex in a frame-invariant strain measure, to account for the finite deformation, and admits three minima corresponding to the tetragonal variants.  
The classical variational treatment of elasticity only identifies stationary points, while of particular interest are the \emph{metastable} solution branches that correspond to local minima of the total energy.  
These metastable branches are to be identified in this work by examining the stability of solutions obtained by numerically solving the boundary value problems arising from the variational formulation.  

Configurations that minimize the total free energy on a given domain with boundaries have been studied in the setting of sharp-interface models by constructing sequences that converge weakly to the minimizer \cite{Ball1987,Chipot1988}.  
Although this approach provides one with good insights to various classes of problems, it allows for arbitrarily fine twin microstructures---a non-physical aspect of the mathematical formulation resulting from the absence of interface energies associated with the martensitic phase boundaries.  
Diffuse-interface models resolve this pathology by including higher-order, strain gradient-dependent terms representing the interfacial energy; the coefficients of the higher-order gradient terms, which control the twin interface thicknesses, give rise to \emph{length-scale parameters}.  
As diffuse-interface models directly taking into account the total free energy, they also make it straightforward to solve general boundary value problems, provided that the high-order gradients in the partial differential equation can be suitably treated, to solve problems that involve energy wells of unequal depths, and to investigate stability/metastability of solution branches via investigation of the second variation of the total energy.  
One-dimensional models have been intensively studied in this context:
Carr et al. \cite{Carr1984} showed, for standard one-dimensional problems with Dirichlet and higher-order Neumann boundary conditions, that mere inclusion of the higher-order strain gradient energy terms only leaves a pair of stable solutions, the global minimizers of the total energy, which have only a single \emph{phase boundary}.  
These solutions, however, do not represent experimentally obtained microstructures that have twin layers separated by multiple phase boundaries.  
This gap in the representation was resolved by Truskinovsky and Zanzotto \cite{Truskinovsky1995,Truskinovsky1996} and by Vainchtein and co-workers \cite{Vainchtein1998,Vainchtein1999} by adding to the model an elastic support that represents the multi-dimensional effect. This allowed the successful recovery of metastable solution branches corresponding to local minimizers of the total energy.  
Healey and Miller \cite{Healey2007} studied an anti-plane shear model of martensite-martensite phase transformations for pure Dirichlet problems in two dimensions, and obtained metastable solution branches over a range of values of the length scale parameter.  
Numerical work in this field has included \emph{branch-tracking} techniques to obtain metastable solution branches and evaluation of the second variation of the total free energy to assess their stability \cite{Vainchtein1998,Vainchtein1999,Healey2007}.  

While one-dimensional problems and some restricted, linearized two-dimensional problems may be partially aided by analysis, the complete, nonlinear, three-dimensional treatment at finite strain with general boundary conditions must be numerical.  
Rudraraju et al. \cite{Rudraraju2014} have adopted spline-based (isogeometric analytic) numerical methods to obtain three-dimensional solutions to general boundary value problems of Toupin's theory of gradient elasticity at finite strain \cite{Toupin1962}.  
This approach makes it possible to study a wide diversity of problems.  
In this communication we present numerical solutions to diffuse-interface problems of phase transformations between martensitic variants under traction loading in three dimensions.  
The fundamental framework to study metastable branches is the same as that used previously in the literature \cite{Vainchtein1998,Vainchtein1999,Healey2007}.  
We obtain numerical solutions to boundary value problems from a starting guess.  
Particular solution branches are tracked as the strain gradient length scale parameter is varied, and stability of a given branch is determined by numerically examining the second variation of the total free energy.  
To the best of our knowledge, this is the first three-dimensional study of twin microstructures and their stability using gradient-regularized non-convex elasticity.  
Of further note is that we apply arbitrary boundary conditions.  
Crucial to our work is the numerical framework derived from the work of Rudraraju et al. \cite{Rudraraju2014}.

In Sec. \ref{S:One-dimensional_example} we present an overview of the fundamental ideas using a simple problem in one dimension.  
Three-dimensional problems are then studied in Sec. \ref{S:Three-dimensional_example} employing virtually the same numerical techniques.  
Conclusions and future studies are discussed in Sec. \ref{S:Conclusion}.

\section{A one-dimensional primer}\label{S:One-dimensional_example}
Branching in three dimensions being our eventual concern, it is instructive to first study a related problem in one dimension. 

\begin{figure}
\begin{center}
        \begin{subfigure}[b]{5.5cm}
            \centering
            \includegraphics[scale=0.25]{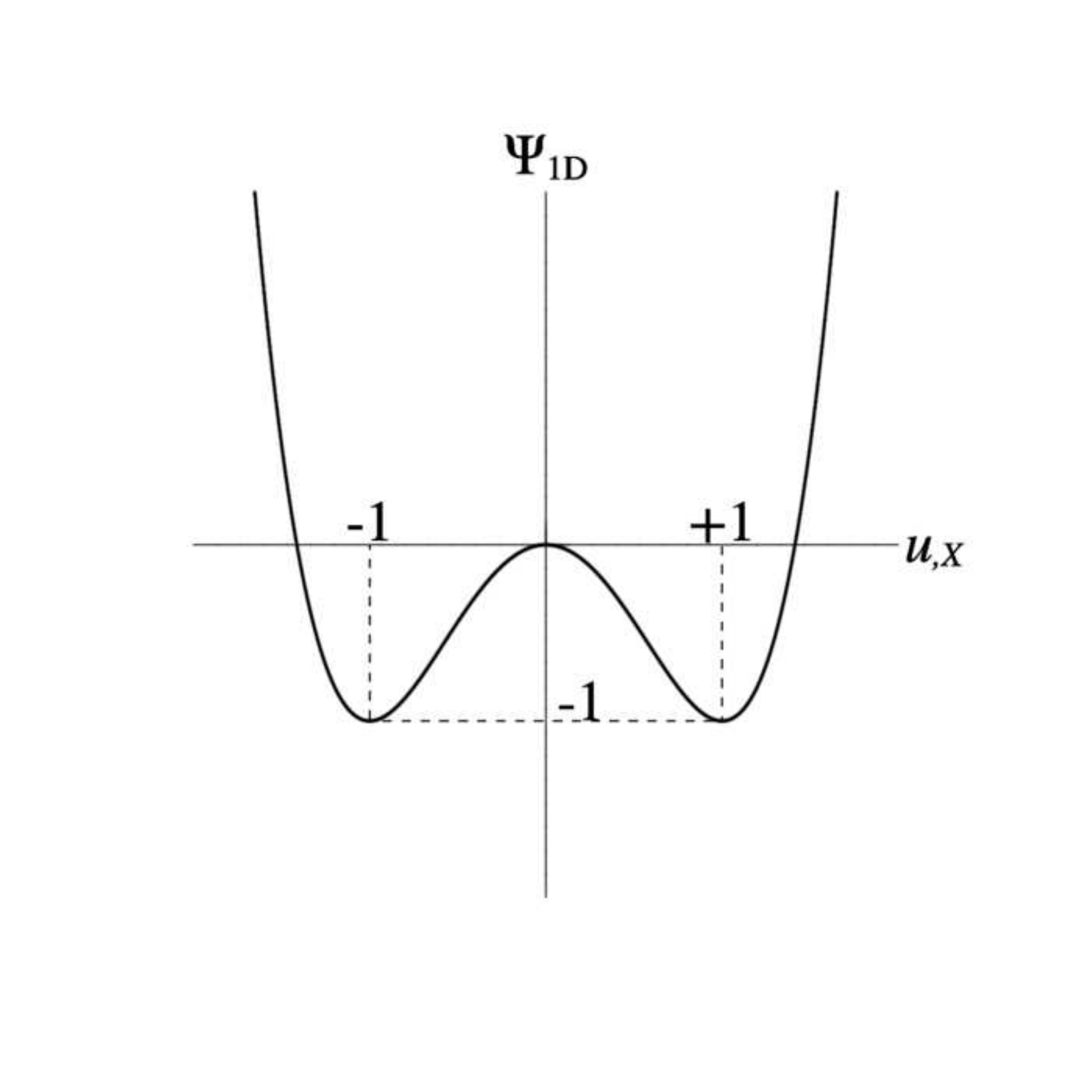}
            \caption{}
            \label{Fi:plot_free-energy-density}
        \end{subfigure}
        ~
        \begin{subfigure}[b]{5.5cm}
            \centering
            \includegraphics[scale=0.25]{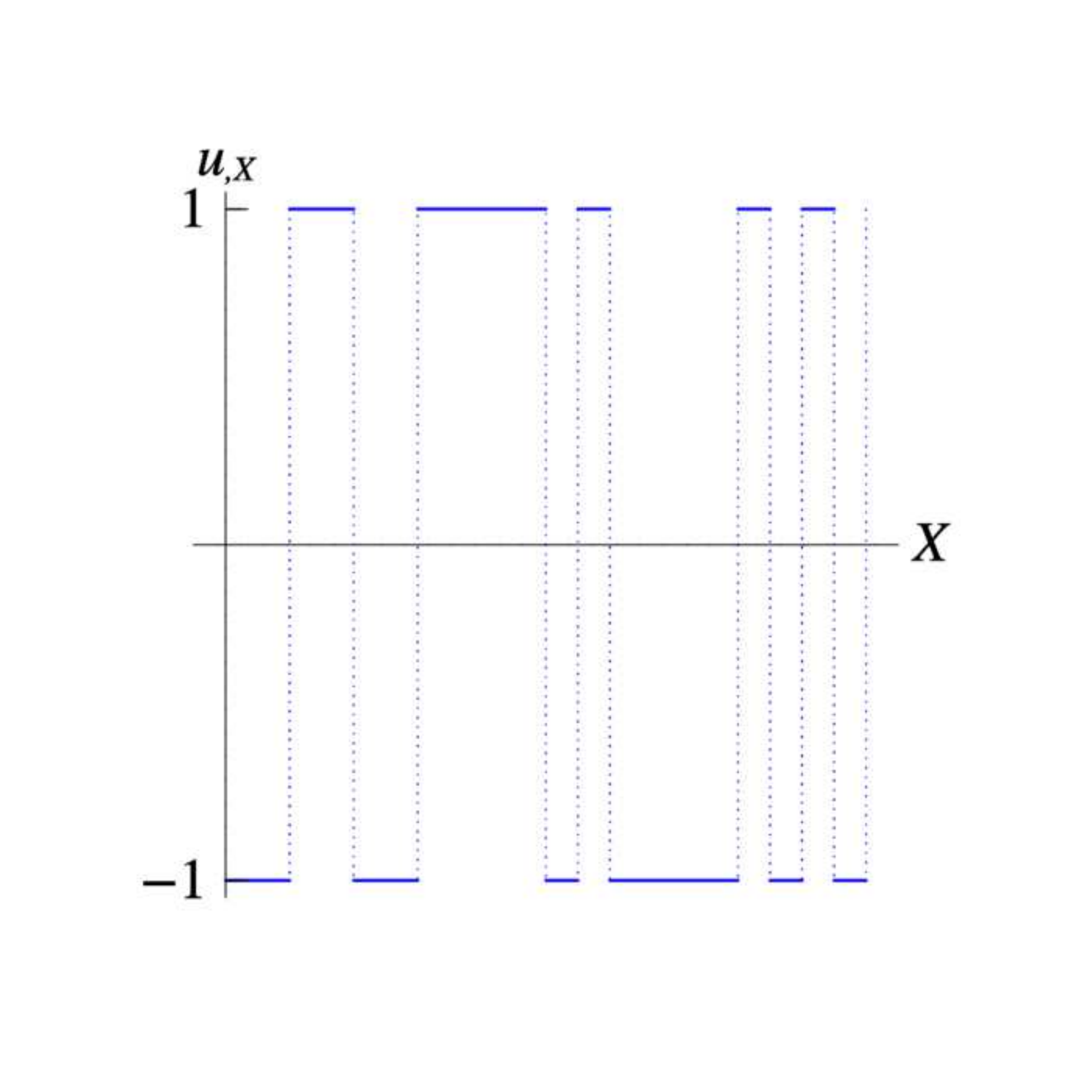}
            \caption{}
            \label{Fi:ux_1d_alt}
        \end{subfigure}
\caption{ Plots of (\subref{Fi:plot_free-energy-density}) free energy density function $\Psi_{\text{1D}}$ projected onto the $\Psi_{\text{1D}}\!\!-\!u_{,X}$ plane and (\subref{Fi:ux_1d_alt}) a typical solution that minimizes the total free energy when $l=0$ in \eqref{E:Psi_1d}.  }
\label{Fi:plot_strain-energy-density_soln}
\end{center}
\end{figure}

The one-dimensional free energy density, $\Psi_{\text{1D}}$, is defined as a function of strain and strain-gradient derived from the solution field $u(X)$ as:\\
\begin{align}
\Psi_{\text{1D}}=(u_{,X}^4-2u_{,X}^2)+l^2 u_{,XX}^2,\label{E:Psi_1d}
\end{align}
where $l$ is the strain gradient length-scale parameter.  
This energy density \eqref{E:Psi_1d} is non-convex with respect to the strain component $u_{,X}$; see Fig. \ref{Fi:plot_free-energy-density}.  
In the absence of the strain-gradient contribution; i.e., with $l=0$ in \eqref{E:Psi_1d}, this non-convex density function characterizes fields $u$ that are composed purely of two variants, one with $u_{,X}=-1$ and the other with $u_{,X}=+1$; see Fig. \ref{Fi:ux_1d_alt}.
In this setting laminae (sub-domains) of these two variants form with arbitrary size, and in principle, infinitely fine microstructures can develop.  The length scale parameter $l$ precludes the existence of such twinned microstructures of infinite fineness by penalizing the interfaces between them. It also introduces a characteristic length scale to the problem.

We seek solution fields $u \in \mathcal{S}_{1D}$ on $\overline{\Omega}_{1D}$, where $\Omega_{1D}=(0,1)$, that satisfy standard and higher-order Dirichlet boundary conditions:
\begin{subequations}
\begin{align}
u=0,\; u_{,X}=0 \quad &\text{ on } X=0,\\
u=d,\; u_{,X}=0 \quad &\text{ on } X=1,
\end{align}
\label{E:bc_u_1D}%
\end{subequations}
where $d=2^{-10}$, and globally/locally minimize the total free energy corresponding to the density function \eqref{E:Psi_1d}:
\begin{align}
\Pi_{\text{1D}}:=\int_{\Omega_{\text{1D}}}\Psi_{\text{1D}}\dif\Omega_{\text{1D}}.
\end{align}
To this end, we define admissible test functions $w \in \mathcal{V}_{1D}$ that satisfy:
\begin{subequations}
\begin{align}
w=0,\; w_{,X}=0 \quad &\text{ on } X=0,\\
w=0,\; w_{,X}=0 \quad &\text{ on } X=1,
\end{align}
\label{E:bc_w_1D}%
\end{subequations}
and solve the following weak form of the boundary value problem derived from the variational argument: Find $u \in \mathcal{S}_{1D}$ such that $\forall w \in\mathcal{V}_{1D}$,
\begin{align}
D\Pi_{\text{1D}}[w]=\int_{\Omega_{\text{1D}}}(w_{,X}P+w_{,XX}B)\dif\Omega_{\text{1D}}=0,\label{E:D_Pi_1D}
\end{align}
where $P$ is the first Piola-Kirchhoff stress and $B$ is the higher-order stress, defined as:
\begin{subequations}
\begin{align}
P&:=\partial\Psi_{\text{1D}}/\partial u_{,X},\label{P1D}\\
B&:=\partial\Psi_{\text{1D}}/\partial u_{,XX}\label{B1D}.
\end{align}
\end{subequations}
Stability of the solutions is then assessed by examining the positive definiteness of the second variation:
\begin{align}
D^2\Pi_{\text{1D}}[w,w]=\int\limits_{\Omega_{\text{1D}}}\left(
w_{,X}\frac{\partial^2\Psi_{\text{1D}}}{\partial u_{,X}\partial u_{,X}}w_{,X}+
w_{,X}\frac{\partial^2\Psi_{\text{1D}}}{\partial u_{,X}\partial u_{,XX}}w_{,XX}+
w_{,XX}\frac{\partial^2\Psi_{\text{1D}}}{\partial u_{,XX}\partial u_{,X}}w_{,X}+
w_{,XX}\frac{\partial^2\Psi_{\text{1D}}}{\partial u_{,XX}\partial u_{,XX}}w_{,XX}
\right)\dif\Omega_{\text{1D}}>0,\label{E:D2_Pi_1D}
\end{align}
where we have explicitly retained the symmetric second and third terms of the integrand for clarity of the development. Note that weak form \eqref{E:D_Pi_1D} and \eqref{E:bc_u_1D} with \eqref{E:bc_w_1D} leads us to the following strong form:
\begin{align*}
-P_{,X}+B_{,XX}=0\quad\text{ on }\Omega_{\text{1D}},\label{E:strong_1d}
\end{align*}
which possesses fourth-order spatial derivatives due to the constitutive relations \eqref{E:Psi_1d}, \eqref{P1D} and \eqref{B1D}.  
This strong form is not investigated further in this work.  

\subsection{The numerical framework}
We seek numerical solutions $u^h\in \mathcal{S}_{1D}^h\subset\mathcal{S}_{1D}$ to a discretized counterpart of the weak from \eqref{E:D_Pi_1D} with the test functions $w^h\in\mathcal{V}_{1D}^h\subset\mathcal{V}_{1D}$ where:
\begin{subequations}
\begin{align}
\mathcal{S}_{1D}^h &= \{v^h \in H^2(\Omega_\text{1D})\vert v^h = 0,\; v^h_{,X} = 0\;\text{on} \; X = 0, \; v^h = d,\; v^h_{,X} = 0\; \text{on} \; X = 1\}, \\
\mathcal{V}_{1D}^h &= \{v^h \in H^2(\Omega_\text{1D})\vert v^h = 0,\; v^h_{,X} = 0\;\text{on} \; X = \{0,1\}\},
\end{align}
\end{subequations}
where $H^2$ represents the standard Sobolev space of integrable functions with integrable first and second derivatives.  
The problem was solved using isogeometric analysis (IGA) with fourth-order B-spline basis functions defined on $1024$ elements of uniform size on $\Omega_{\text{1D}}$; see Cottrell et al. \cite{Cottrell2009} for a comprehensive treatment of  IGA, and Rudraraju et al. \cite{Rudraraju2014} for its application to the current problem framework, but with weak enforcement of Dirichlet boundary conditions.  
We use the quartic-precision floating-point format and solve for solutions up to an absolute tolerance of $10^{-25}$ in the Euclidean norm of the residual corresponding to the discretized version of \eqref{E:D_Pi_1D}.  The higher than typical precision and more stringent tolerance are important to verify convergence to extrema/saddle points of the rapidly fluctuating free energy functional that governs this problem. This rapid fluctuation underlies the existence of families of stable/unstable solutions, which is the crux of this work.
The second variation in \eqref{E:D2_Pi_1D} is discretized on the same B-spline basis, and stability is assessed by examining the positive definiteness of the resulting symmetric Hessian matrix using the eigenvalue solver \textsf{FEAST v3.0} \cite{Polizzi2009} with the relative accuracy of $O(10^{-8})$.  
Figures were produced using \textsf{mathgl 2.3.0}.  

\subsection{Solution branches and branch-tracking}
\begin{figure}
    \begin{center}
        \begin{tabular}{rp{16.5cm}}
            \parbox[t]{0.5cm}{ }&
            \begin{tabular}{p{3.5cm}p{3.5cm}p{3.5cm}p{3.5cm}p{2.5cm}}
                \hspace{1.2cm}$l\!=\!0.01$ & 
                \hspace{1.2cm}$l\!=\!0.02$ &
                \hspace{1.2cm}$l\!=\!0.04$ &
                \hspace{1.2cm}$l\!=\!0.08$ &
                \vspace{0.5\baselineskip}
            \end{tabular}  \\
            \parbox[t]{0.5cm}{ $A$ } &
            \begin{tabular}{p{3.5cm}p{3.5cm}p{3.5cm}p{3.5cm}p{2.5cm}}
                \includegraphics[scale=0.15]{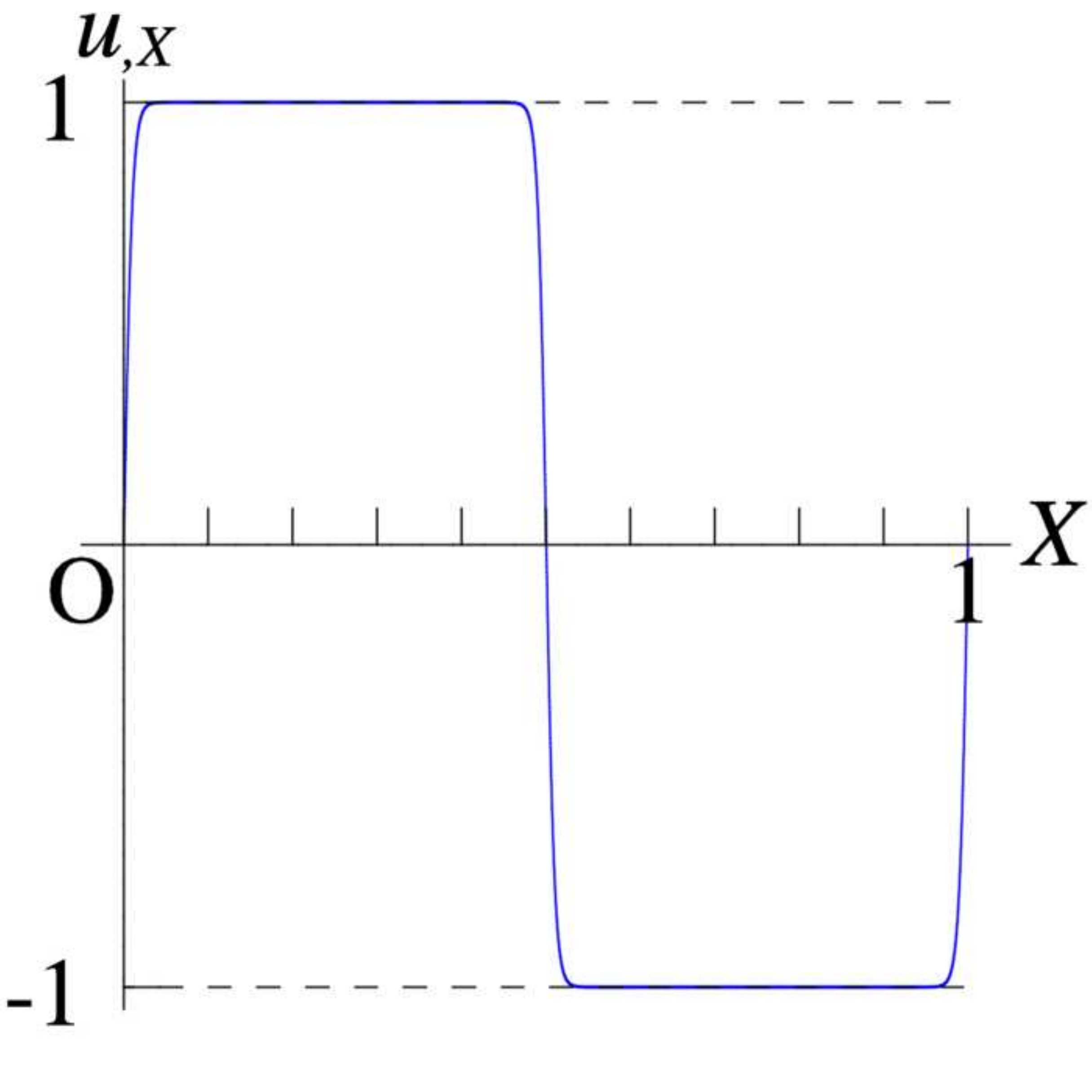} & 
                \includegraphics[scale=0.15]{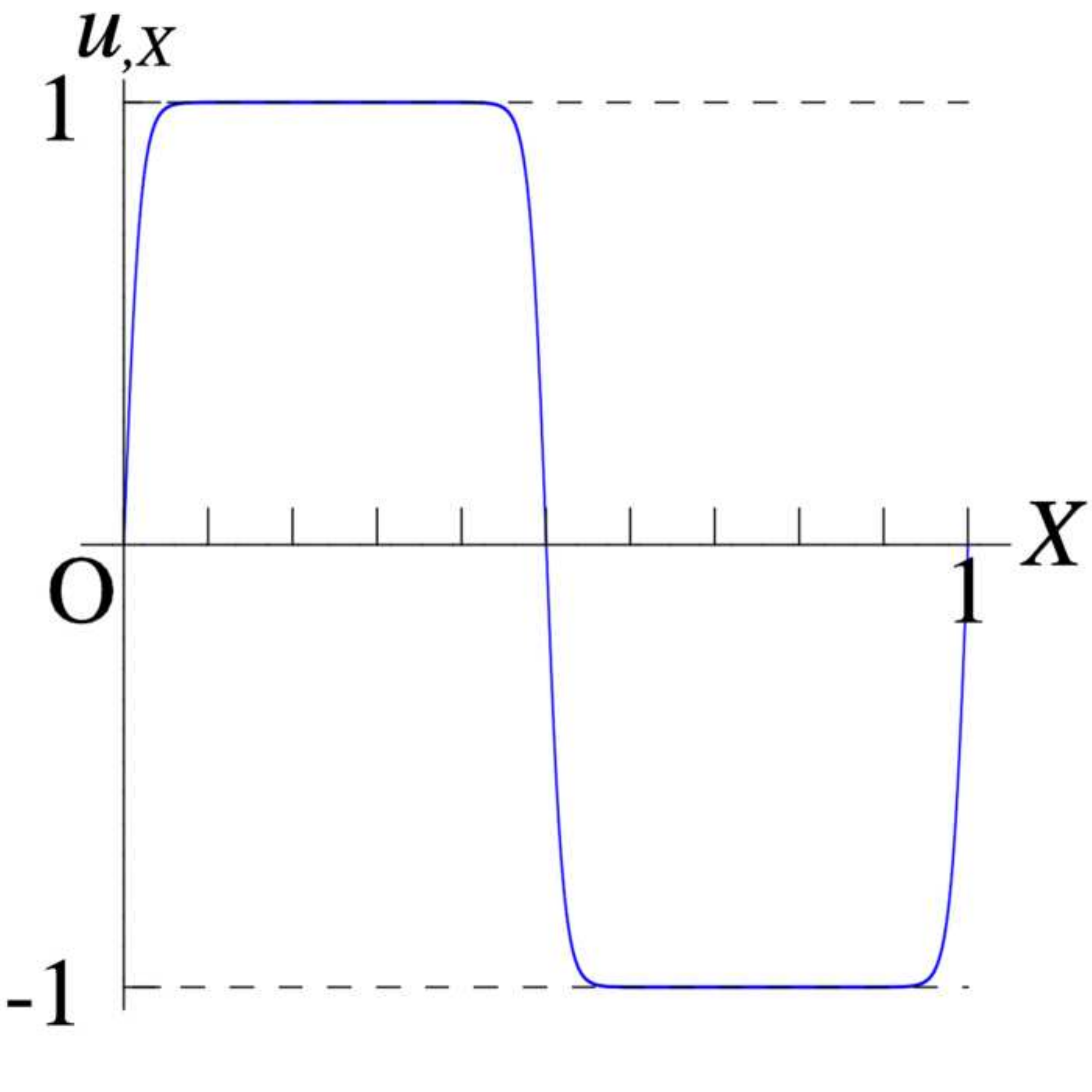} & 
                \includegraphics[scale=0.15]{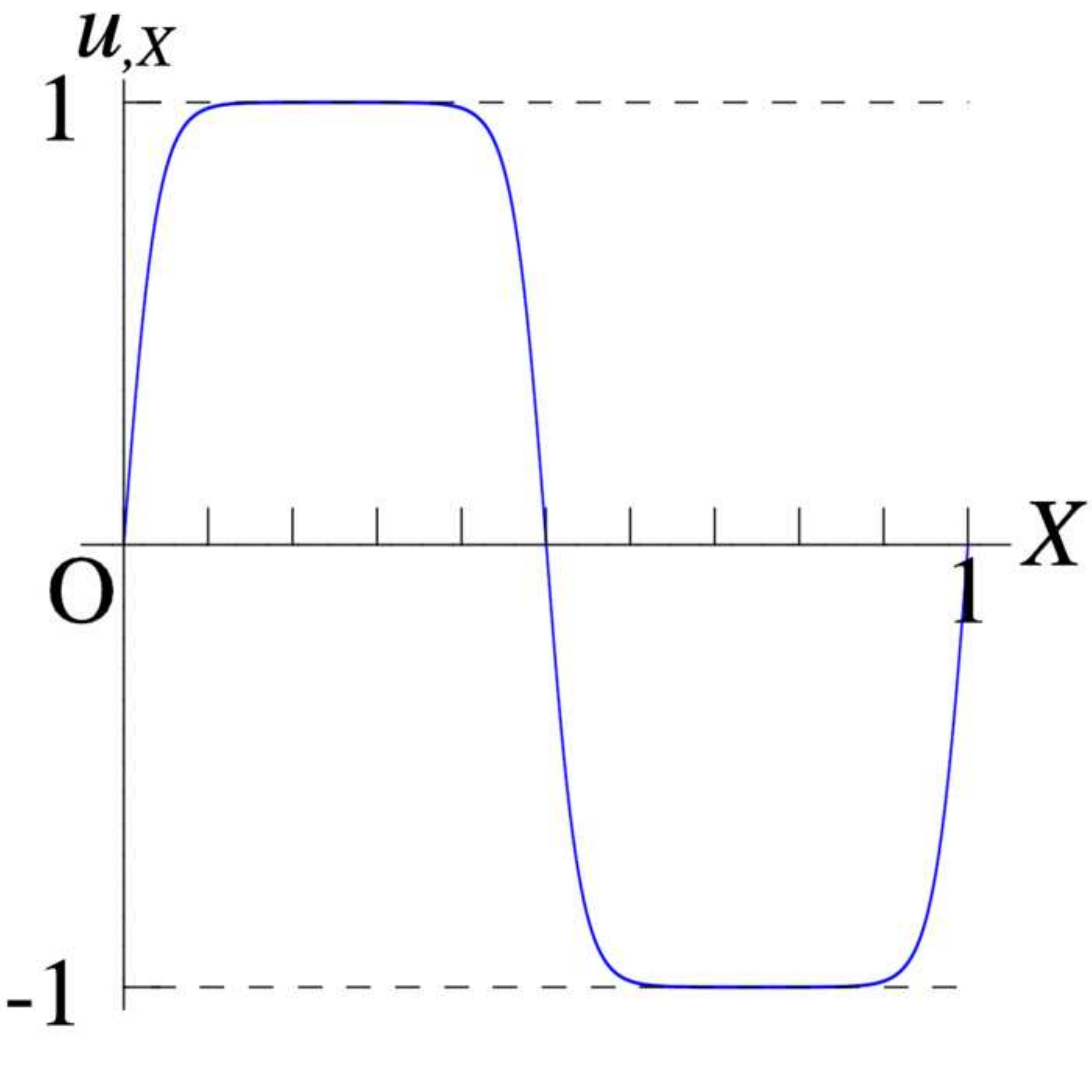} & 
                \includegraphics[scale=0.15]{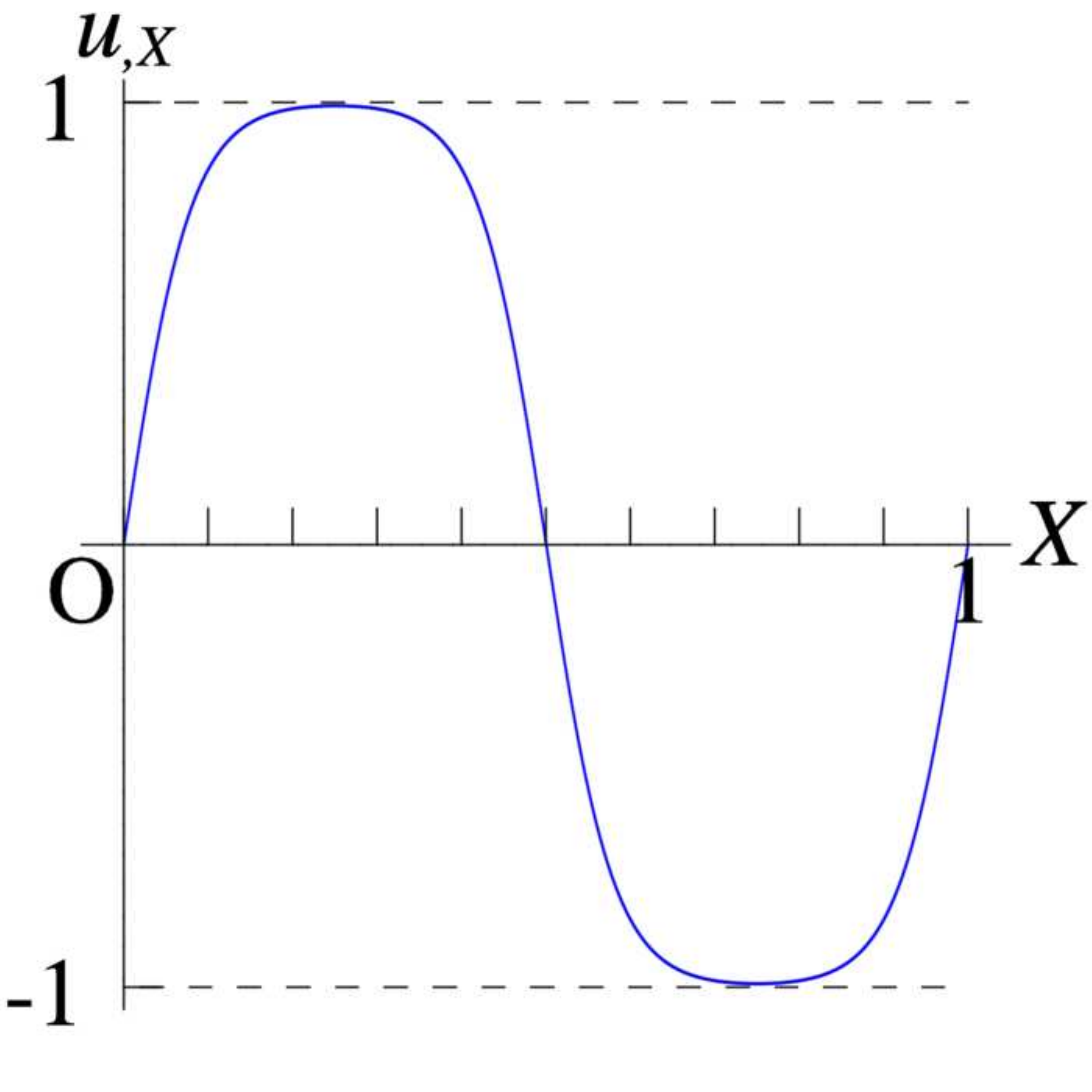} & 
            \end{tabular}  \\
            \parbox[t]{0.5cm}{$B$} &
            \begin{tabular}{p{3.5cm}p{3.5cm}p{3.5cm}p{3.5cm}p{2.5cm}}
                \includegraphics[scale=0.15]{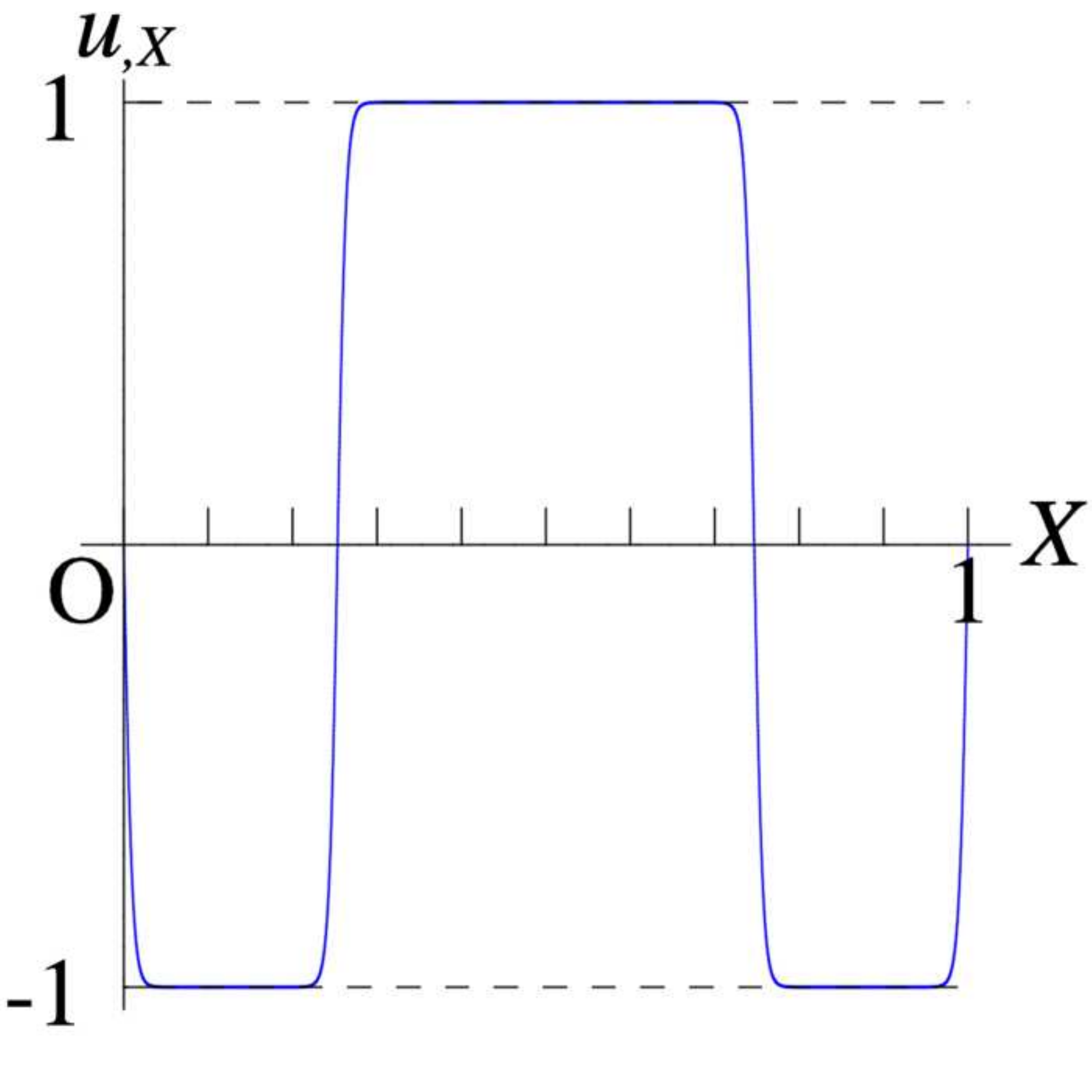} & 
                \includegraphics[scale=0.15]{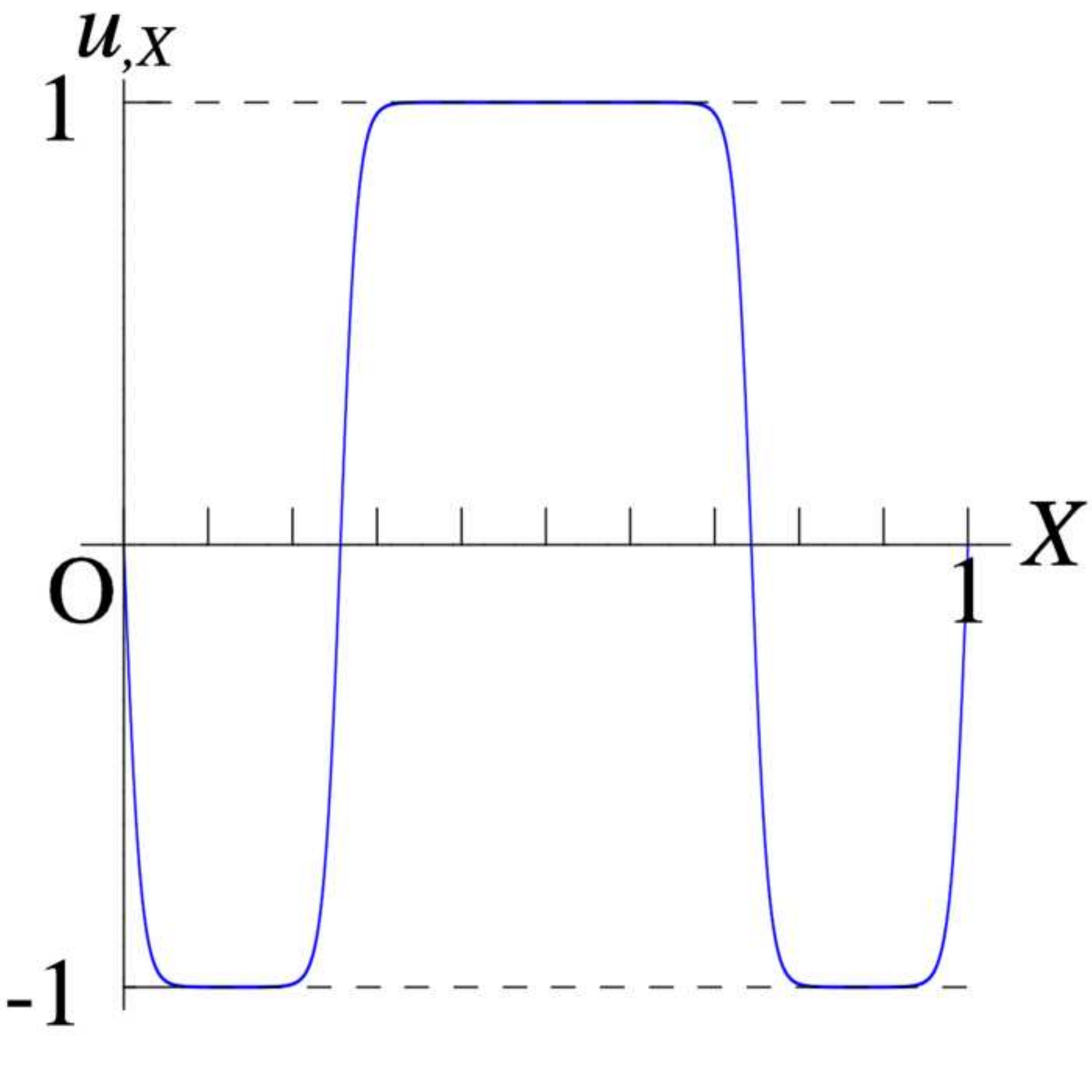} & 
                \includegraphics[scale=0.15]{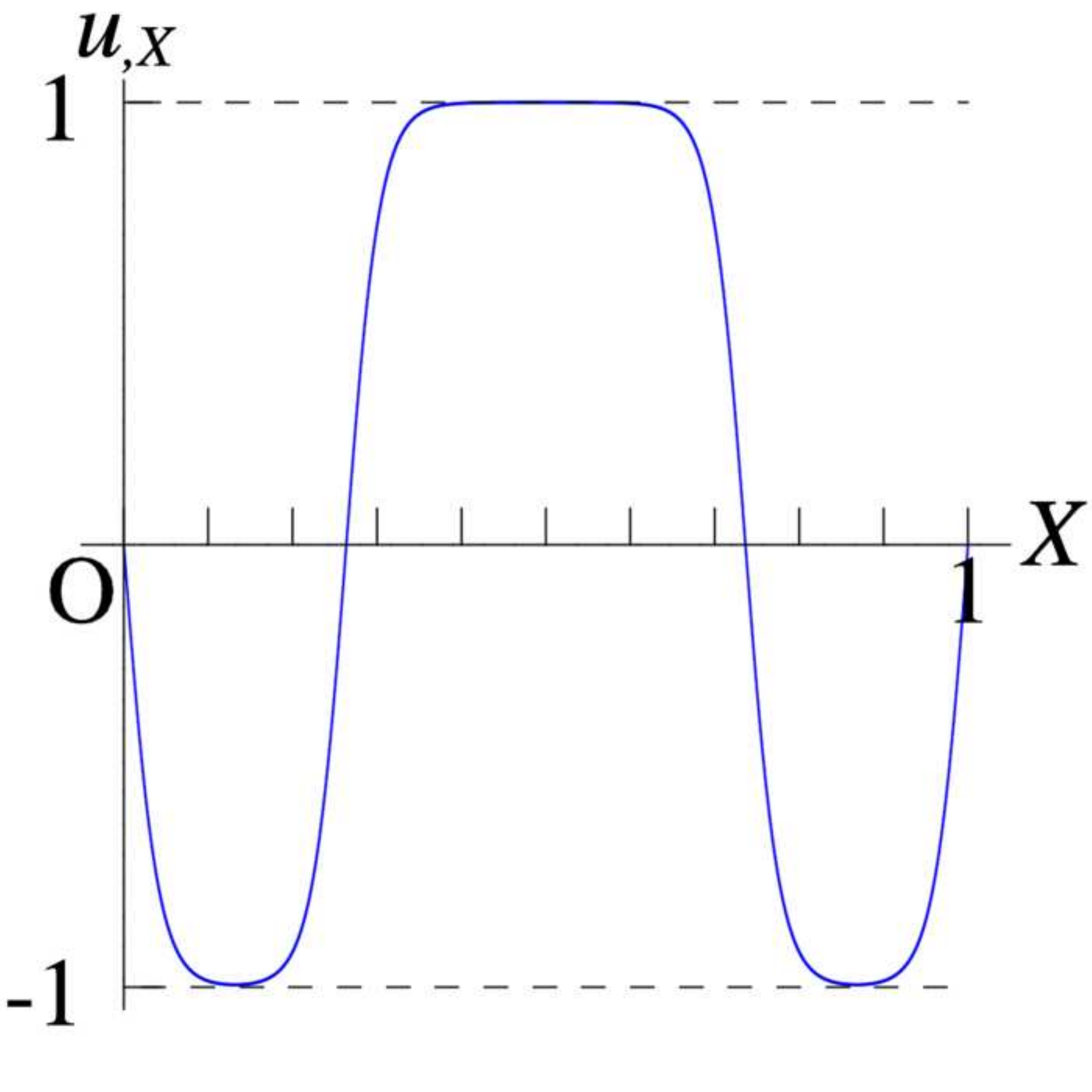} & 
                \includegraphics[scale=0.15]{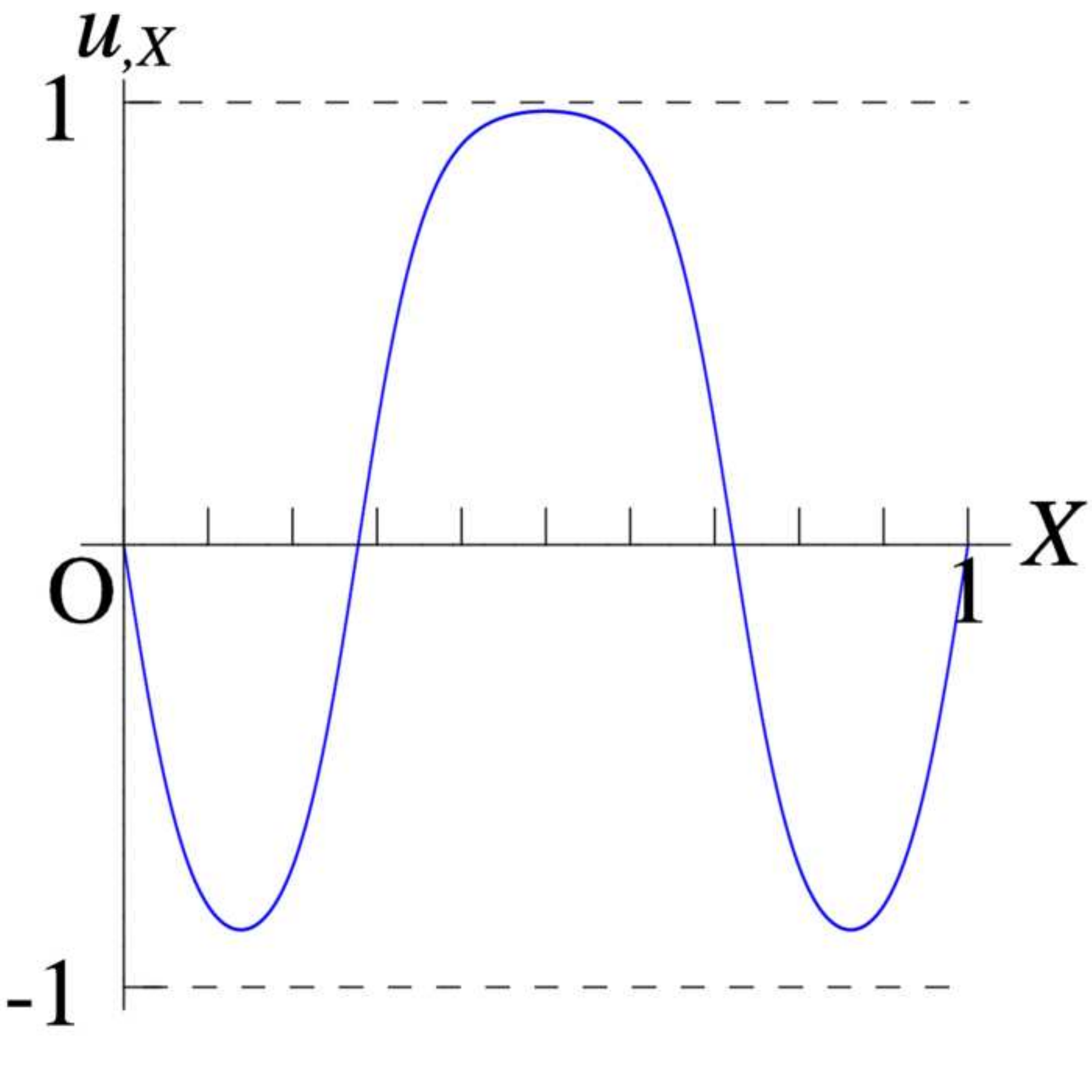} & 
            \end{tabular}  \\
            \parbox[t]{0.5cm}{ $C$ } &
            \begin{tabular}{p{3.5cm}p{3.5cm}p{3.5cm}p{3.5cm}p{2.5cm}}
                \includegraphics[scale=0.15]{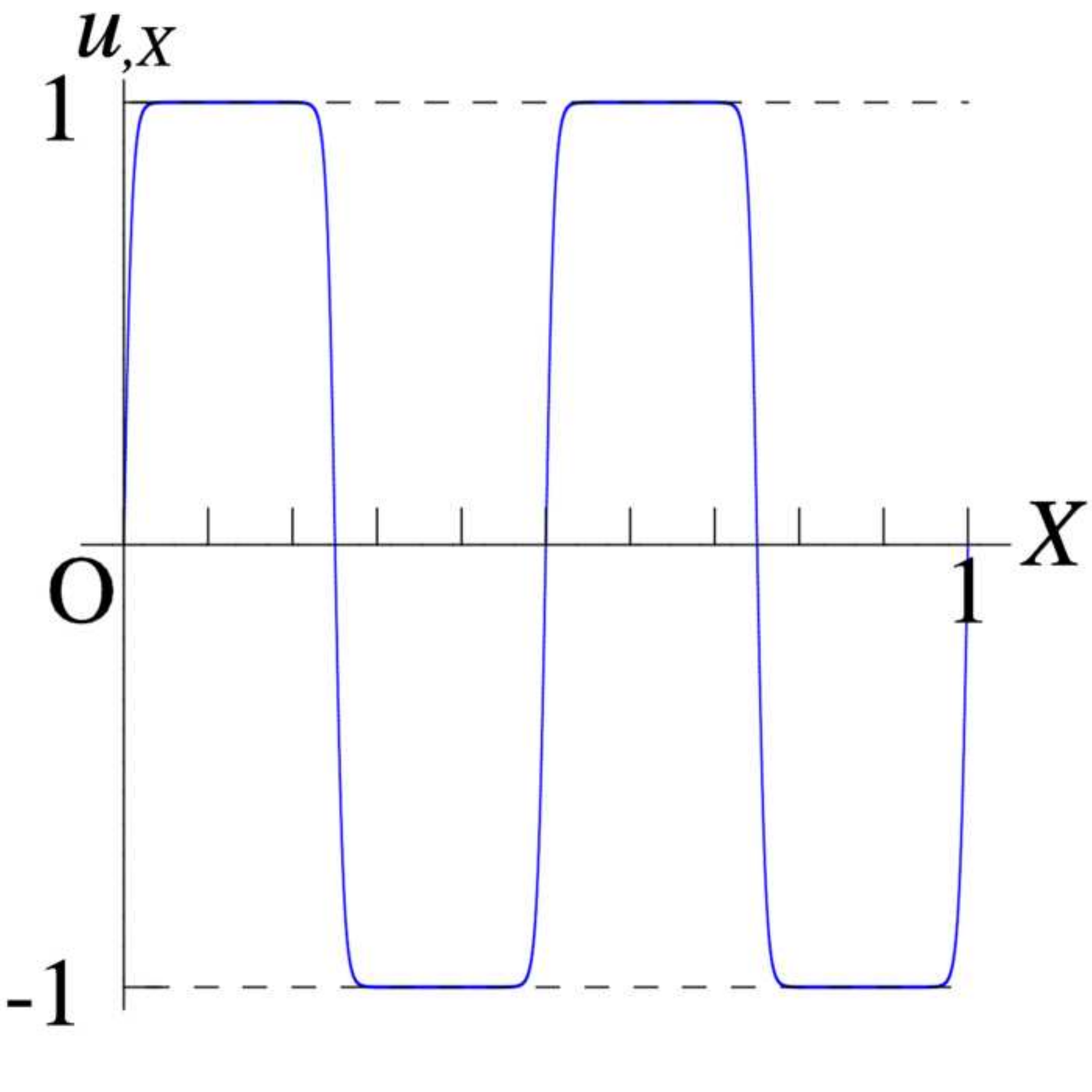} & 
                \includegraphics[scale=0.15]{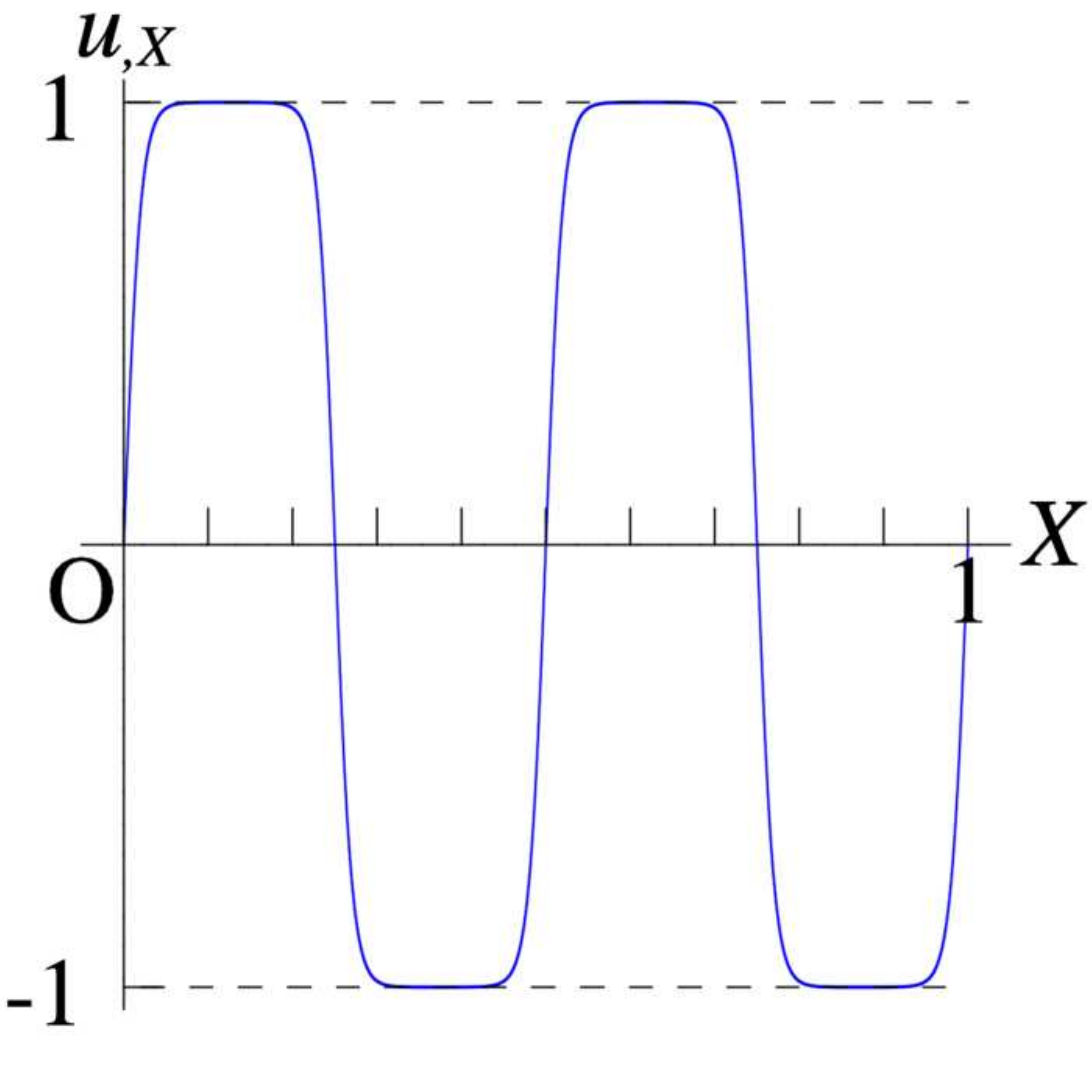} & 
                \includegraphics[scale=0.15]{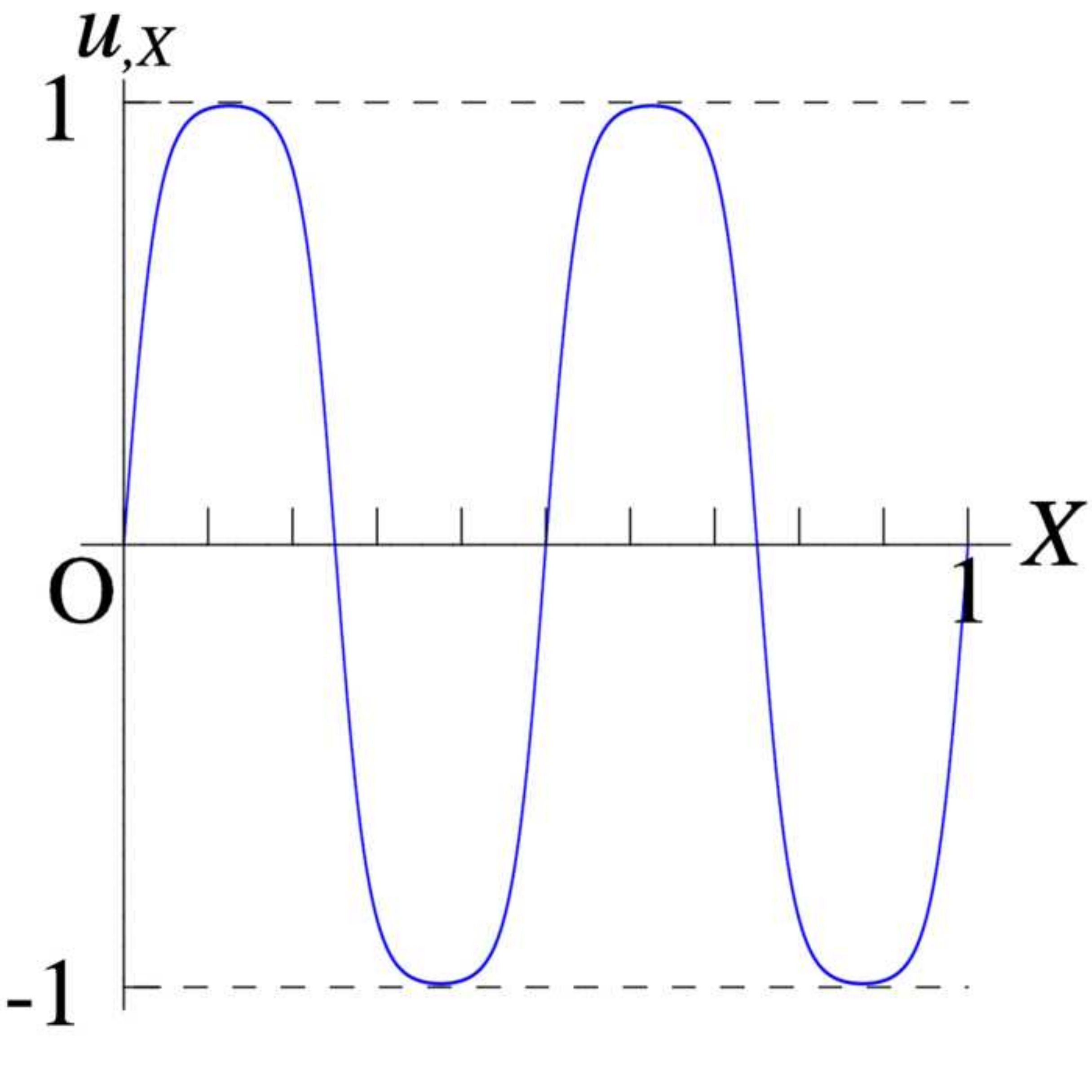} & 
                \includegraphics[scale=0.15]{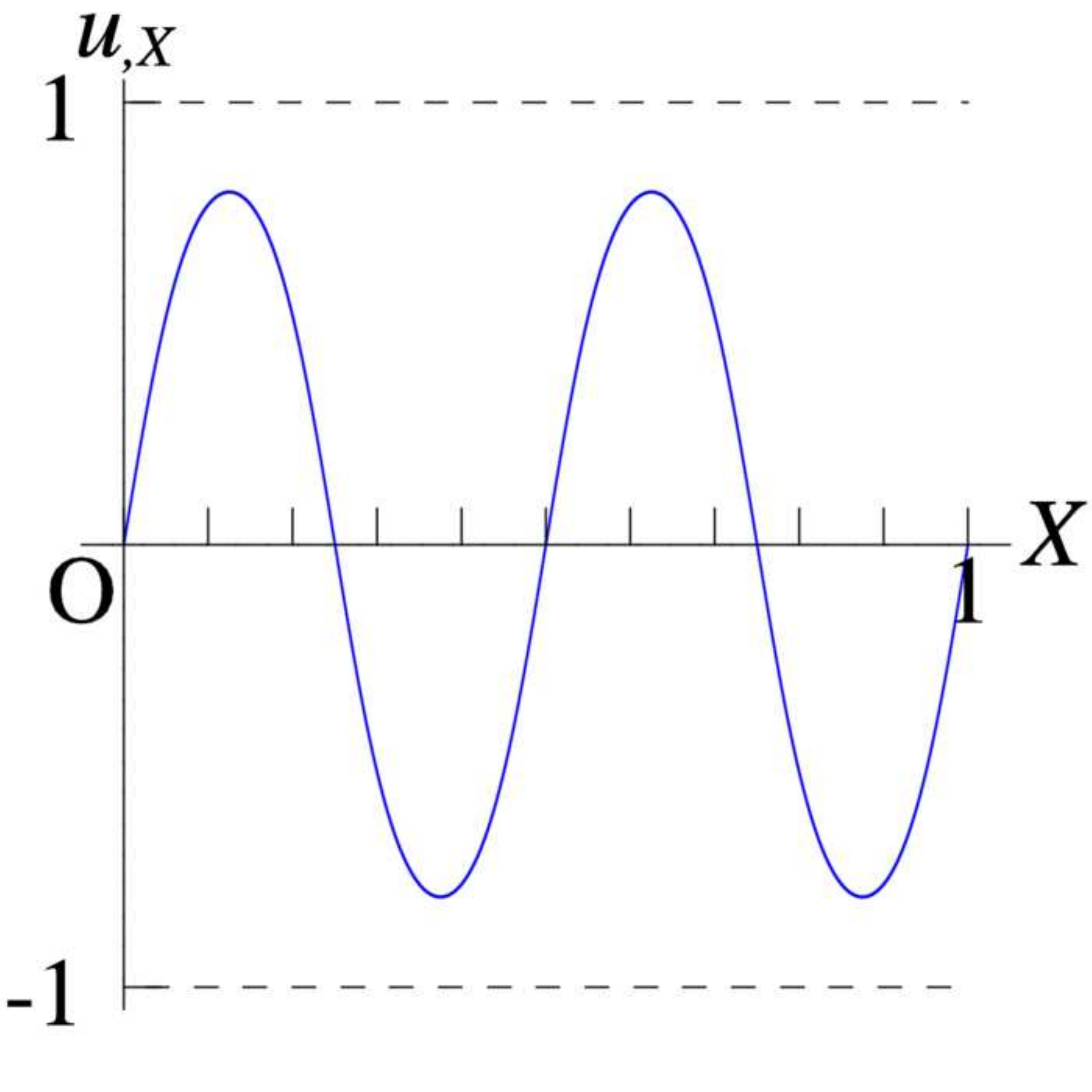} & 
            \end{tabular}  \\
            \parbox[t]{0.5cm}{ $D$ } &
            \begin{tabular}{p{3.5cm}p{3.5cm}p{3.5cm}p{3.5cm}p{2.5cm}}
                \includegraphics[scale=0.15]{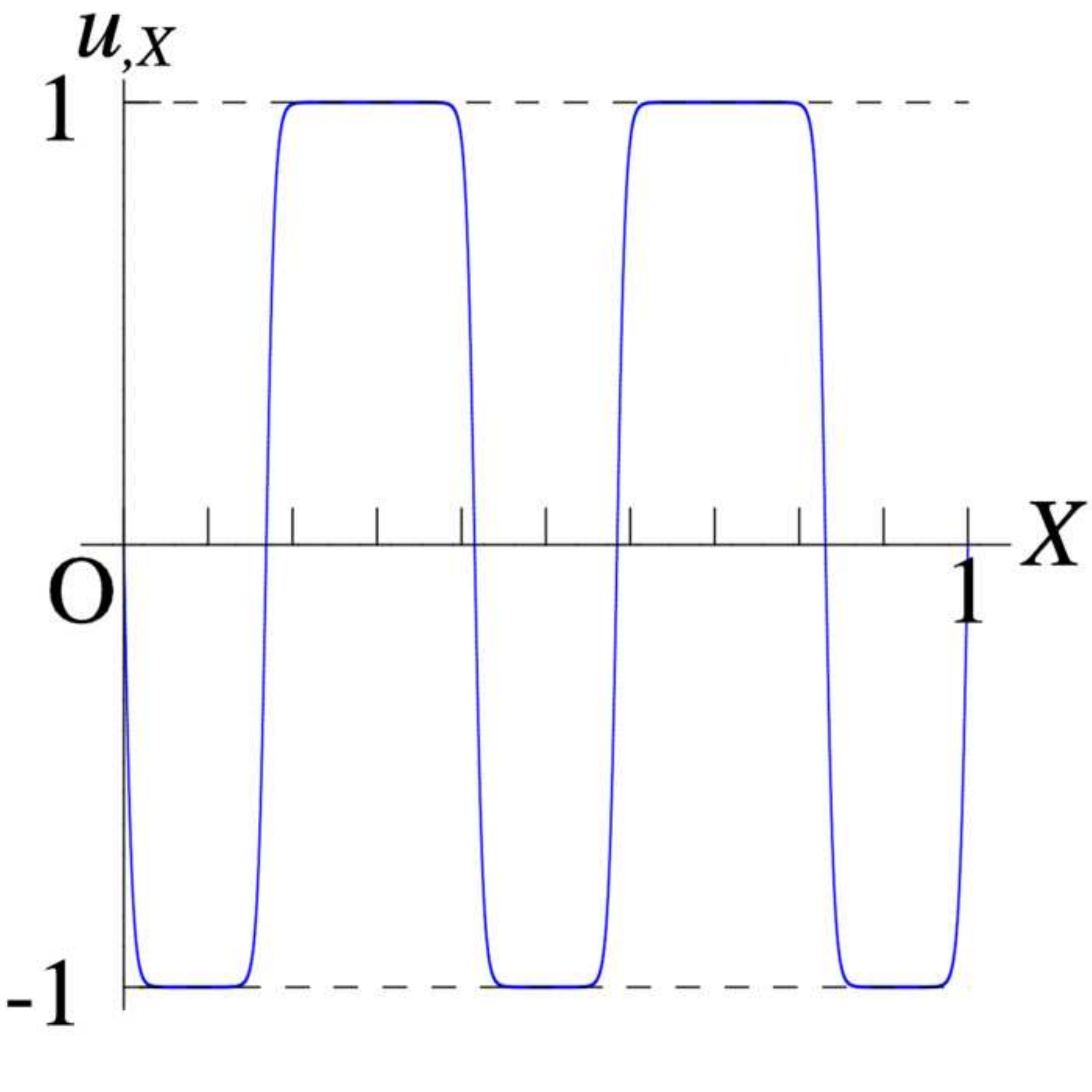} & 
                \includegraphics[scale=0.15]{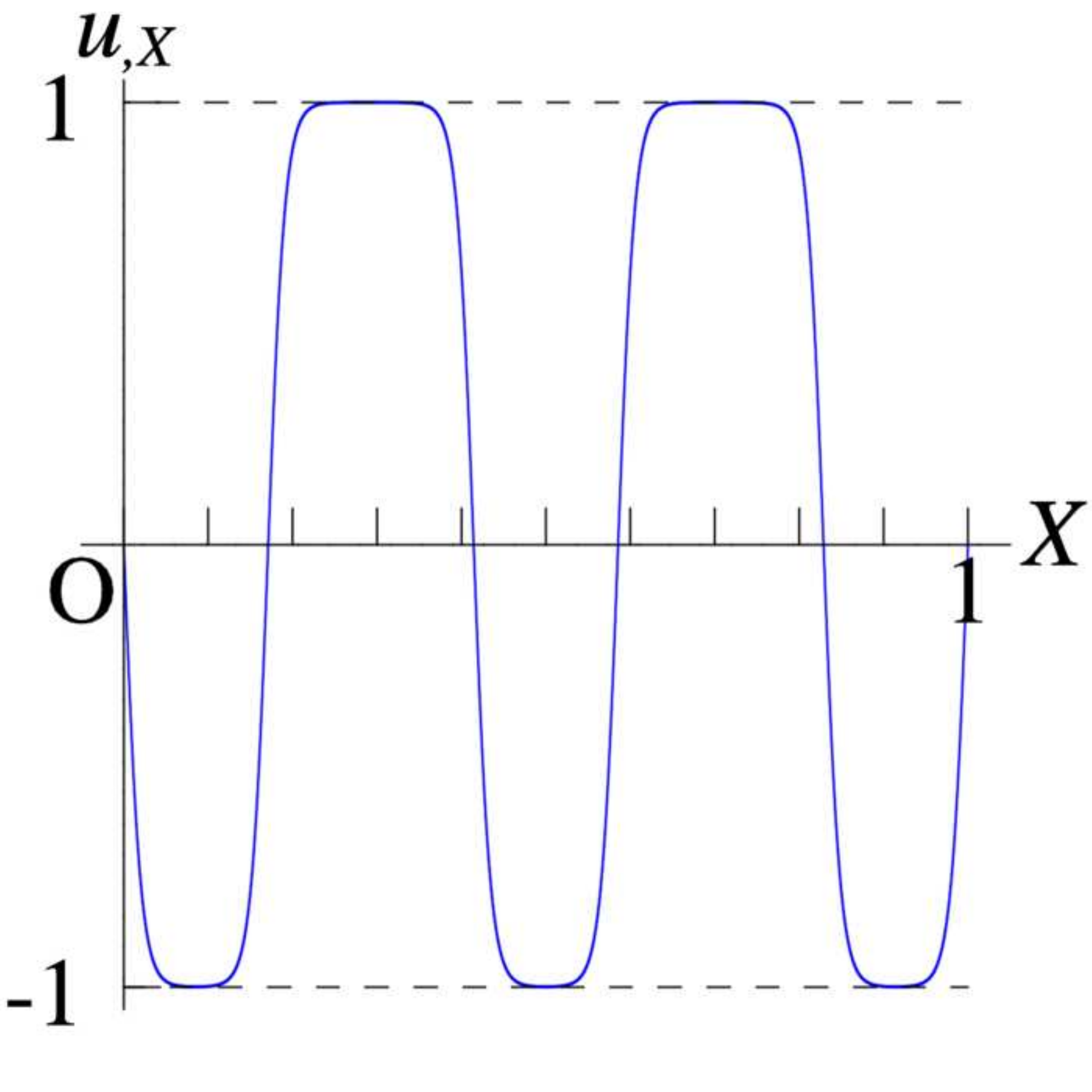} & 
                \includegraphics[scale=0.15]{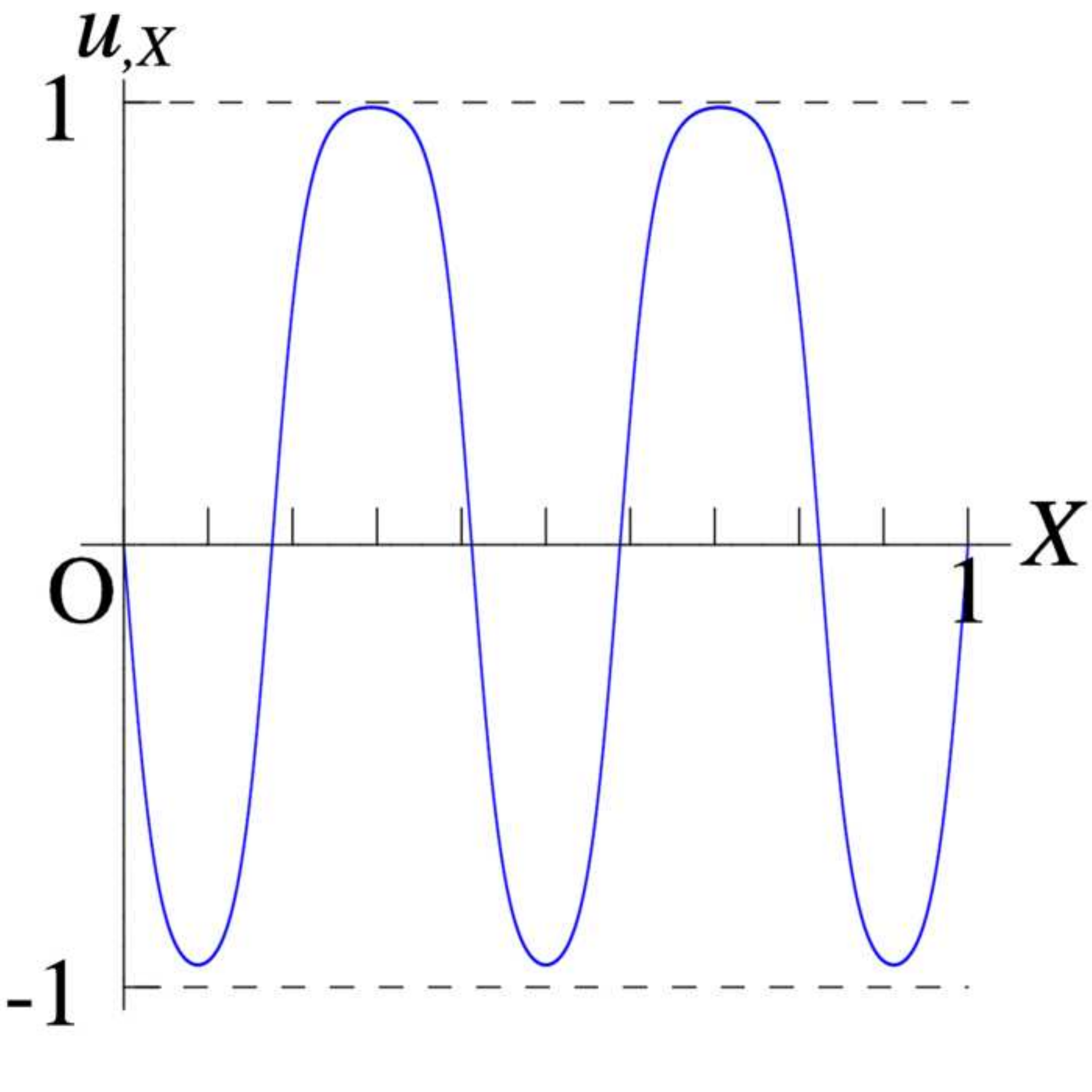} & 
                \includegraphics[scale=0.15]{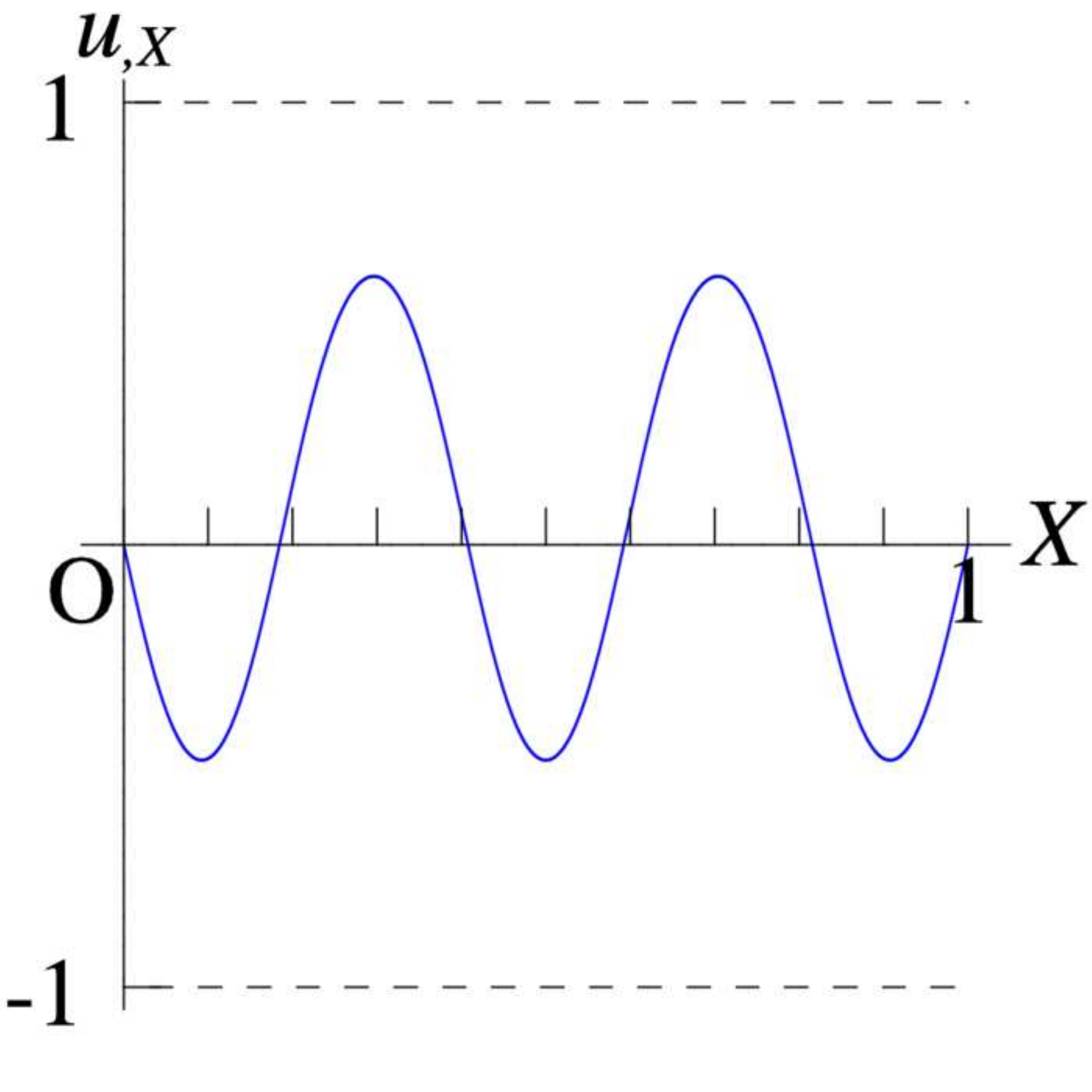} & 
            \end{tabular}  \\
            \parbox[t]{0.5cm}{ $E$ } &
            \begin{tabular}{p{3.5cm}p{3.5cm}p{3.5cm}p{3.5cm}p{2.5cm}}
                \includegraphics[scale=0.15]{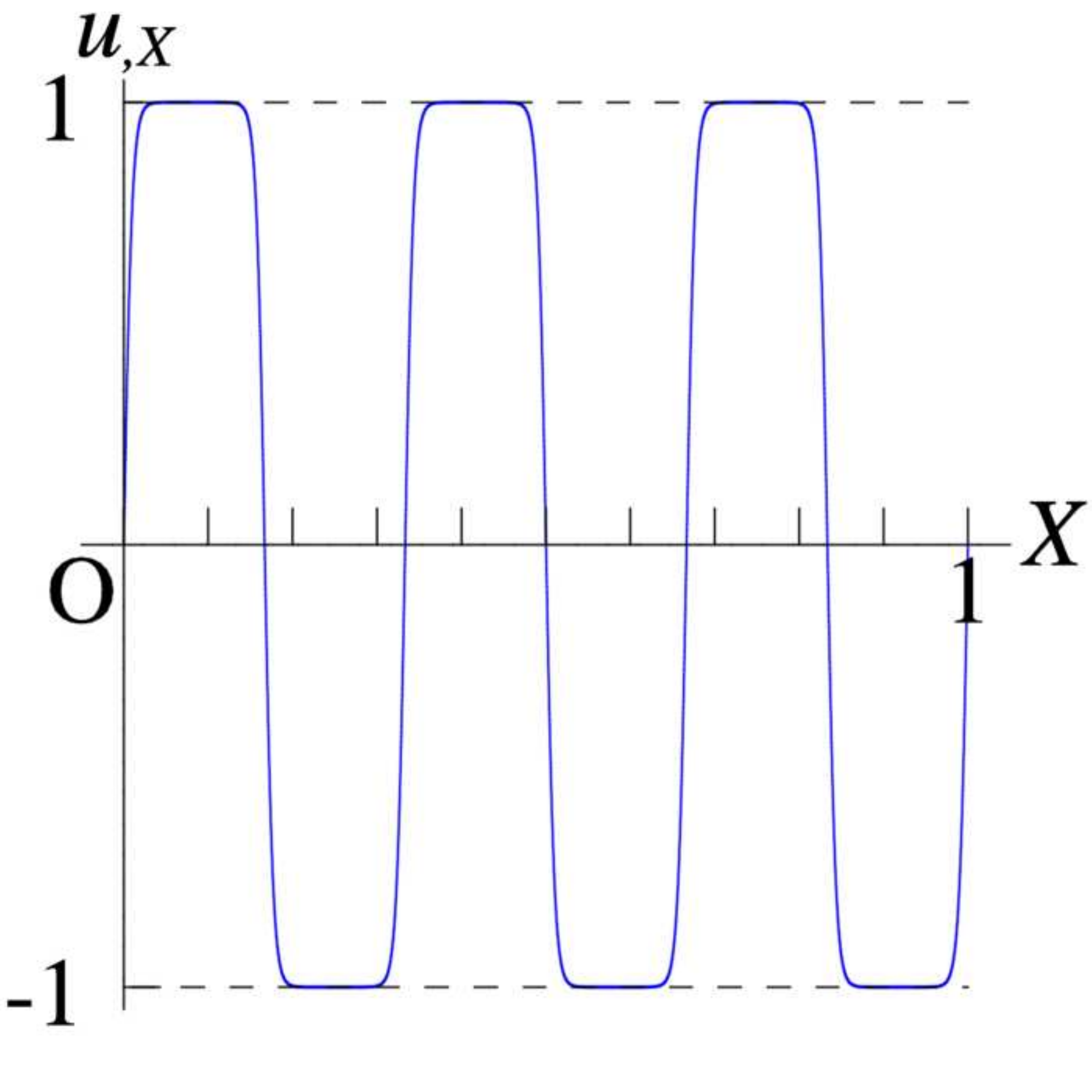} & 
                \includegraphics[scale=0.15]{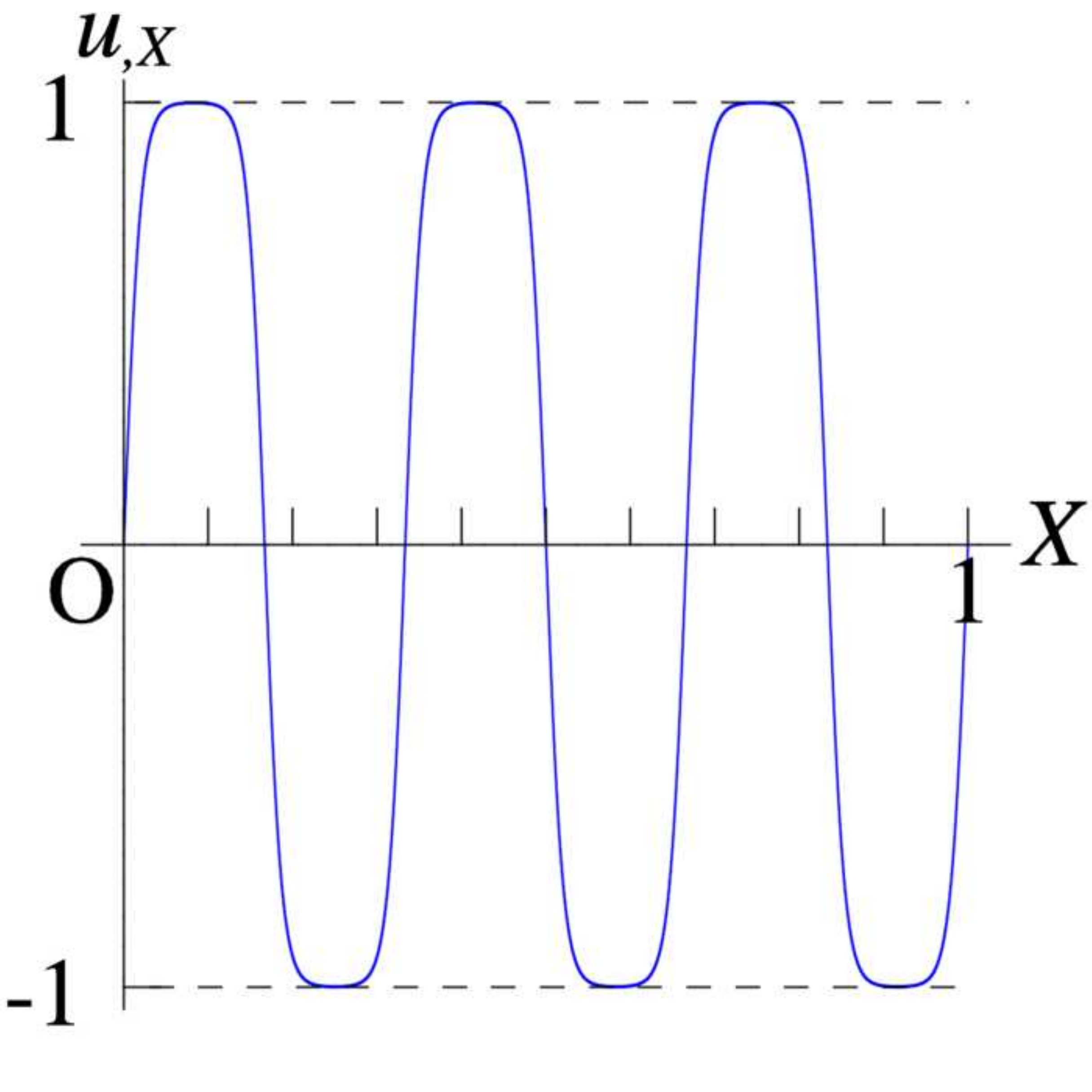} & 
                \includegraphics[scale=0.15]{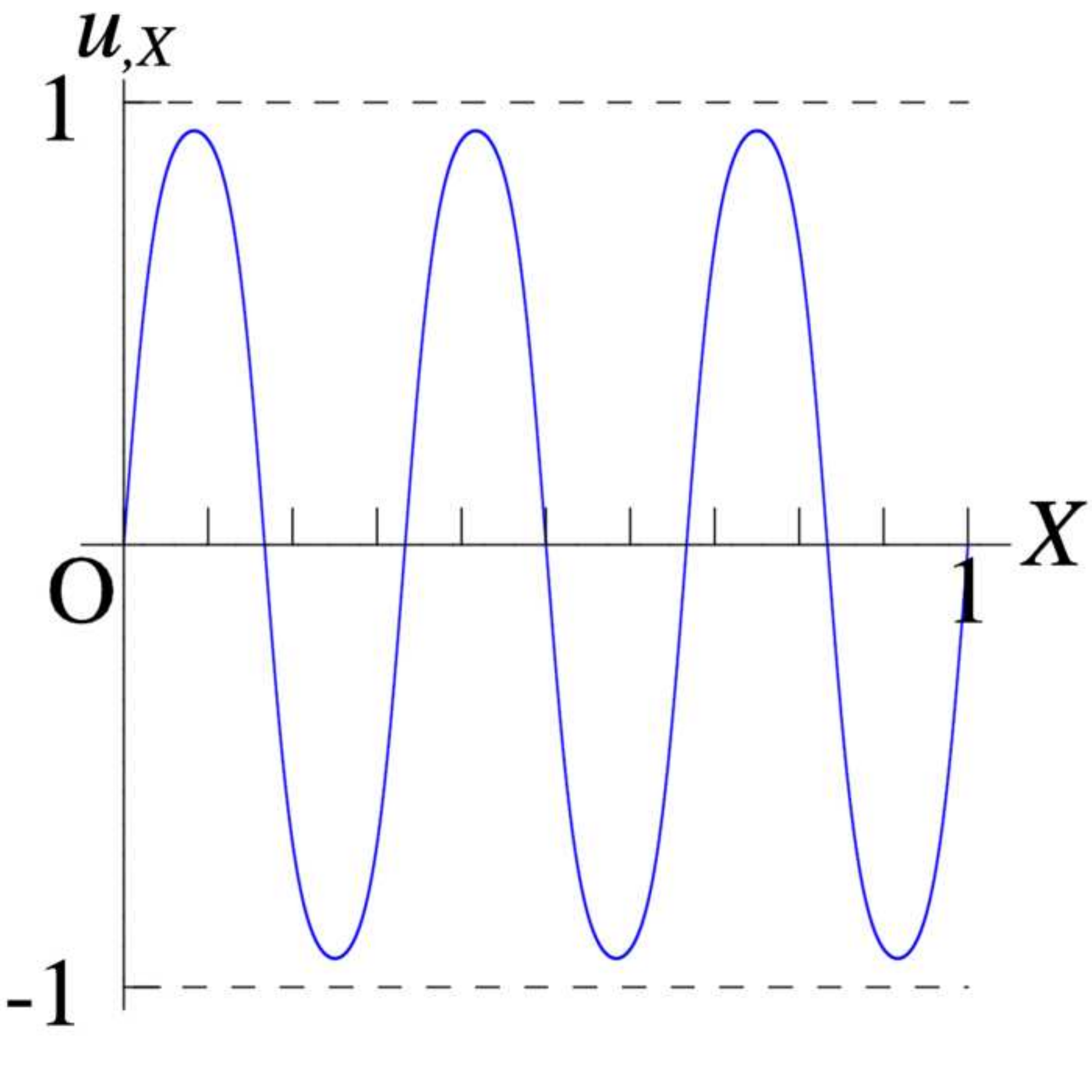} & 
                & 
            \end{tabular}  \\
            \parbox[t]{0.5cm}{ $F$ } &
            \begin{tabular}{p{3.5cm}p{3.5cm}p{3.5cm}p{3.5cm}p{2.5cm}}
                \includegraphics[scale=0.15]{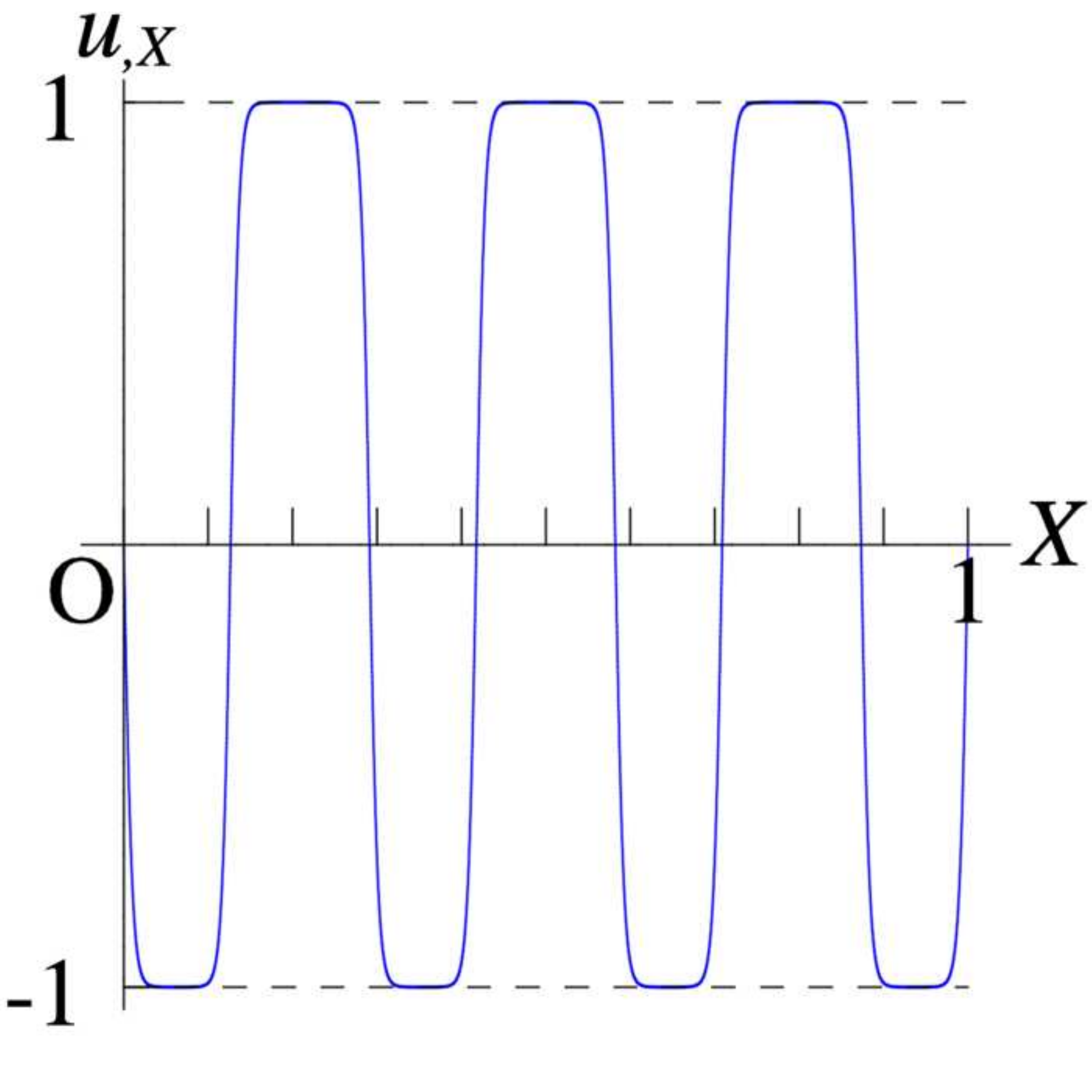} & 
                \includegraphics[scale=0.15]{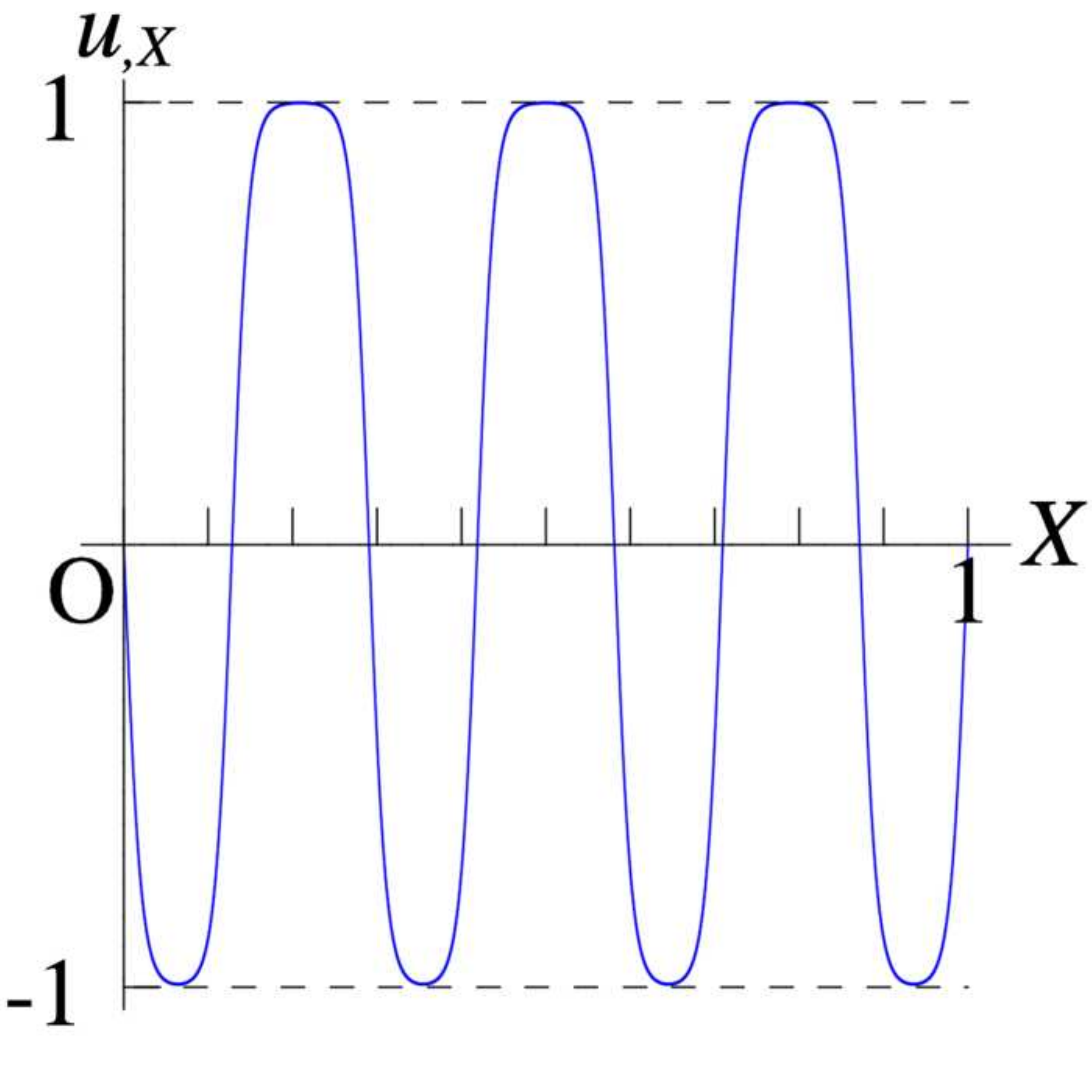} & 
                \includegraphics[scale=0.15]{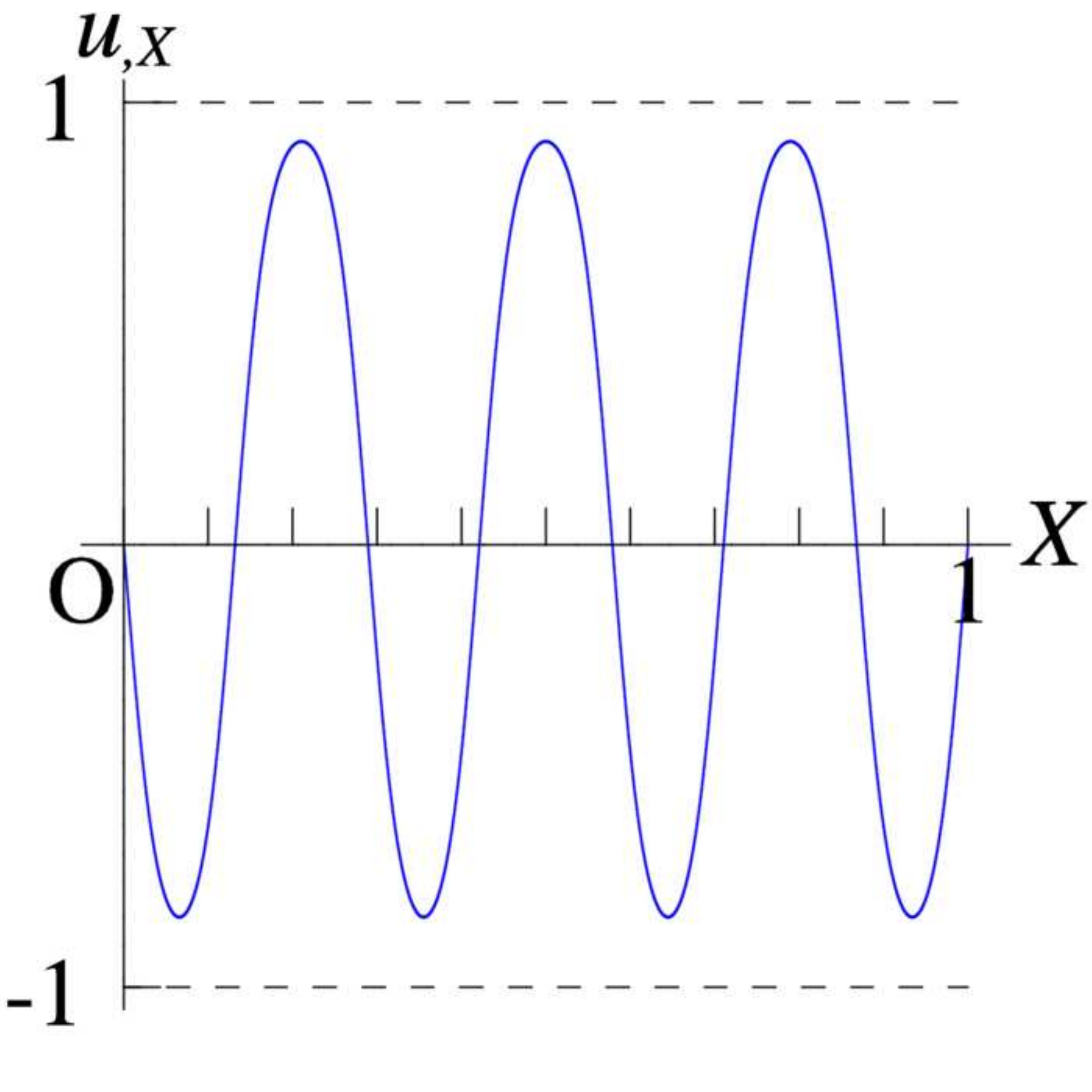} & 
                \includegraphics[scale=0.15]{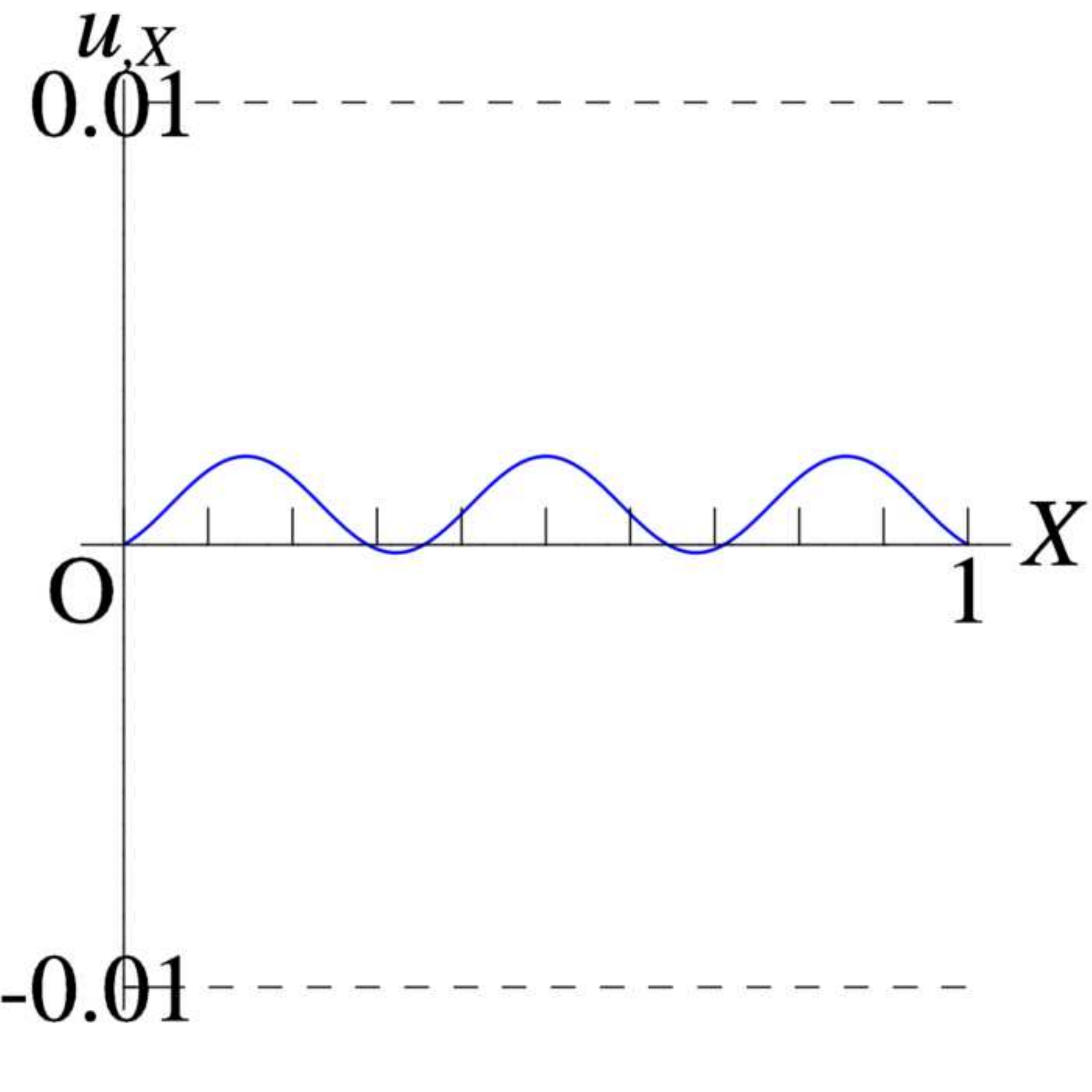} & 
            \end{tabular}
        \end{tabular}
    \end{center}
    \caption{Plots of computed strain $u_{,X}$ for branches A - F at selected values of $l$.  }
    \label{Fi:uX-X}
\end{figure}

Fig. \ref{Fi:PI-l_1d_orig} shows the total free energy, $\Pi_{\text{1D}}$, of solutions to the boundary value problem, plotted against the length scale parameter $l$, where six representative branches are labeled as A - F.  
Numerically computed strains, $u_{,X}$, are plotted against $X$ in Fig. \ref{Fi:uX-X} for branches A-F at selected values of $l$.  

Here, we outline the procedure that we used to obtain the branches shown in Fig. \ref{Fi:PI-l_1d_orig}.  
In solving the nonlinear boundary value problem \eqref{E:D_Pi_1D} and \eqref{E:bc_u_1D} with \eqref{E:bc_w_1D}, we observed that the homogeneous initial guess $u_{\text{init}}=0$ always captures the branch of highest energy at each $l$ as seen in Fig. \ref{Fi:PI-l_1d}, where
these solutions are represented by green squares; discontinuities present in the sequence of green squares are good indicators of the existence of multiple brunches.  
The blue solid lines in Fig. \ref{Fi:PI-l_1d}, on the other hand, are obtained by a simple \emph{branch-tracking} technique.  
We first chose a starting value of $l$ and an initial guess for the solution $u_{\text{init}}$, and solved the problem for $\bar{u}(l)$.  
We then incremented/decremented $l$ by a small amount $\Delta l$ and solved this updated problem for $\bar{u}(l+\Delta l)$ using $\bar{u}(l)$ as the initial guess.  
We repeated this process of using the previous solution as the initial guess for the updated problem to extend the smooth energy curves shown in Fig. \ref{Fi:PI-l_1d}.  
This method helped us to \emph{stay} on the branch that the very first solution happened to fall onto.   
In our numerical experiments the first solutions were obtained in two different ways: using the homogeneous initial guess and using random initial guesses.  
Specifically, branches B, D, and F were first solved using the homogeneous initial guess at $l=0.20$, $l=0.10$, and $l=0.08$, followed by incrementation/decrementation of $l$.  
Branches B, D, and F obtained in this way respect \emph{geometric symmetry} (cf. \cite{Truskinovsky1996}) of the boundary value problem at least for large enough $l$; contrary to the case considered in Sec. 3.4 of \cite{Truskinovsky1996}, elastic supports are absent in our problem and one can see that the translated solution fields $\tilde{u}(X):=u(X)-d/2$ of these branches satisfy $\tilde{u}(X)=-\tilde{u}(1-X)$, and thus $\tilde{u}_{,X}(X)=\tilde{u}_{,X}(1-X)$.  
Branches A, C, and E, on the other hand, are asymmetric. Those branches were obtained using random initial guesses at $l=0.10$, $l=0.10$, and $l=0.05$, and then branch-tracking.  

\subsection{Stability}
\begin{figure}
\begin{center}
        \begin{subfigure}[b]{5.5cm}
            \centering
            \includegraphics[scale=0.25]{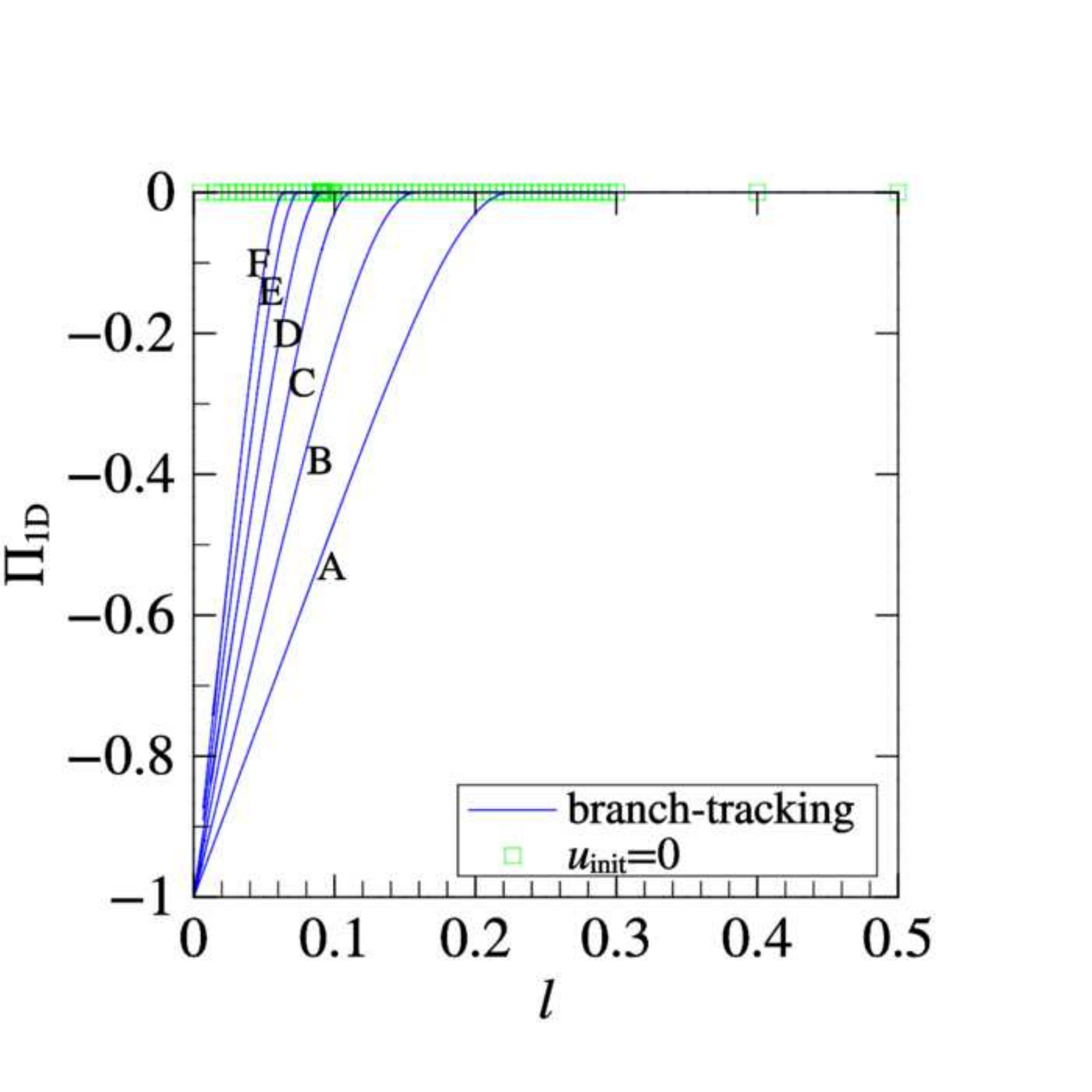}
            \caption{}
            \label{Fi:PI-l_1d_orig}
        \end{subfigure}
        ~
        \begin{subfigure}[b]{5.5cm}
            \centering
            \includegraphics[scale=0.25]{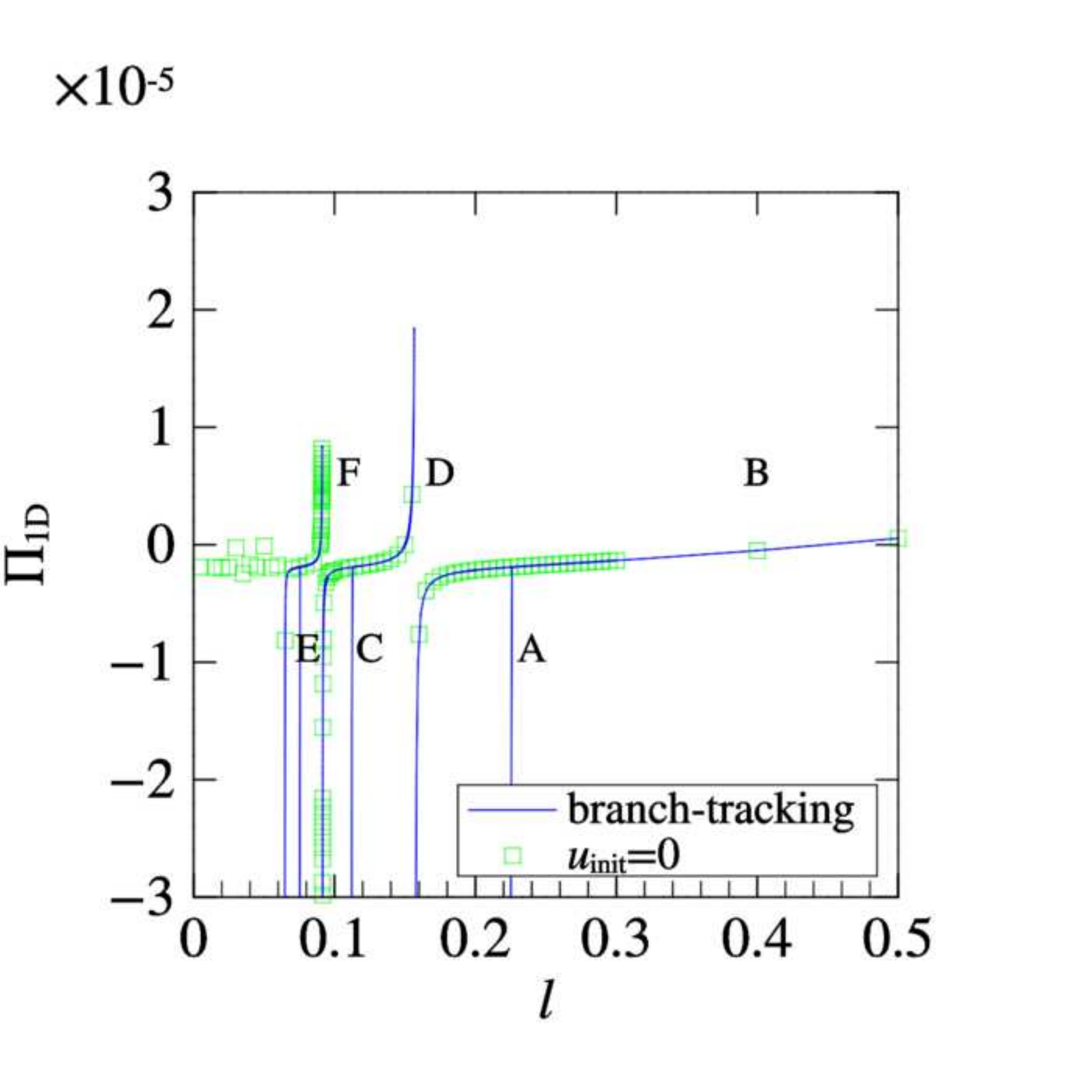}
            \caption{}
            \label{Fi:PI-l_1d_zoom1}
        \end{subfigure}
\caption{ Plots of total free energy $\Pi_{\text{1D}}$ v.s. length scale parameter $l$ on (\subref{Fi:PI-l_1d_orig}) $(l,\Pi_{\text{1D}}) \in [0,0.5]\times[-1,3\times 10^{-5}]$ and (\subref{Fi:PI-l_1d_zoom1}) $(l,\Pi_{\text{1D}}) \in[0,0.5]\times[-3\times 10^{-5},3\times10^{-5}]$ for selected branches.  
Green squares were computed using the homogeneous initial guess, $u_{\text{init}}=0$, while blue solid curves were obtained by branch-tracking.  
Branches A - F are labeled.  }
\label{Fi:PI-l_1d}
\end{center}
\end{figure}

We performed a numerical stability analysis by evaluating the positive definiteness of the Hessians corresponding to the second variation \eqref{E:D2_Pi_1D} of the continuous problem for branches A - F, of increasing total free energy.  We recall that if we find multiple solution branches corresponding to local minima of the total free energy, these are classified as \emph{metastable}.
Our analysis showed that the lowest branch, A, is stable up to $l=0.225$, at which value it meets branch B in the $\Pi_{\text{1D}}-l$ space, branch B is stable for $l>0.225$, and there exists no metastable solution branch at any $l$; see Figs. \ref{Fi:PI-l_1d_orig} and \ref{Fi:PI-l_1d_zoom1}. 
In \cite{Carr1984} it was shown that, if higher-order Neumann boundary conditions $u_{,XX}=0$ are applied at both ends instead of the higher-order Dirichlet conditions $u_{,X}=0$ as in \eqref{E:bc_u_1D}, the only stable branch is the one of lowest energy and no metastable branch exists for the type of energy density defined in \eqref{E:Psi_1d}.  
Although the boundary conditions employed here are different, leading to the particular $\Pi_{\text{1D}}-l$ free energy landscapes in Fig. \ref{Fi:PI-l_1d_orig}, our observation is essentially consistent with the analysis given in \cite{Carr1984}.

\section{Branching and stability of solutions in three-dimensional, non-convex elasticity}\label{S:Three-dimensional_example}
We now turn our attention to the main focus of this communication: branching in three-dimensional problems.  
In related two dimensional work \cite{Healey2007} the branching of solutions was studied for two-phase elastic solids with pure Dirichlet boundaries using a non-convex free energy density function also regularized by strain-gradient terms.  
There, the associated Euler-Lagrange equation admitted the trivial solution, allowing for a local bifurcation analysis. The equations were first linearized to find \emph{bifurcation points}, i.e. values of the length-scale parameter at which solution branches bifurcate from the trivial solution.  
The local bifurcation analysis was then followed by a global bifurcation analysis, where solution branches were continued along the length-scale parameter from those bifurcation points, solving the original nonlinear equation using a branch-tracking technique.  
The stability of those branches was then assessed by numerically checking the positive definiteness of the second variation of the total free energy.
The same technique was also used by Vainchtein and co-workers \cite{Vainchtein1999} for one-dimensional problems.  

Here, we carry out a numerical, three-dimensional study of an elastic solid that undergoes phase transformations between three tetragonal variants under traction loads to form branches.  
We chose to work on a body subject to traction because one of our ultimate goals is the simulation of shape-memory alloys, where such traction boundary conditions naturally arise.  
For the boundary value problem governed by non-convex elasticity, regularized by Toupin's theory of gradient elasticity at finite strain, only pure numerical approaches are feasible to compute branches, using the methods described in Sec. \ref{S:One-dimensional_example}, instead of the ones based on local bifurcation analysis.  The stability of each solution is assessed as described in Sec. \ref{S:One-dimensional_example} using the second variation of the total free energy.

\begin{figure}
\begin{center}
        \begin{subfigure}[b]{6.5cm}
            \centering
            \includegraphics[scale=0.15]{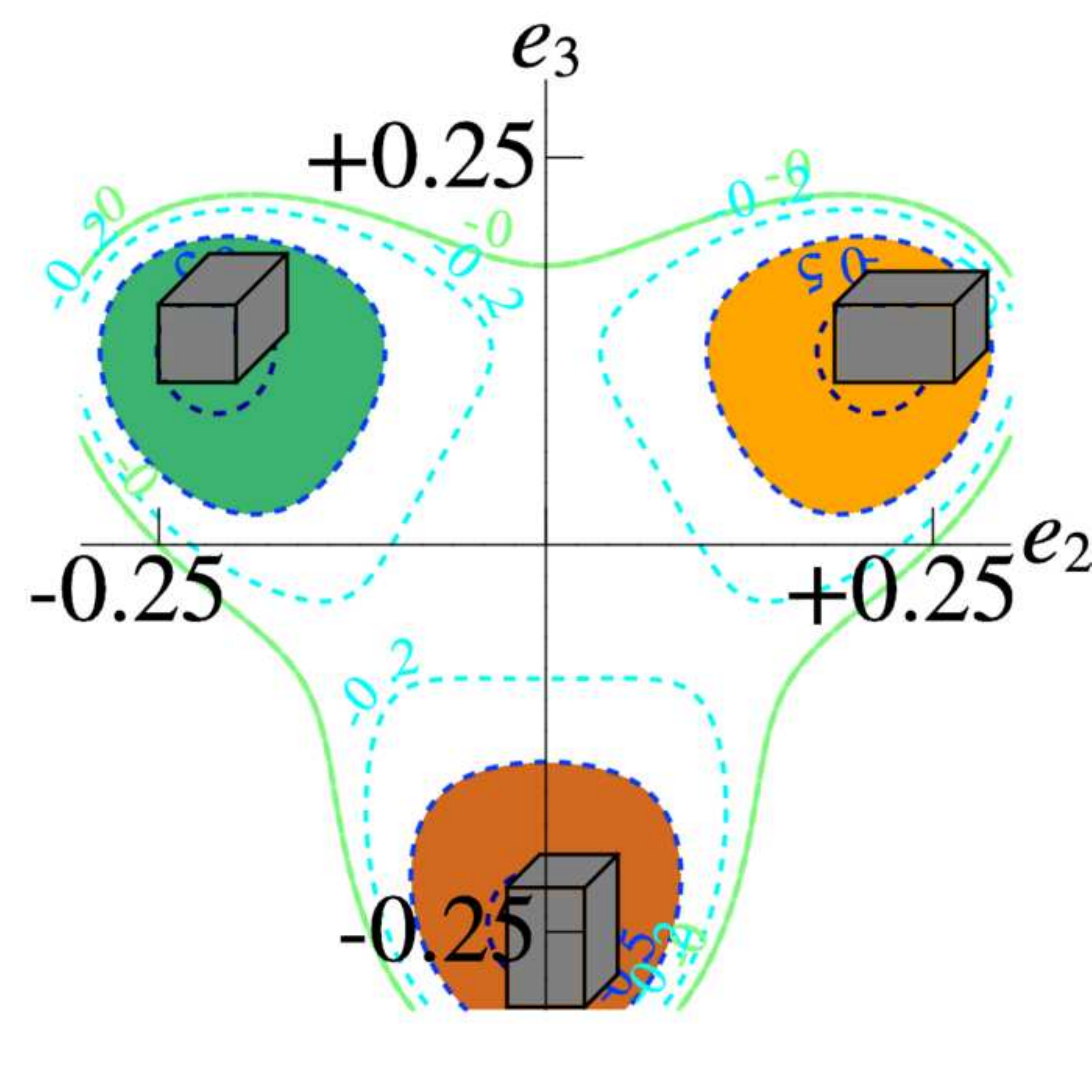}
            \includegraphics[scale=0.15]{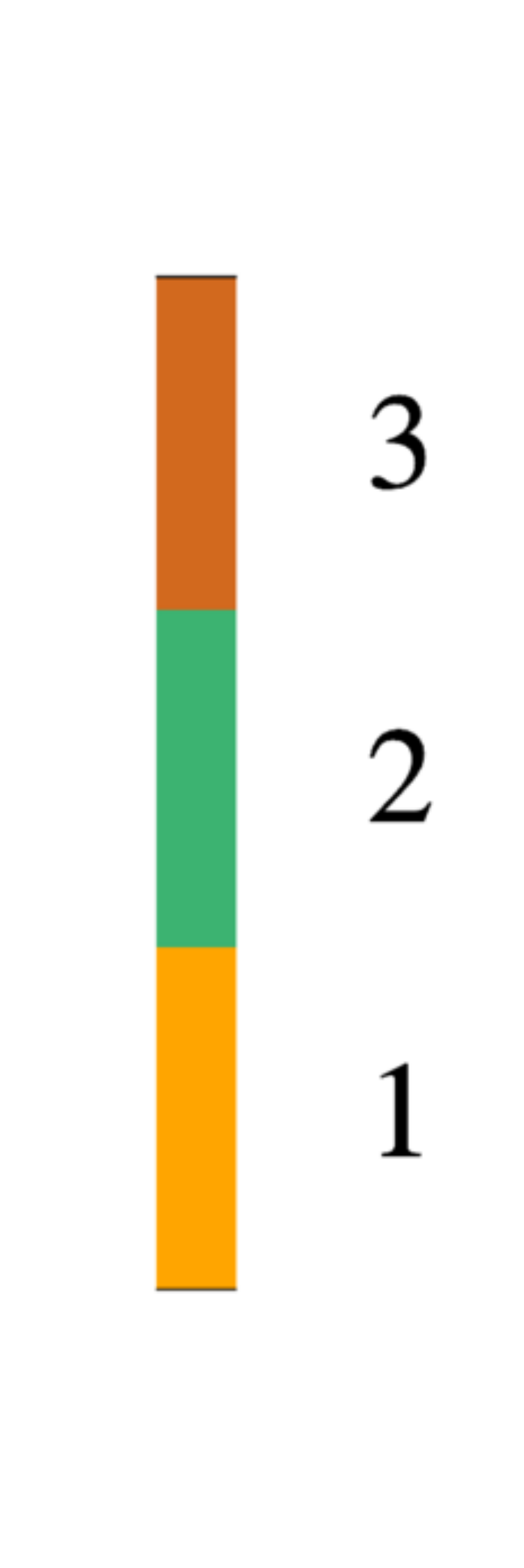}
            \caption{}
            \label{Fi:3well_diagram}
        \end{subfigure}
        ~
        \begin{subfigure}[b]{6.5cm}
            \centering
            \includegraphics[scale=0.15]{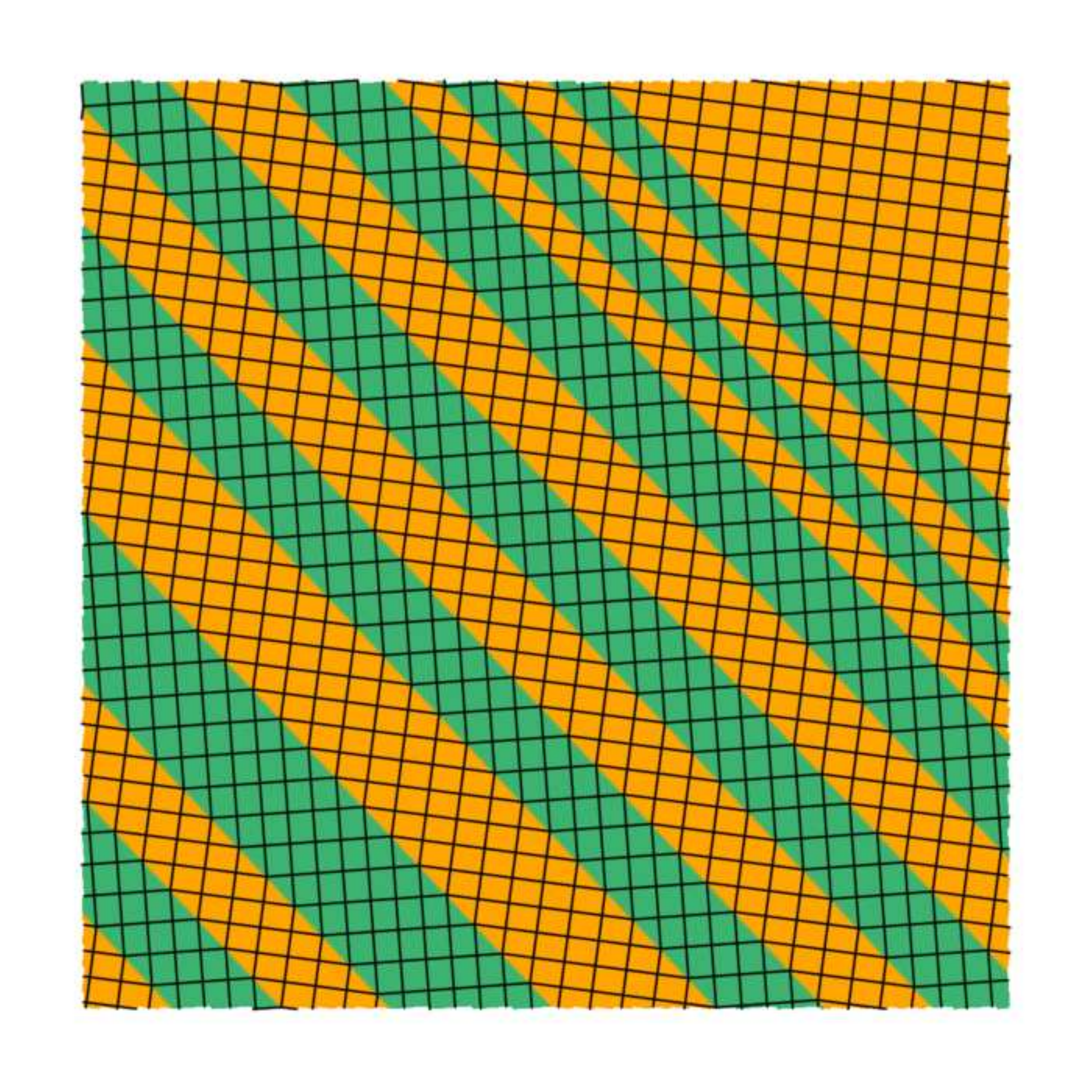}
            \includegraphics[scale=0.15]{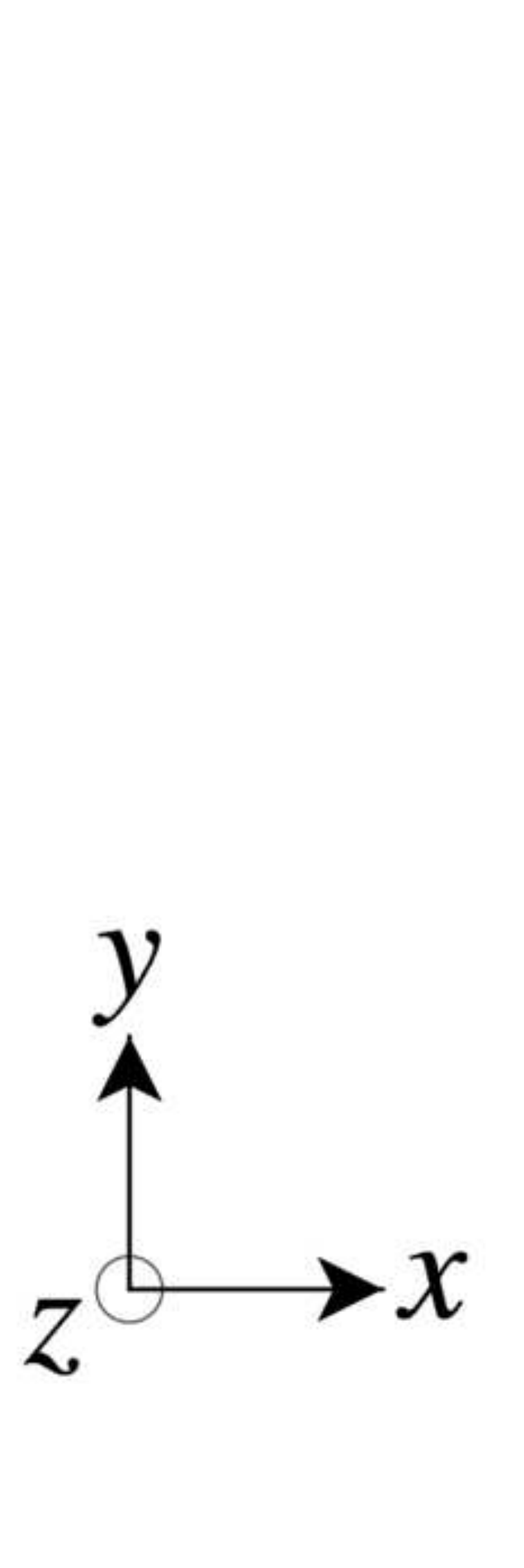}
            \caption{}
            \label{Fi:twin}
        \end{subfigure}
\caption{ (\subref{Fi:3well_diagram}) Contour plots of $\Psi$ projected onto the $e_2-e_3$ plane; contour curves of $\Psi=0,-0.2,-0.5,-0.9$ are shown.  
Three tetragonal variants are colored/numbered.  
A tetragonal variant is regarded as present at a point of the body if the energy density on this $e_2-e_3$ plane is less than $-0.5$ at that point.  
(\subref{Fi:twin}) A two-dimensional slice of a typical microstructure in three dimensions that minimizes the total free energy when $l=0$ in \eqref{E:Psi}.  Twinning is formed between Variants 1 and 2.  }
\label{Fi:plot_contour_3variants}
\end{center}
\end{figure}

The non-dimensionalized free-energy density function is defined in terms of gradients and strain gradients of the displacement field $\bs{u}(\bs{X})$ as:
\begin{align}
\Psi&:
=B_1e_1^2
+B_2(e_2^2+e_3^2)
+B_3e_3(e_3^2-3e_2^2)
+B_4(e_2^2+e_3^2)^2
+B_5(e_4^2+e_5^2+e_6^2)\notag\\
&\phantom{:}+l^2(e_{2,1}^2+e_{2,2}^2+e_{2,3}^2+e_{3,1}^2+e_{3,2}^2+e_{3,3}^2),
\label{E:Psi}
\end{align}
where $B_1,...,B_5$ are constants and the following reparameterized strain measures were used:
\begin{subequations}
\begin{alignat}{3}
&e_1=\frac{E_{11}+E_{22}+E_{33}}{\sqrt{3}},\qquad &&e_2 =\frac{E_{11}-E_{22}}{\sqrt{2}},\qquad &&e_3 =\frac{E_{11}+E_{22}-2E_{33}}{\sqrt{6}} \label{E:strains1-3}\\ 
&e_4=E_{23}=E_{32},\qquad &&e_5 =E_{13}=E_{31},\qquad &&e_6 =E_{12}=E_{21},\label{E:strains4-6}
\end{alignat}
\end{subequations}
where $E_{IJ}=1/2(F_{kI}F_{kJ}-\delta_{IJ})$ are components of the Green-Lagrange strain tensor, $F_{iJ}=\delta_{iJ}+u_{i,J}$ being components of the deformation gradient tensor.  
Here as elsewhere ${(\hspace{1pt}\cdot\hspace{1pt})_{,J}}$ denotes spatial derivatives with respect to the reference rectangular Cartesian coordinate variable $X_J$ ($J=1,2,3$).  
Throughout this work we set $B_5=180$, $B_1=3.25B_5$, $B_2=-1.5/r^2$, $B_3=1.0/r^3$, and $B_4=1.5/r^4$, where $r=0.25$, unless otherwise noted.  
At $e_2,e_3 = 0$, corresponding to deformations that reduce to volumetric dilatations in the infinitesimal strain limit, the free energy density function \eqref{E:Psi}  possesses a local maximum in $e_2-e_3$ space and represents the cubic austenite crystal structure.  The reference, unstrained state is also in the cubic austenite structure. Thus defined,  $\Psi$ is non-convex with respect to the strain variables $e_2$ and $e_3$ and possesses three minima, or \emph{energy wells}, of unit depth located at a distance of $0.25$ from the origin on the $e_2-e_3$ plane; see Fig. \ref{Fi:3well_diagram}.  
These three energy wells represent three martensitic variants of symmetrically equivalent tetragonal crystal structures, elongated in the $X_1$-, $X_2$-, and $X_3$-directions, respectively, that are colored/numbered in Fig. \ref{Fi:3well_diagram}.
As an example, a characteristic configuration that achieves minimum energy density of $-1$ almost everywhere for $l=0$ in \eqref{E:Psi} appears in Fig. \ref{Fi:twin}, showing laminae of variant 1 and variant 2.  
In our numerical example a tetragonal variant is regarded as present at a point of the body if the energy density on this $e_2-e_3$ plane is less than $-0.5$ at that point.

We are interested in solution fields $\bs{u}(\bs{X})$ on a unit cube $\overline{\Omega}$, where $\Omega=(0,1)^3$, that satisfy the following Dirichlet boundary conditions:
\begin{subequations}
\begin{align}
u_i=0,\; u_{i,1}=0 \quad &\text{ on } X_1=0,\quad i=1,2,3,\\
u_1=0,\; u_{1,1}=0 \quad &\text{ on } X_1=1,
\end{align}
\label{E:bc_u}
\end{subequations}
and globally/locally minimize the total free energy corresponding to the density function \eqref{E:Psi}:
\begin{align}
\Pi:=\int_{\Omega}\Psi\dif\Omega-\int_{\Gamma_{X_1=1}}(u_2T_2+u_3T_3)\dif\Gamma_{X_1=1},\label{E:Pi}
\end{align}
where $T_2=T_3=0.01$ are the standard tractions on the reference boundary $X_1=1$ denoted by $\Gamma_{X_1=1}$. The standard traction vanishes on the boundaries $X_2 = \{0,1\}$, $X_3 = \{0,1\}$, and higher-order tractions vanish wherever higher-order Dirichlet conditions are not prescribed \cite{Toupin1962,Toupin1964}. 
We define admissible test functions $\bs{w} \in \mathcal{V}$ that satisfy:
\begin{subequations}
\begin{align}
w_i=0,\; w_{i,1}=0 \quad &\text{ on } X_1=0,\quad i=1,2,3,\\
w_1=0,\; w_{1,1}=0 \quad &\text{ on } X_1=1,
\end{align}
\label{E:bc_w}
\end{subequations}
and solve the following weak form of the boundary value problem derived from variational arguments \cite{Toupin1962,Toupin1964,Rudraraju2014}:
\begin{align}
D\Pi[\bs{w}]=\int_{\Omega}(w_{i,J}P_{iJ}+w_{i,JK}B_{iJK})\dif\Omega-\int_{\Gamma_{X_1=1}}(w_2T_2+w_3T_3)\dif\Gamma_{X_1=1}=0,\quad\forall\bs{w}\in \mathcal{V}\label{E:D_Pi}
\end{align}
where the first Piola-Kirchhoff stress tensor and the higher-order stress tensor in component form are:
\begin{subequations}
\begin{align}
P_{iJ}&:=\partial\Psi/\partial F_{iJ},\label{P3D}\\
B_{iJK}&:=\partial\Psi/\partial F_{iJ,K}.\label{B3D}
\end{align}
\end{subequations}
We then assess stability of each solution by numerically checking the positive definiteness of the second variation as:
\begin{align}
D^2\Pi[\bs{w},\bs{w}]&=\int\limits_{\Omega}\left(
w_{i,I}\frac{\partial^2\Psi}{\partial F_{iI}\partial F_{jJ}}w_{j,J}
+w_{i,I}\frac{\partial^2\Psi}{\partial F_{iI}\partial F_{jJ,K}}w_{j,JK}
+w_{i,IL}\frac{\partial^2\Psi}{\partial F_{iI,L}\partial F_{jJ}}w_{j,J}
+w_{i,IL}\frac{\partial^2\Psi}{\partial F_{iI,L}\partial F_{jJ,K}}w_{j,JK}
\right)\dif\Omega>0,\label{E:D2_PI}
\end{align}
where, as in the one-dimensional case, we have explicitly retained the symmetric second and third terms of the integrand for clarity of the development.  Note that, following standard variational arguments \cite{Toupin1962,Toupin1964}, one can derive the strong form of the boundary value problem corresponding to the weak form \eqref{E:D_Pi} and \eqref{E:bc_u} with \eqref{E:bc_w} as:
\begin{align}
-P_{iJ,J}+B_{iJK,JK}=0\quad\text{ on }\Omega,\notag
\end{align}
along with Neumann/higher-order Neumann conditions:
\begin{subequations}
\begin{alignat}{2}
P_{i1}-B_{i11,1}-2B_{i12,2}-2B_{i13,3}&=T_i,\; B_{i11}&&=0,\; (i=2,3) \;\text{on}\; X_1=1,\notag\\
P_{iL}-2(B_{iL1,1}+B_{iL2,2}+B_{iL3,3})+B_{iLL,L}&=0,\;B_{iLL}&&=0,\;(i=1,2,3,\; \text{no sum on } L)\; \text{on}\; X_L=\{0,1\}\;(L=2,3).\notag  
\end{alignat}
\end{subequations}
A more detailed treatment of these boundary conditions is found in the works of Toupin \cite{Toupin1962,Toupin1964}.  

\subsection{Numerics}
We seek numerical solutions $\bs{u}^h\in \mathcal{S}^h\subset\mathcal{S}$ to a finite-dimensional counterpart of the weak form \eqref{E:D_Pi} defined for $\bs{w}^h\in\mathcal{V}^h\subset\mathcal{V}$, where:
\begin{subequations}
\begin{align}
\mathcal{S}^h &= \{\bs{v}^h \in H^2(\Omega)\vert v_i^h=0,\; v^h_{i,1}=0\;\text{on} \; X_1=0,\;\text{for}\; i=1,2,3,\;v^h_1=0,\;v^h_{1,1}=0\;\text{on}\; X_1 = 1\}, \\
\mathcal{V}^h &= \{\bs{v}^h \in H^2(\Omega)\vert v_i^h=0,\; v^h_{i,1}=0\;\text{on} \; X_1=0,\;\text{for}\; i=1,2,3,\;v^h_1=0,\;v^h_{1,1}=0\;\text{on}\; X_1 = 1\}.
\end{align}
\end{subequations}
The problem was solved using IGA.  
The finite-dimensional subspaces $\mathcal{S}^h$ and $\mathcal{V}^h$ were constructed using a second-order, $C^1$-continuous, B-spline basis defined in three dimensions on $64^3$, $128^3$, and $256^3$ elements of uniform size, which enforce the Dirichlet/higher-order Dirichlet conditions \emph{strongly}.  
IGA was previously employed to solve a range of boundary value problems with Toupin's theory of gradient elasticity at finite strain by Rudraraju et al. \cite{Rudraraju2014}, with higher-order Dirichlet conditions applied \emph{weakly}.  
Our code \cite{IGAP4GradElast} is written in \textsf{C}.  
We use \textsf{Mathematica 10} to symbolically produce elementwise residual/tangent evaluation routines, \textsf{PETSc 3.7.4} \cite{petsc-web-page,petsc-user-ref,petsc-efficient} for iterative linear/nonlinear solvers, 
\textsf{SLEPc 3.7.3} \cite{slepc1,slepc2,slepc3} for an eigenvalue problem solver, and \textsf{mathgl 2.3.0} for plots.  
Specifically, \textsf{MINRES} with Jacobi preconditioner and a backtracking line search method with cubic-order approximation were chosen for iterative solvers.  
We used the double-precision floating-point format with absolute tolerance of $10^{-12}$ on the residual of the discretized, matrix-vector weak form.  The floating point precision and residual tolerance were relaxed relative to the one-dimensional problem for numerical efficiency. In practice, the lower precision and less stringent tolerance were found to be adequate after using the more demanding thresholds in the one-dimensional case.  
The second variation \eqref{E:D2_PI} was discretized on the same B-spline basis, and the stability of each solution was assessed by extracting the lower end of the spectrum of eigenvalues of the corresponding symmetric Hessian matrices.  
For the eigenvalue problem, we used an absolute convergence error tolerance of $1\times 10^{-6}$.

\begin{figure}
\begin{center}
    \begin{subfigure}[b]{6.5cm}
        \centering
        \includegraphics[scale=0.25]{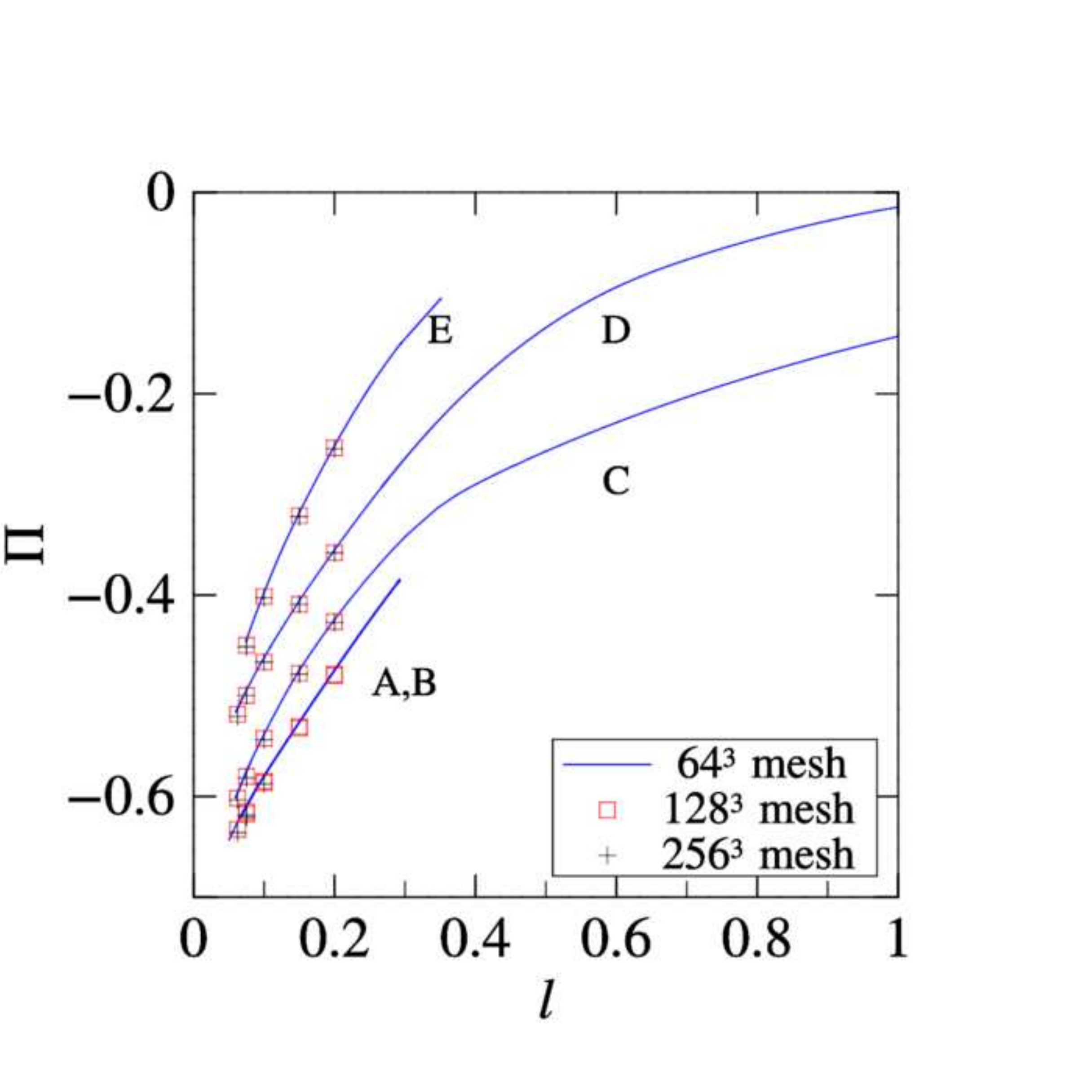}
        \caption{}
        \label{Fi:PI-l_orig}
    \end{subfigure}
    ~
    \begin{subfigure}[b]{6.5cm}
        \centering
        \includegraphics[scale=0.25]{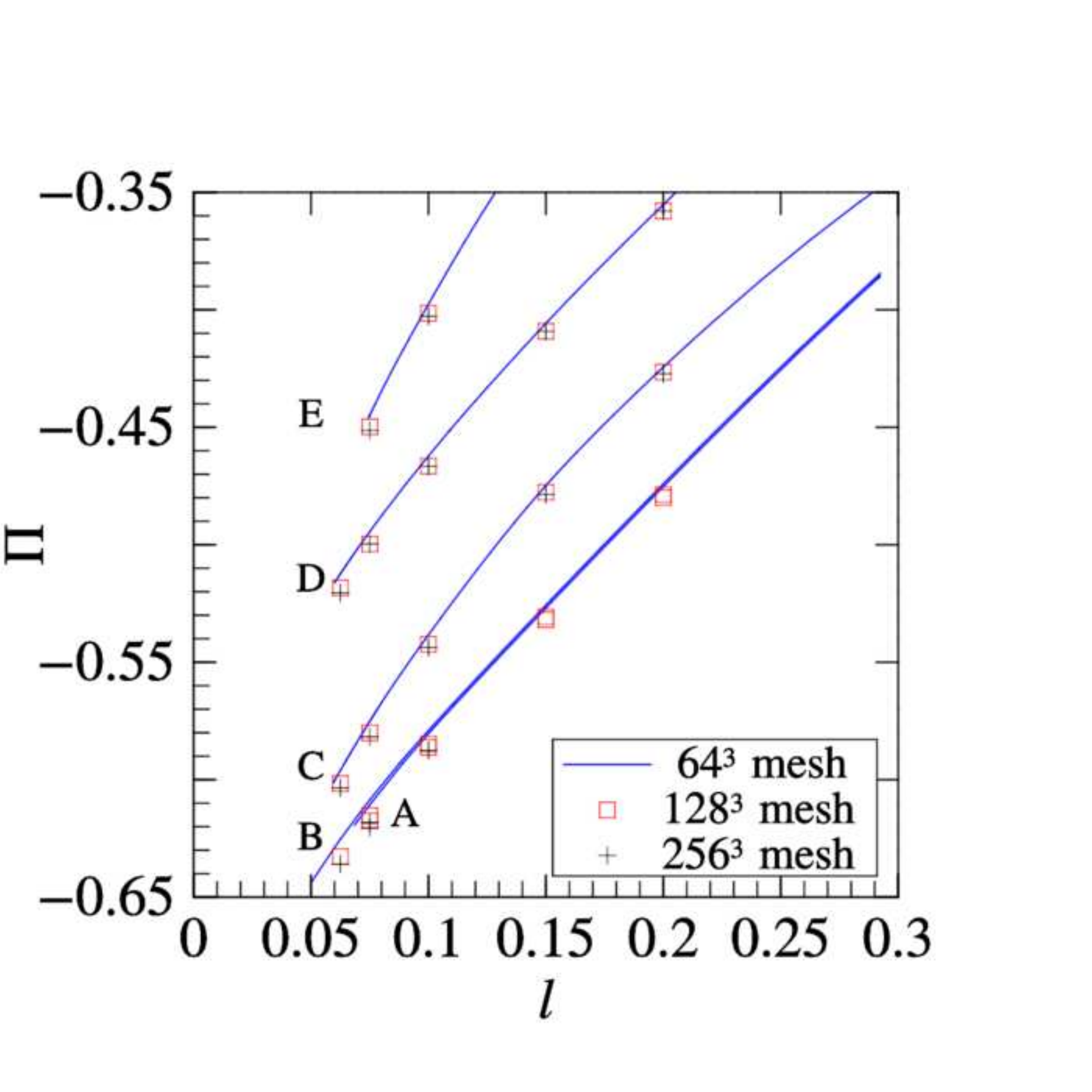}
        \caption{}
        \label{Fi:PI-l_zoom}
    \end{subfigure}
    \caption{Plots of total free energy $\Pi$ v.s. length scale parameter $l$ for branches A - E on (\subref{Fi:PI-l_orig}) $(l,\Pi) \in[0,1]\times[-0.70,0]$ and (\subref{Fi:PI-l_zoom}) $(l,\Pi) \in[0,0.3]\times[-0.65,-0.35]$. Note that branches A and B very nearly overlap in these plots.
    The blue curves were obtained by the branch-tracking technique on a $64^3$ mesh.  
    The data points represented by red squares were obtained by refining these solutions on a $128^3$ mesh at selected values of $l$ and those represented by black '+' signs were obtained by further refinement on a $256^3$ mesh.  }
\label{Fi:PI-l}
\end{center}
\end{figure}

\subsection{Solution branches and branch-tracking}
\begin{figure}
    \begin{center}
        \begin{tabular}{rp{16.5cm}}
            \parbox[t]{0.5cm}{ }&
            \begin{tabular}{p{2.8cm}p{2.8cm}p{2.8cm}p{2.8cm}p{2.8cm}p{2.5cm}}
                \hspace{0.8cm}$l\!=\!0.0625$ & 
                \hspace{0.8cm}$l\!=\!0.0750$ & 
                \hspace{0.8cm}$l\!=\!0.1000$ &
                \hspace{0.8cm}$l\!=\!0.1500$ &
                \hspace{0.8cm}$l\!=\!0.2000$ &
                \vspace{0.5\baselineskip}
            \end{tabular}  \\
            \parbox[t]{0.5cm}{ $E$ } &
            \begin{tabular}{p{2.8cm}p{2.8cm}p{2.8cm}p{2.8cm}p{2.8cm}p{2.5cm}}
                &
                \includegraphics[scale=0.12]{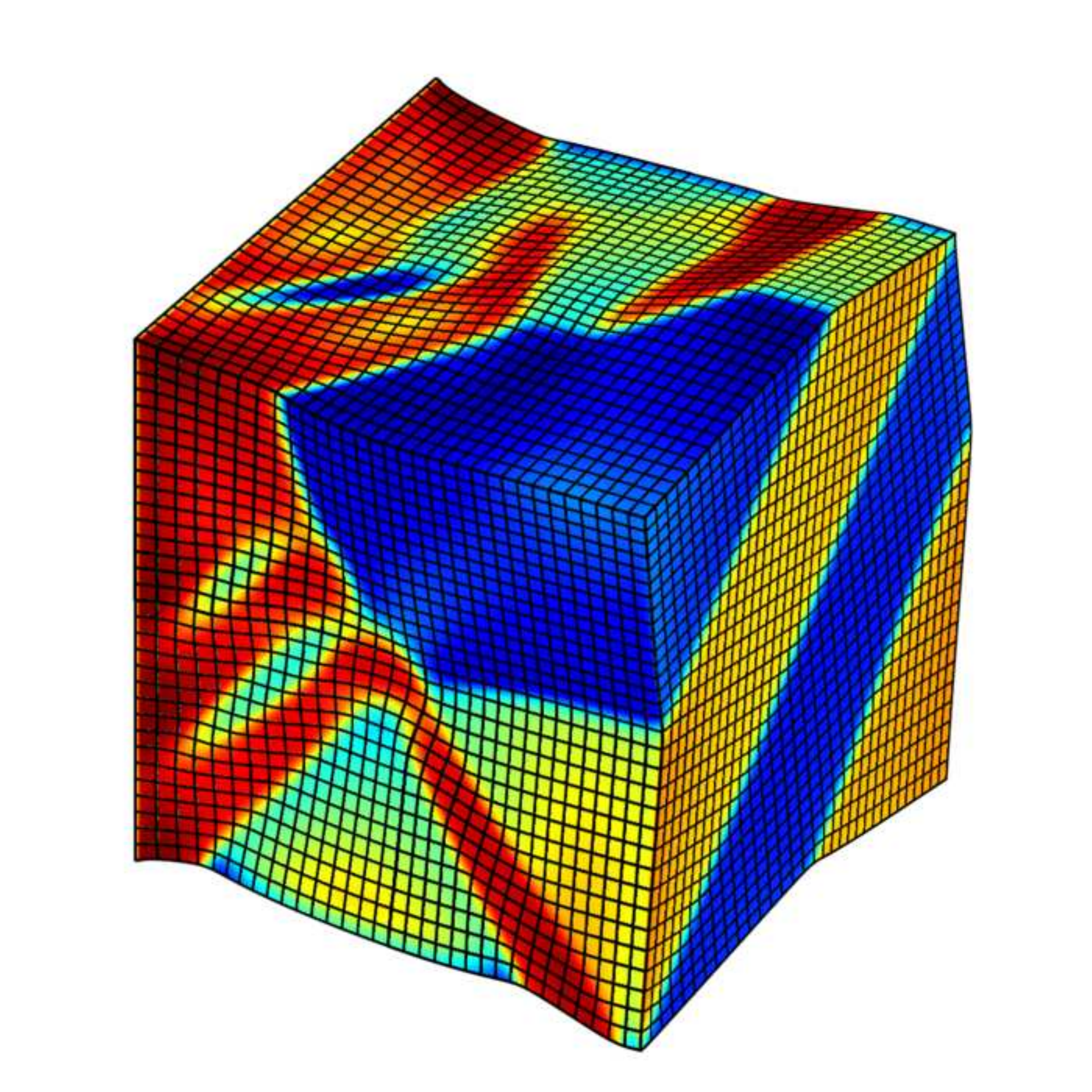} & 
                \includegraphics[scale=0.12]{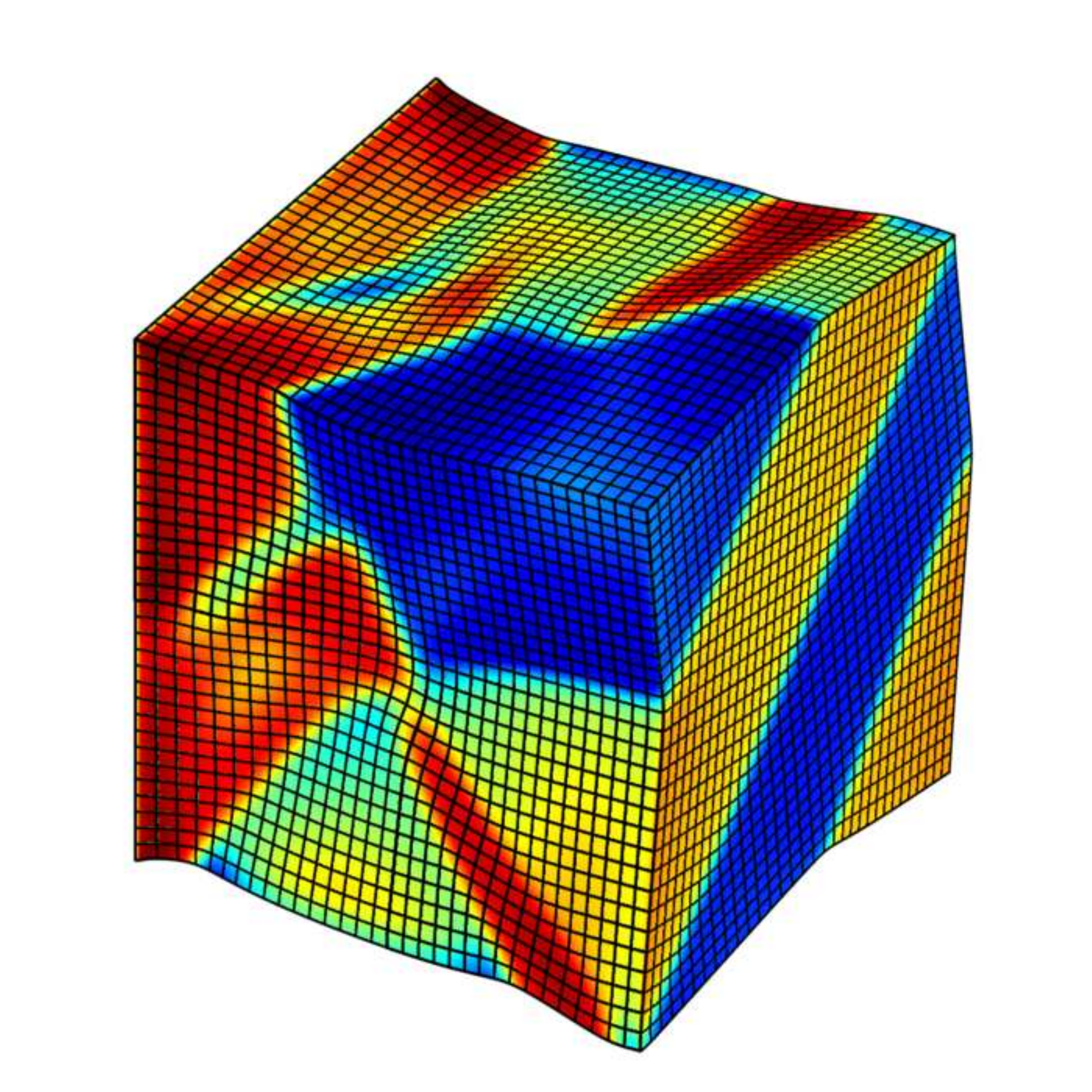} &
                \includegraphics[scale=0.12]{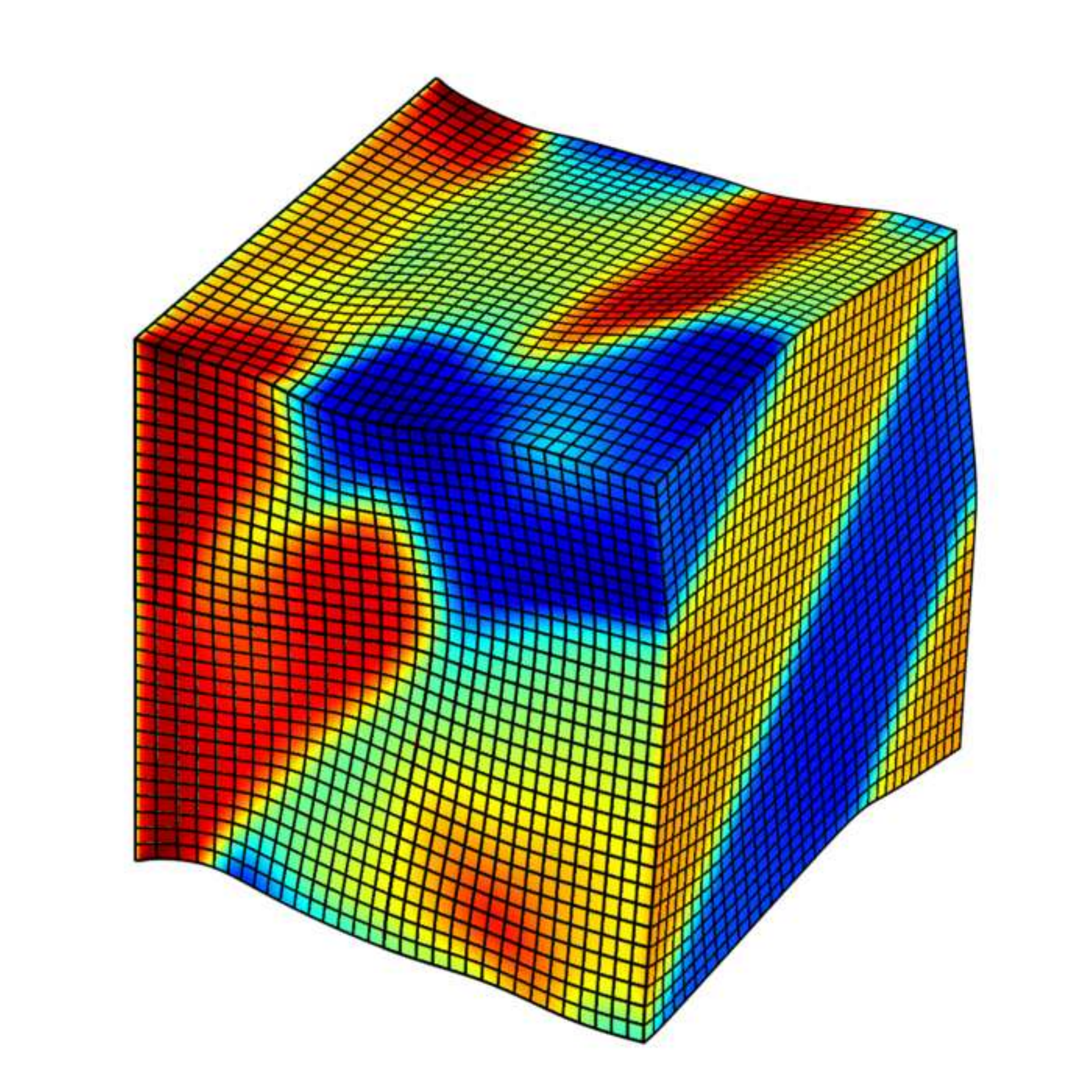} &
                \includegraphics[scale=0.12]{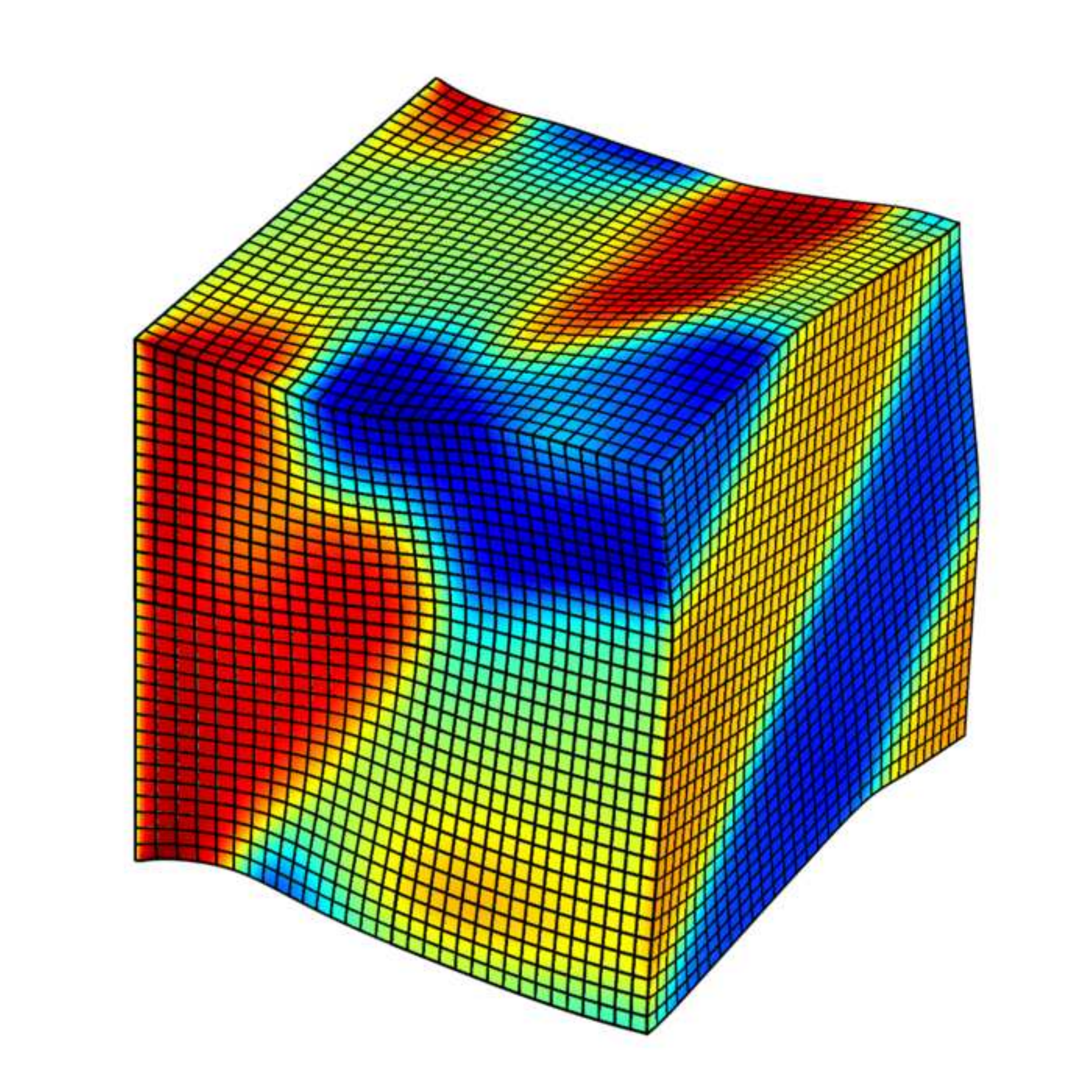} &
            \end{tabular}  \\
            \parbox[t]{0.5cm}{ $D$ } &
            \begin{tabular}{p{2.8cm}p{2.8cm}p{2.8cm}p{2.8cm}p{2.8cm}p{2.5cm}}
                \includegraphics[scale=0.12]{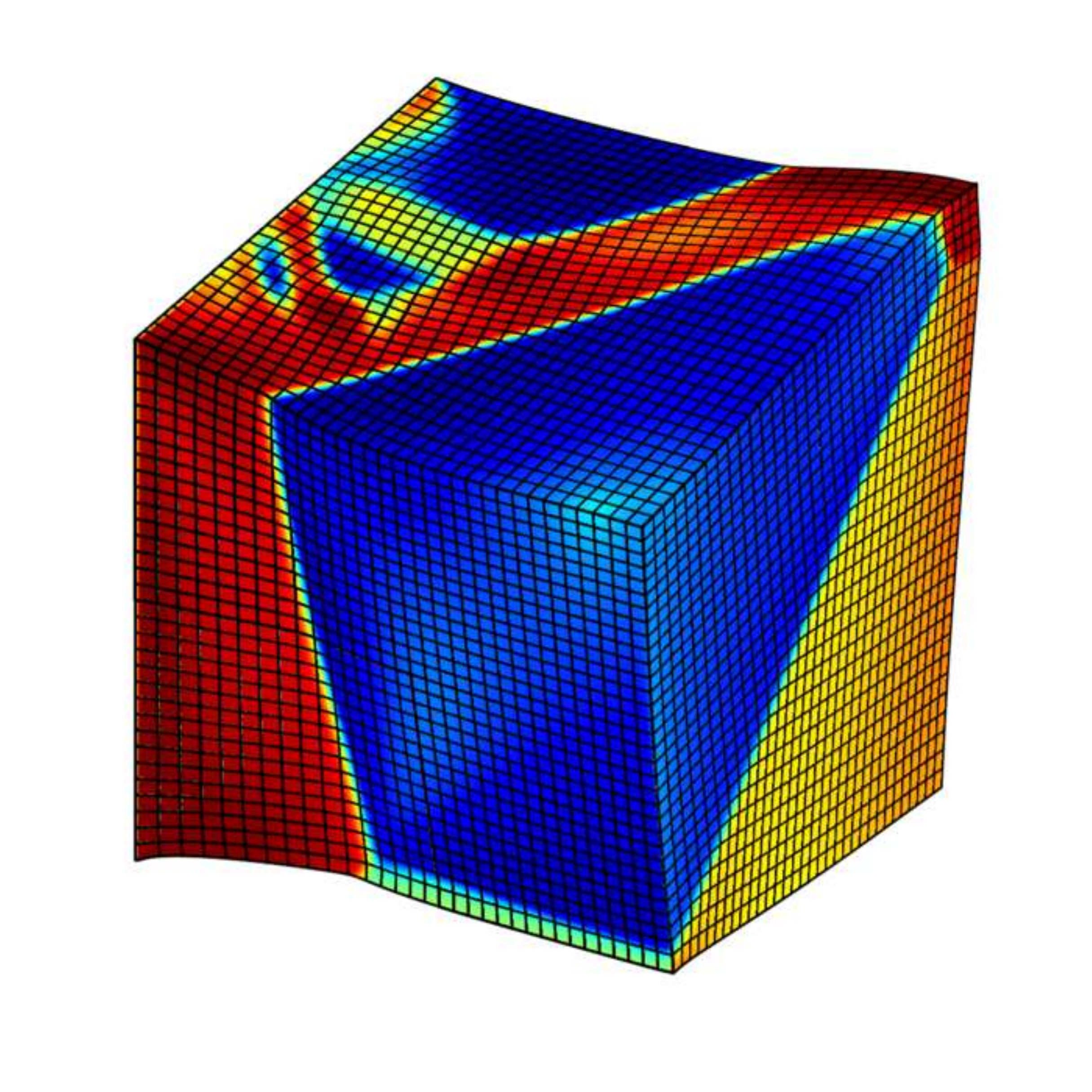} &
                \includegraphics[scale=0.12]{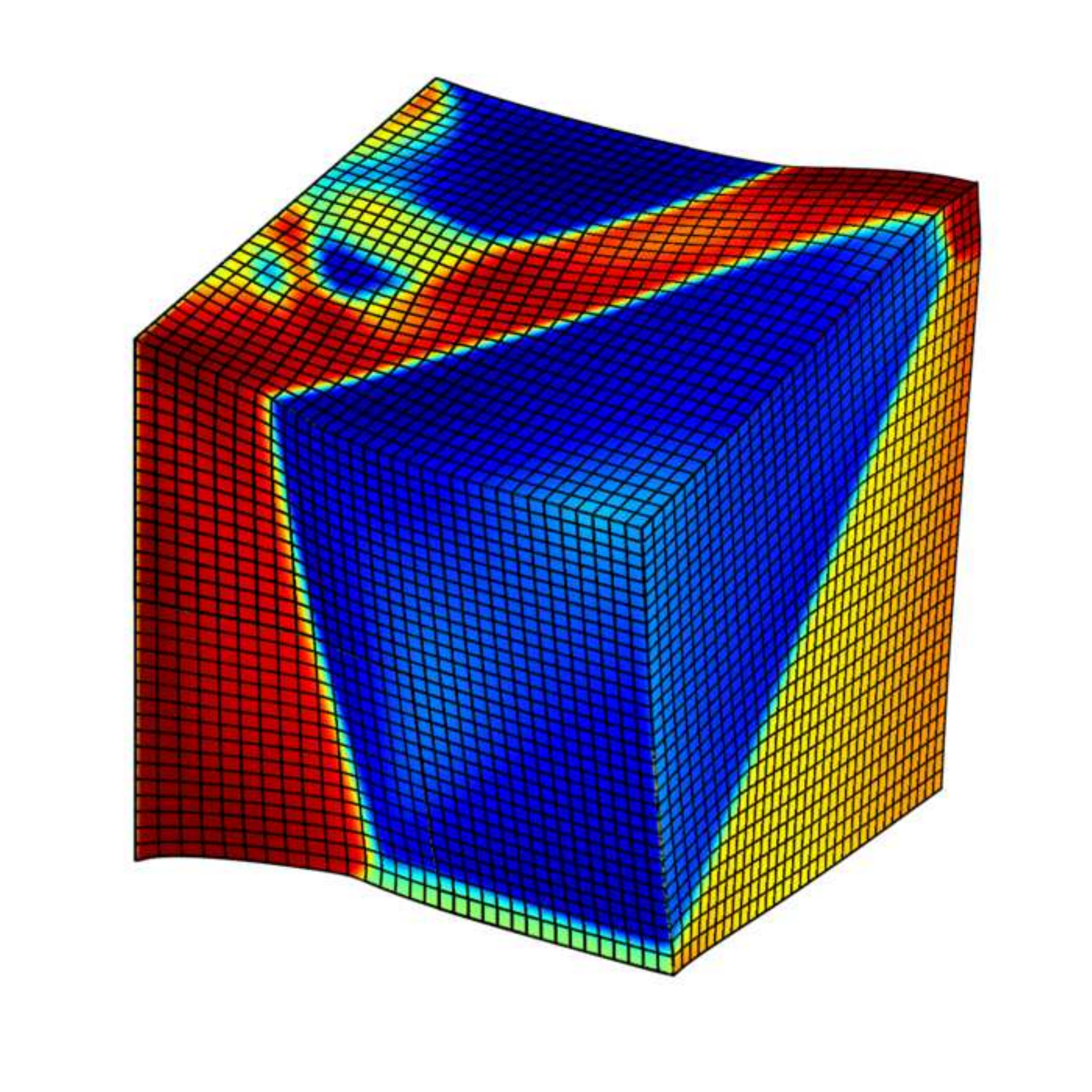} & 
                \includegraphics[scale=0.12]{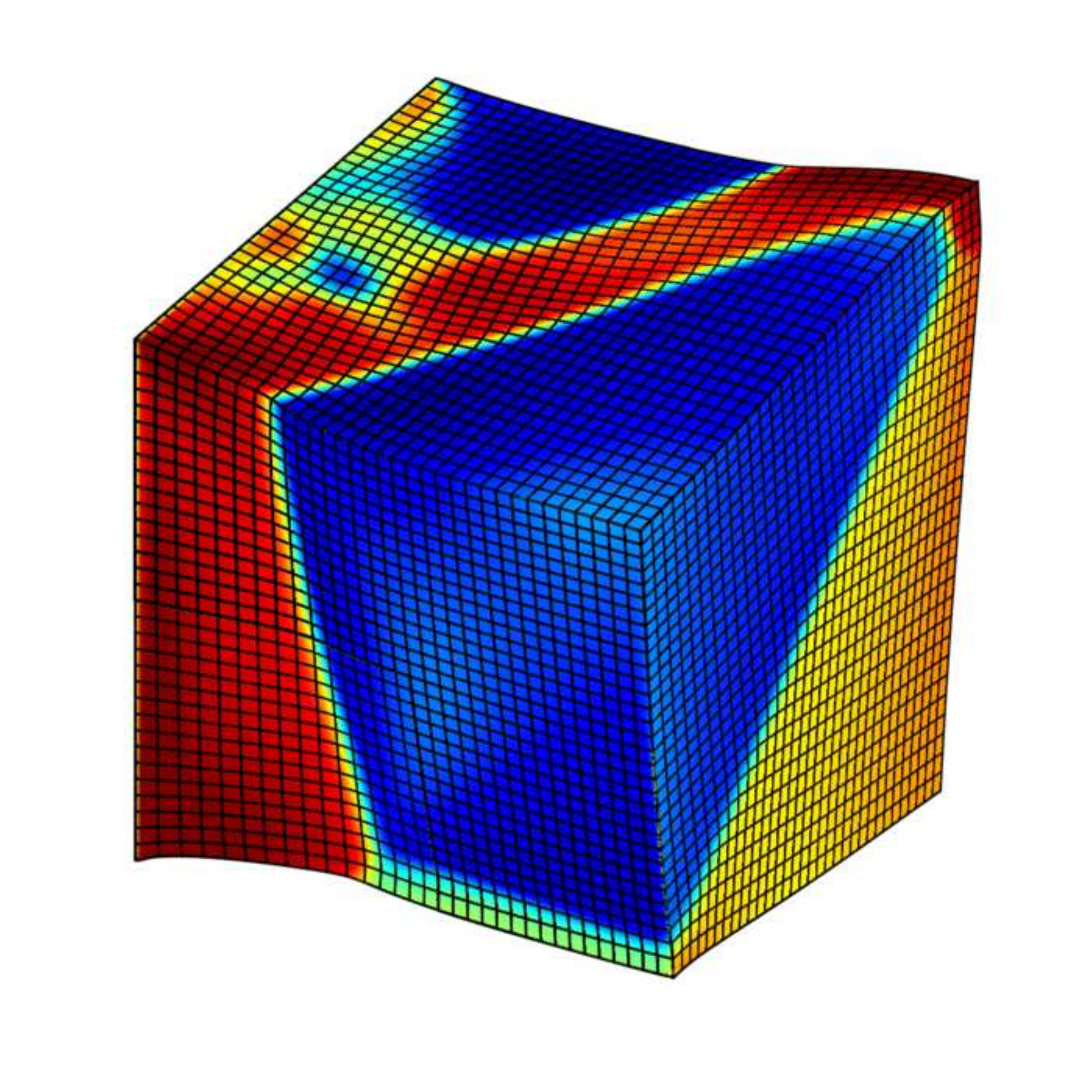} &
                \includegraphics[scale=0.12]{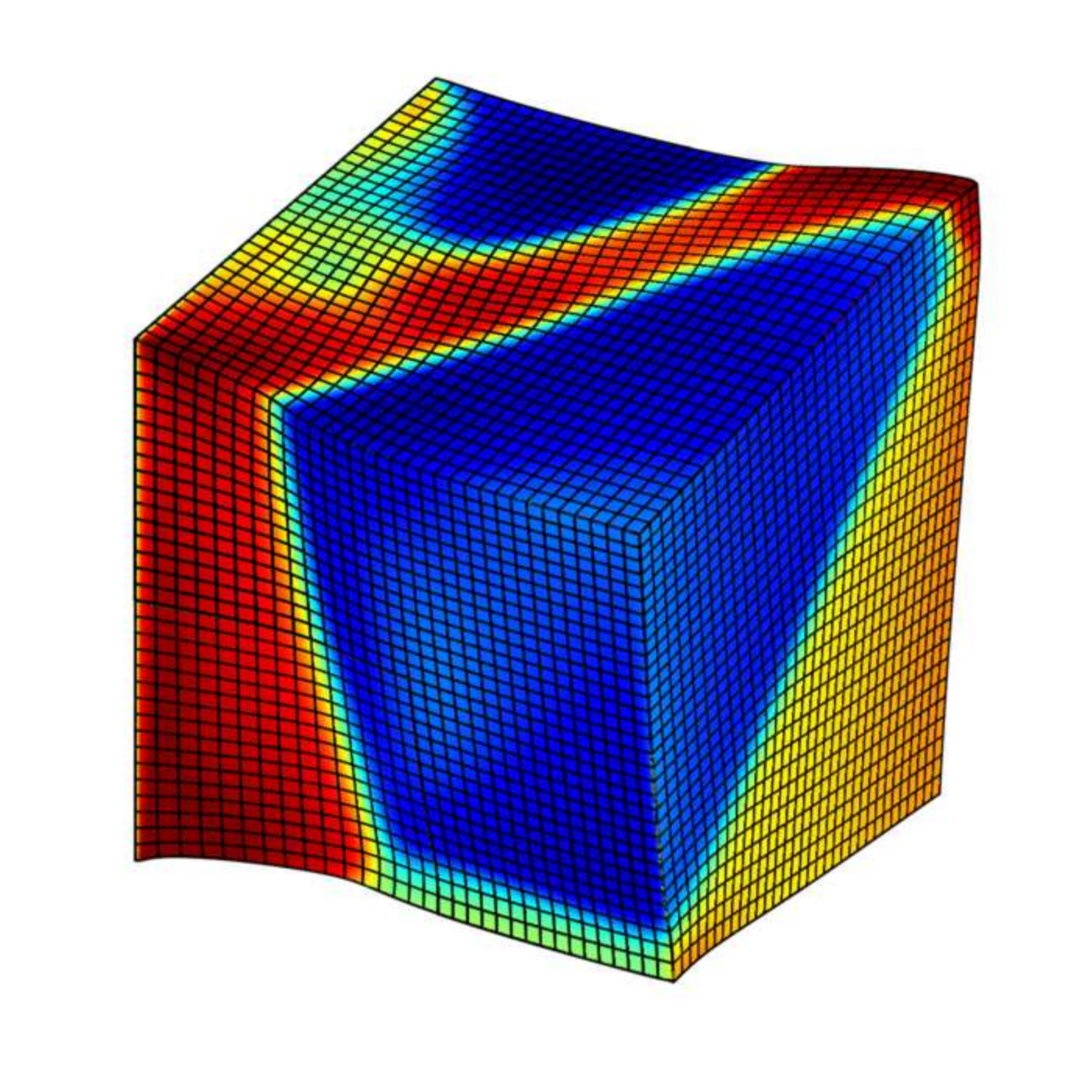} &
                \includegraphics[scale=0.12]{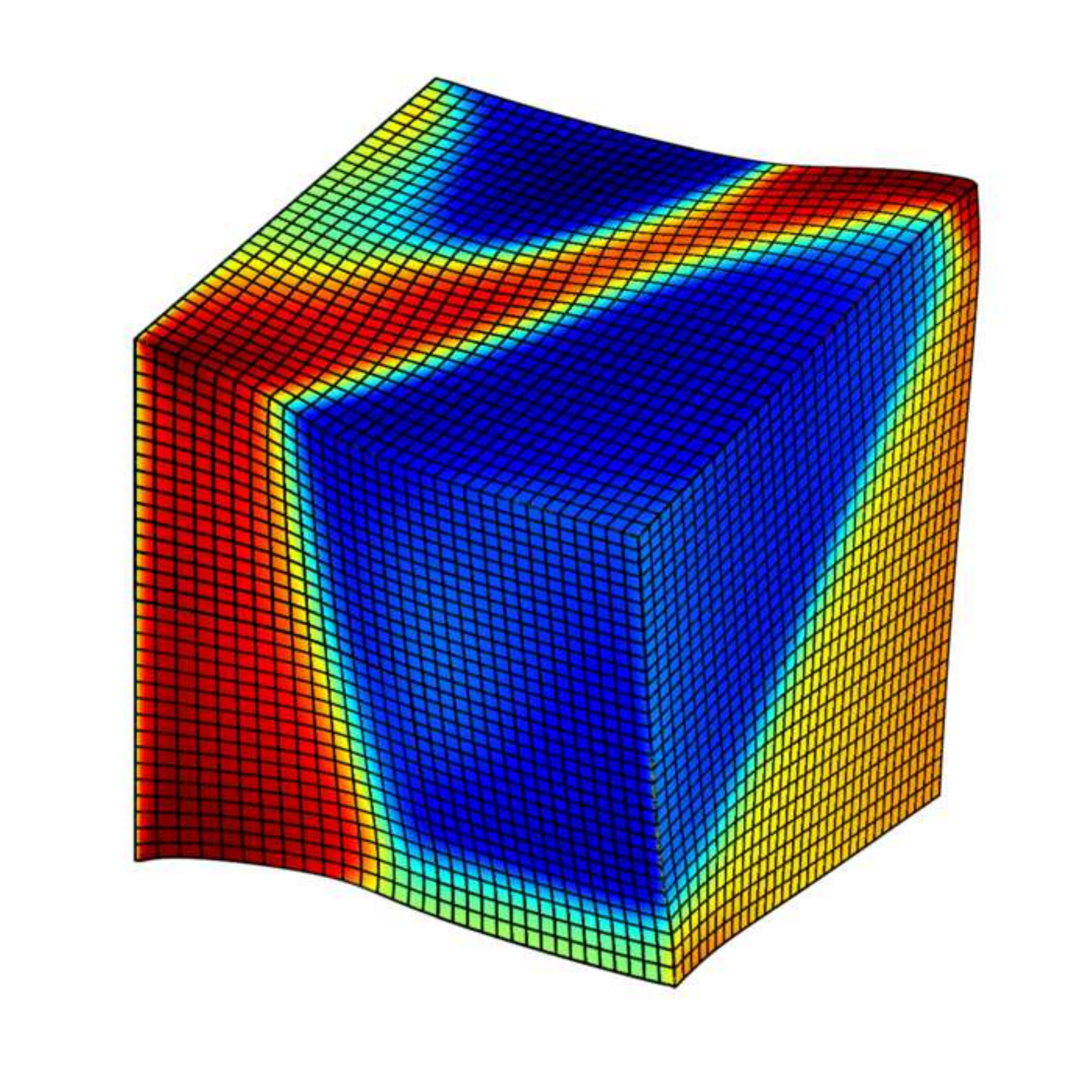} &
            \end{tabular}  \\
            \parbox[t]{0.5cm}{ $C$ } &
            \begin{tabular}{p{2.8cm}p{2.8cm}p{2.8cm}p{2.8cm}p{2.8cm}p{2.5cm}}
                \includegraphics[scale=0.12]{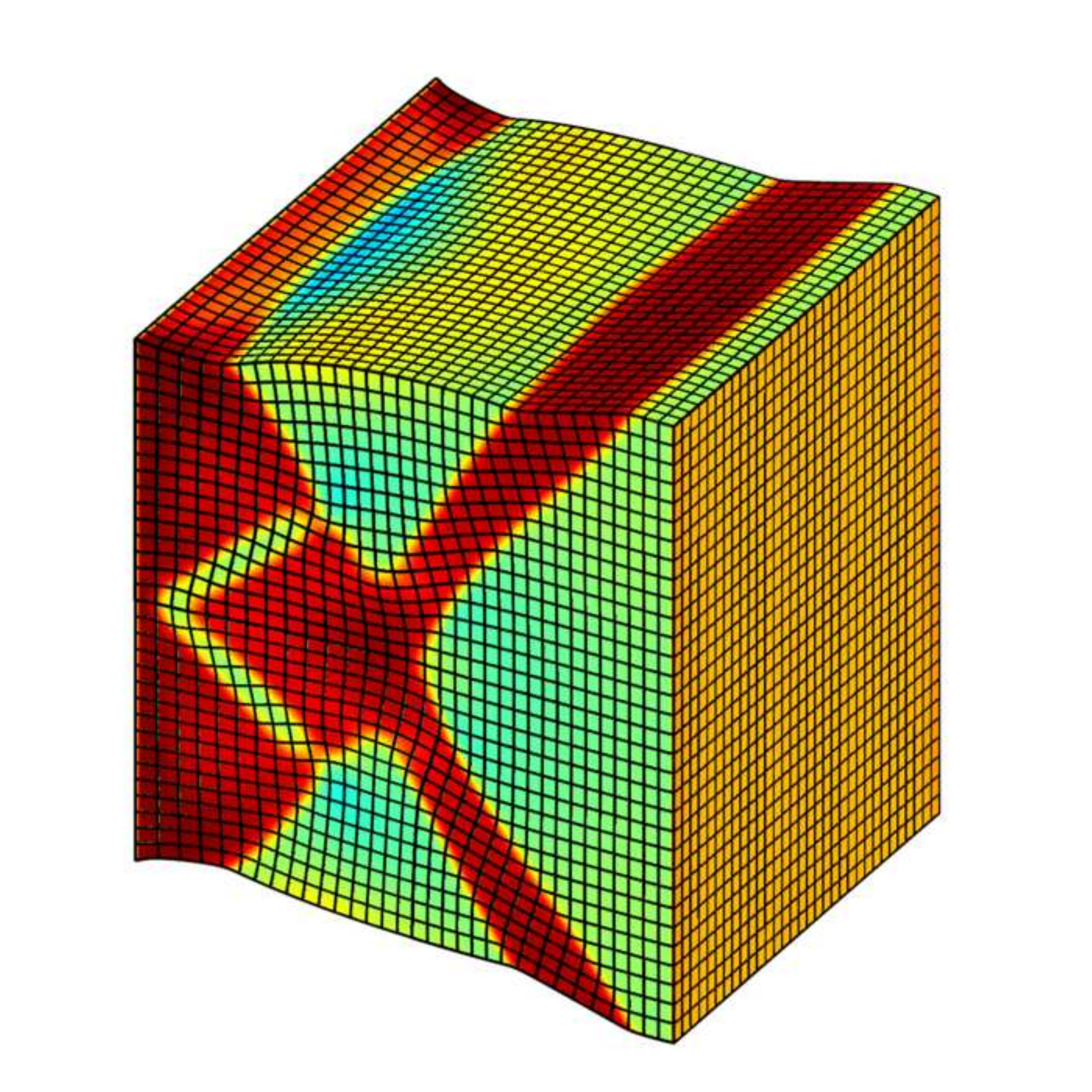} &
                \includegraphics[scale=0.12]{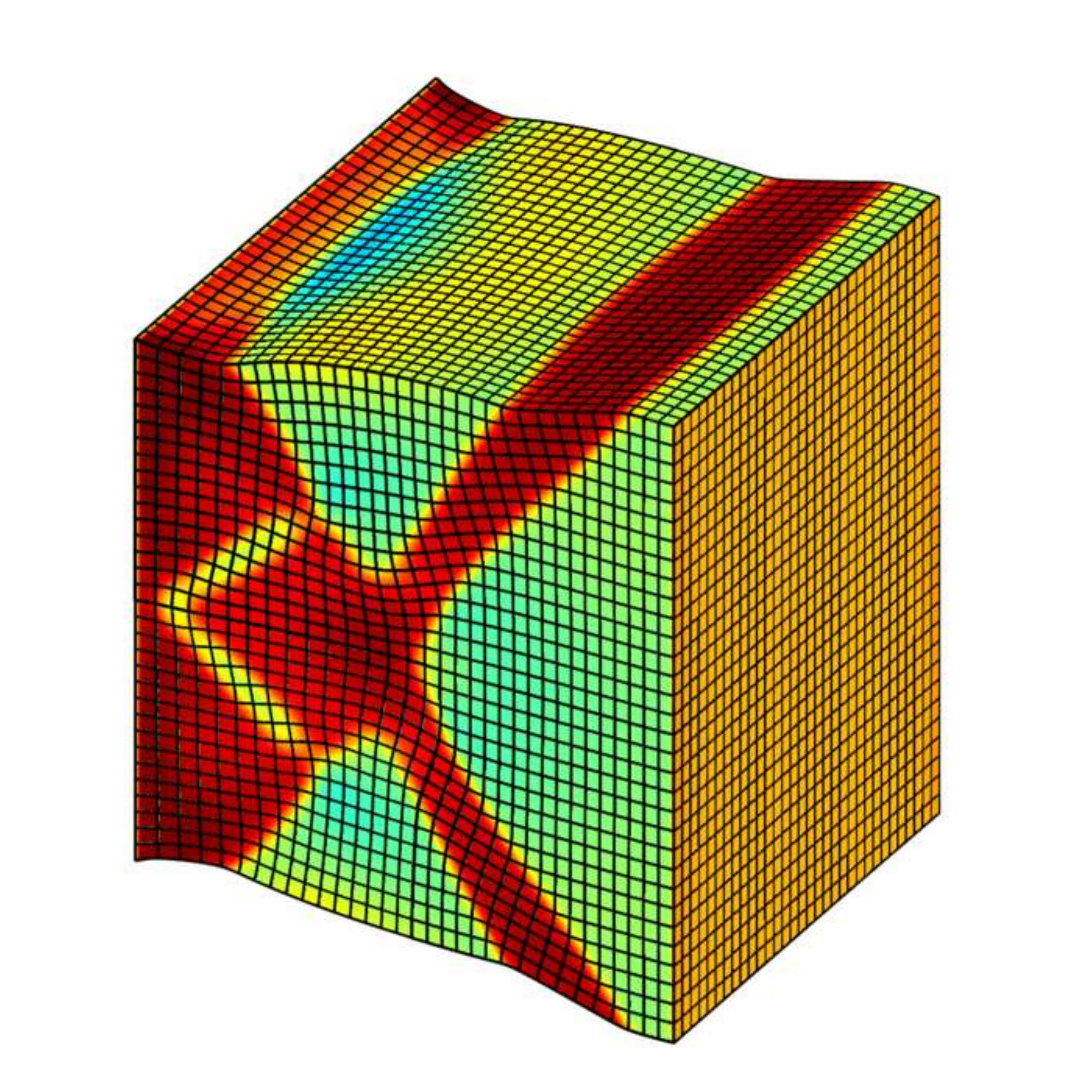} & 
                \includegraphics[scale=0.12]{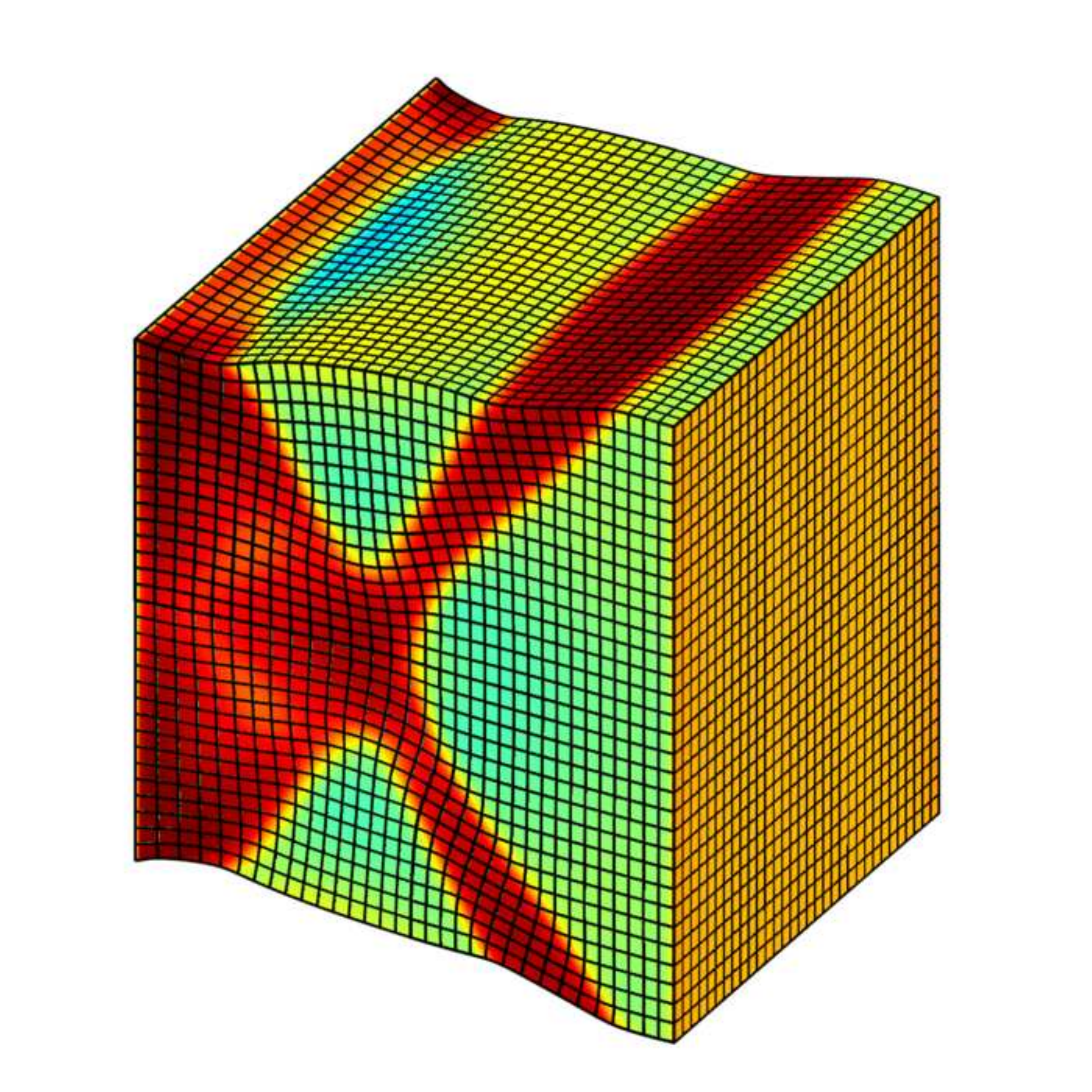} &
                \includegraphics[scale=0.12]{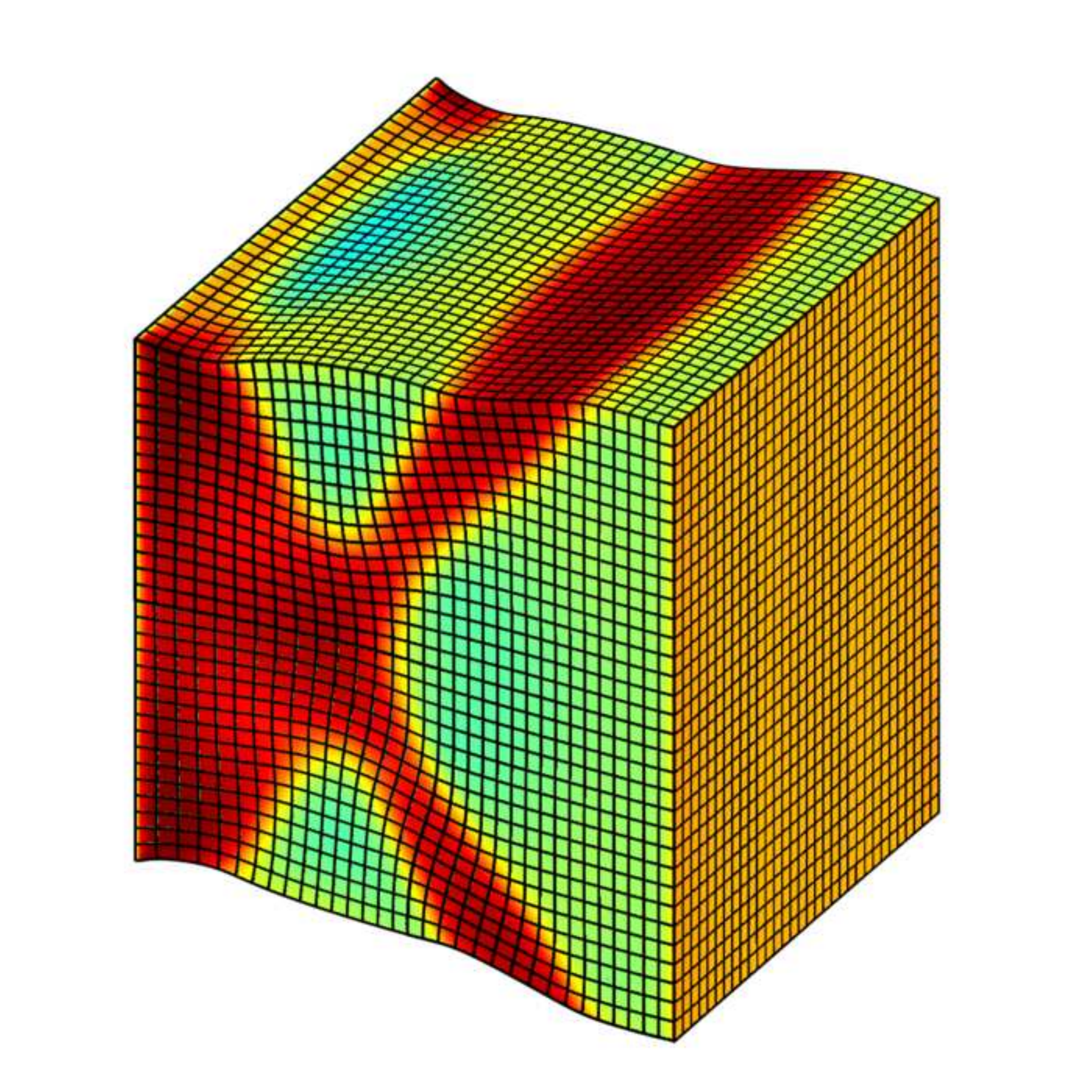} &
                \includegraphics[scale=0.12]{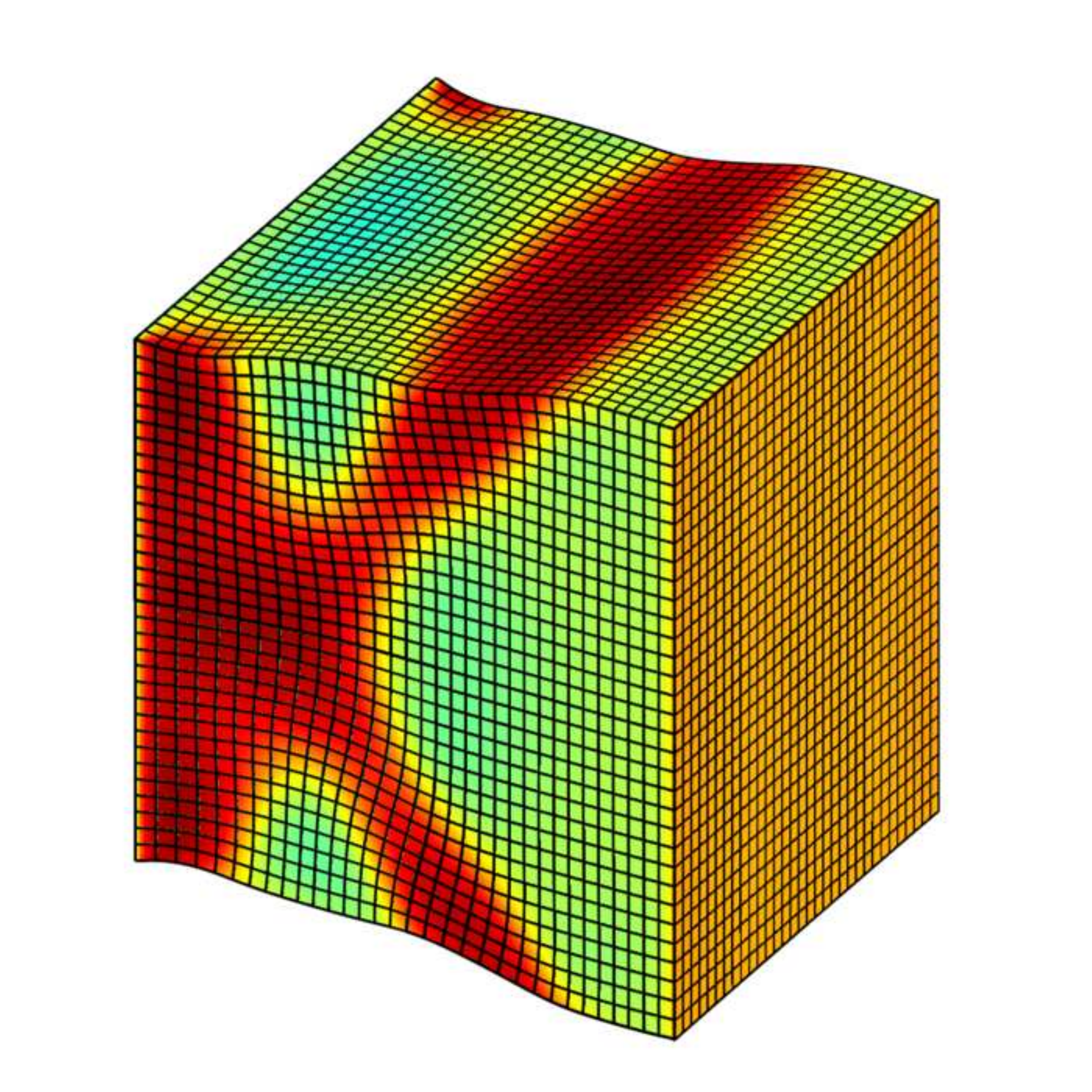} &
            \end{tabular}  \\
            \parbox[t]{0.5cm}{ $B$ } &
            \begin{tabular}{p{2.8cm}p{2.8cm}p{2.8cm}p{2.8cm}p{2.8cm}p{2.5cm}}
                \includegraphics[scale=0.12]{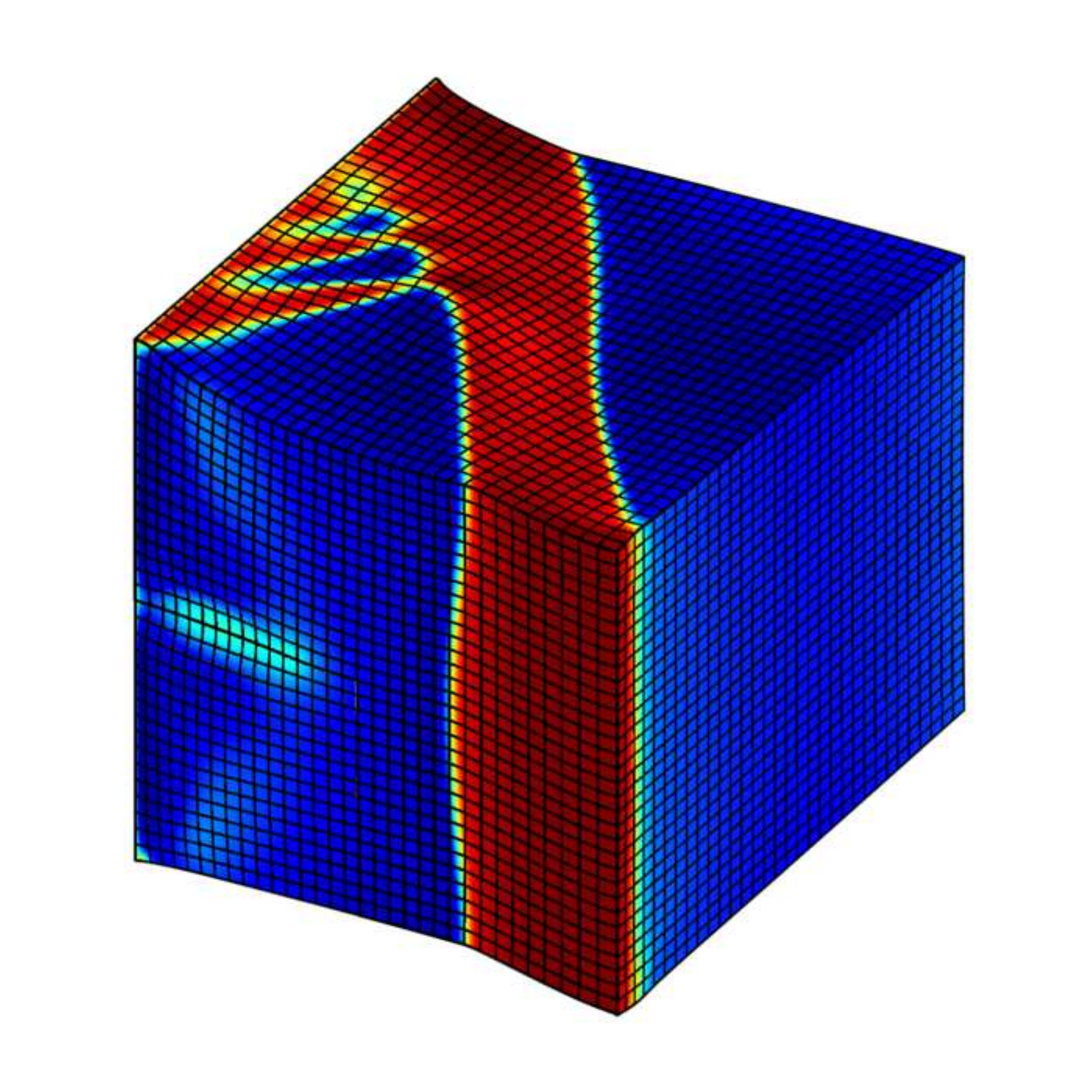} &
                \includegraphics[scale=0.12]{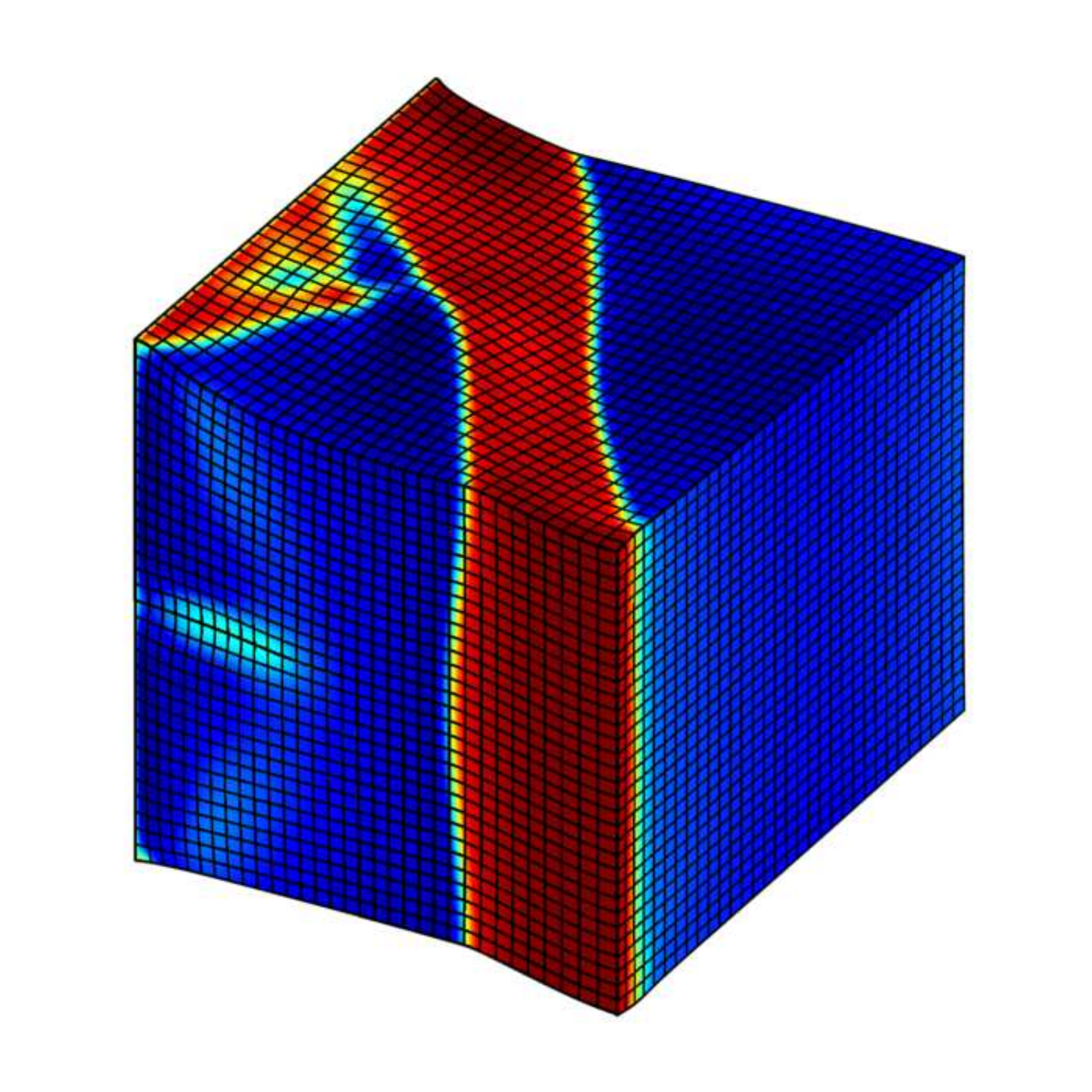} & 
                \includegraphics[scale=0.12]{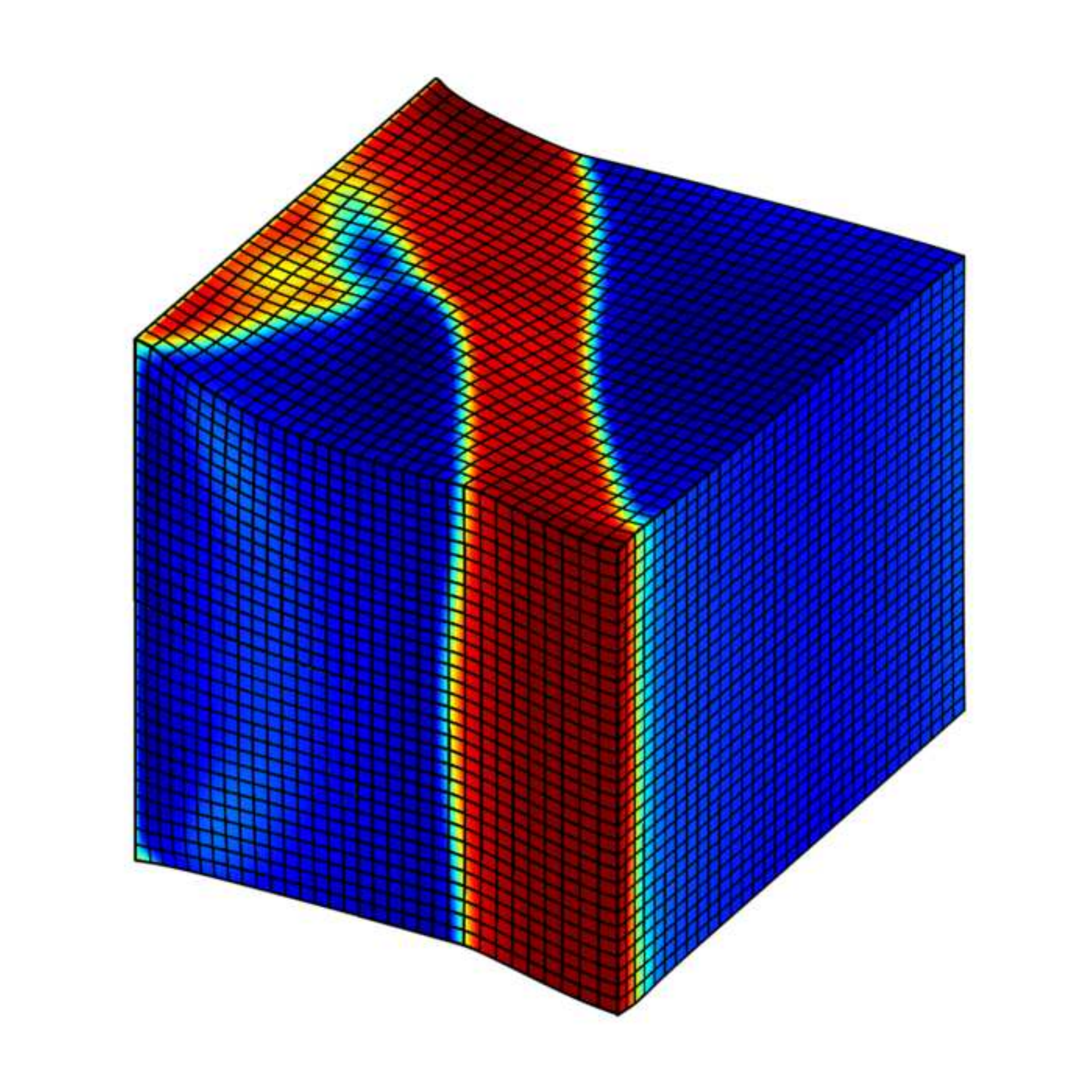} &
                \includegraphics[scale=0.12]{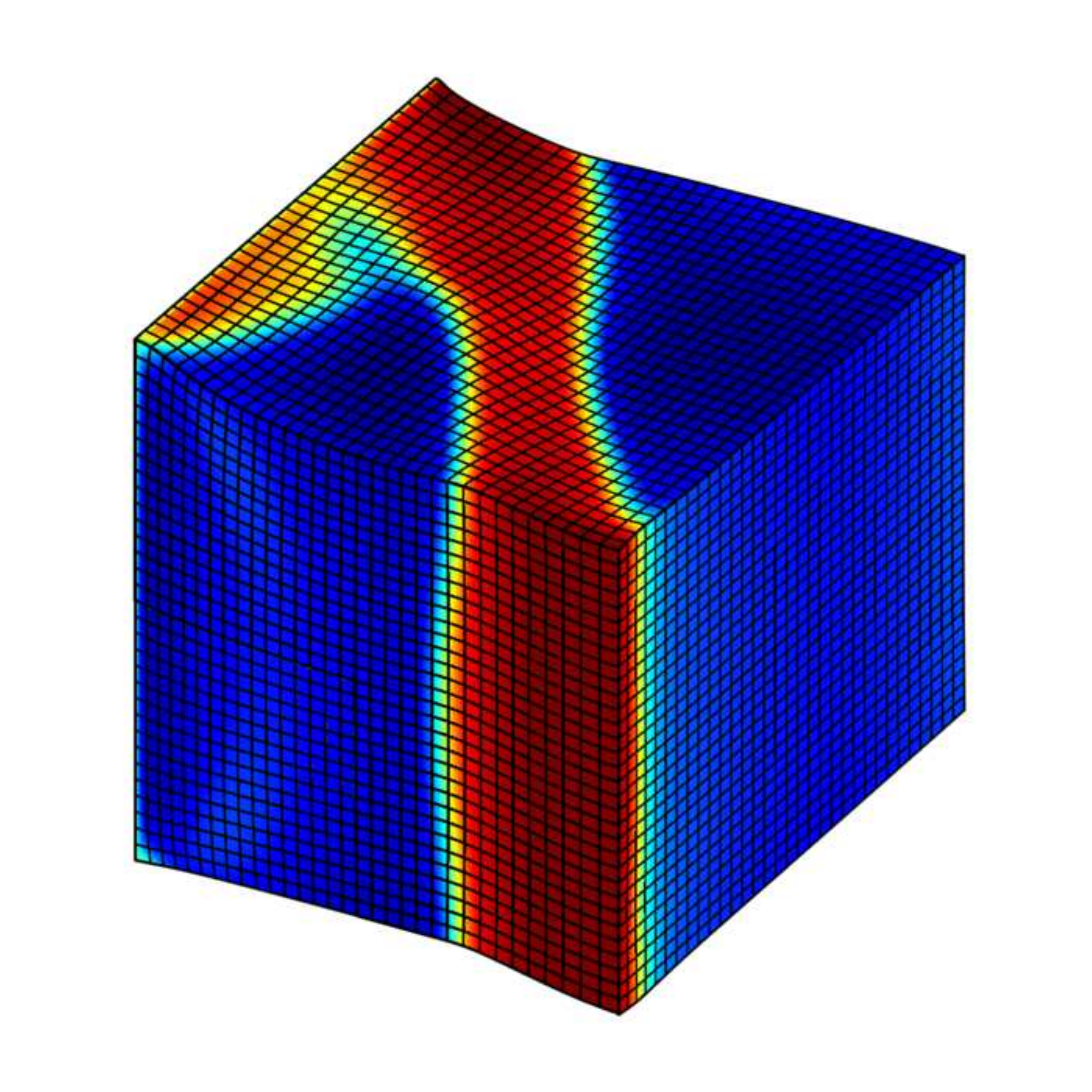} &
                \includegraphics[scale=0.12]{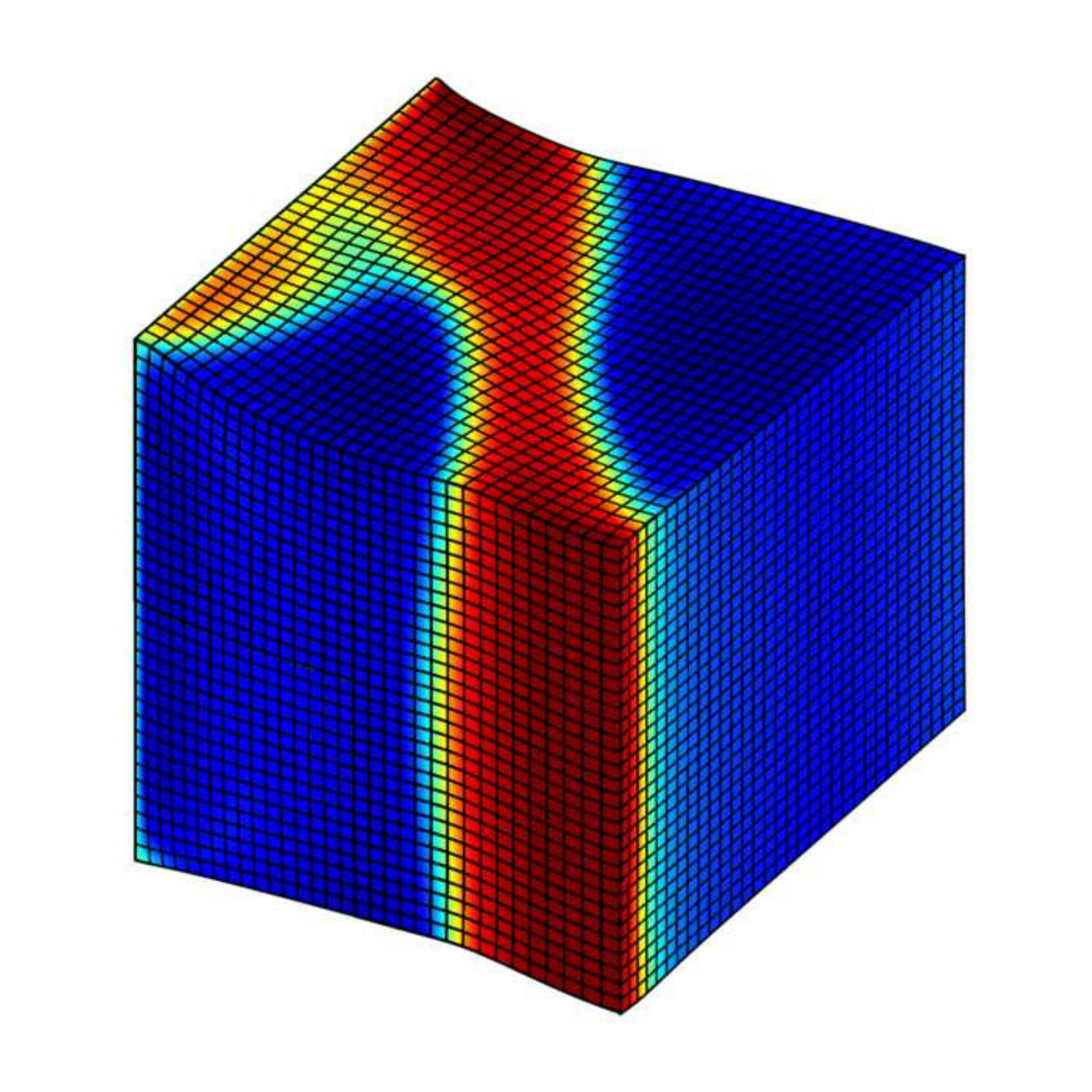} &
                \includegraphics[scale=0.12]{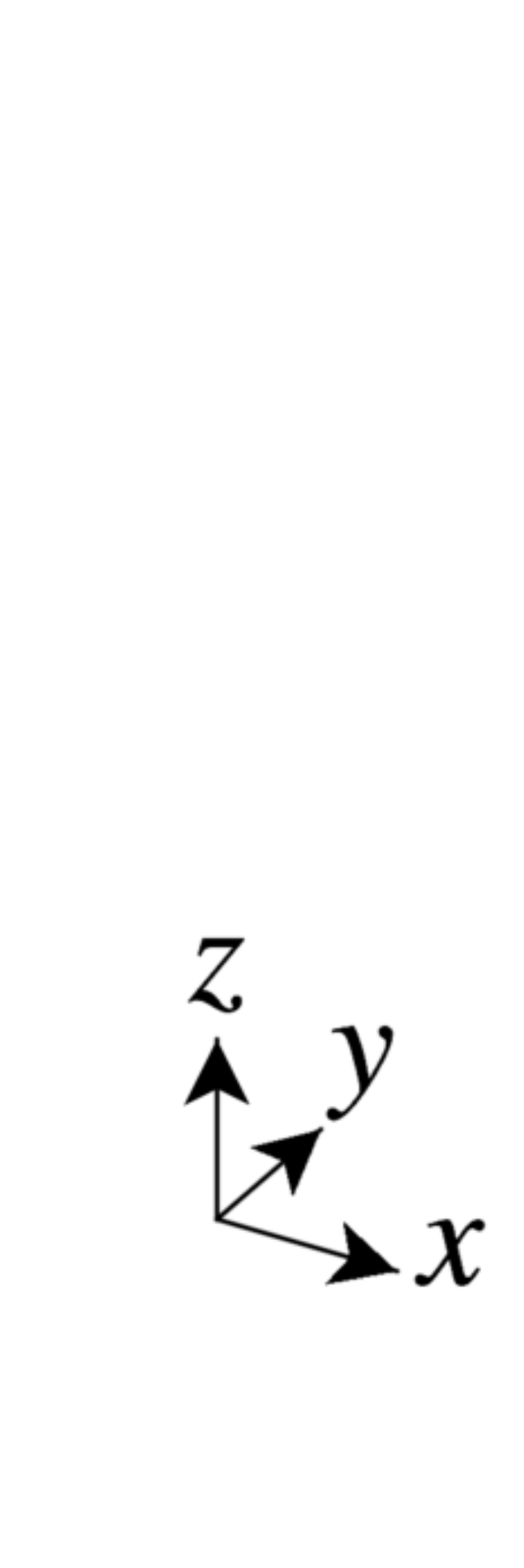}
            \end{tabular}  \\
            \parbox[t]{0.5cm}{ $A$ } &
            \begin{tabular}{p{2.8cm}p{2.8cm}p{2.8cm}p{2.8cm}p{2.8cm}p{2.5cm}}
                &
                \includegraphics[scale=0.12]{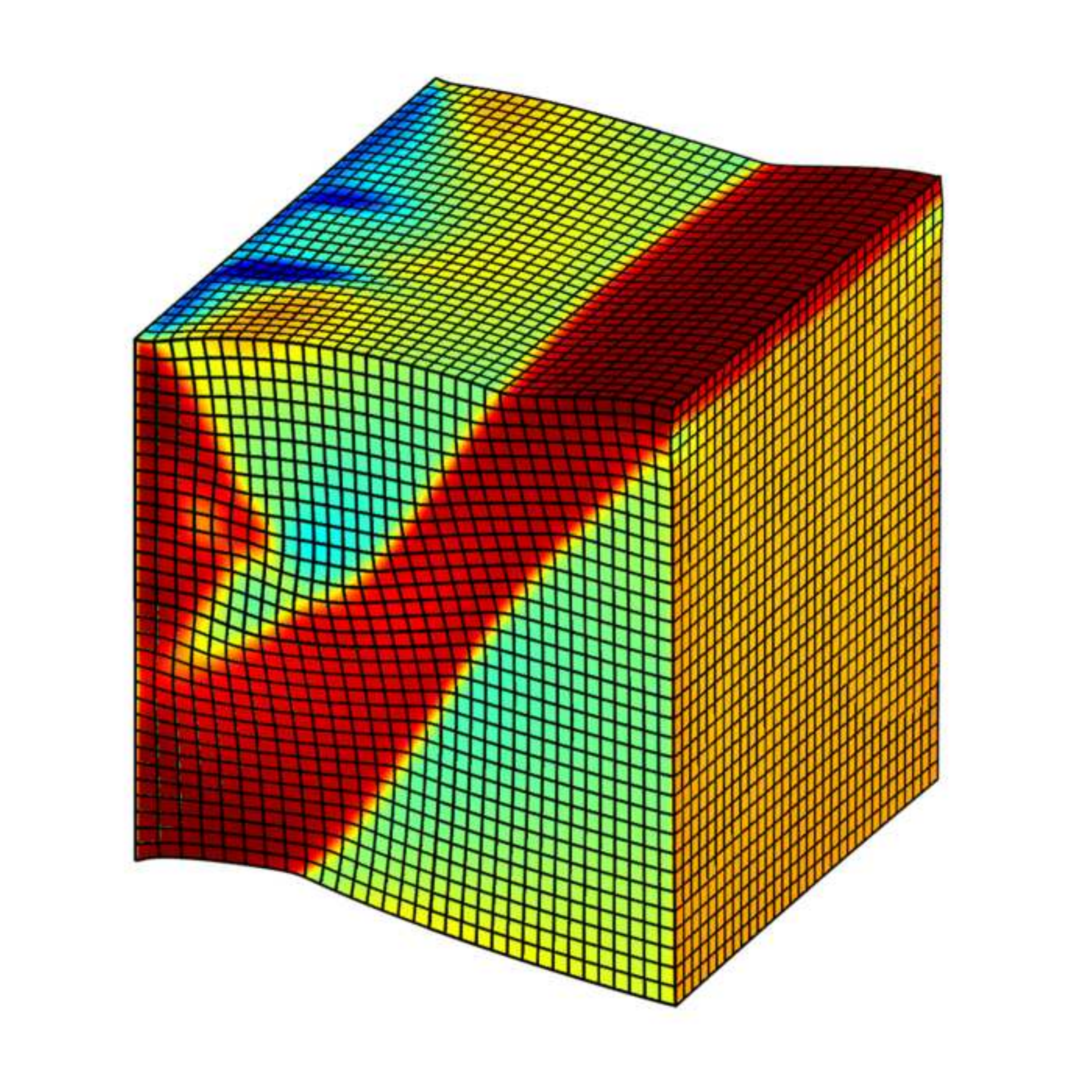} & 
                \includegraphics[scale=0.12]{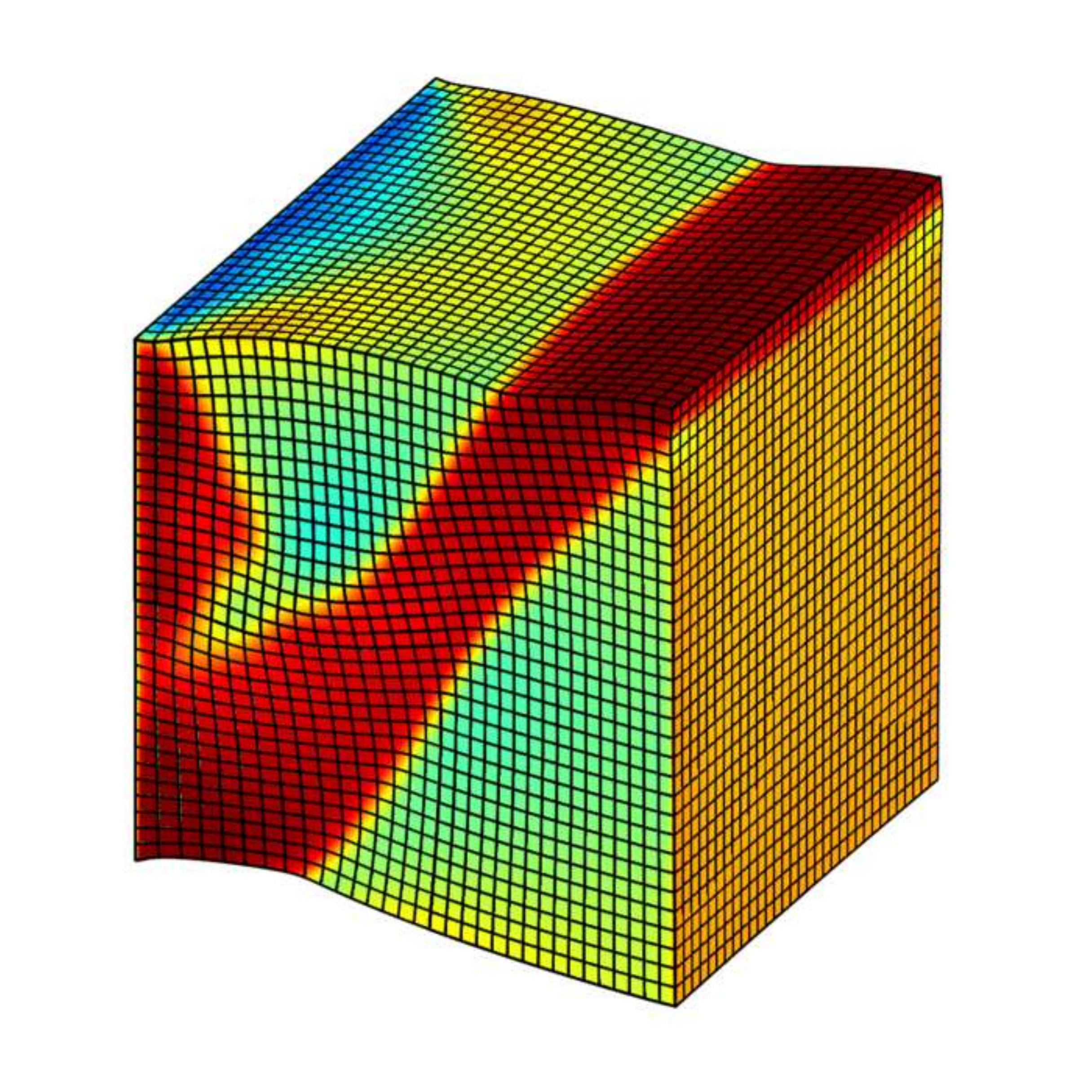} &
                \includegraphics[scale=0.12]{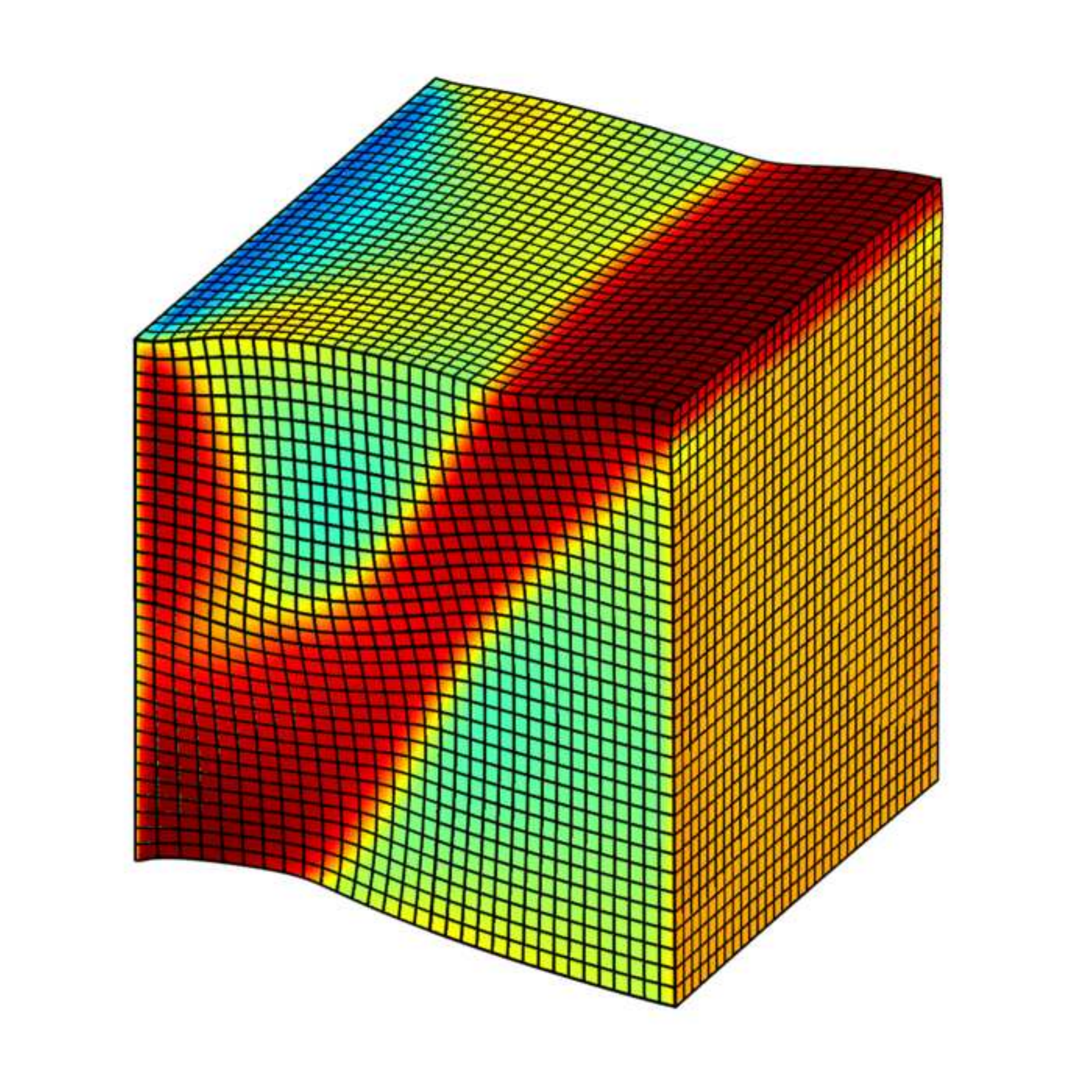} &
                \includegraphics[scale=0.12]{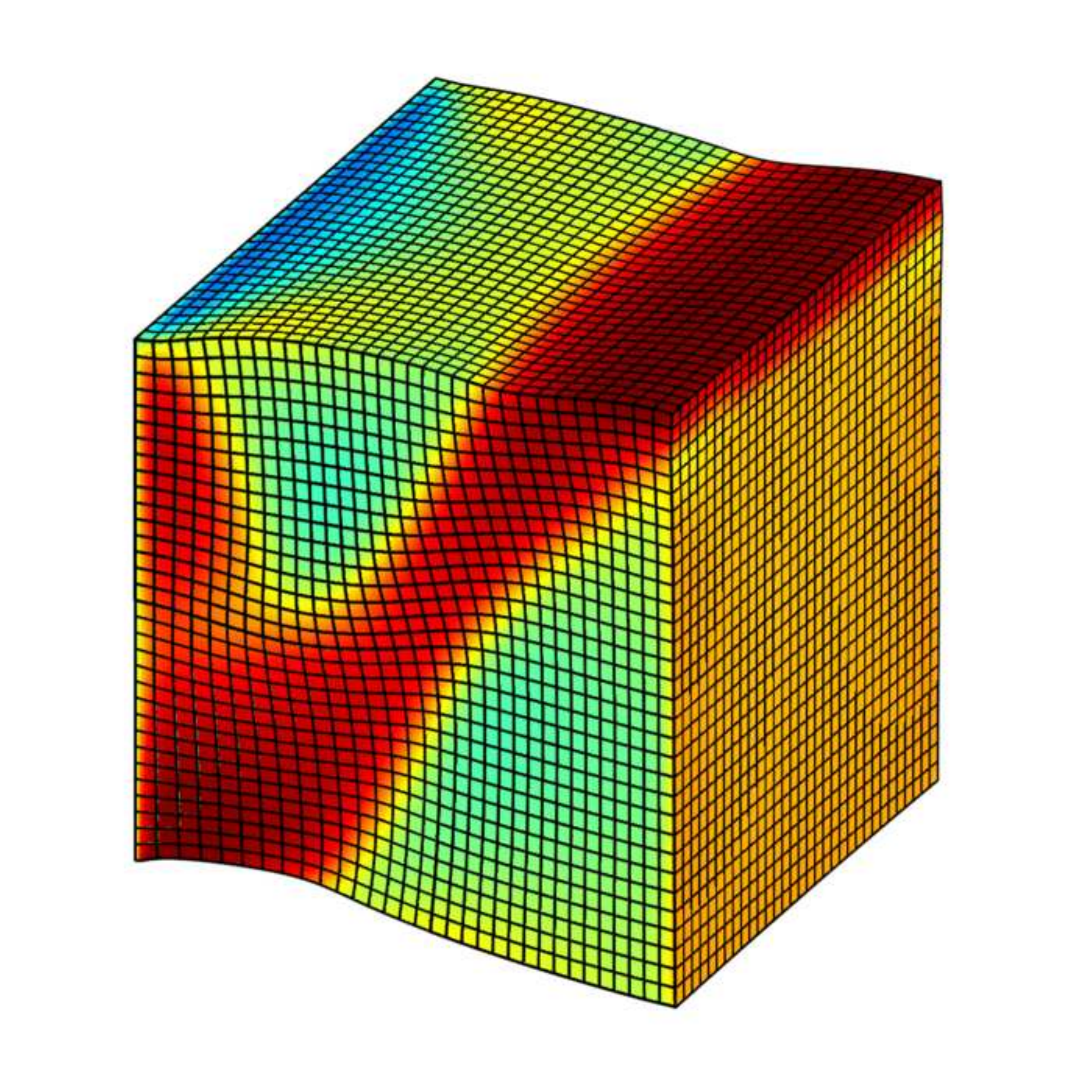} &
                \includegraphics[scale=0.12]{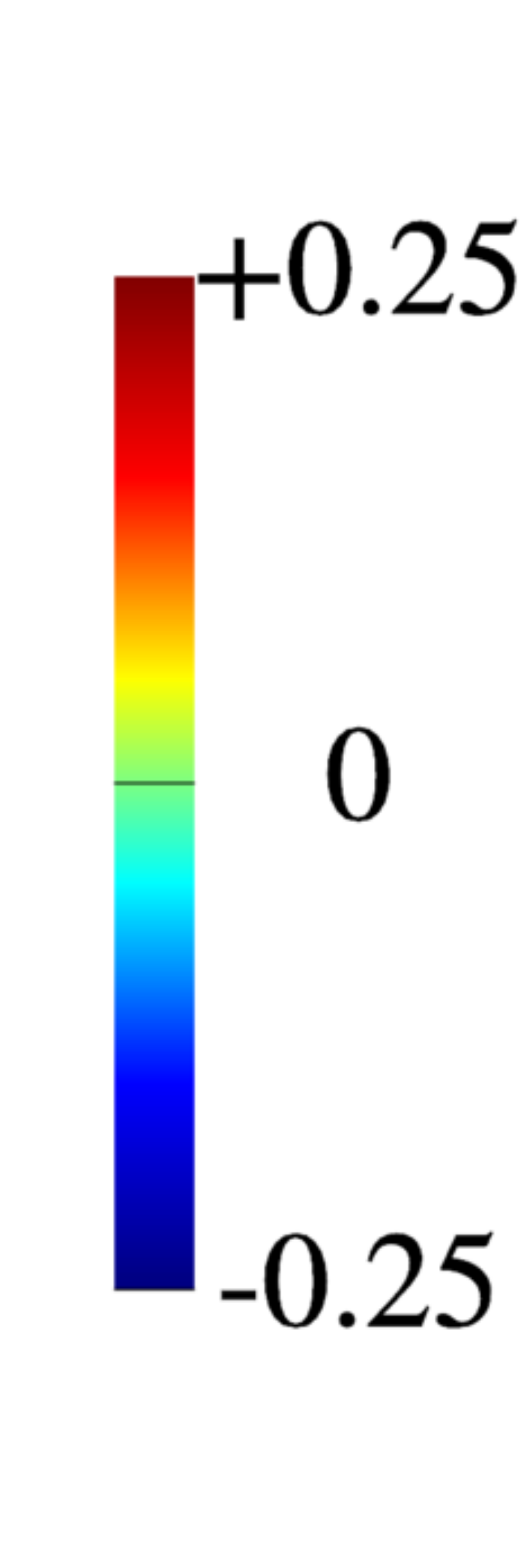}
            \end{tabular}
        \end{tabular}
    \end{center}
    \caption{Field values of $e_2$ for branches A - E on deformed configurations for selected values of $l$.  
    Solutions on the $128^3$ mesh have been overlaid with a $32^3$ plotting mesh.  }
    \label{Fi:e1}
\end{figure}

\begin{figure}
    \begin{center}
        \begin{tabular}{rp{16.5cm}}
            \parbox[t]{0.5cm}{ }&
            \begin{tabular}{p{2.8cm}p{2.8cm}p{2.8cm}p{2.8cm}p{2.8cm}p{2.5cm}}
                \hspace{0.8cm}$l\!=\!0.0625$ & 
                \hspace{0.8cm}$l\!=\!0.0750$ & 
                \hspace{0.8cm}$l\!=\!0.1000$ &
                \hspace{0.8cm}$l\!=\!0.1500$ &
                \hspace{0.8cm}$l\!=\!0.2000$ &
                \vspace{0.5\baselineskip}
            \end{tabular}  \\
            \parbox[t]{0.5cm}{ $E$ } &
            \begin{tabular}{p{2.8cm}p{2.8cm}p{2.8cm}p{2.8cm}p{2.8cm}p{2.5cm}}
                &
                \includegraphics[scale=0.12]{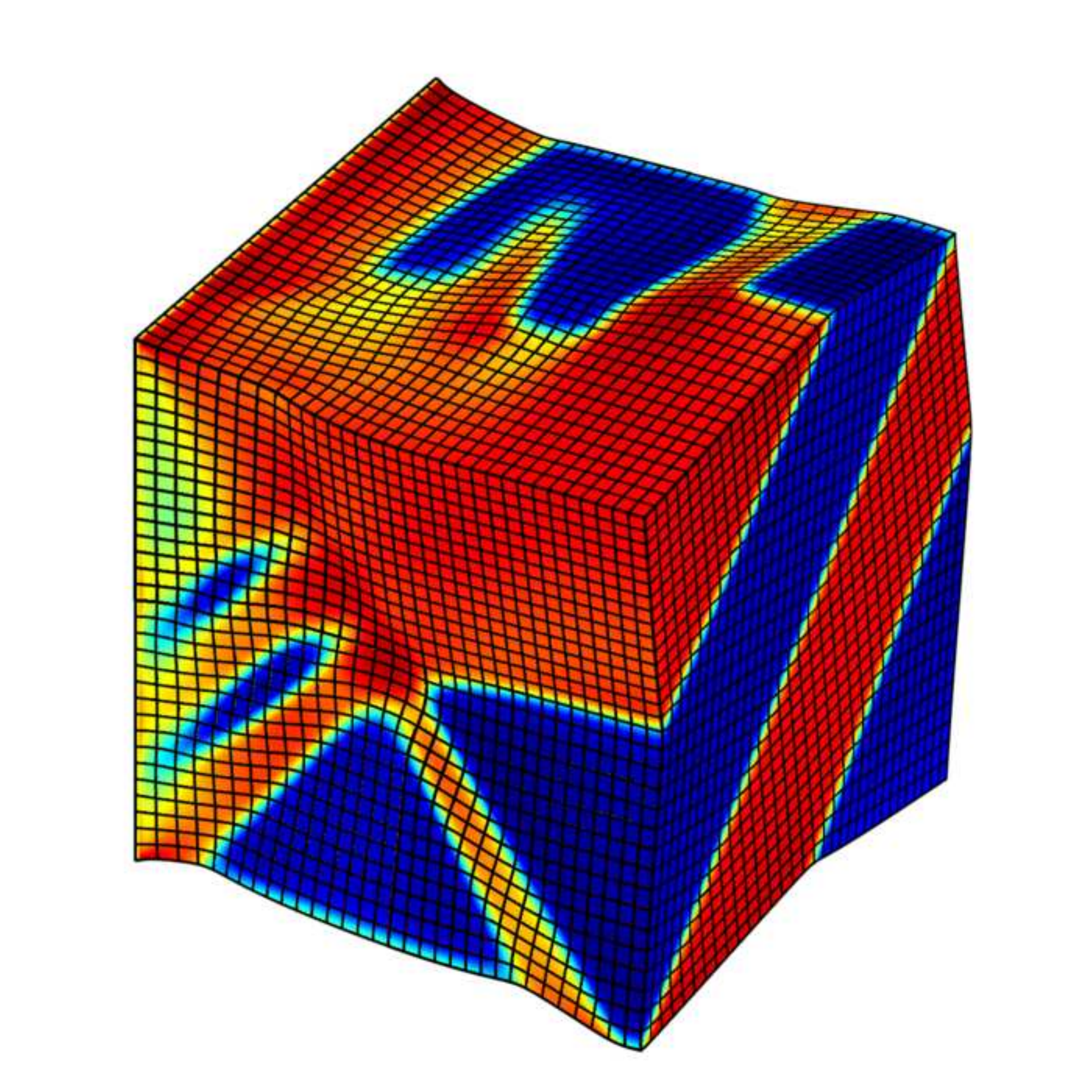} & 
                \includegraphics[scale=0.12]{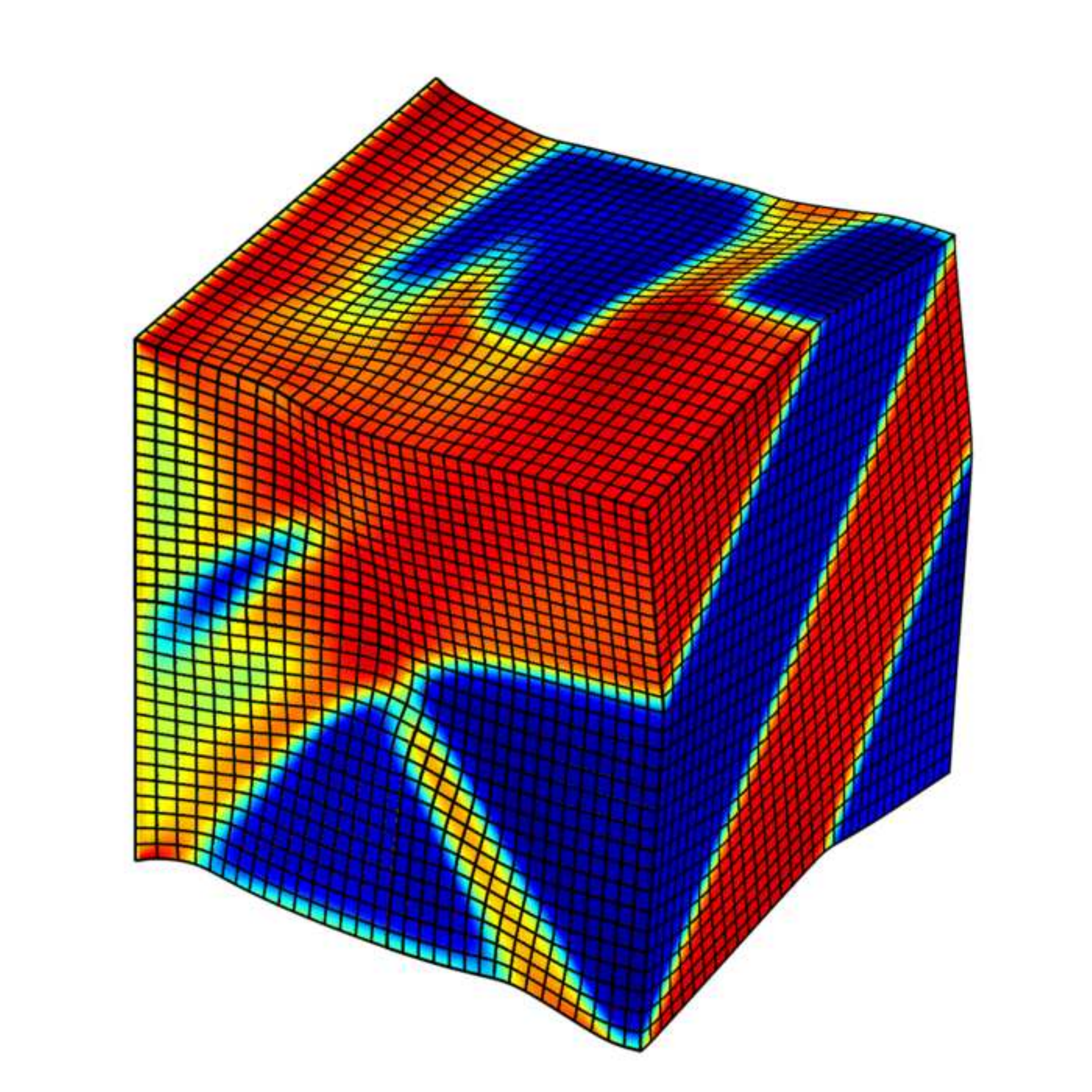} &
                \includegraphics[scale=0.12]{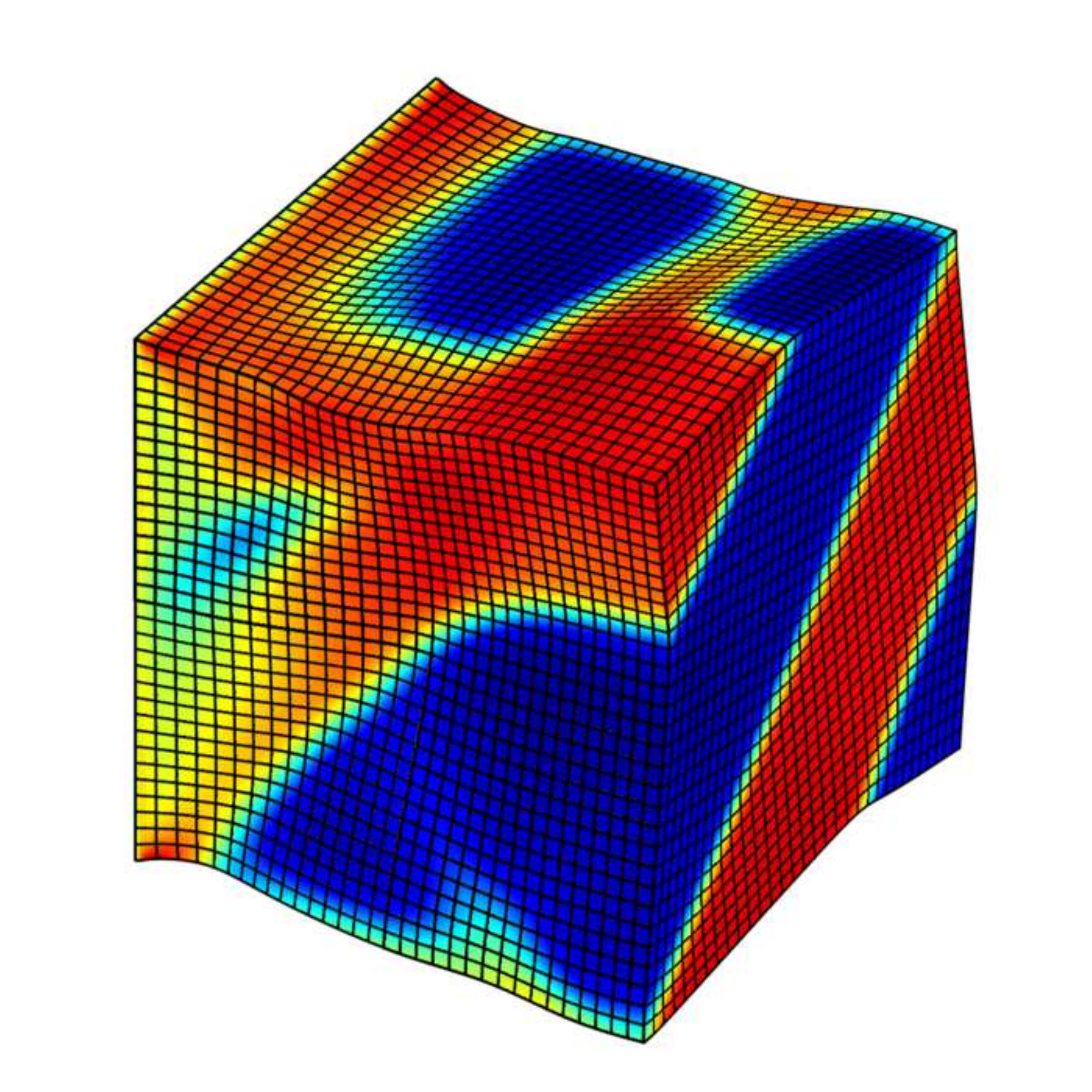} &
                \includegraphics[scale=0.12]{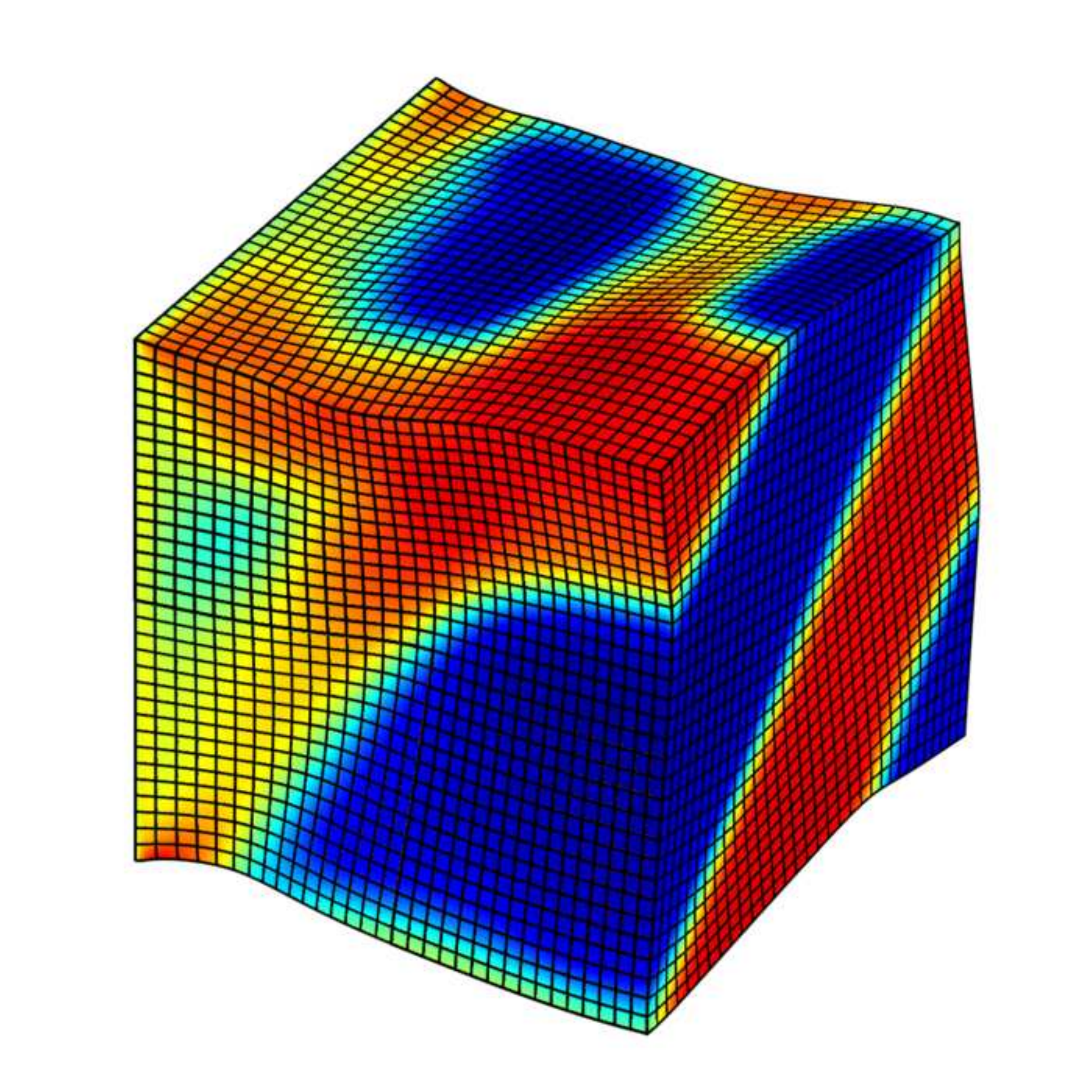} &
            \end{tabular}  \\
            \parbox[t]{0.5cm}{ $D$ } &
            \begin{tabular}{p{2.8cm}p{2.8cm}p{2.8cm}p{2.8cm}p{2.8cm}p{2.5cm}}
                \includegraphics[scale=0.12]{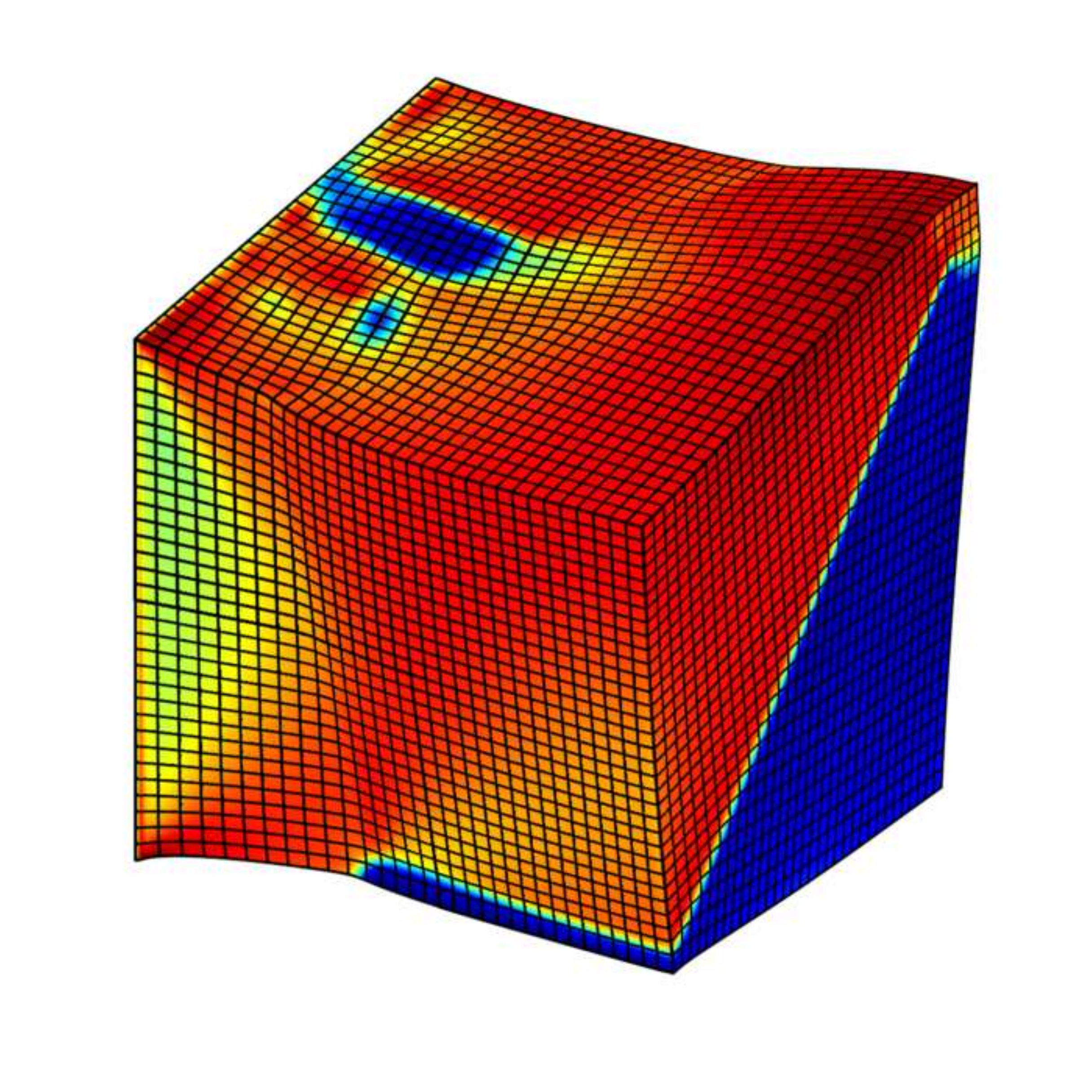} &
                \includegraphics[scale=0.12]{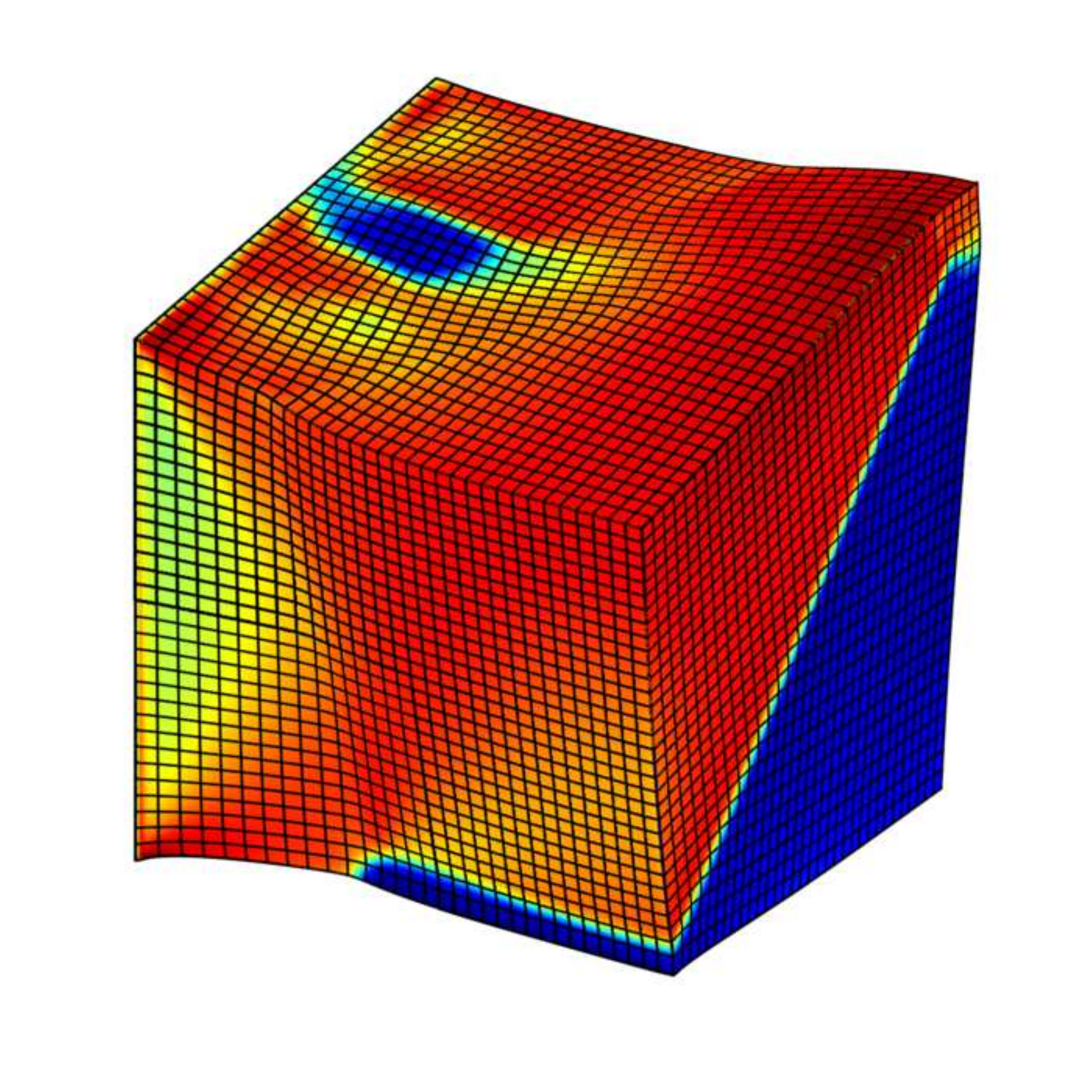} & 
                \includegraphics[scale=0.12]{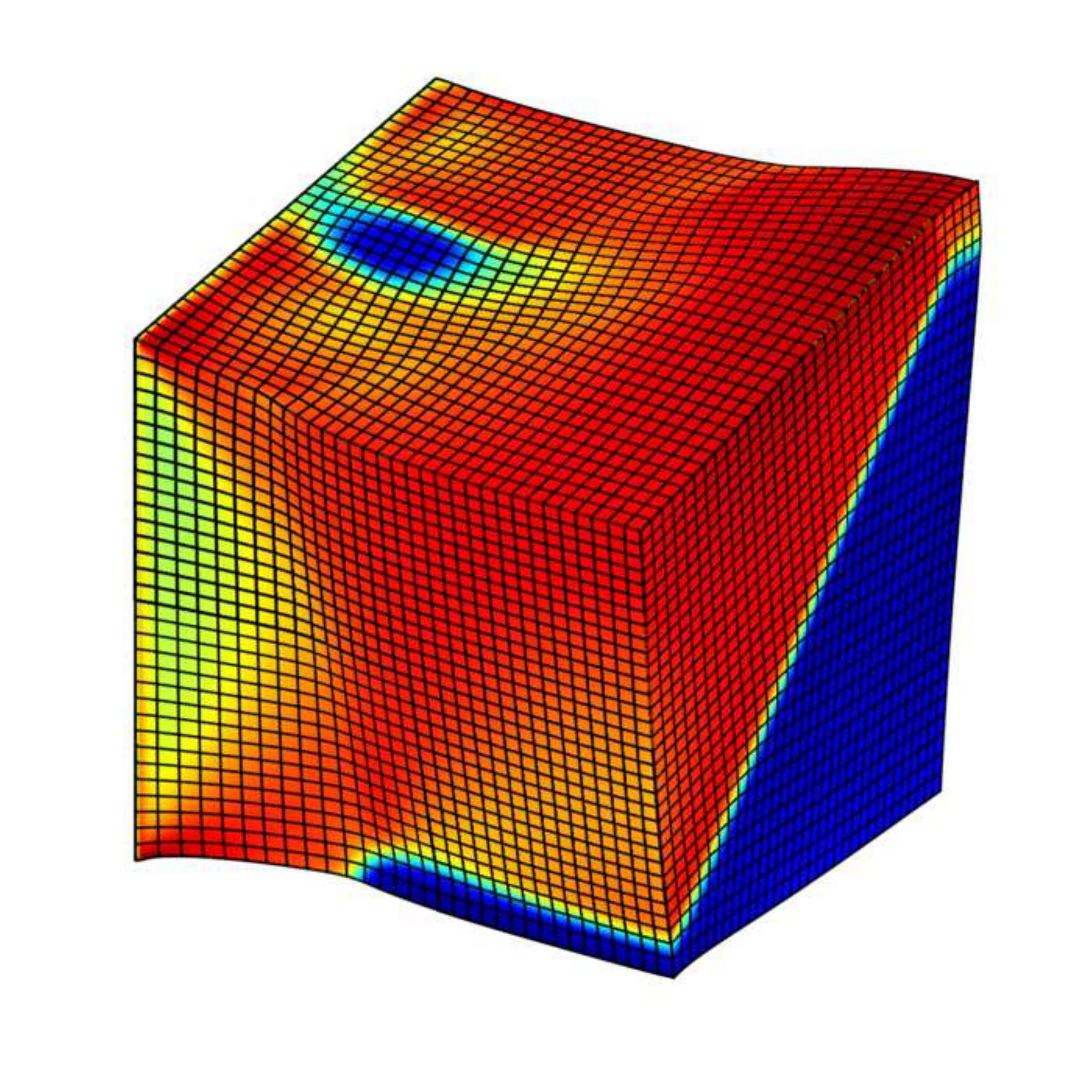} &
                \includegraphics[scale=0.12]{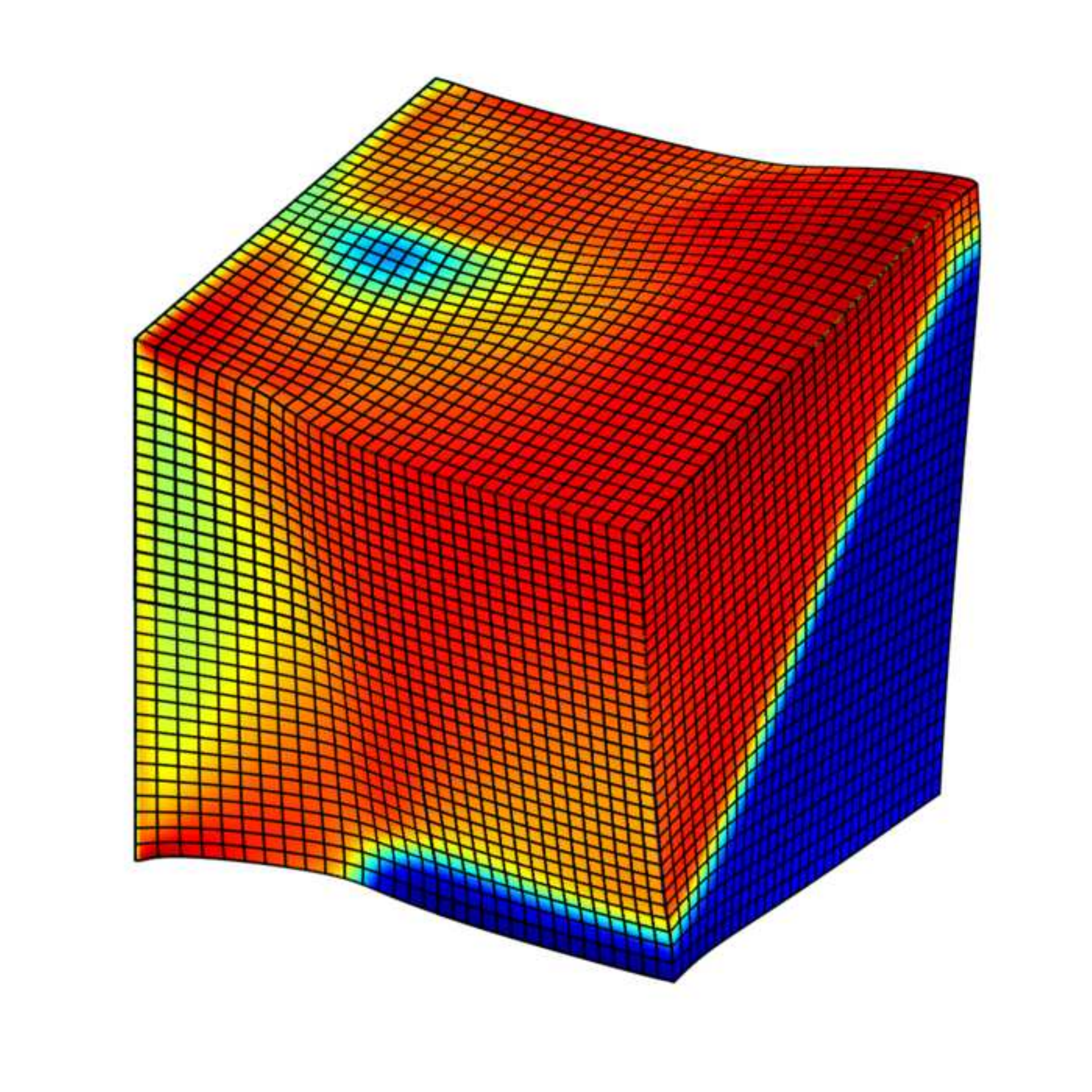} &
                \includegraphics[scale=0.12]{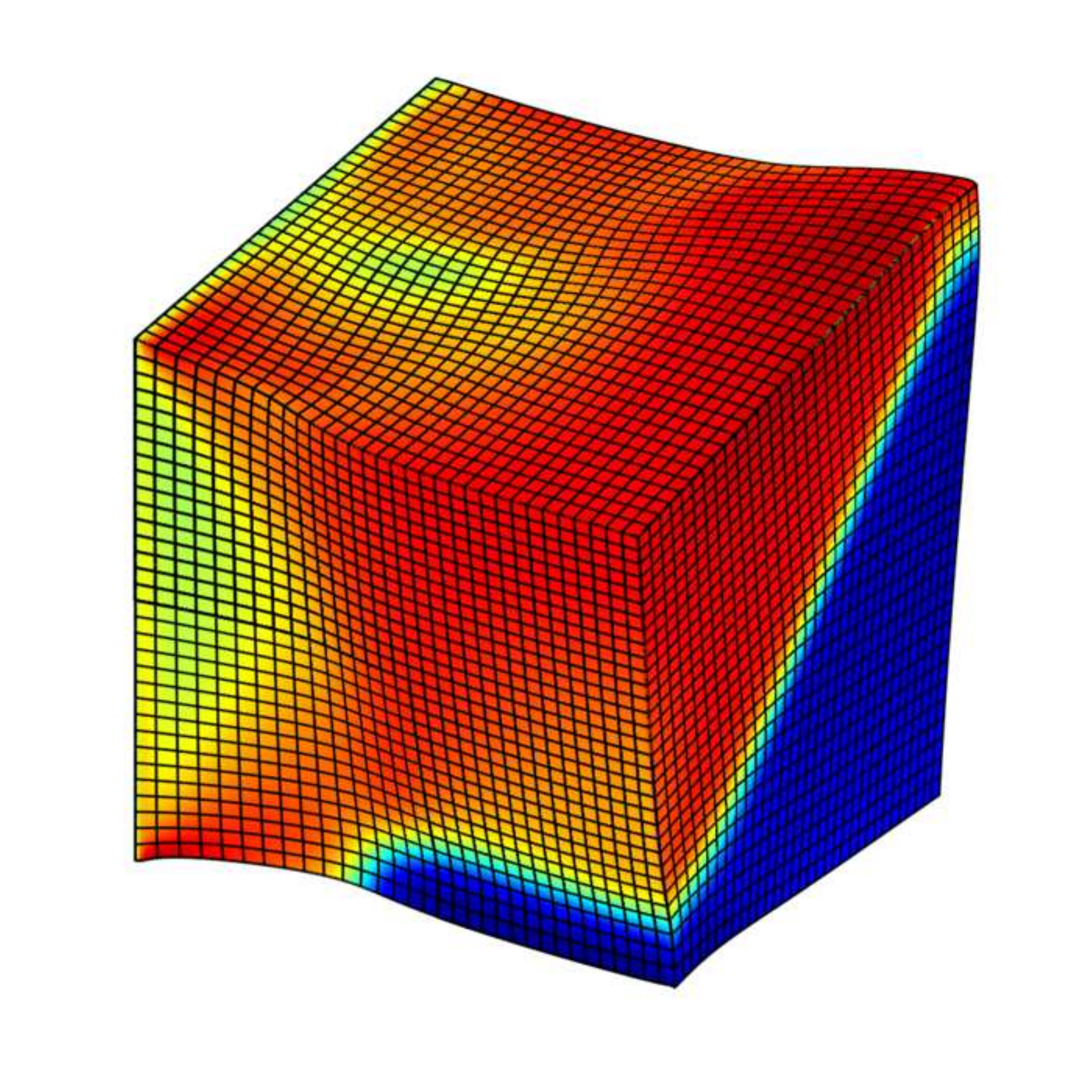} &
            \end{tabular}  \\
            \parbox[t]{0.5cm}{ $C$ } &
            \begin{tabular}{p{2.8cm}p{2.8cm}p{2.8cm}p{2.8cm}p{2.8cm}p{2.5cm}}
                \includegraphics[scale=0.12]{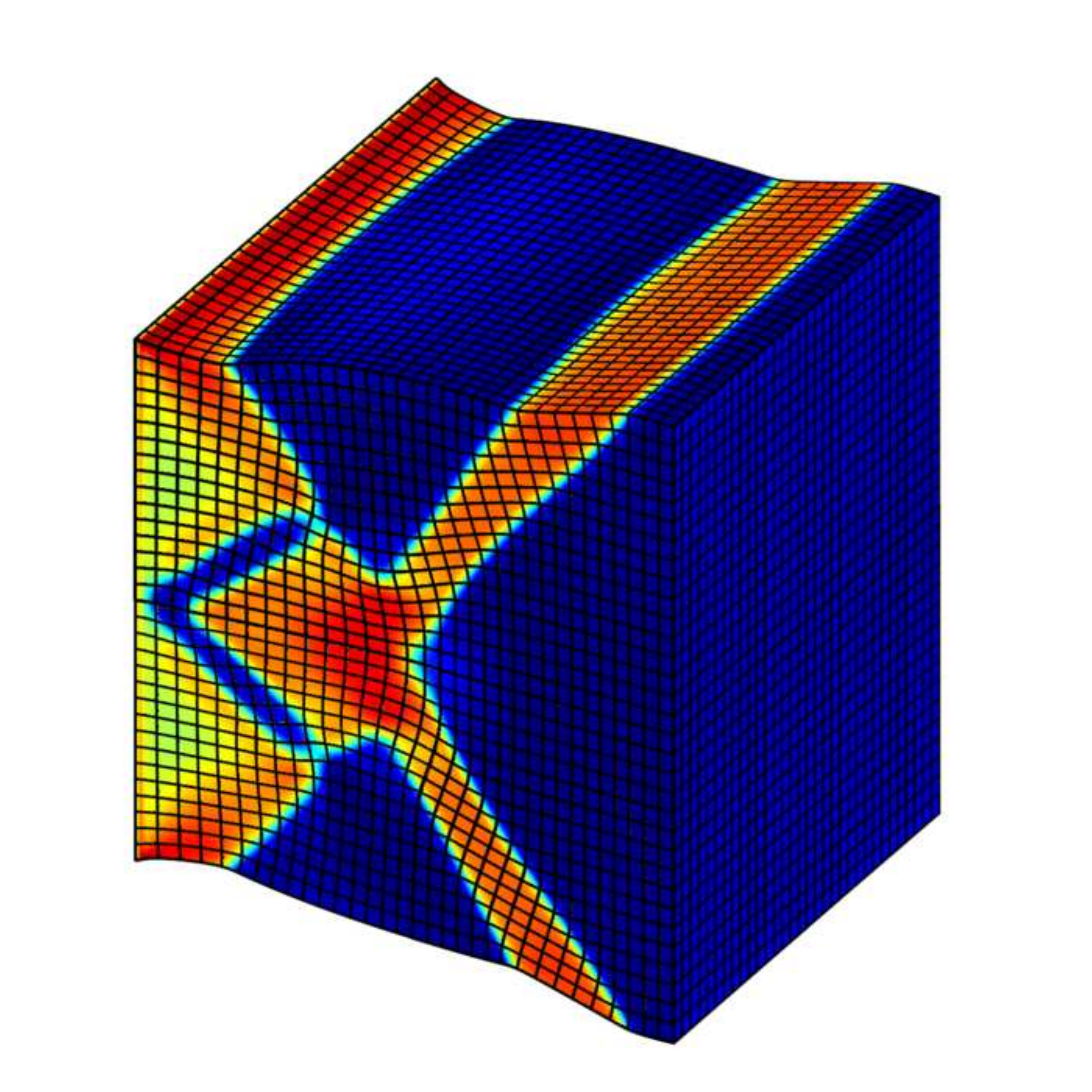} &
                \includegraphics[scale=0.12]{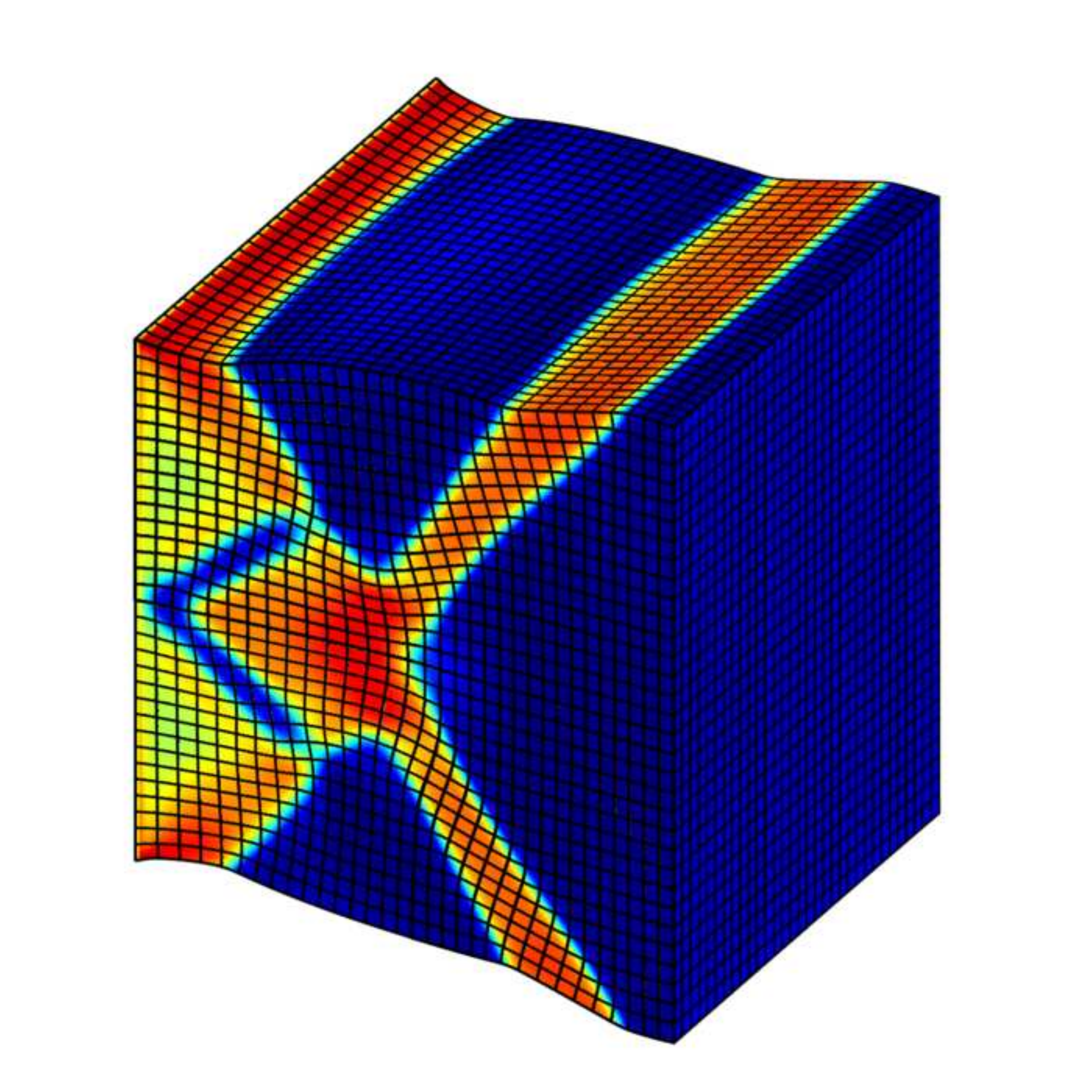} & 
                \includegraphics[scale=0.12]{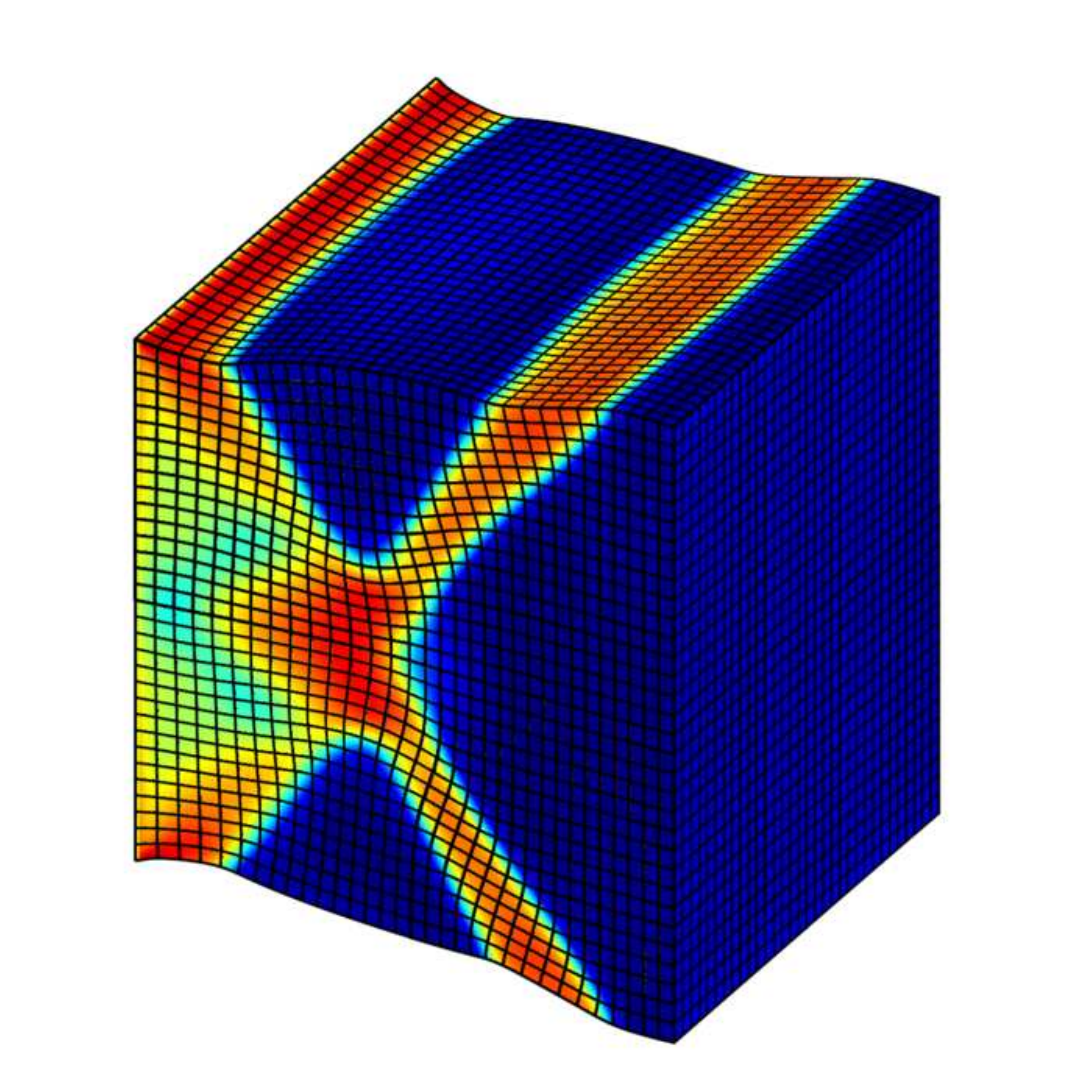} &
                \includegraphics[scale=0.12]{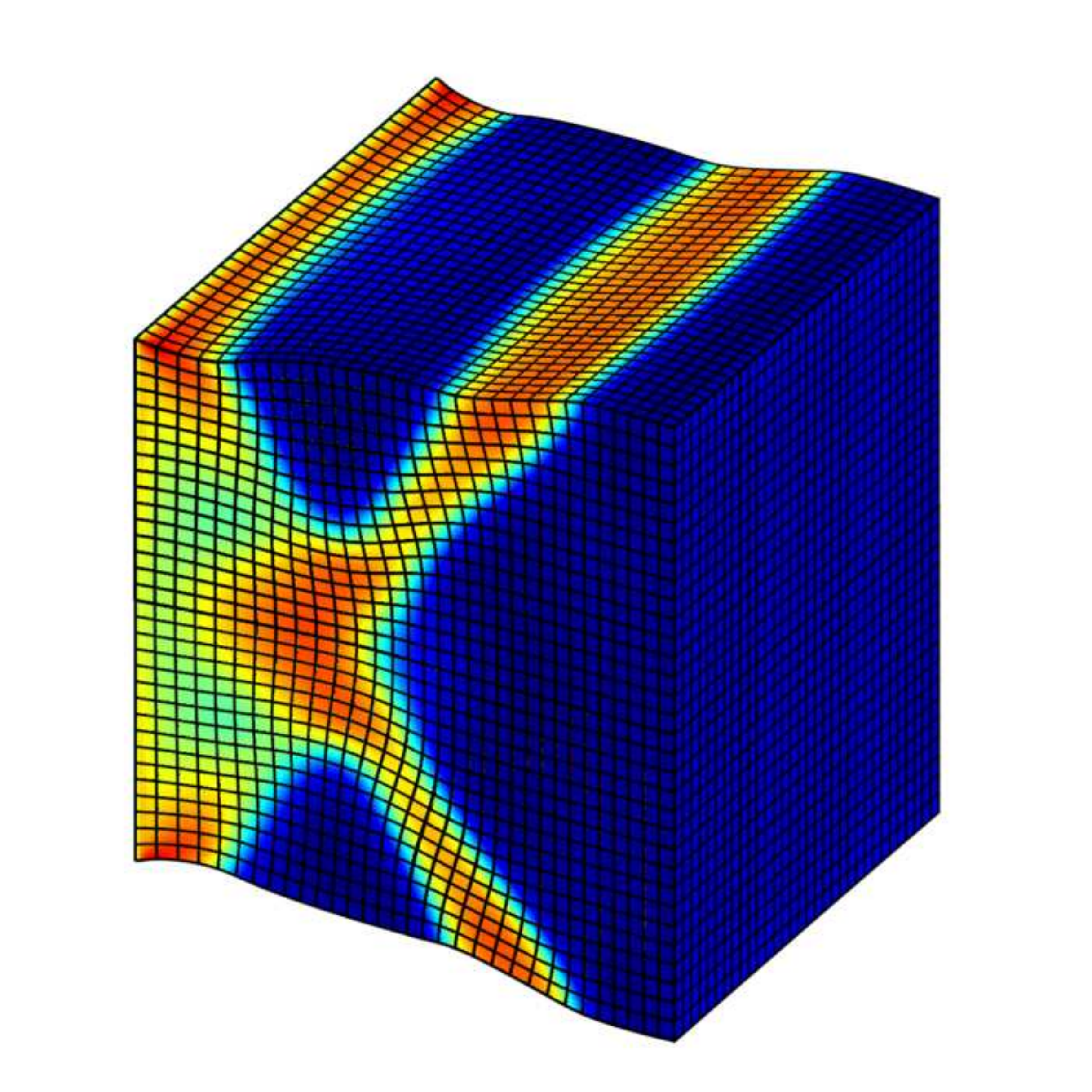} &
                \includegraphics[scale=0.12]{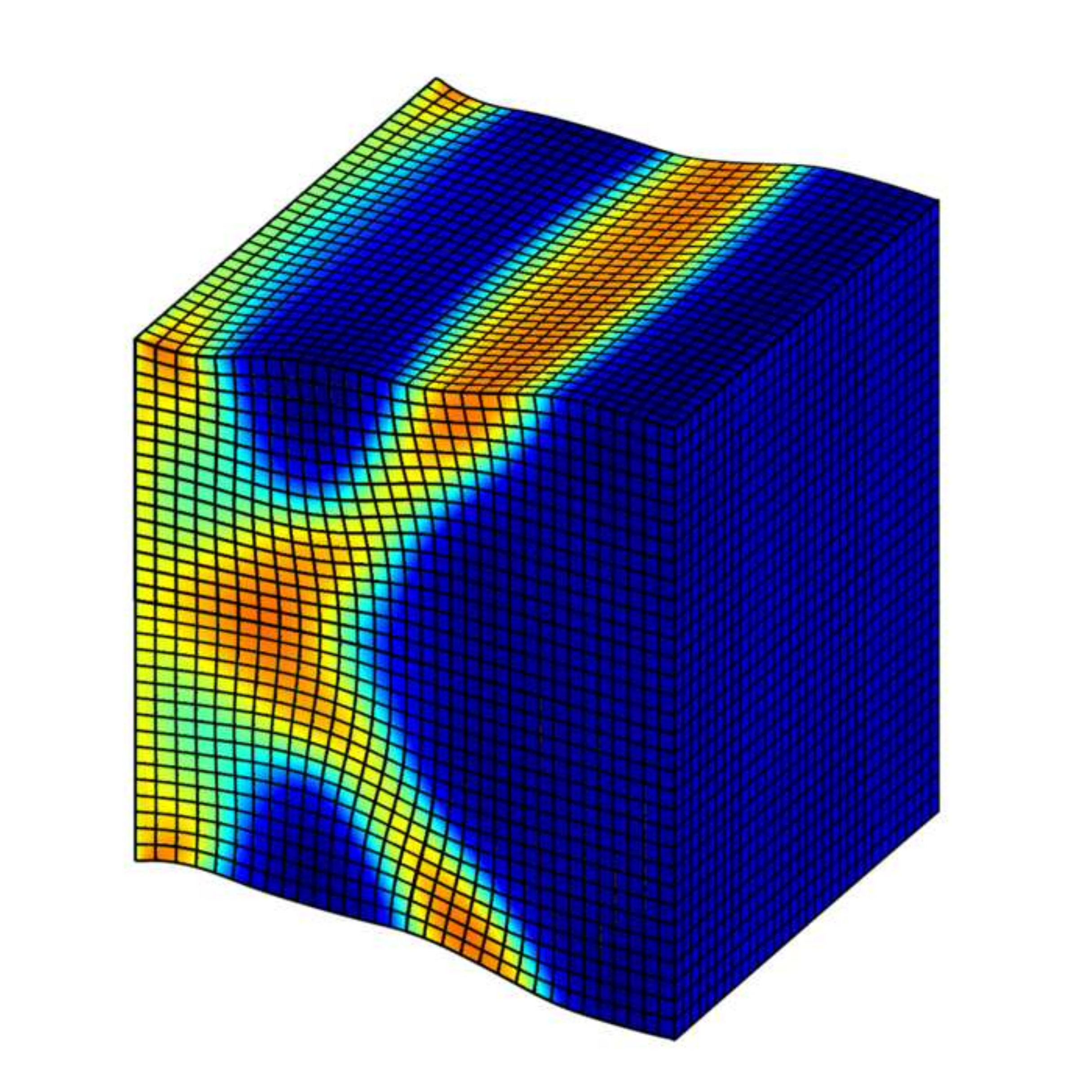} &
            \end{tabular}  \\
            \parbox[t]{0.5cm}{ $B$ } &
            \begin{tabular}{p{2.8cm}p{2.8cm}p{2.8cm}p{2.8cm}p{2.8cm}p{2.5cm}}
                \includegraphics[scale=0.12]{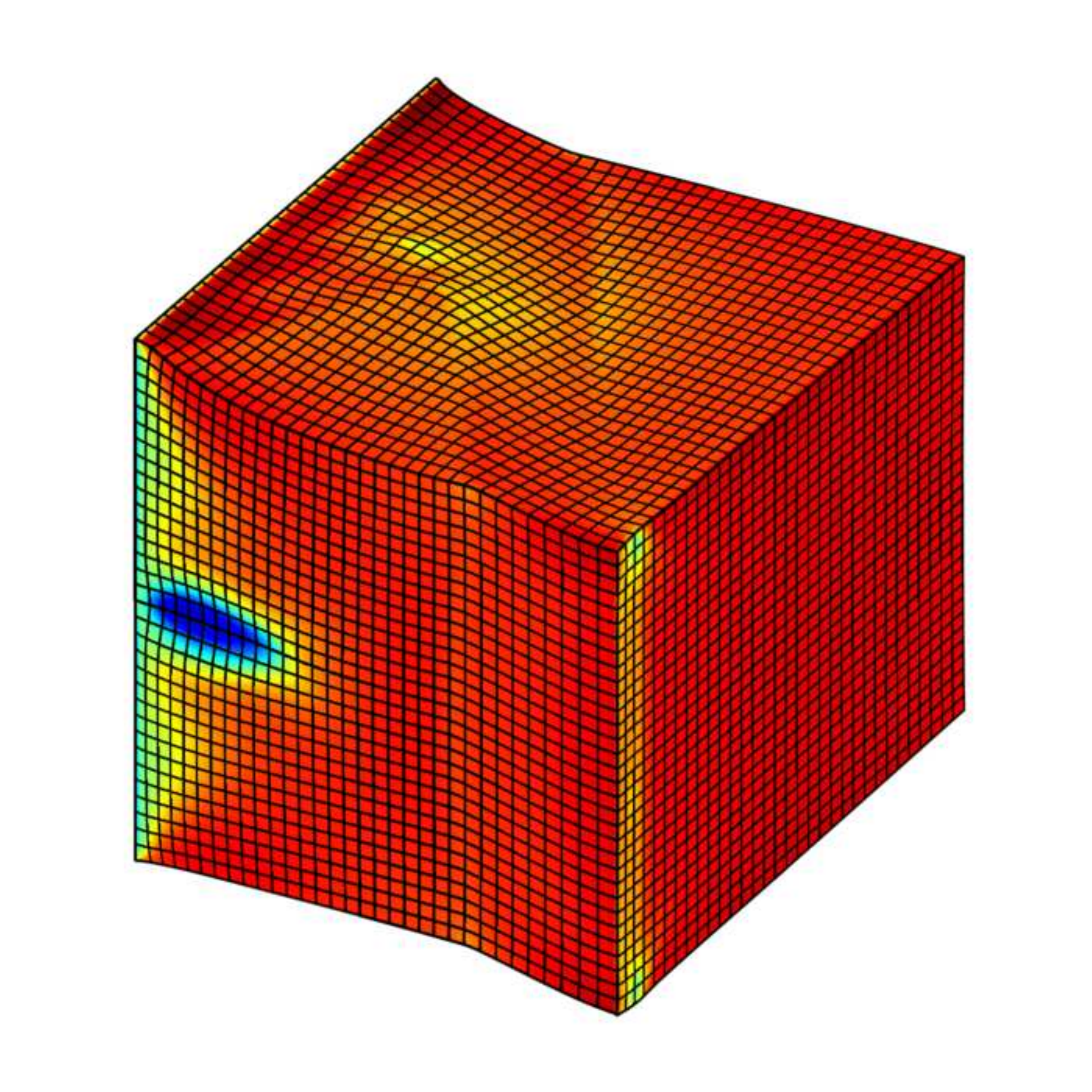} &
                \includegraphics[scale=0.12]{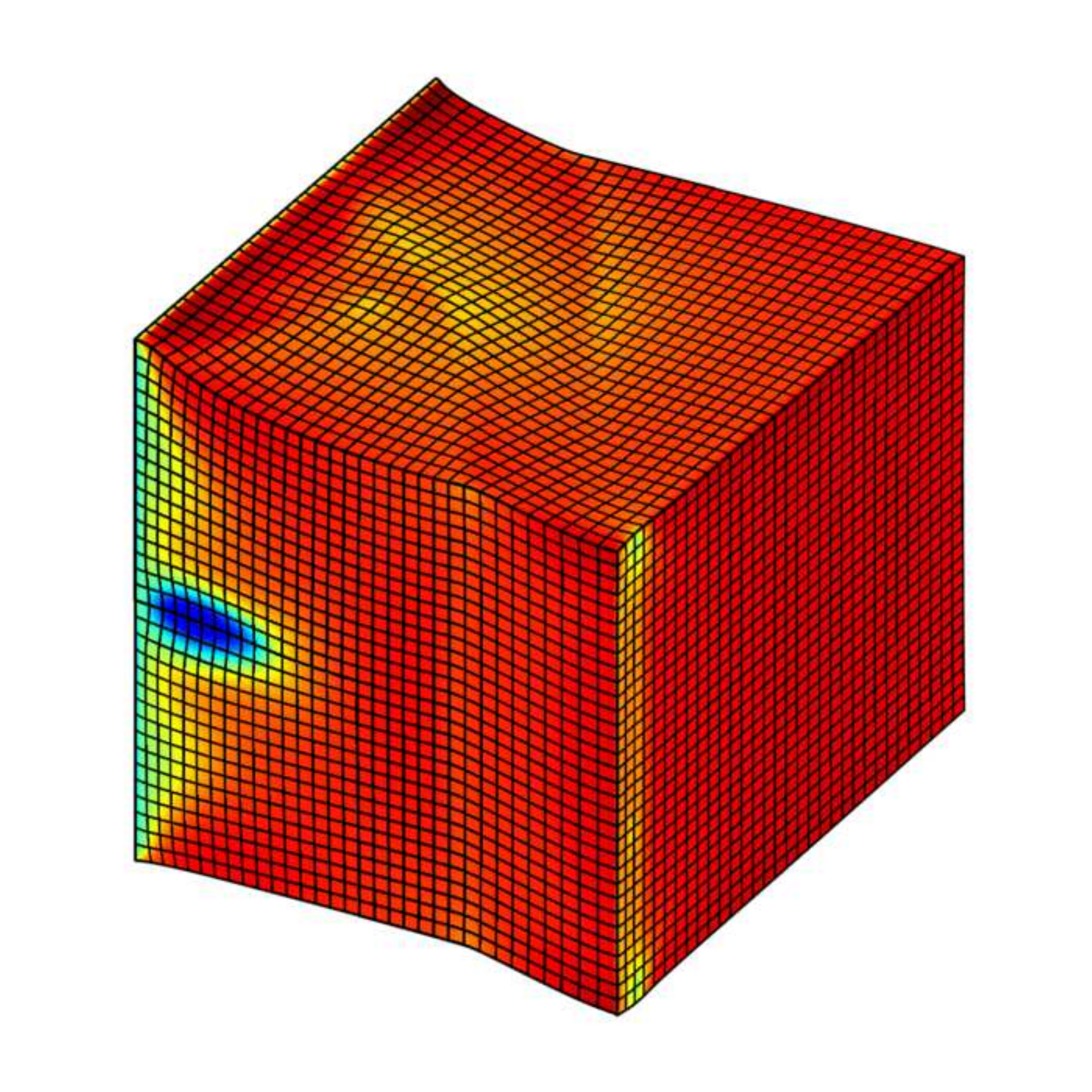} & 
                \includegraphics[scale=0.12]{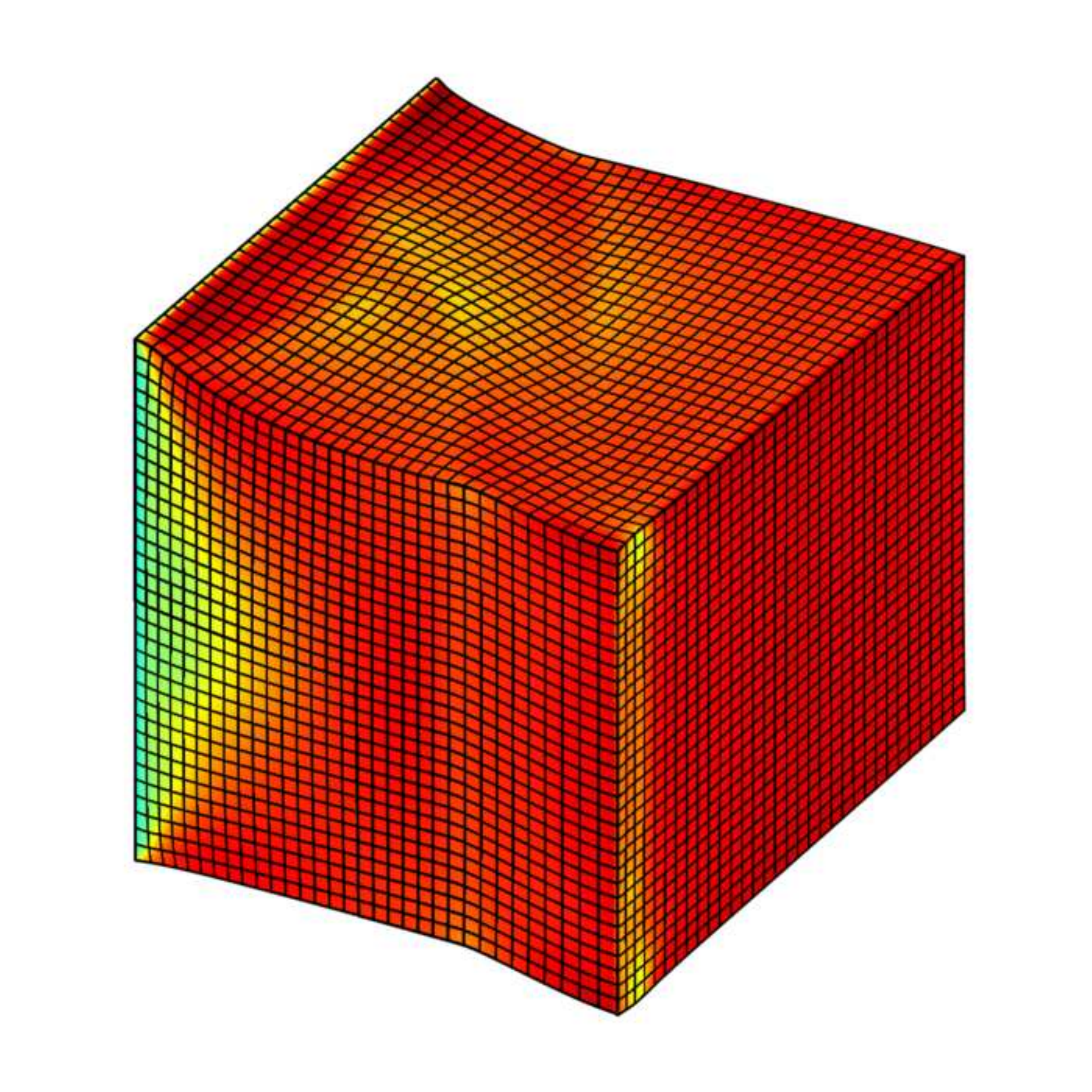} &
                \includegraphics[scale=0.12]{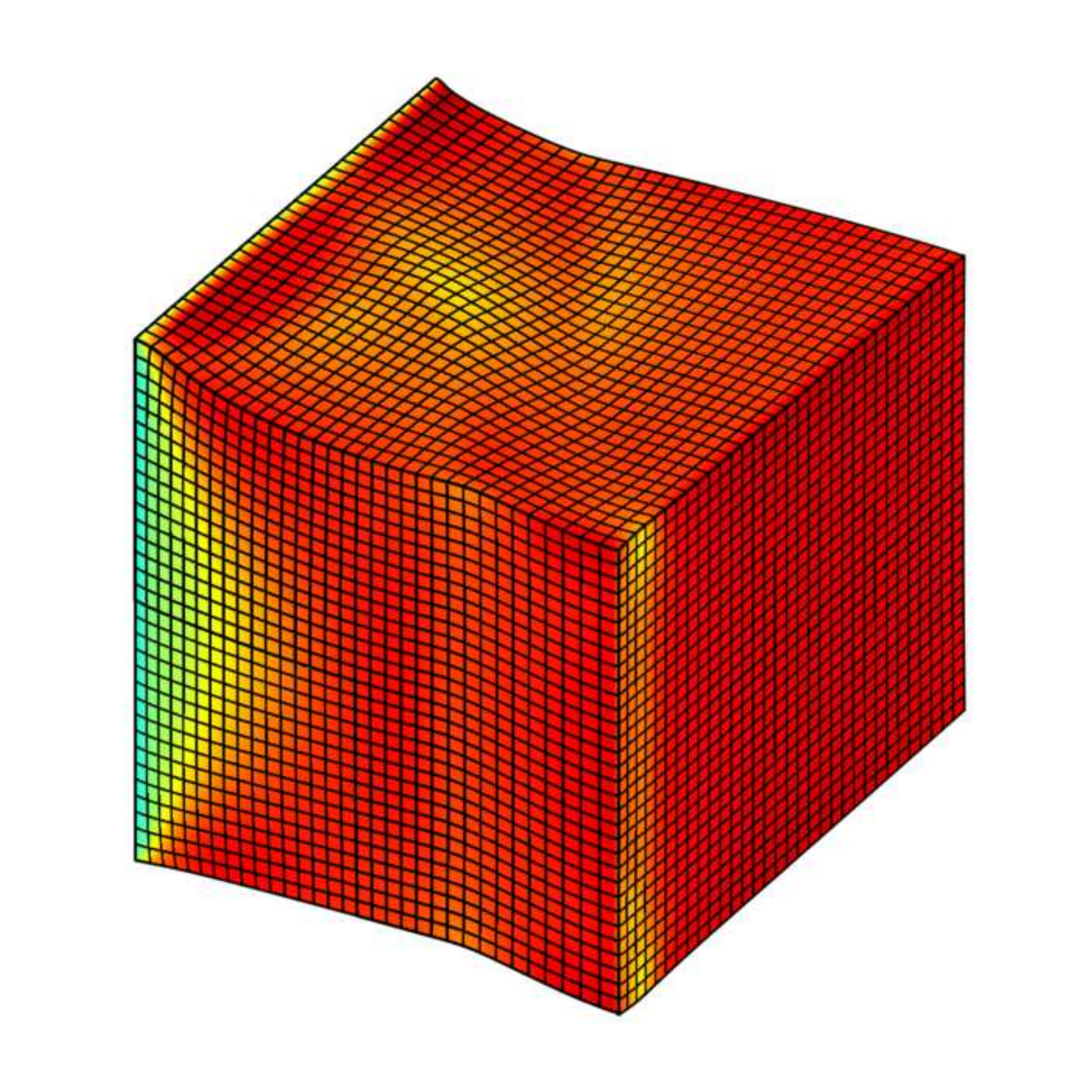} &
                \includegraphics[scale=0.12]{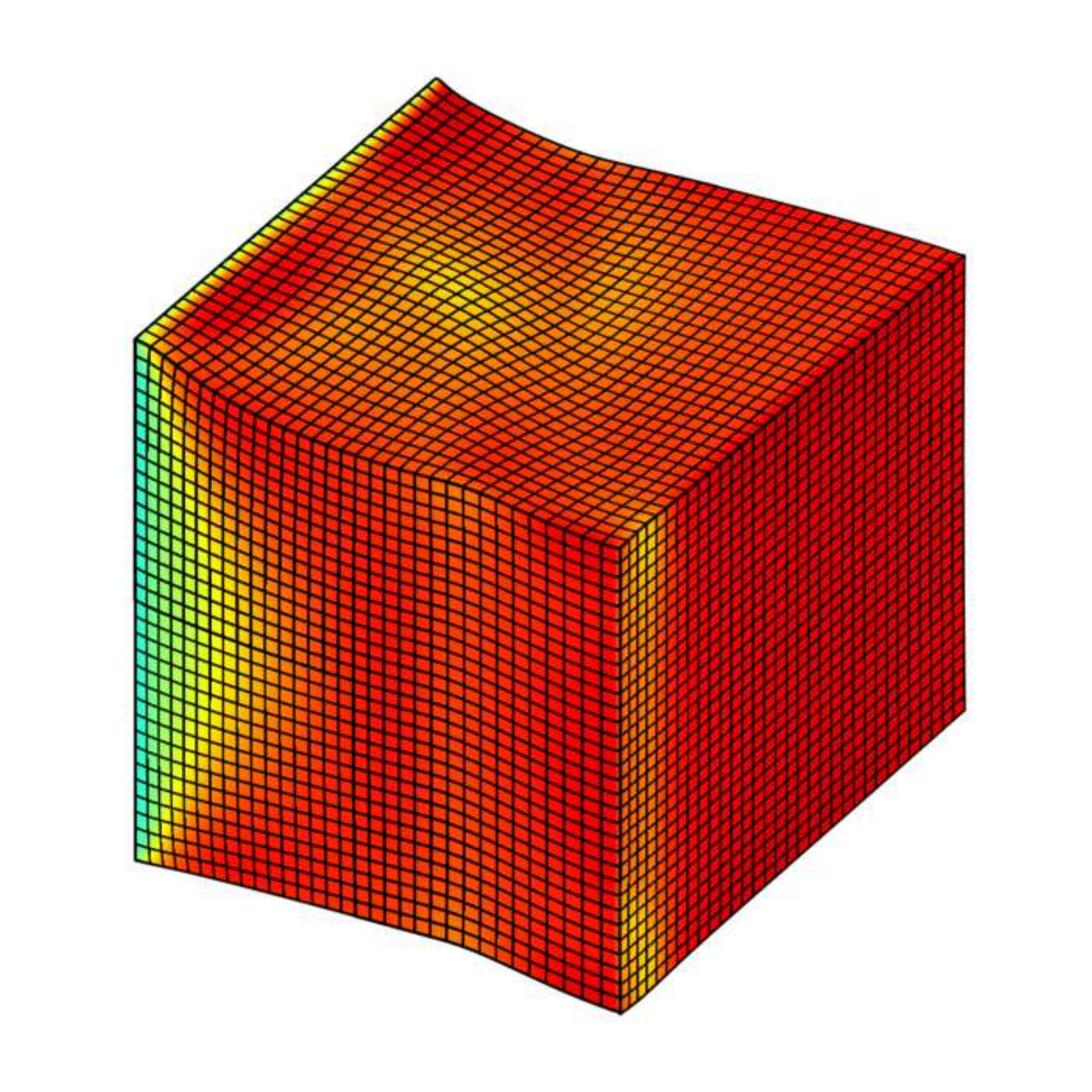} &
                \includegraphics[scale=0.12]{plot/plot_xyz.pdf}
            \end{tabular}  \\
            \parbox[t]{0.5cm}{ $A$ } &
            \begin{tabular}{p{2.8cm}p{2.8cm}p{2.8cm}p{2.8cm}p{2.8cm}p{2.5cm}}
                &
                \includegraphics[scale=0.12]{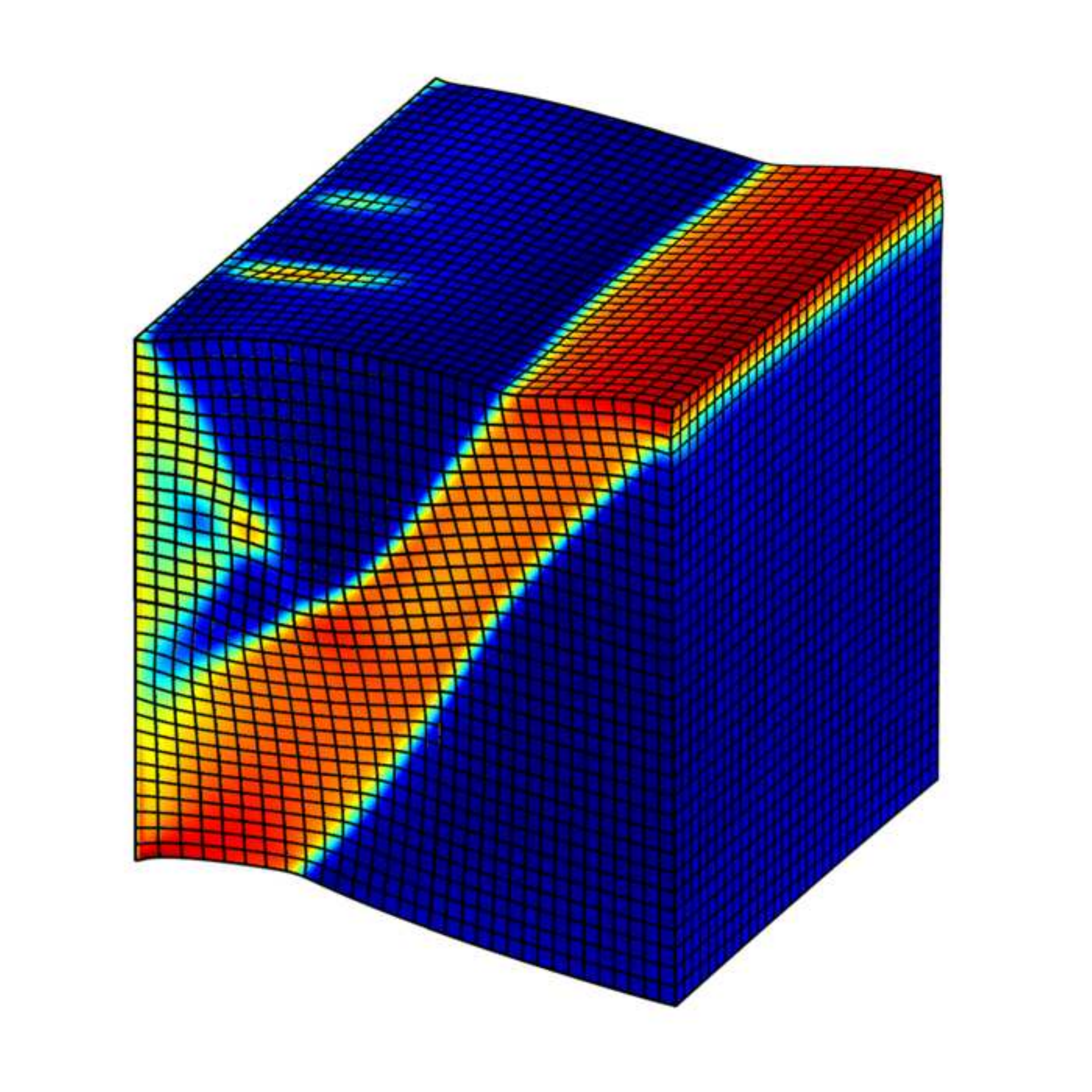} & 
                \includegraphics[scale=0.12]{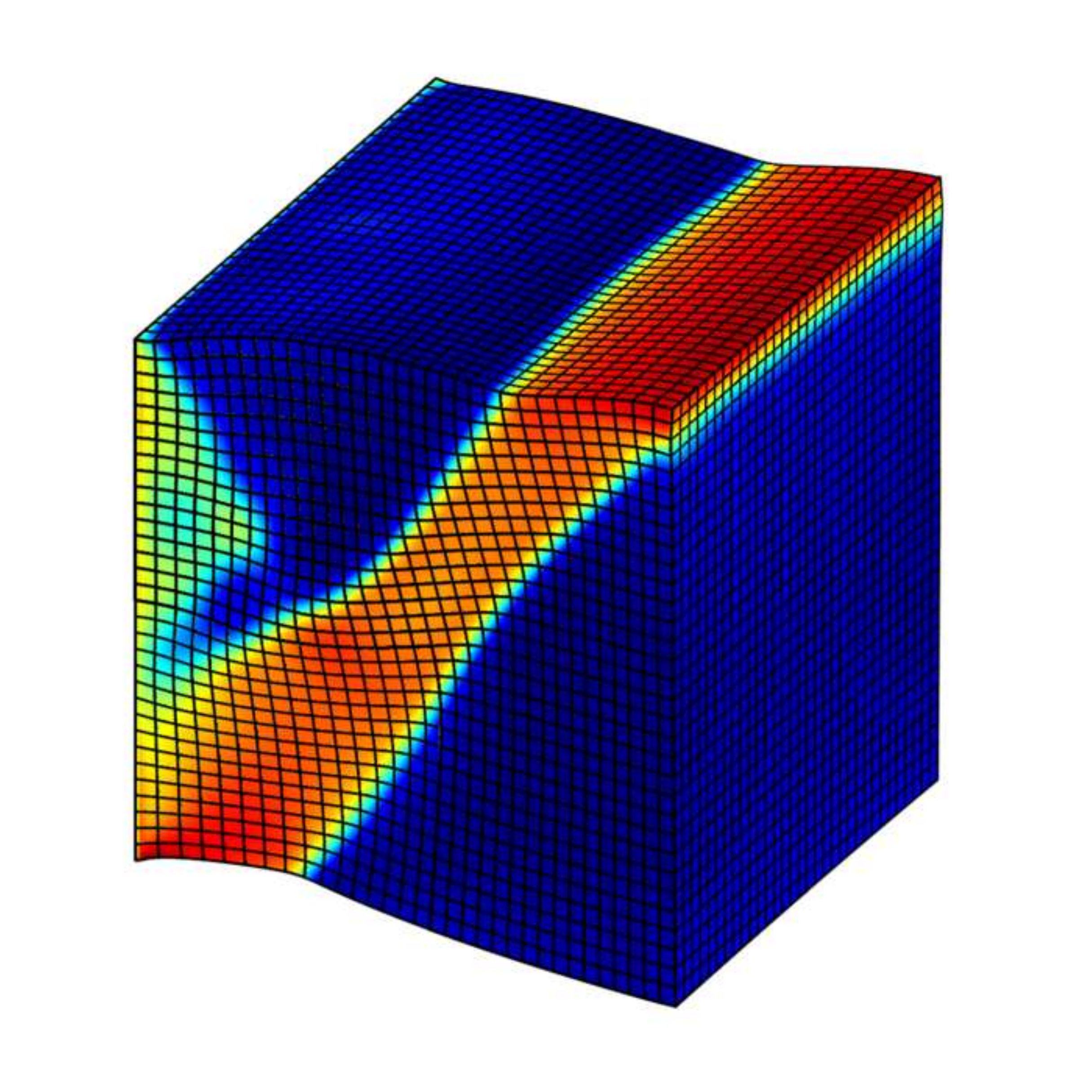} &
                \includegraphics[scale=0.12]{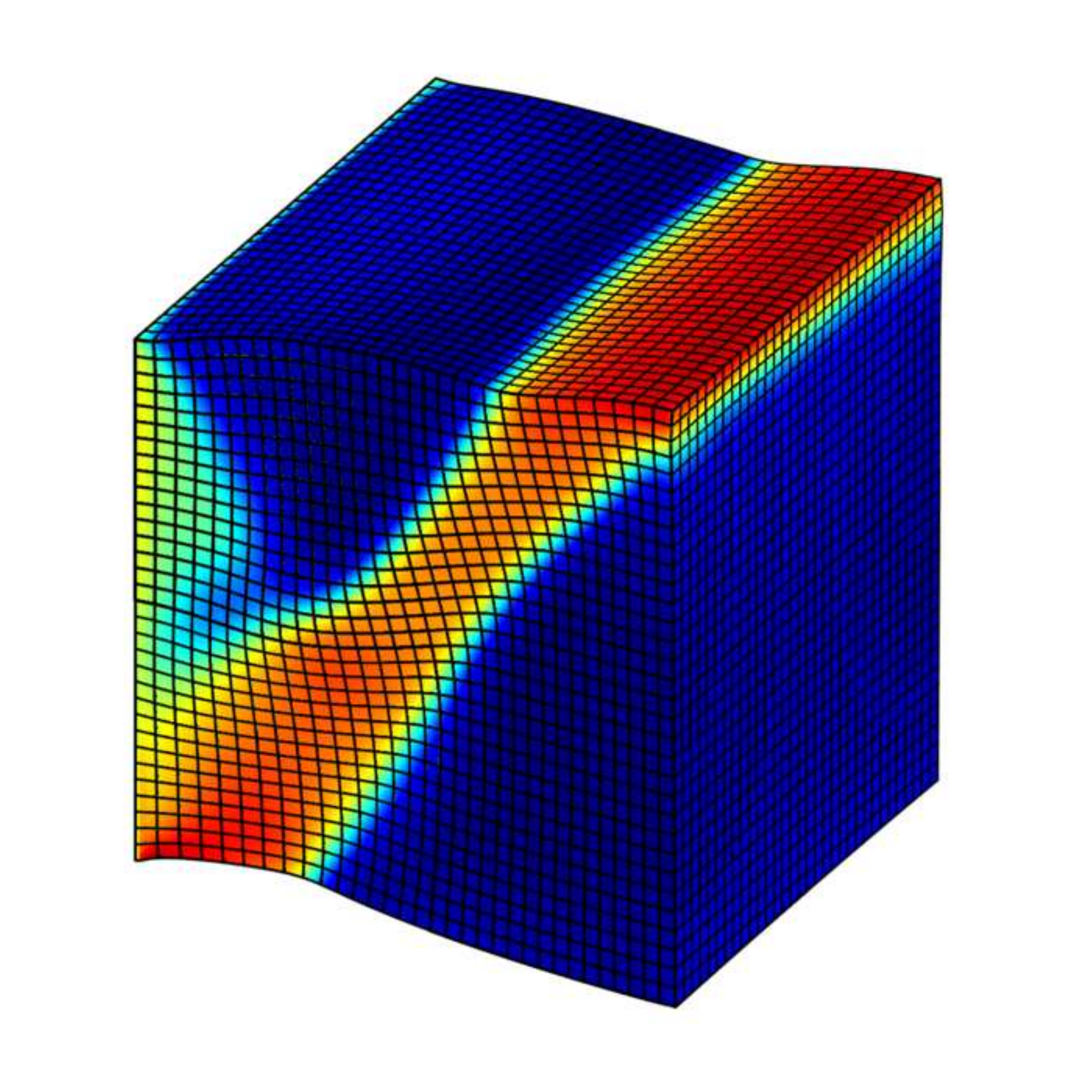} &
                \includegraphics[scale=0.12]{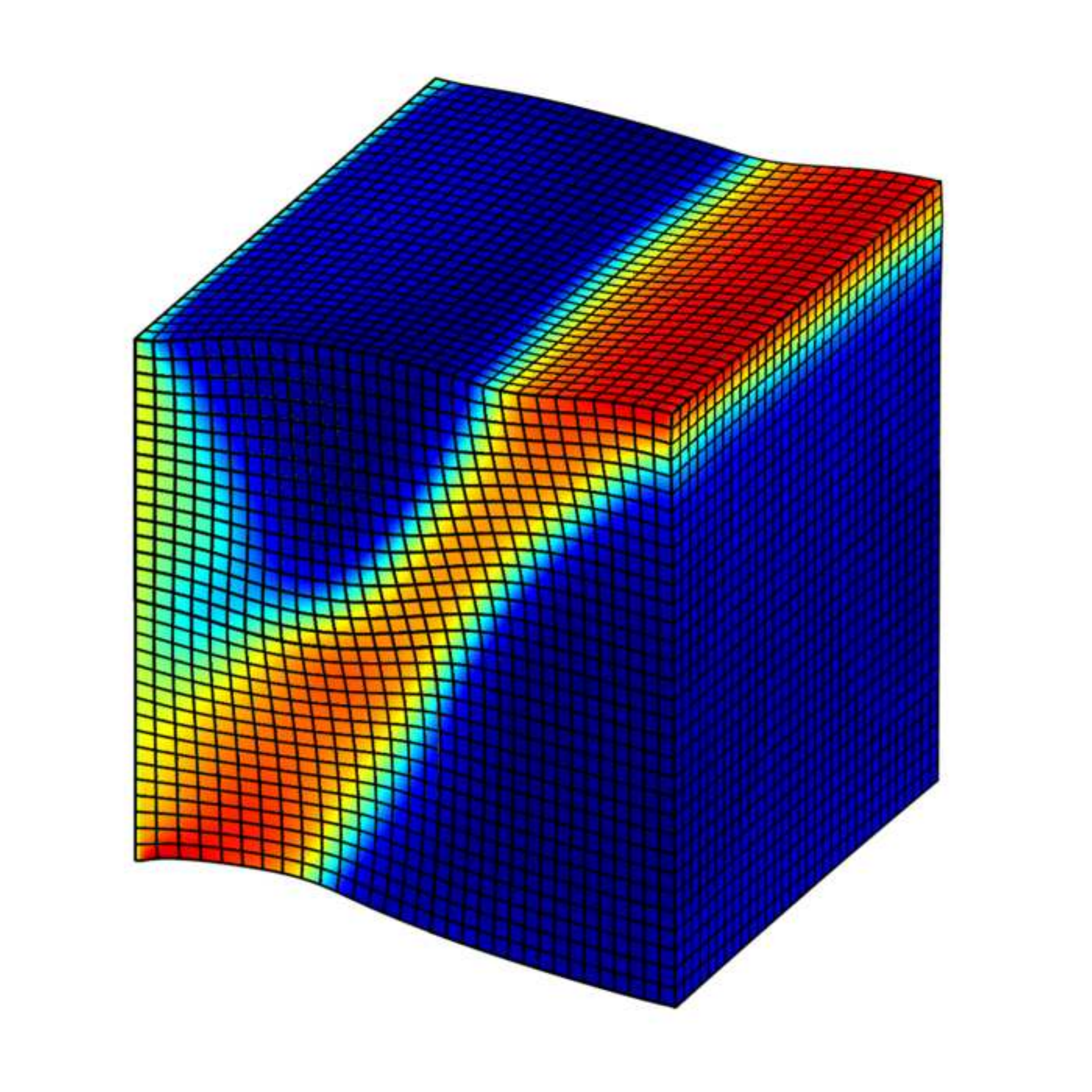} &
                \includegraphics[scale=0.12]{plot/plot_colorbar.pdf}
            \end{tabular}
        \end{tabular}
    \end{center}
    \caption{Field values of $e_3$ for branches A - E on deformed configurations for selected values of $l$.  
    Solutions on the $128^3$ mesh have been overlaid with a $32^3$ plotting mesh.  }
    \label{Fi:e2}
\end{figure}

We solved the boundary value problem on the $64^3$ mesh, obtained five different branches, denoted by A - E, and computed the total free energy for these solutions over ranges of values of $l$; these computed values appear as solid curves in Fig. \ref{Fi:PI-l}.  
At selected values of $l$, $l=0.0625$, $0.0750$, $0.1000$, $0.1500$, and/or $0.2000$, solutions were refined on a $128^3$ mesh and a $256^3$ mesh, and are plotted in Fig. \ref{Fi:PI-l} by red squares and black '+' signs, respectively.  
Reparameterized strains, $e_2$ and $e_3$ \eqref{E:strains1-3}, obtained on the $128^3$ mesh are plotted in Figs. \ref{Fi:e1} and \ref{Fi:e2} for each branch at selected values of $l$.  
In addition, values of $e_2$ and $e_3$ were computed at $33^3$ uniformly spaced points in the body and were plotted on the $e_2-e_3$ space in Fig. \ref{Fi:e1-e2}, superposed on the three-well diagram presented in Fig. \ref{Fi:3well_diagram}.  The strain states corresponding to tetragonal variants 1, 2, and 3 appear in the orange, green, and brown wells, respectively.  
Fig. \ref{Fi:phase} shows three-variant plots that delineate sub-domains of the variants; material points that lie in the interfaces between any two variants are colored in dark gray.  Note that while branches A and B very nearly overlap in the $\Pi-l$ space of Fig.  \ref{Fi:PI-l}, the strains, microstructures and energy landscapes of these branches are actually vastly different as seen in Figs. \ref{Fi:e1}-\ref{Fi:phase}.

The curves in Fig. \ref{Fi:PI-l} were obtained using a simple branch-tracking technique employed in the one-dimensional example in Sec. \ref{S:One-dimensional_example} with slight modification.  
For the purpose of demonstration, we aimed to obtain \emph{moderately low-energy} microstructures for which tetragonal variants are well developed in the body for the free energy coefficient $B_5=180$ and in the vicinity of $l=0.10$.  
Since such low-energy microstructures are under large strains, direct computation at $(B_5,l)=(180,0.10)$ is not practical as we in general do not have \emph{good} initial guesses with which we can obtain converged solutions in three dimensions.  
We thus first computed high-energy solutions under small strains at relatively large values of $(B_5,l)$, and then employed branch-trackings in decreasing $B_5$- and $l$-directions down to $(B_5,l)=(180,0.10)$, which would develop lower-energy microstructures under larger strains.  
Small-strain, high-energy solutions were obtained using either the homogeneous initial guess or random initial guesses of small magnitude.  
The former was used for branch E and the latter was used for branches A-D.

To eventually obtain solutions well resolved on the $64^3$ mesh, we first computed the initial high-energy solutions on an $8^3$ mesh or on a $16^3$ mesh and successively refined them before branch-trackings.
Finer initial meshes, say $32^3$, would produce finer and more interesting microstructures at the end of the iterative process, but at the expense of greater computational complexity.  Our goal being demonstration of the series of techniques, however, we do not pursue this approach here.  
For branch E, for instance, a solution was first computed with $(B_5,l)=(500,0.54)$ on the $16^3$ mesh using the homogeneous initial guess and successively refined onto the $64^3$ mesh, which was then followed by branch-trackings as $(B_5,l)=(500,0.54)\rightarrow(500,0.15)\rightarrow(180,0.15)\rightarrow(180,0.10)$.  
Branch D, on the other hand, was first computed with $(B_5,l)=(500,0.50)$ using a random initial guess of $\vert\bs{u}^h\vert \sim O(10^{-2})$ on the $8^3$ mesh, successively refined onto the $64^3$ mesh, and subject to branch-trackings as $(B_5,l)=(500,0.54)\rightarrow(500,0.27)\rightarrow(180,0.27)\rightarrow(180,0.10)$.

\subsection{Numerical convergence}
For added confidence that the branches obtained in this study are not numerical artifacts, we studied convergence of these solutions with mesh refinement; 
at selected values: $l=0.0625$, $0.0750$, $0.1000$, $0.1500$, and/or $0.2000$, solutions were refined on a $128^3$ mesh and a $256^3$ mesh, and corresponding energy values were plotted in Fig. \ref{Fi:PI-l} as red squares and black '+' signs, respectively.  
Fig. \ref{Fi:PI-l_orig} implies that, in general, solutions computed on the $64^3$ mesh are energetically well converged and Fig. \ref{Fi:PI-l_zoom} indicates that solutions are better resolved for larger values of $l$, where the interface thickness is wider and microstructures are coarser.
Refinement of branches A and B was more challenging especially for larger values of $l$, and refinement on the $256^3$ mesh was only performed for $l\leq 0.1000$ for these branches.  
These branches also experience slightly larger deviation when refined as seen in Fig. \ref{Fi:PI-l_zoom}.  
We did not pursue these problems further as these branches are of little practical interest.  

\begin{figure}
    \begin{center}
        \begin{tabular}{rp{16.6cm}}
            \parbox[t]{0.5cm}{ }&
            \begin{tabular}{p{4.7cm}p{4.7cm}p{4.7cm}p{2.5cm}}
                \hspace{2.3cm}$64^3$&
                \hspace{2.3cm}$128^3$&
                \hspace{2.3cm}$256^3$&
                \vspace{0.5\baselineskip}
            \end{tabular}  \\
            \parbox[t]{0.5cm}{ $D$ } &
            \begin{tabular}{p{4.7cm}p{4.7cm}p{4.7cm}p{2.5cm}}
                \includegraphics[scale=0.20]{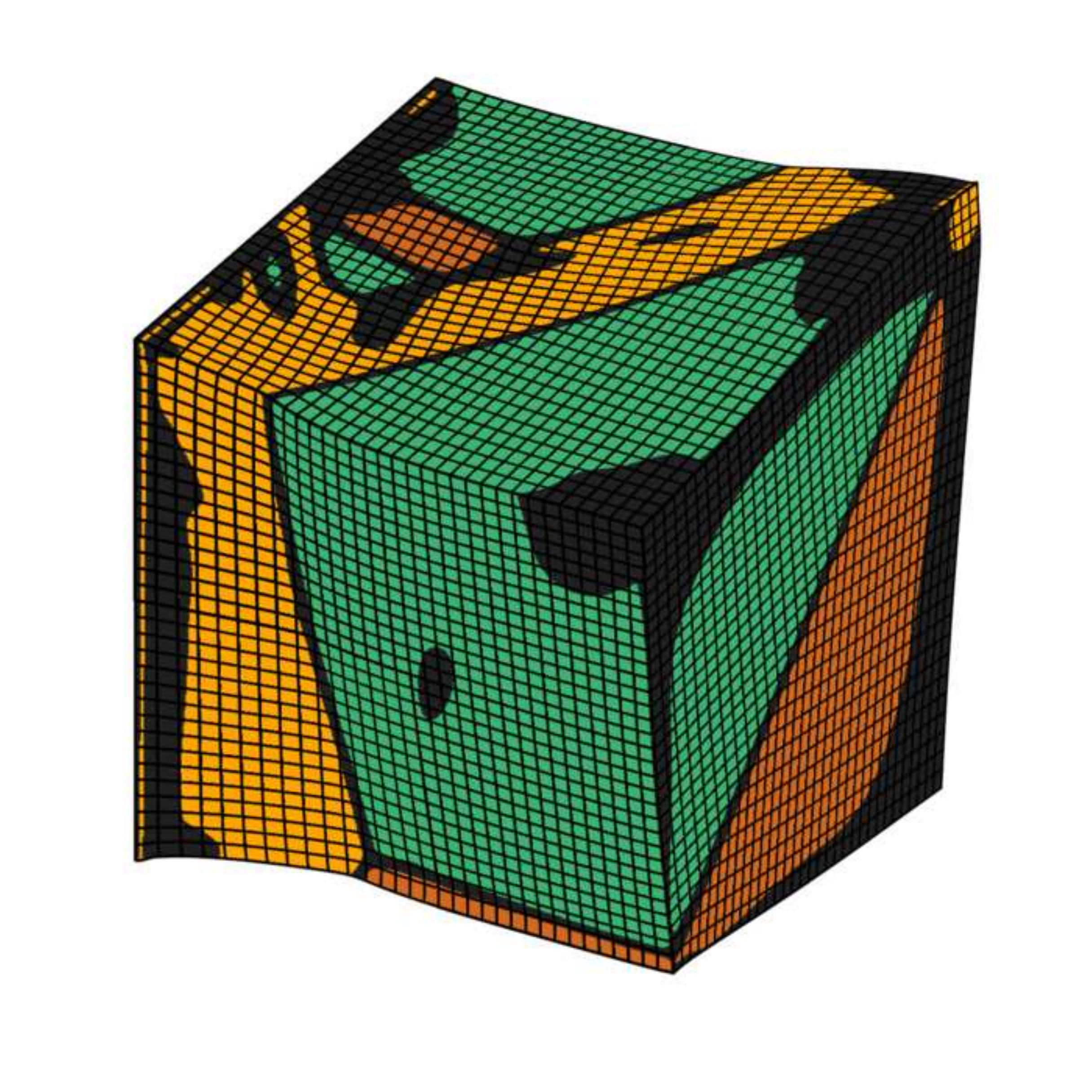} &
                \includegraphics[scale=0.20]{plot/phase_cube_7_2_21_0_25000000_0_50000000_180_00000000_3_25000000_0_06250000.pdf} &
                \includegraphics[scale=0.20]{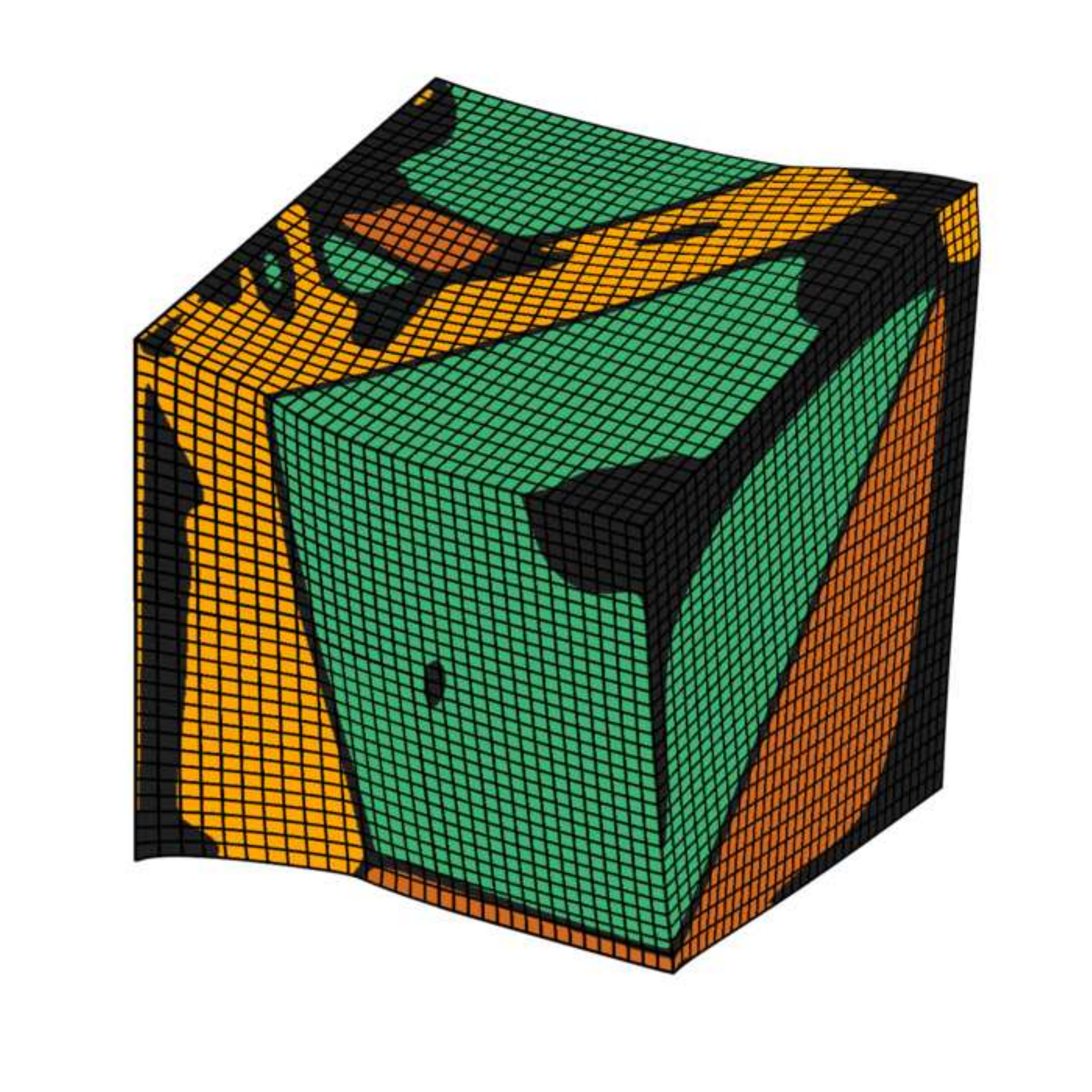} &
                
            \end{tabular}  \\
            \parbox[t]{0.5cm}{ $C$ } &
            \begin{tabular}{p{4.7cm}p{4.7cm}p{4.7cm}p{2.5cm}}
                \includegraphics[scale=0.20]{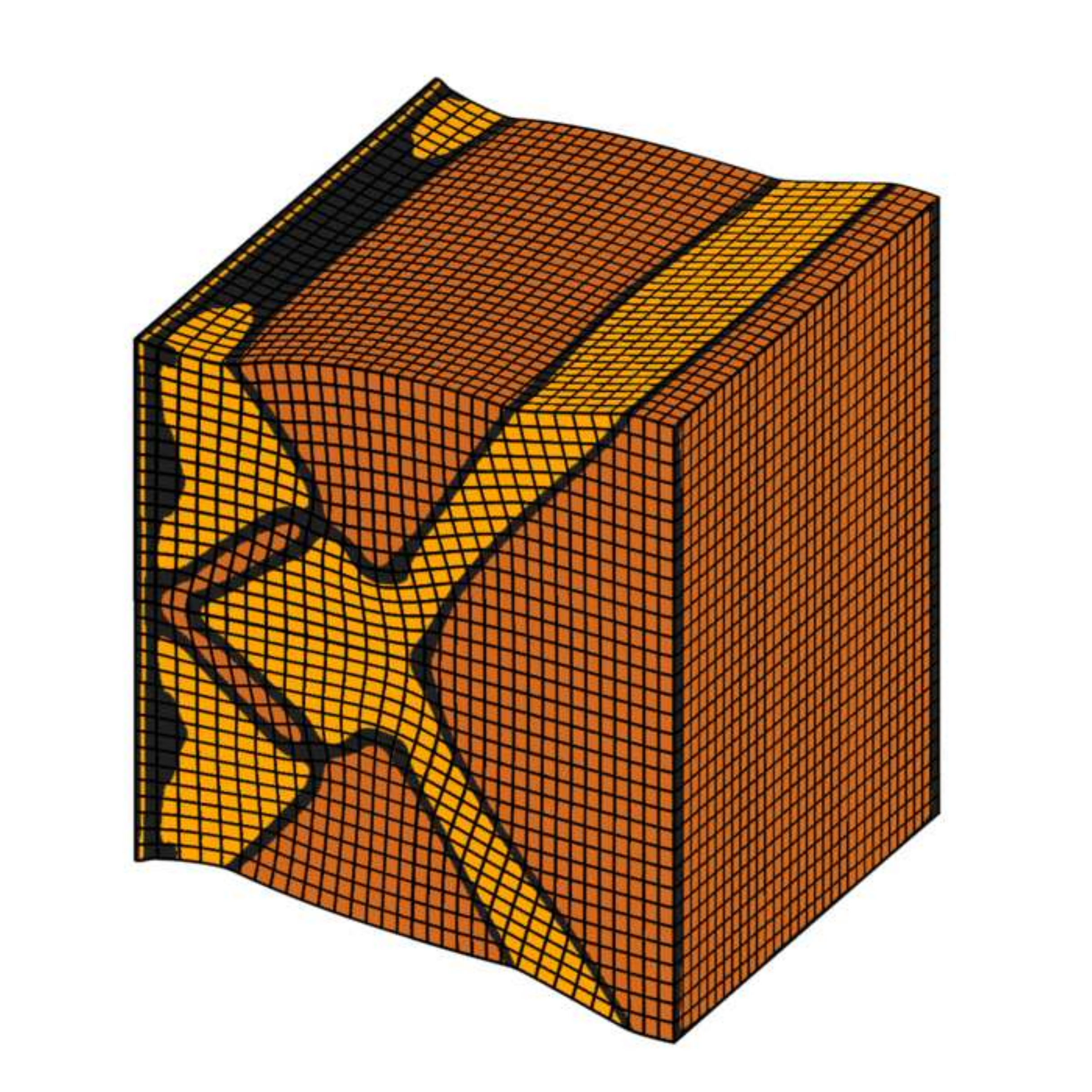} &
                \includegraphics[scale=0.20]{plot/phase_cube_7_2_33_0_25000000_0_50000000_180_00000000_3_25000000_0_06250000.pdf} &
                \includegraphics[scale=0.20]{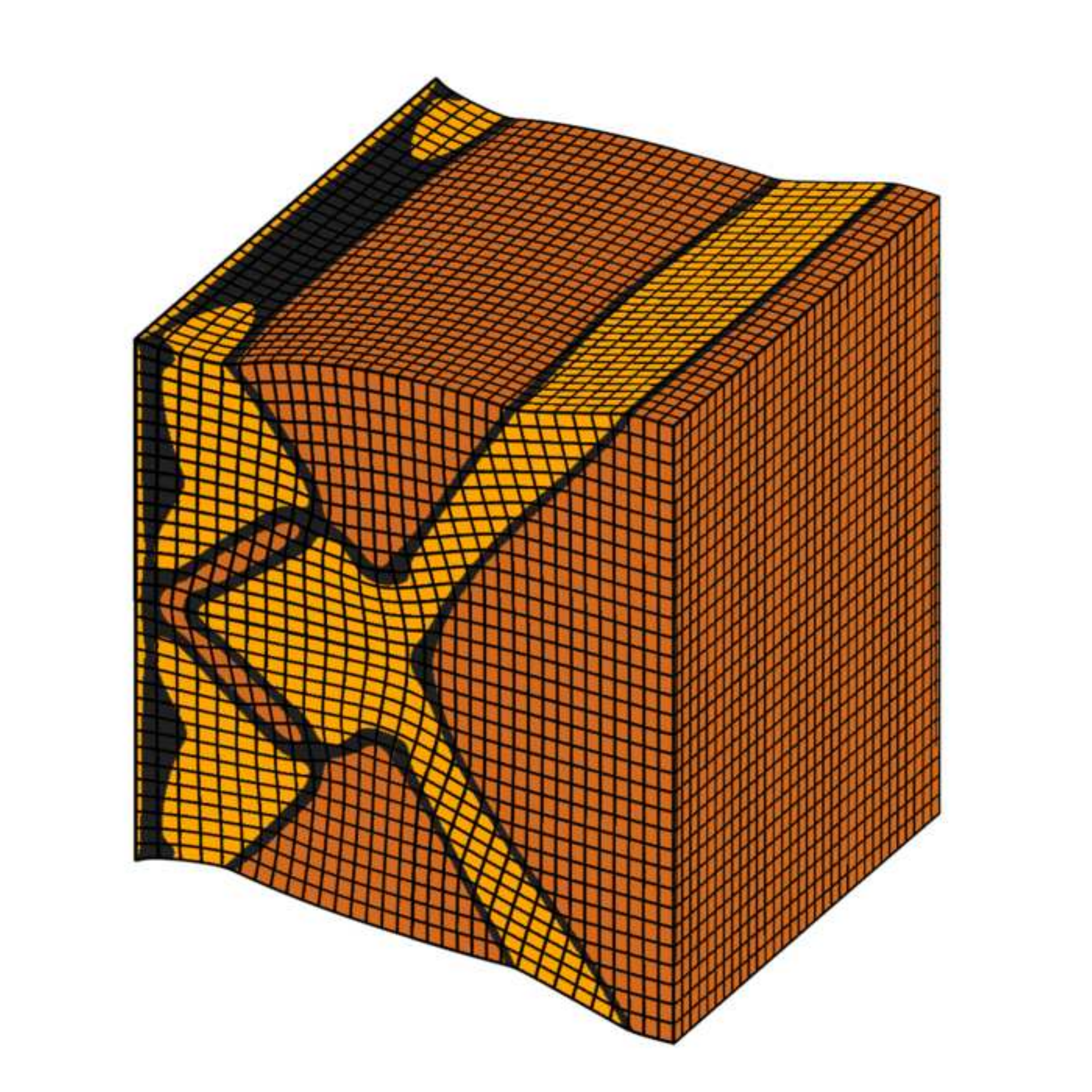} &
                \includegraphics[scale=0.15]{plot/plot_xyz.pdf}
            \end{tabular}  \\
            \parbox[t]{0.5cm}{ $B$ } &
            \begin{tabular}{p{4.7cm}p{4.7cm}p{4.7cm}p{2.5cm}}
                \includegraphics[scale=0.20]{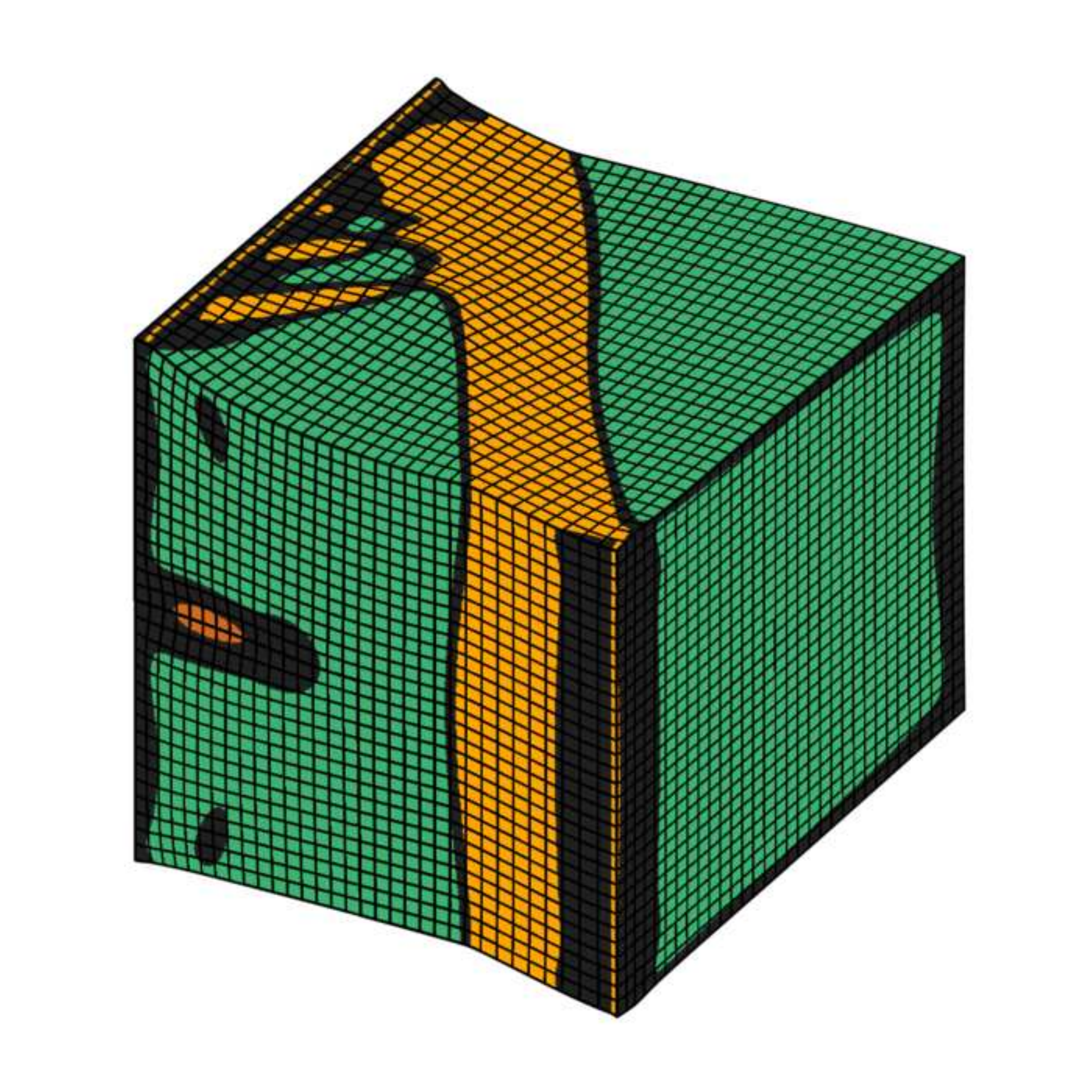} &
                \includegraphics[scale=0.20]{plot/phase_cube_7_2_34_0_25000000_0_50000000_180_00000000_3_25000000_0_06250000.pdf} &
                \includegraphics[scale=0.20]{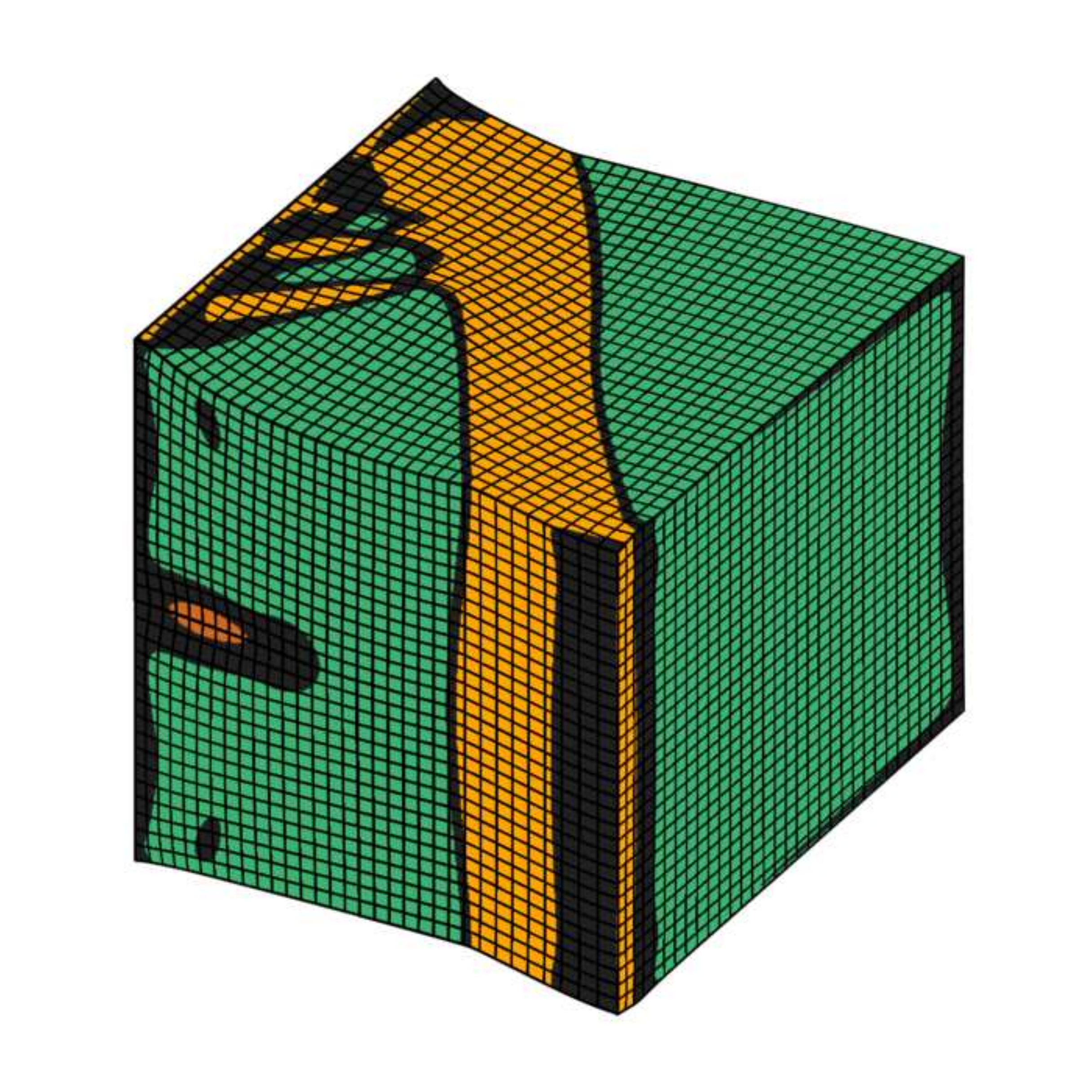} &
                \includegraphics[scale=0.15]{plot/plot_colorbar_phase.pdf}
            \end{tabular}
        \end{tabular}
    \end{center}
    \caption{Convergence study of the distribution of tetragonal variants in the microstructure for branches B, C, and D.  Solutions were computed at $l=0.0625$ on the $64^3$, $128^3$, and $256^3$ meshes.  All plots were overlaid with a $32^3$ plotting mesh.  }
    \label{Fi:convergence}
\end{figure}

Fig.\ref{Fi:convergence}, on the other hand, shows the distribution of the three variants making up the microstructure at $l=0.0625$ for the $64^3$, $128^3$, and $256^3$ meshes.  
Attention is drawn to the near complete convergence of solutions on the $64^3$ mesh in the sense of microstructure.

\subsection{Stability/metastability}
The numerical study of the boundary value problem was followed by stability analysis, where positive definiteness of the Hessian derived from the second variation \eqref{E:D2_PI} was numerically checked for each branch.  
Table \ref{Ta:eig} shows the smallest eigenvalues of the Hessians for branches A - E on different refinement levels at selected values of $l$; positive values therefore imply stability/metastability.  
Branch A, though the smallest eigenvalue is negative at $l=0.0750$ probably due to poor resolution, is most likely to be stable, which is consistent with that branch A is most likely to be the lowest-energy branch.  
Branch B, on the other hand, seems to gain stability somewhere between $l=0.1000$ and $l=0.1500$.  
Portions of branch C and branch D would also be good candidates for metastable portions of branches, considering the convergence behavior of the eigenvalues with mesh refinement.  The stability/metastability behavior of the branches with decreasing $l$ points to the existence of as yet undiscovered branches at increasingly finer microstructure. Their resolution is only limited by the numerical expense of ever finer meshes.

\subsection{Parametric dependence of twin microstructures}
\begin{figure}
    \begin{center}
        \begin{tabular}{rp{16.5cm}}
            \parbox[t]{0.5cm}{ }&
            \begin{tabular}{p{2.8cm}p{2.8cm}p{2.8cm}p{2.8cm}p{2.8cm}p{2.5cm}}
                \hspace{0.8cm}$l\!=\!0.0625$ & 
                \hspace{0.8cm}$l\!=\!0.0750$ & 
                \hspace{0.8cm}$l\!=\!0.1000$ &
                \hspace{0.8cm}$l\!=\!0.1500$ &
                \hspace{0.8cm}$l\!=\!0.2000$ &
                \vspace{0.5\baselineskip}
            \end{tabular}  \\
            \parbox[t]{0.5cm}{ $E$ } &
            \begin{tabular}{p{2.8cm}p{2.8cm}p{2.8cm}p{2.8cm}p{2.8cm}p{2.5cm}}
                &
                \includegraphics[scale=0.12]{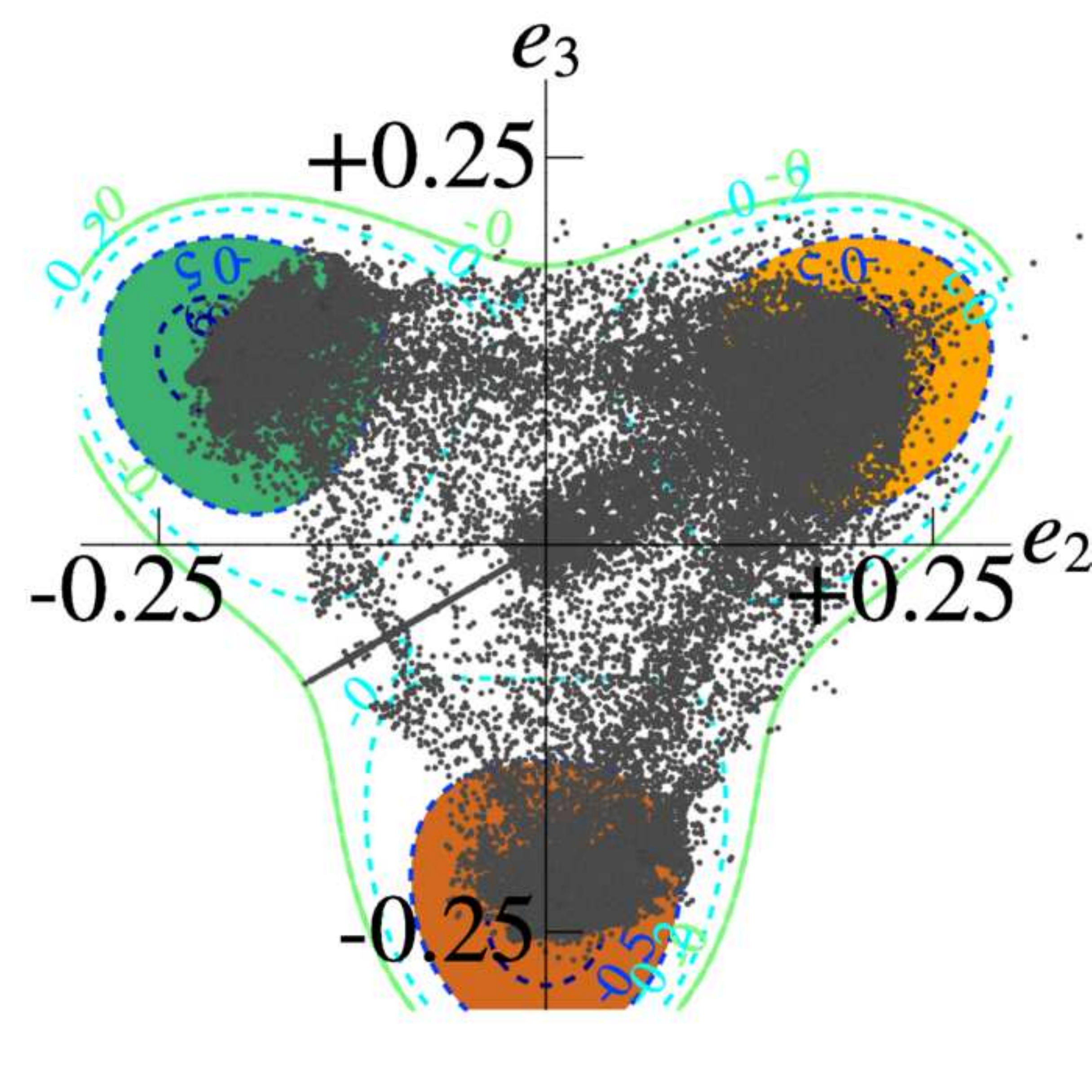} & 
                \includegraphics[scale=0.12]{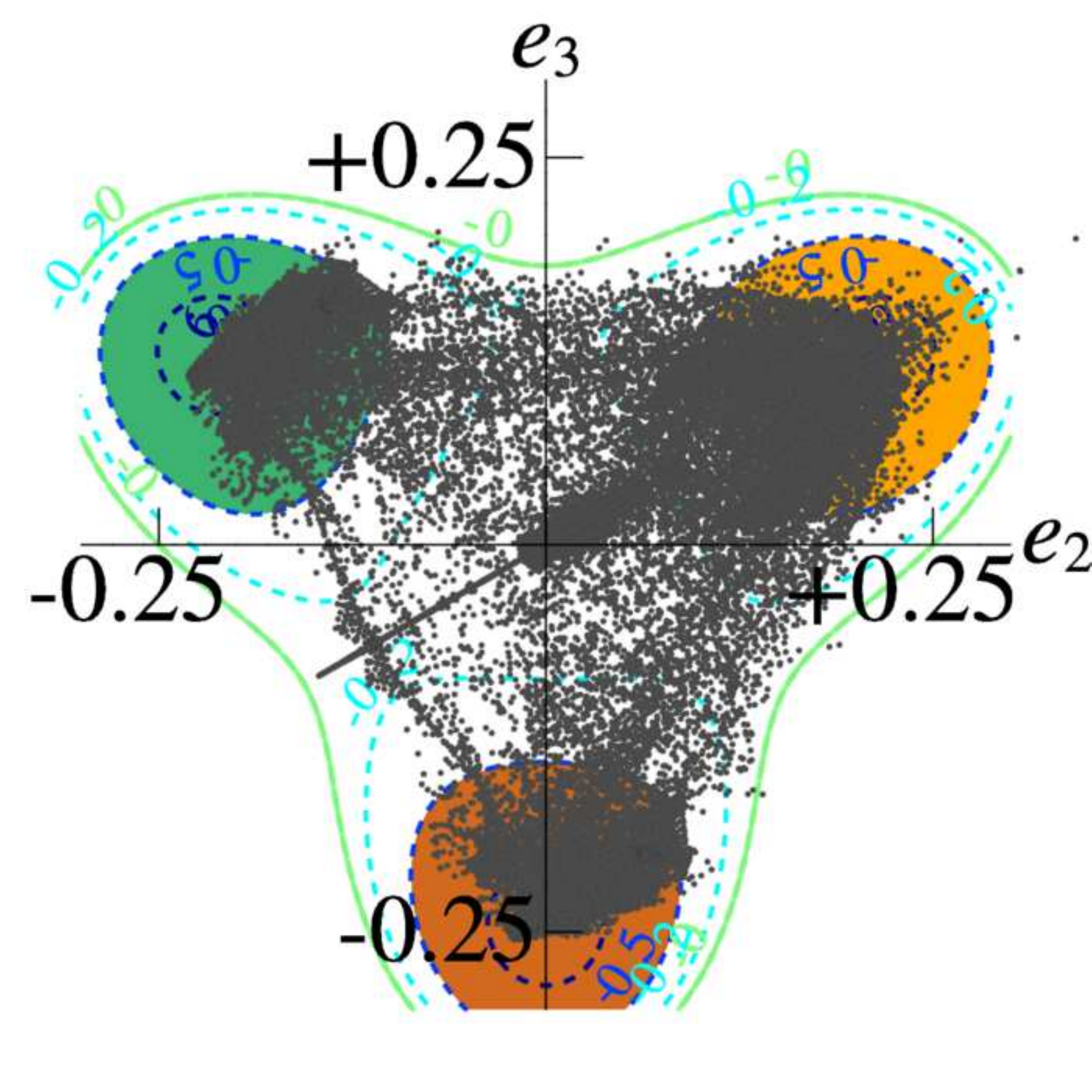} &
                \includegraphics[scale=0.12]{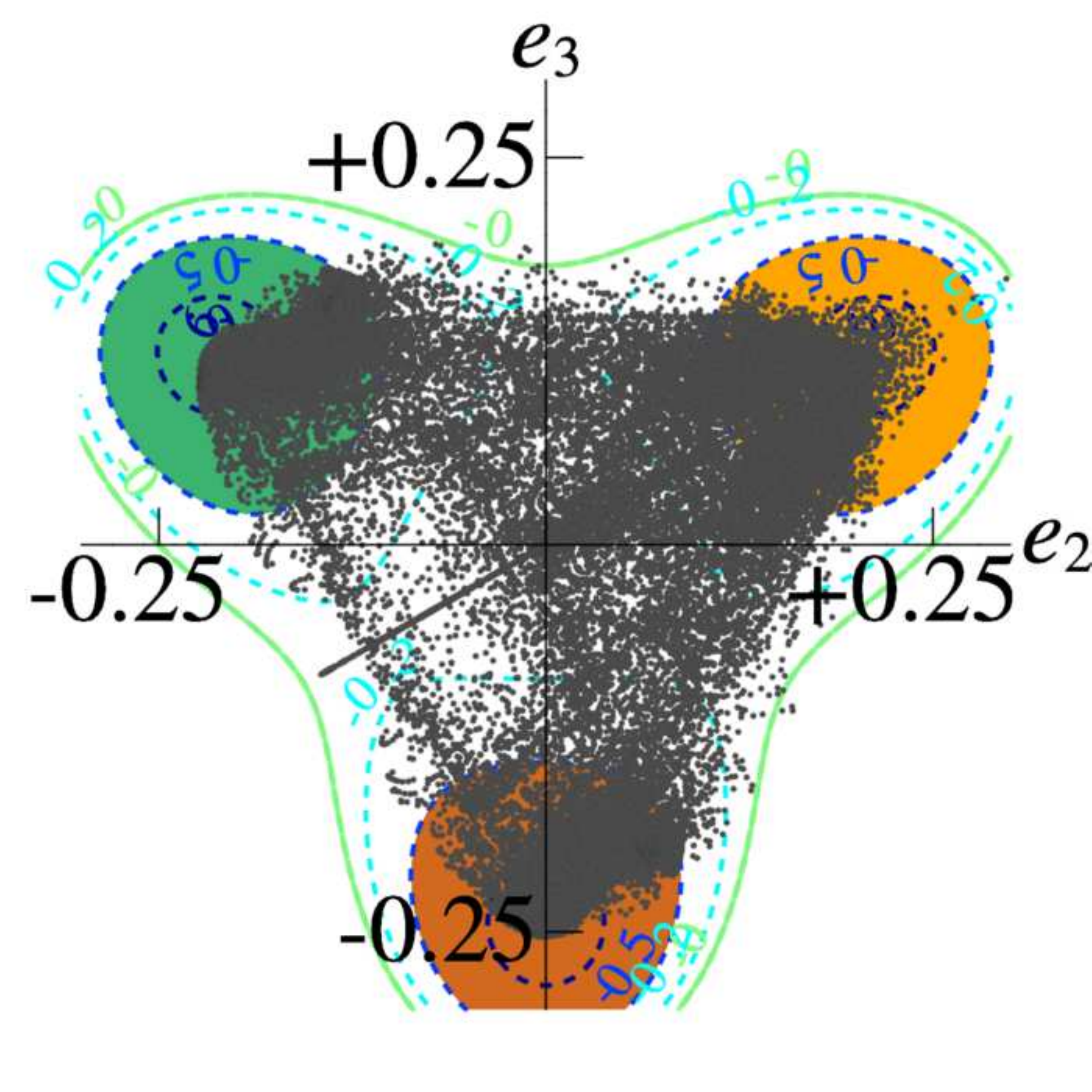} &
                \includegraphics[scale=0.12]{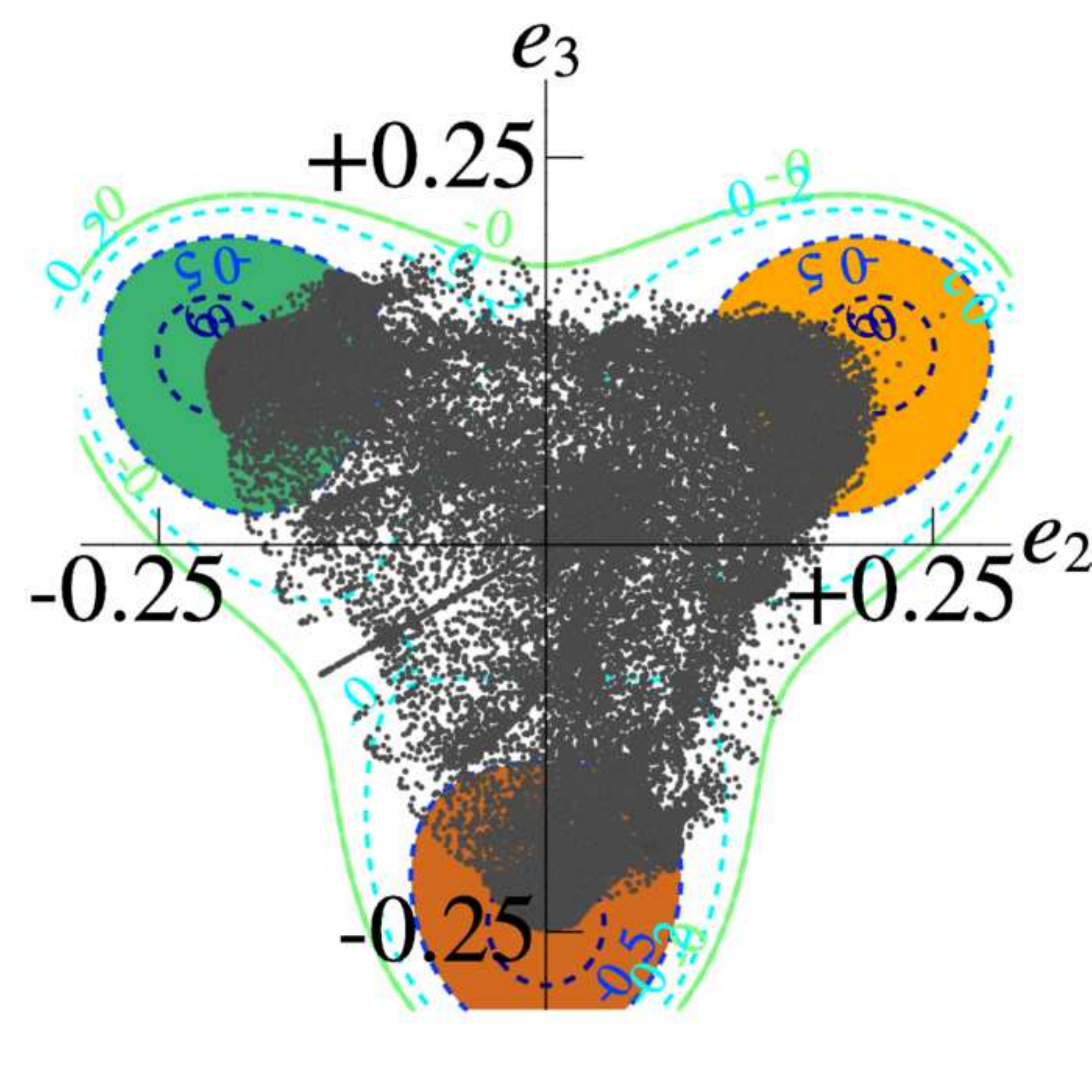} &
            \end{tabular}  \\
            \parbox[t]{0.5cm}{ $D$ } &
            \begin{tabular}{p{2.8cm}p{2.8cm}p{2.8cm}p{2.8cm}p{2.8cm}p{2.5cm}}
                \includegraphics[scale=0.12]{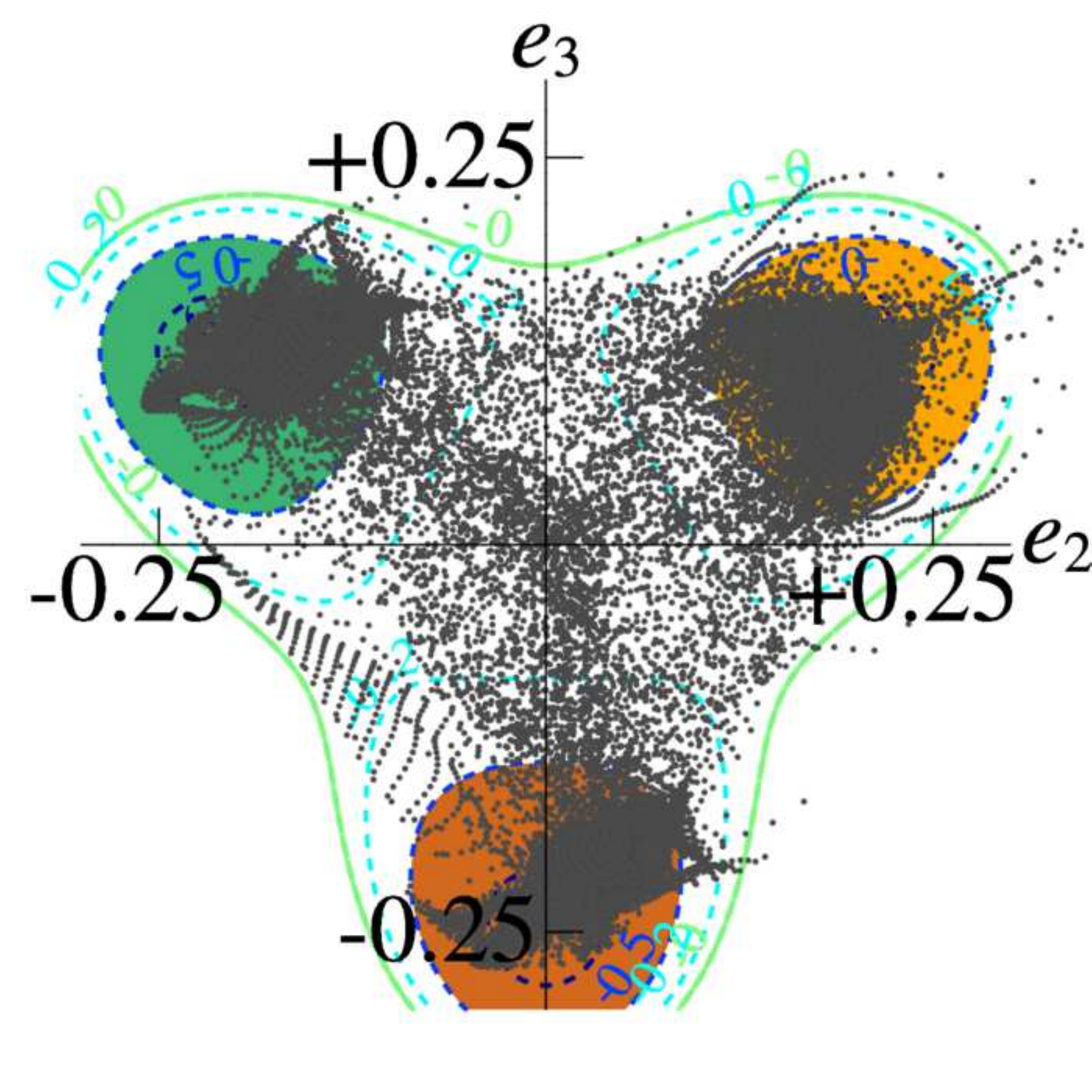} &
                \includegraphics[scale=0.12]{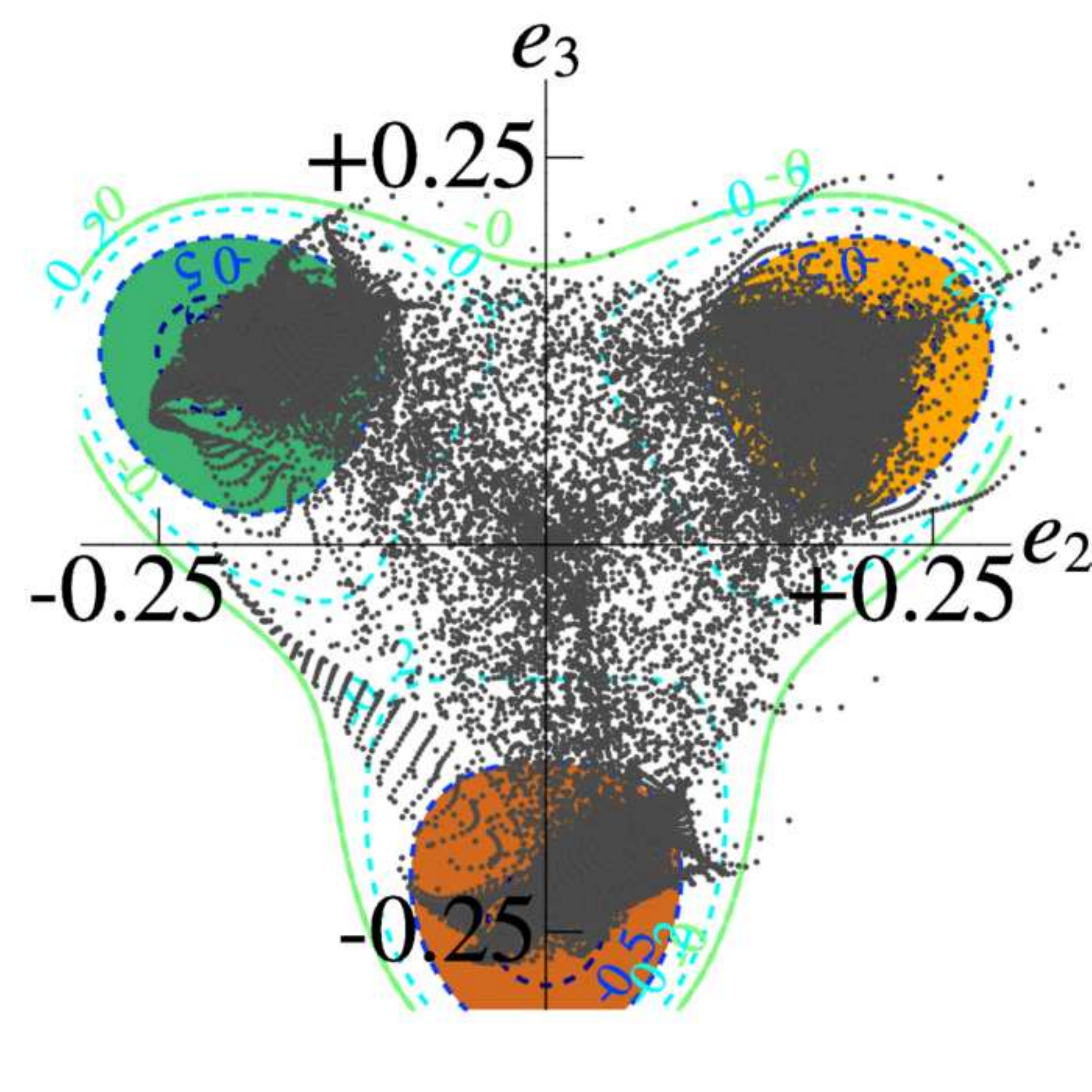} & 
                \includegraphics[scale=0.12]{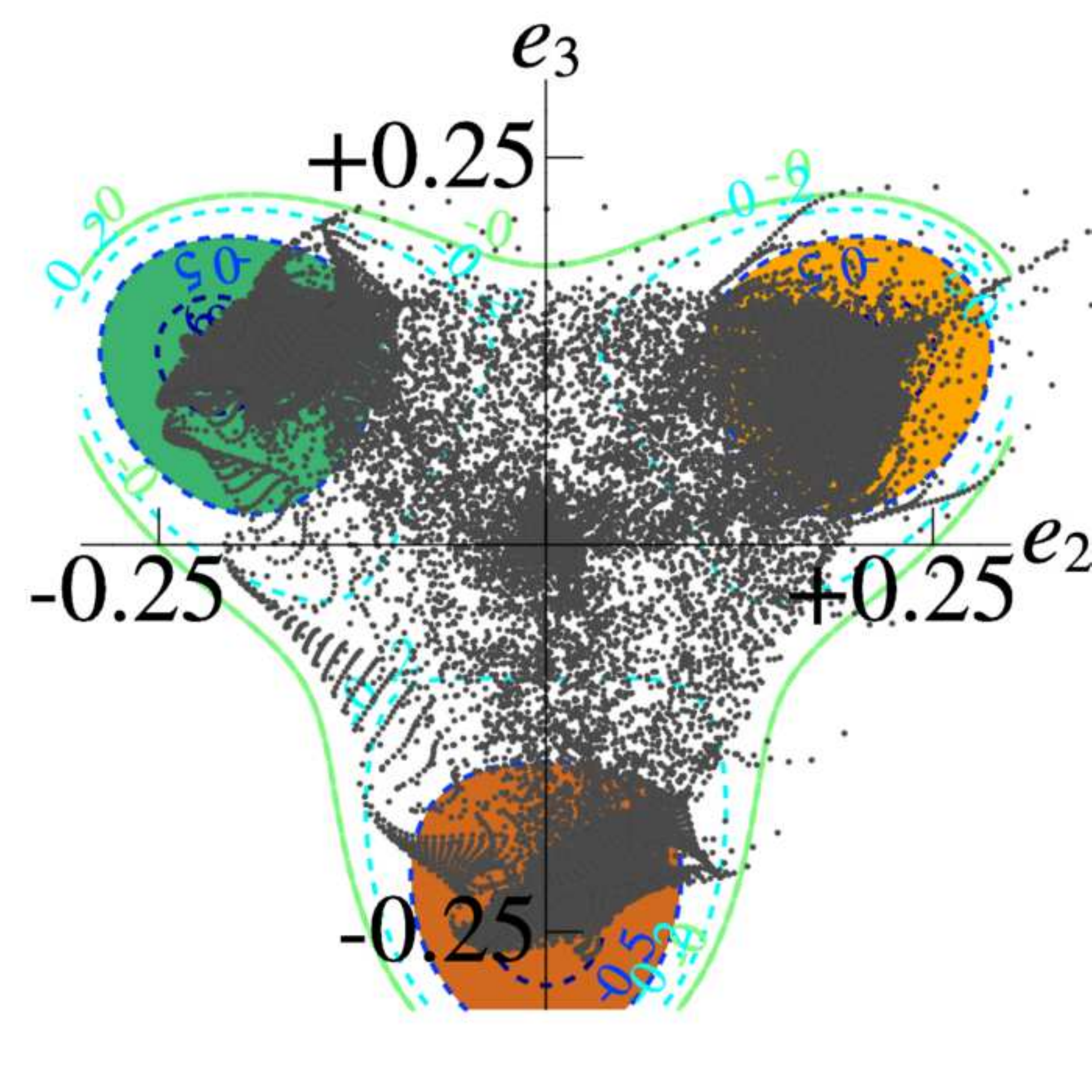} &
                \includegraphics[scale=0.12]{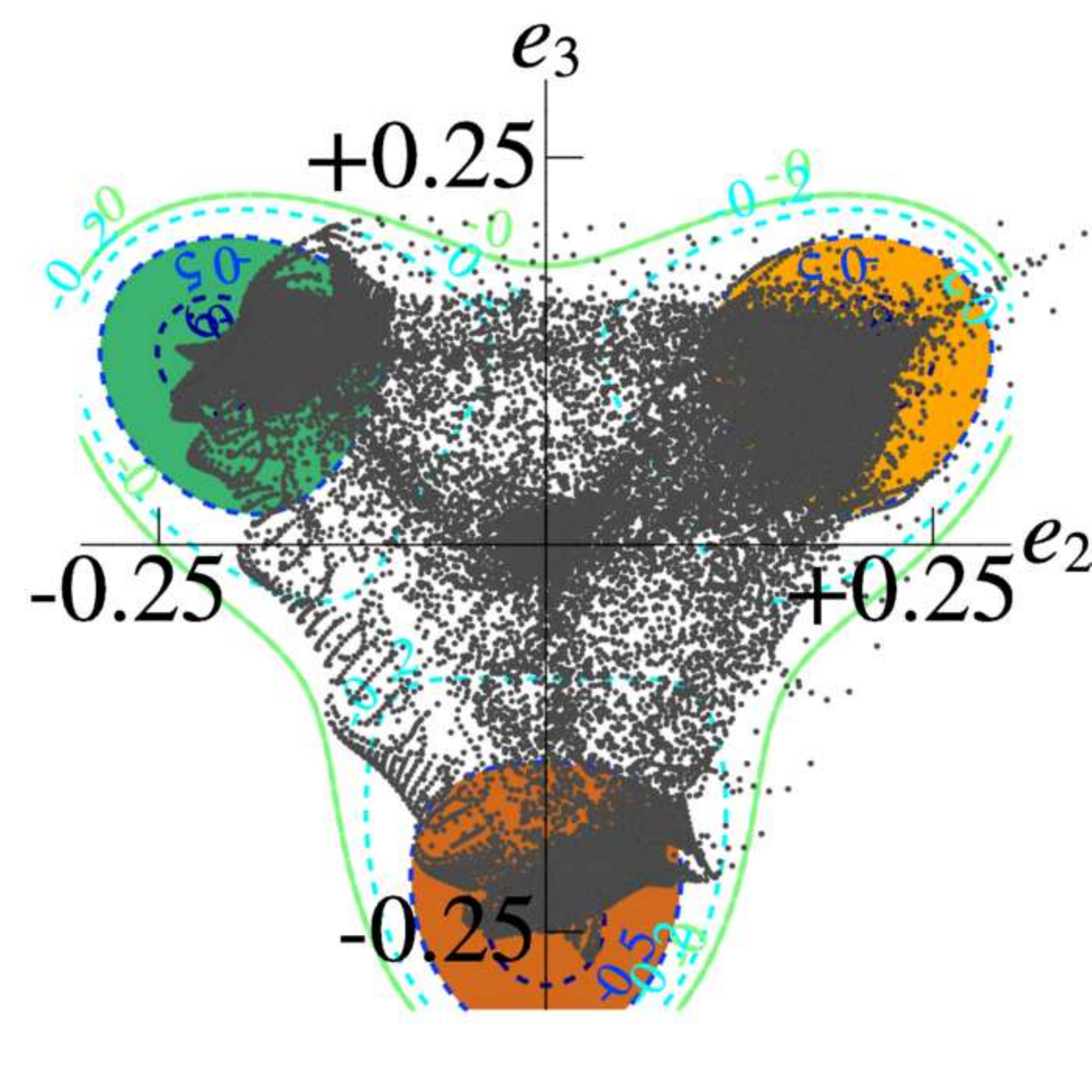} &
                \includegraphics[scale=0.12]{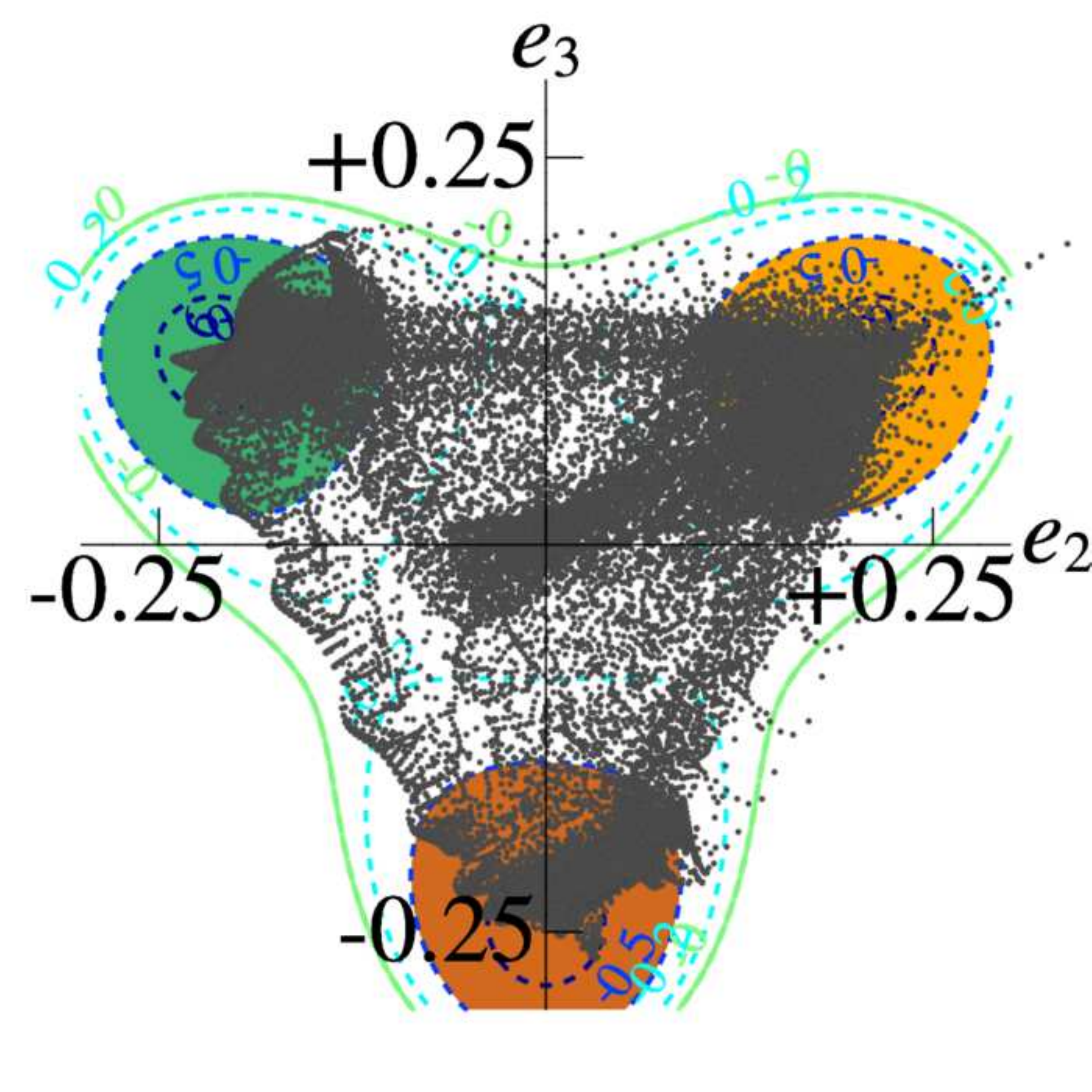} &
            \end{tabular}  \\
            \parbox[t]{0.5cm}{$C$} &
            \begin{tabular}{p{2.8cm}p{2.8cm}p{2.8cm}p{2.8cm}p{2.8cm}p{2.5cm}}
                \includegraphics[scale=0.12]{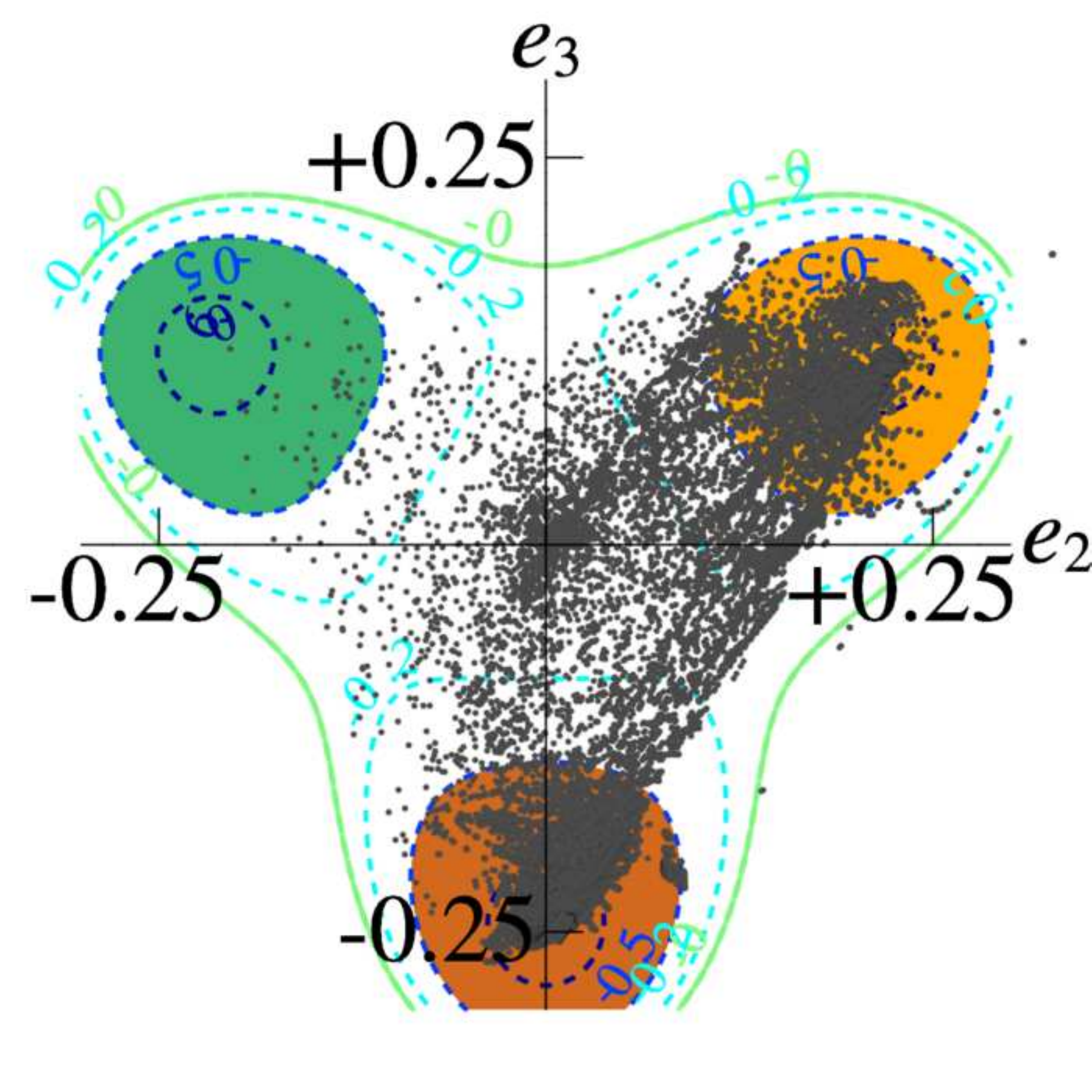} &
                \includegraphics[scale=0.12]{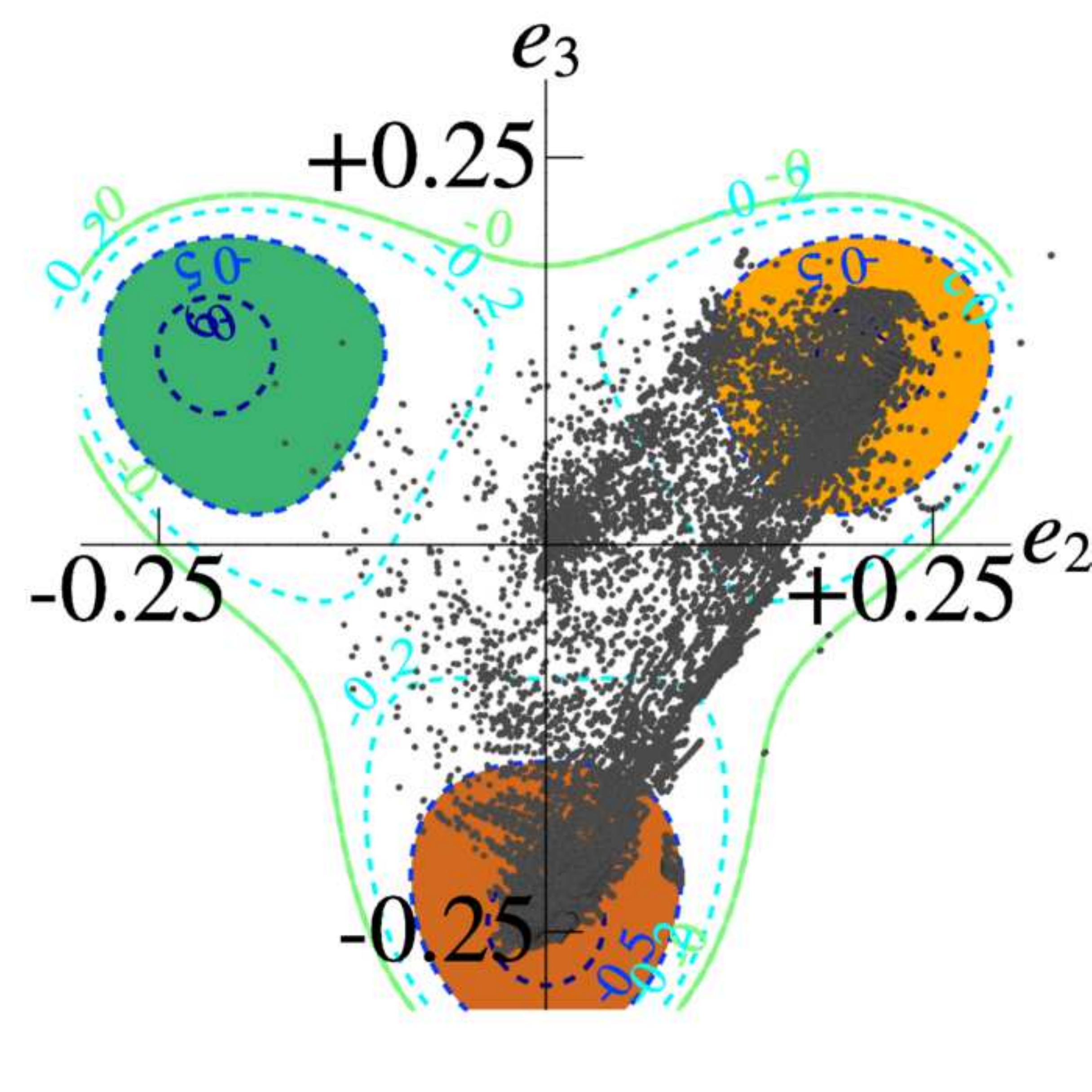} & 
                \includegraphics[scale=0.12]{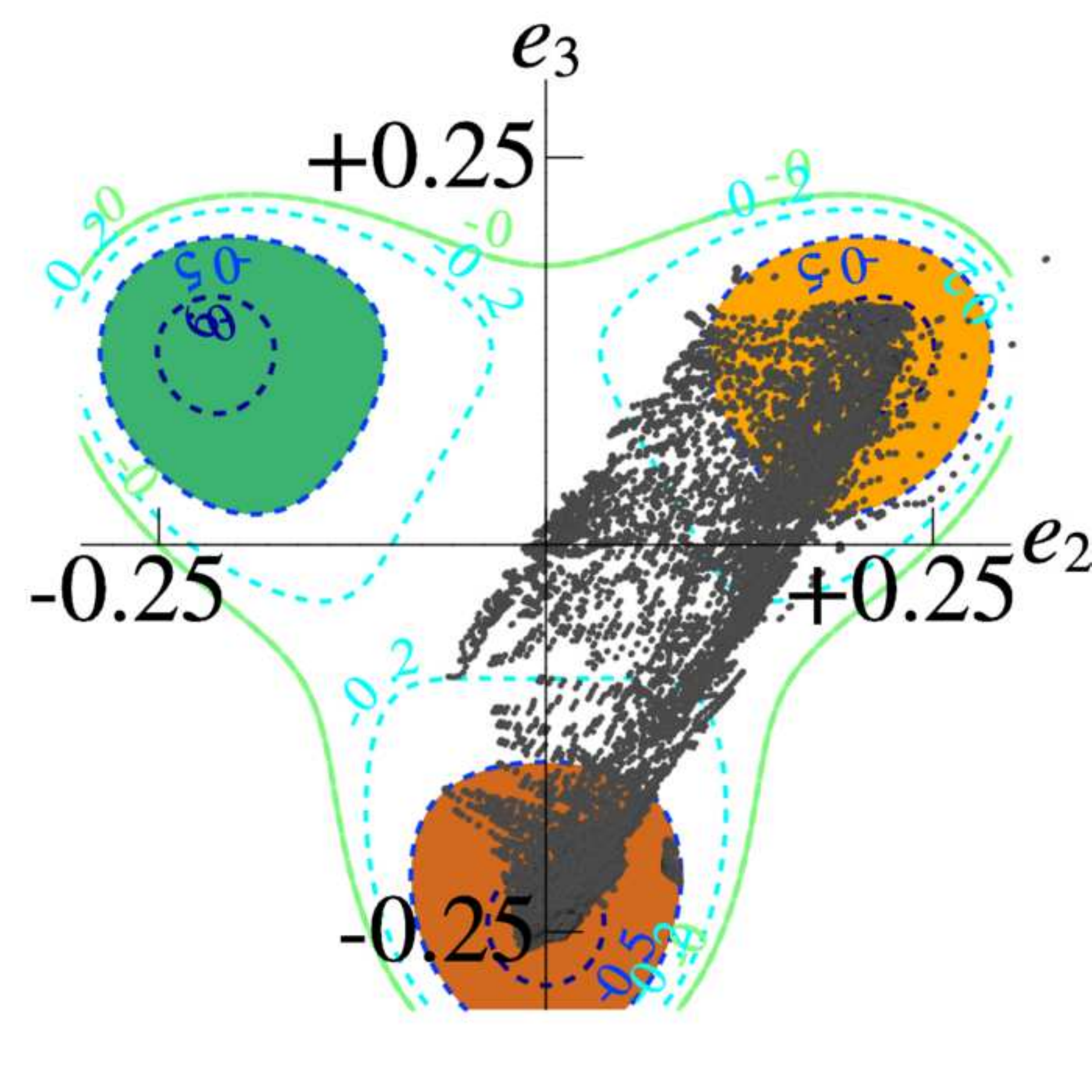} &
                \includegraphics[scale=0.12]{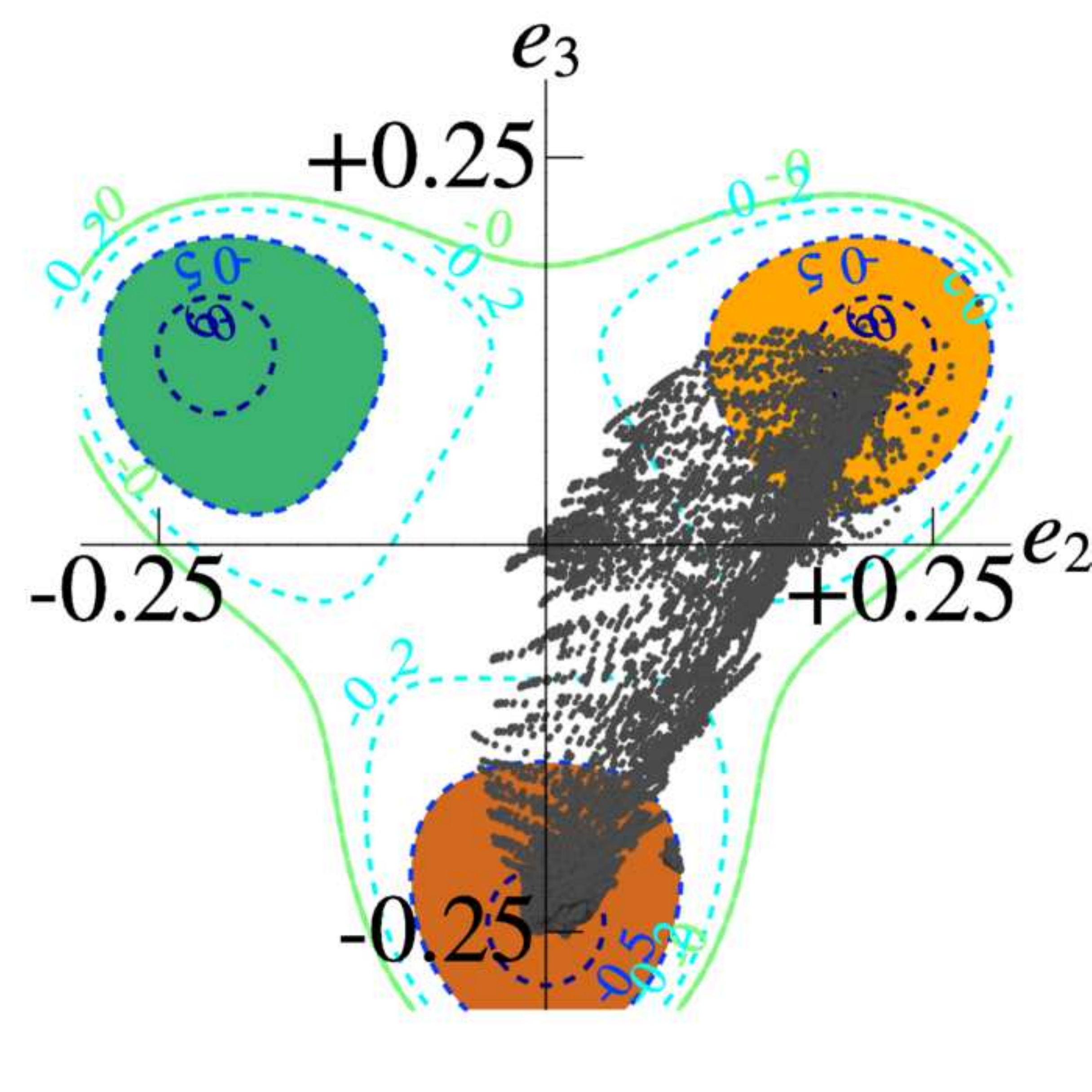} &
                \includegraphics[scale=0.12]{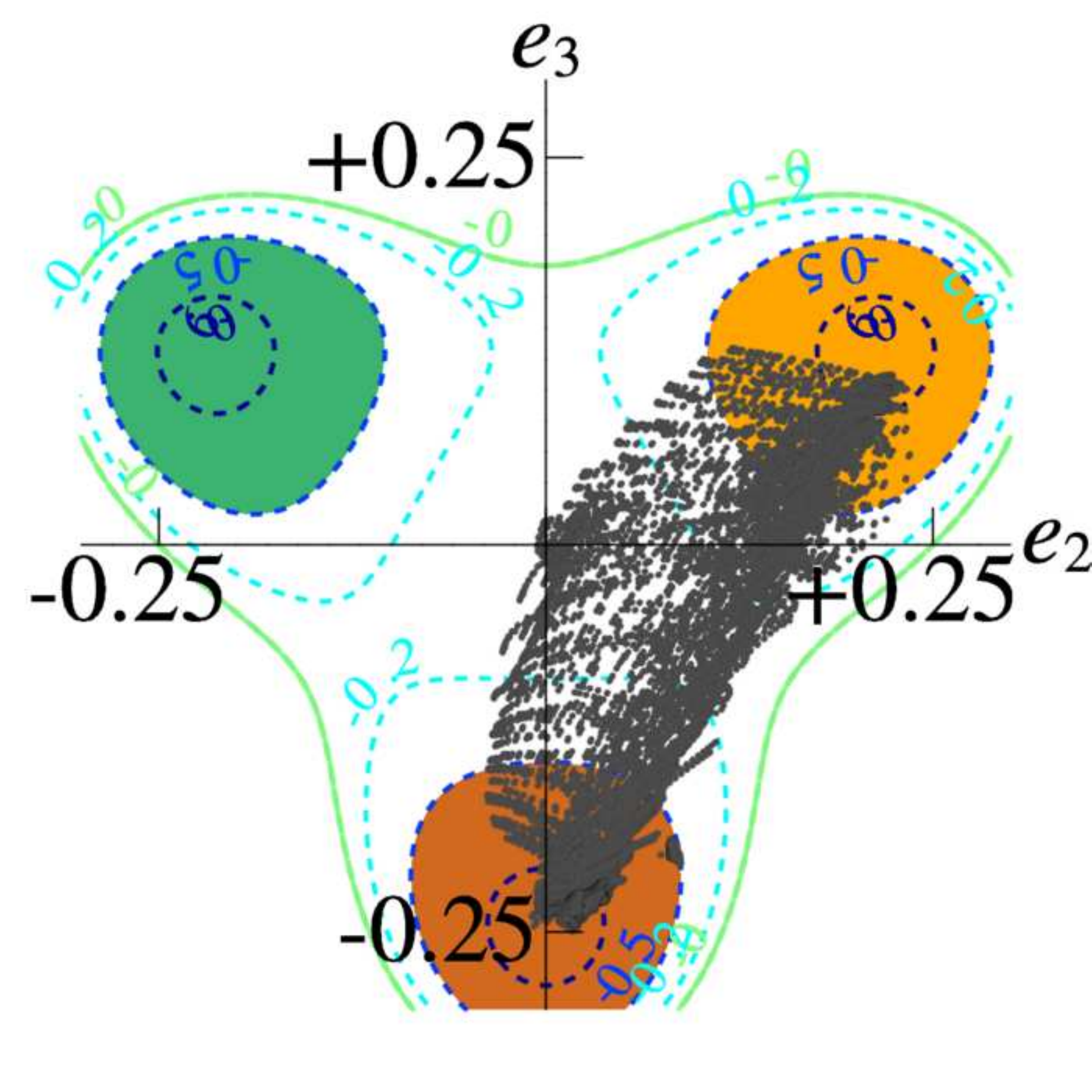} &
            \end{tabular}  \\
            \parbox[t]{0.5cm}{ $B$ } &
            \begin{tabular}{p{2.8cm}p{2.8cm}p{2.8cm}p{2.8cm}p{2.8cm}p{2.5cm}}
                \includegraphics[scale=0.12]{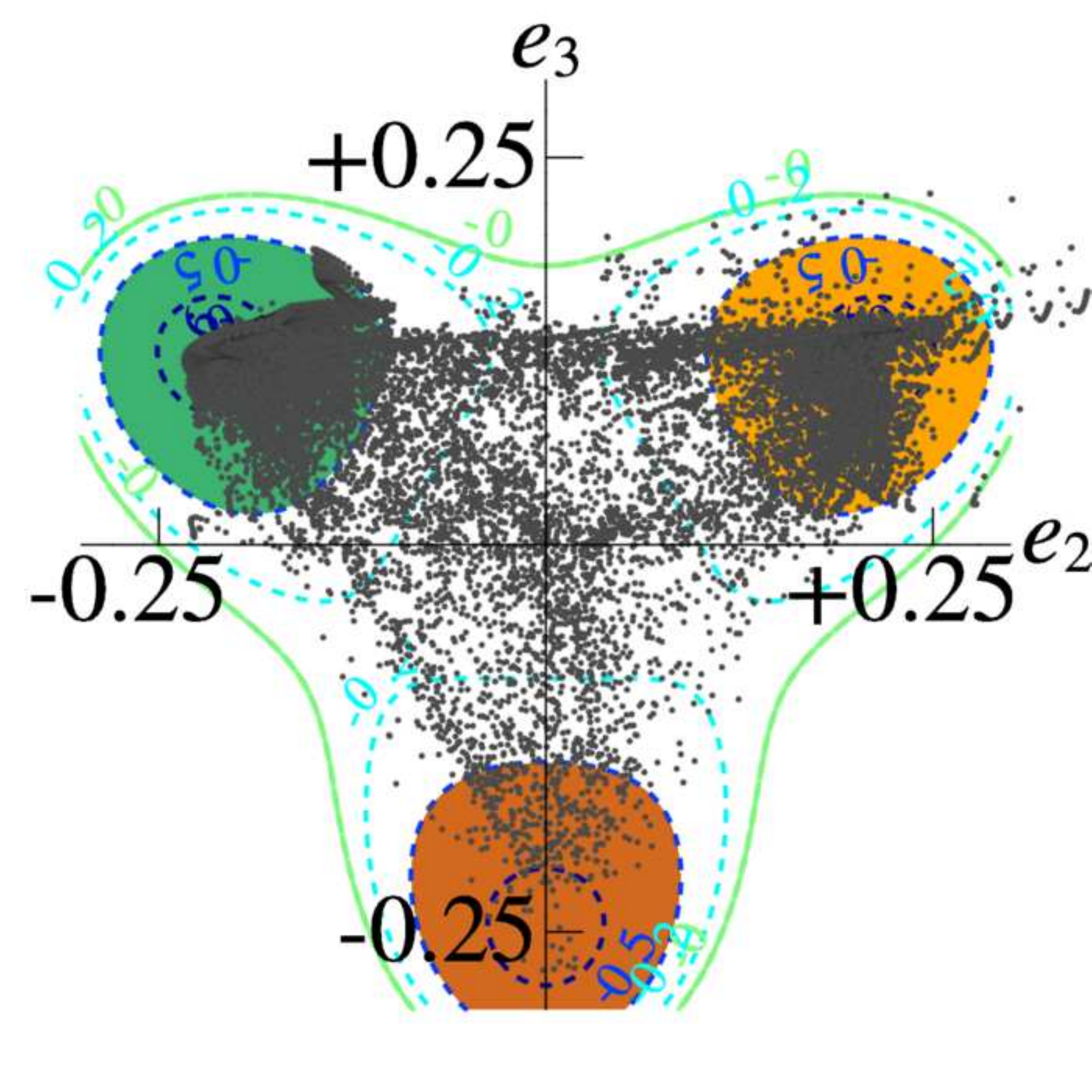} &
                \includegraphics[scale=0.12]{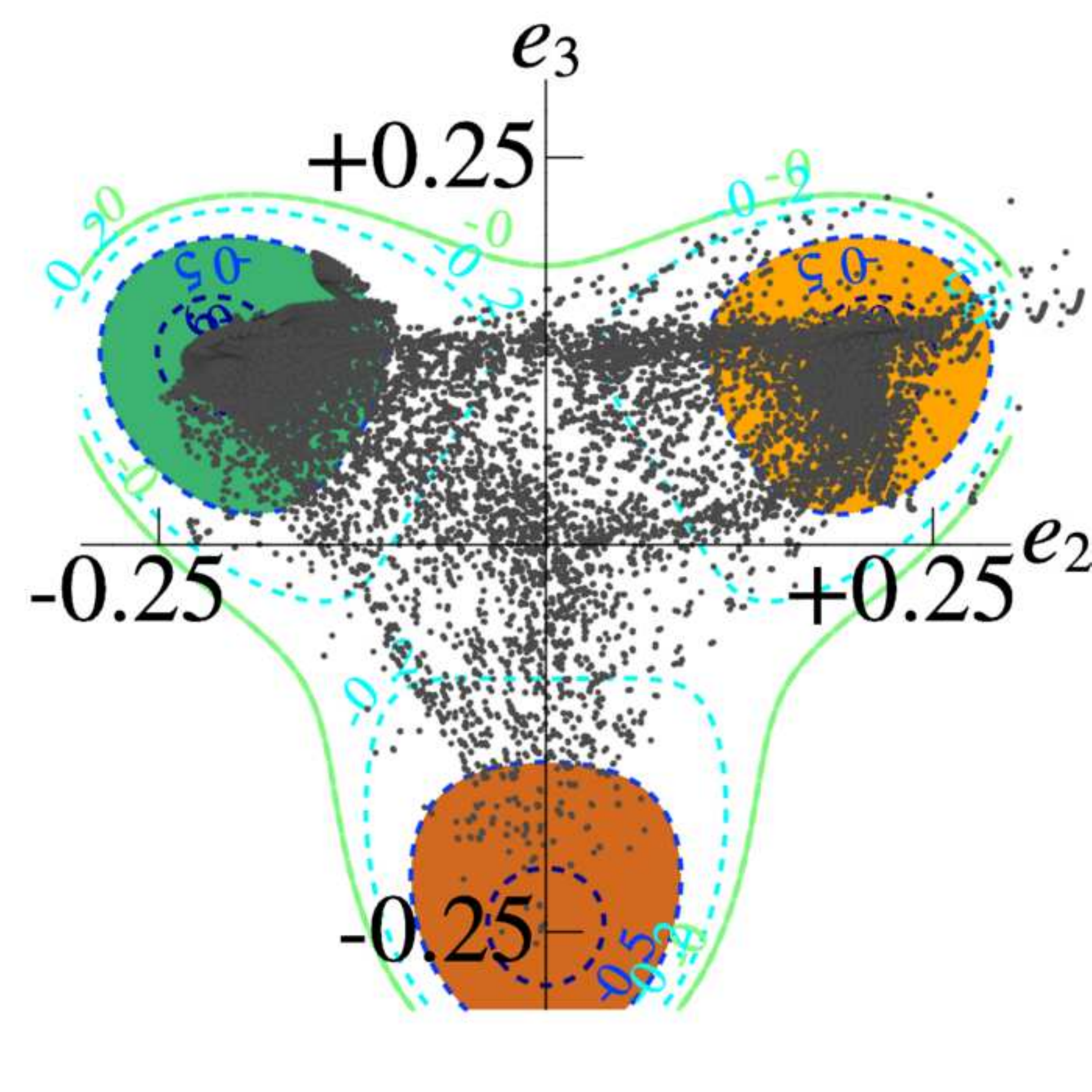} & 
                \includegraphics[scale=0.12]{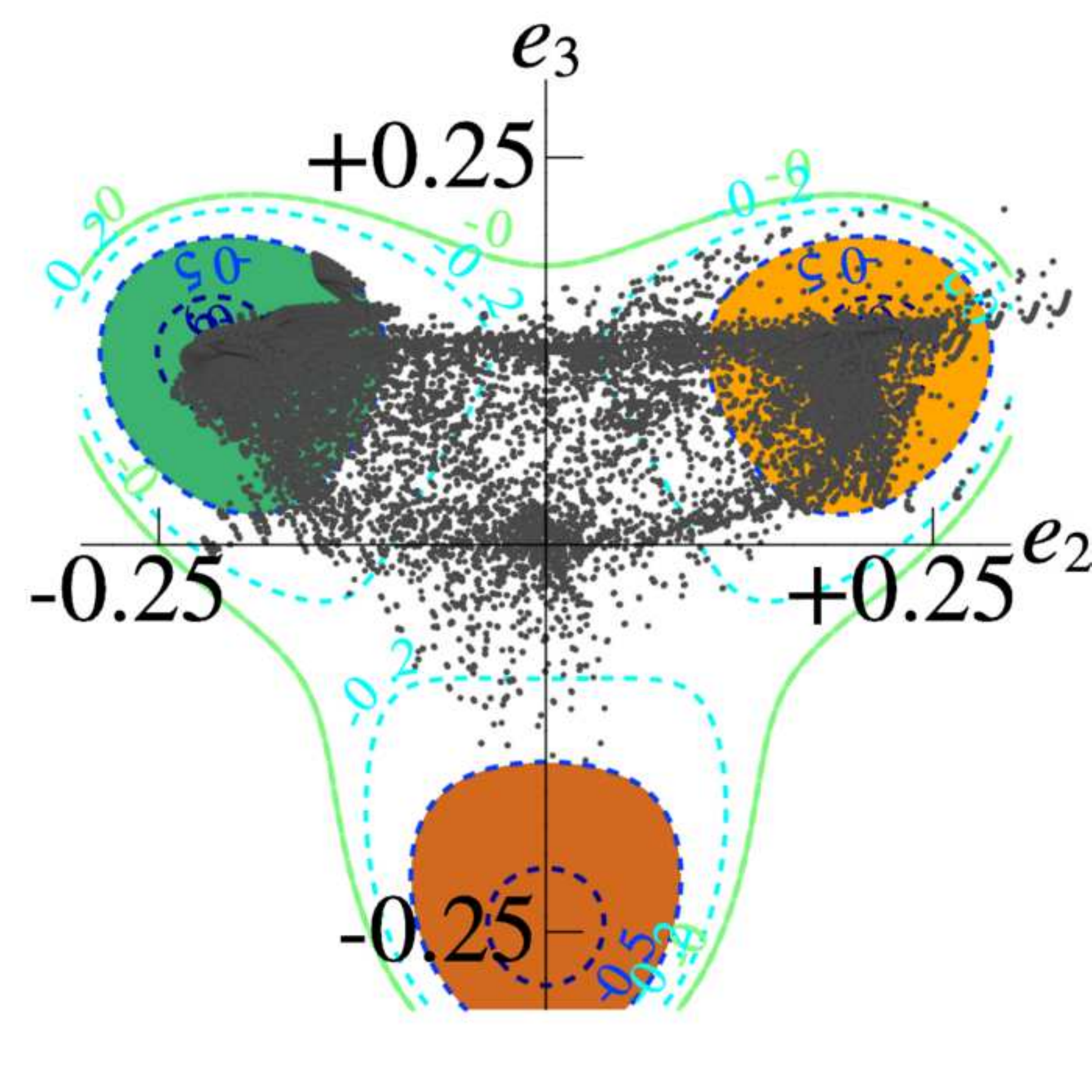} &
                \includegraphics[scale=0.12]{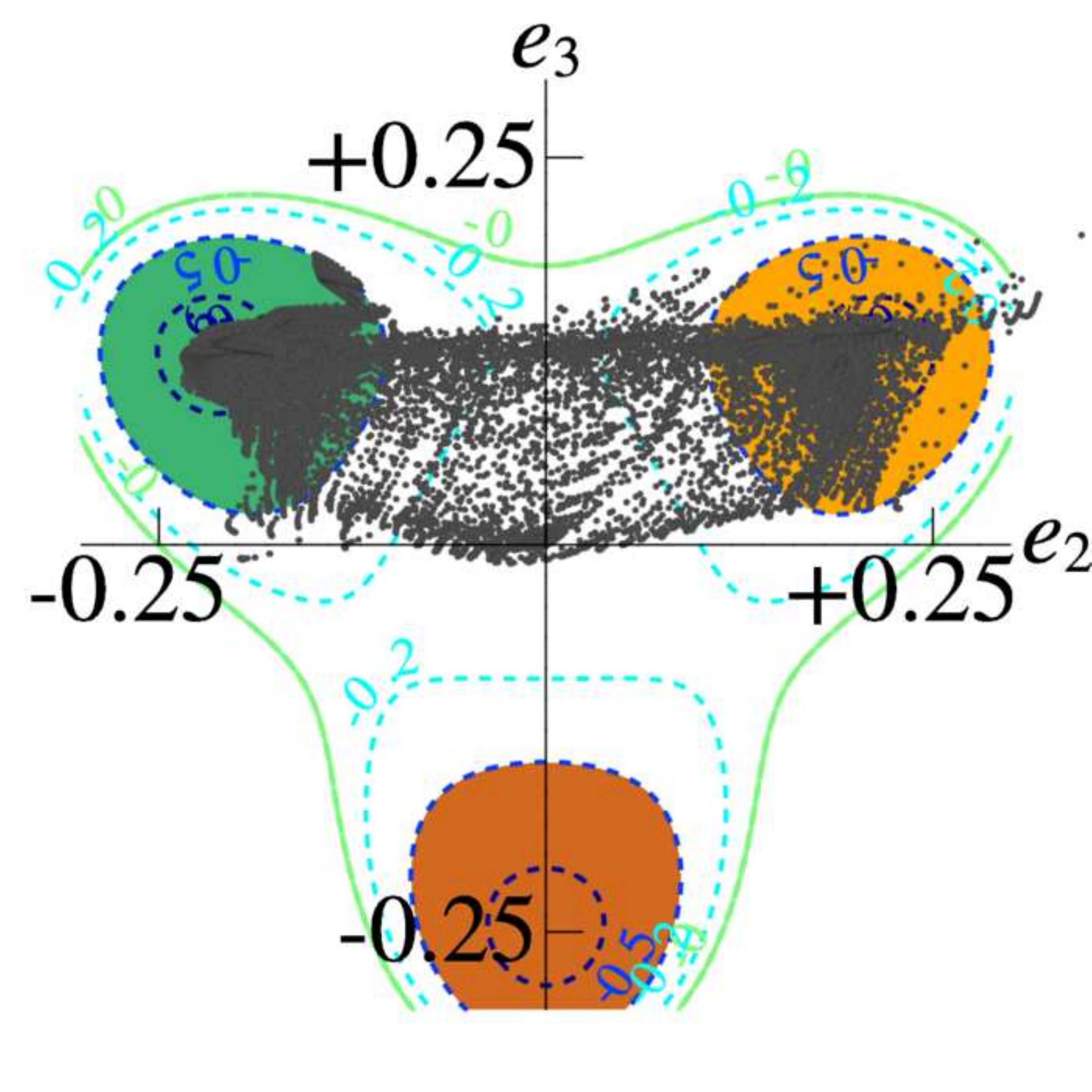} &
                \includegraphics[scale=0.12]{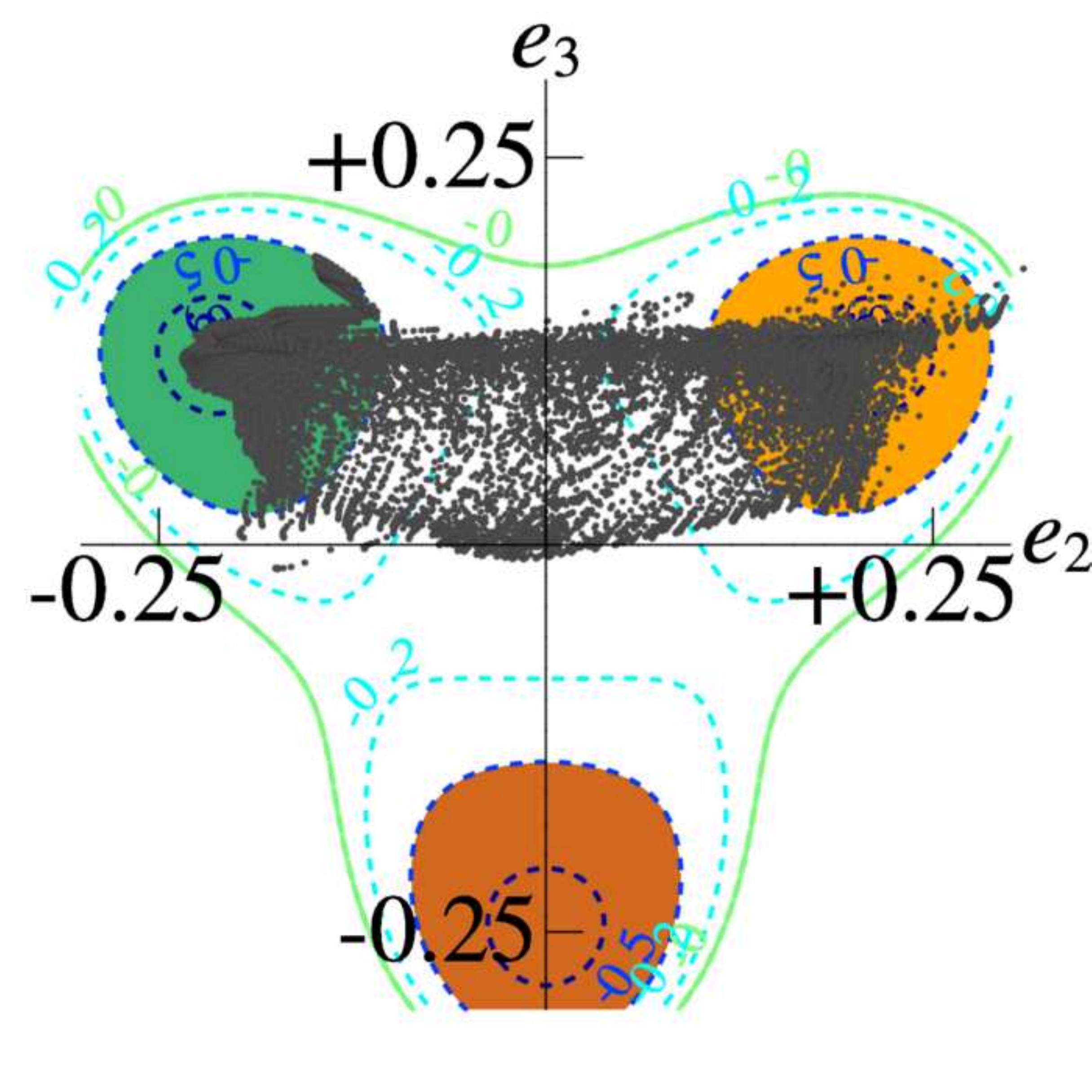} &
            \end{tabular}  \\
            \parbox[t]{0.5cm}{ $A$ } &
            \begin{tabular}{p{2.8cm}p{2.8cm}p{2.8cm}p{2.8cm}p{2.8cm}p{2.5cm}}
                &
                \includegraphics[scale=0.12]{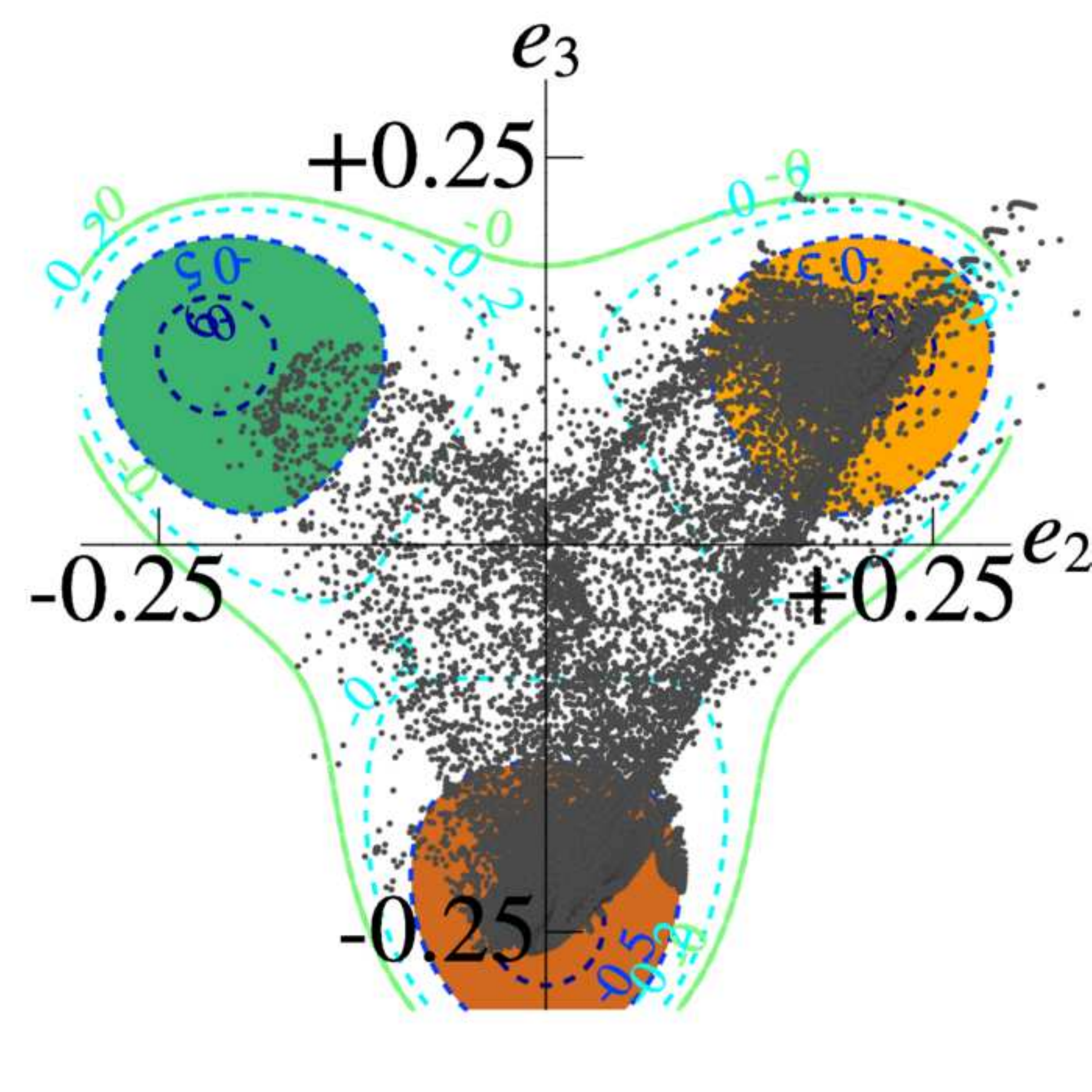} & 
                \includegraphics[scale=0.12]{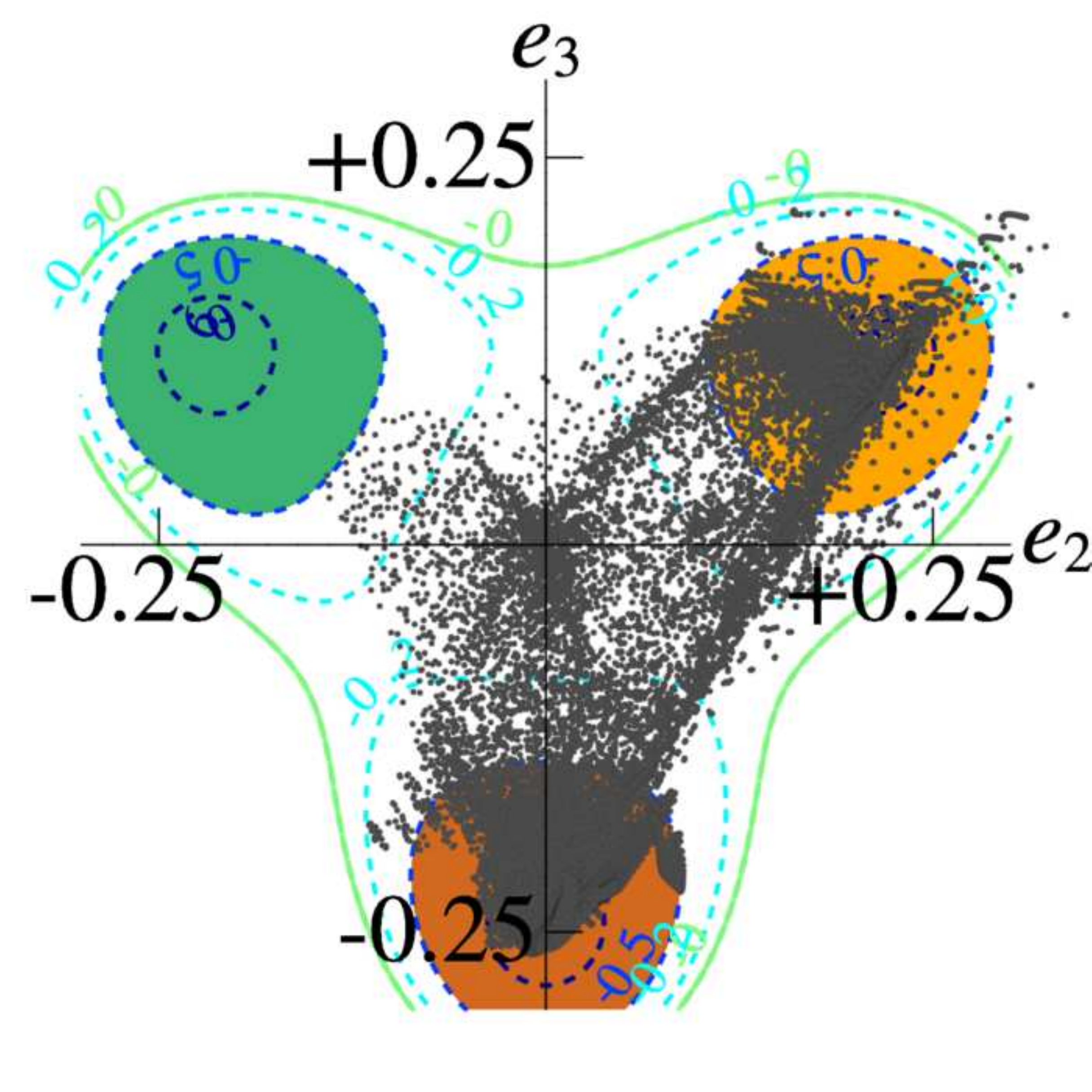} &
                \includegraphics[scale=0.12]{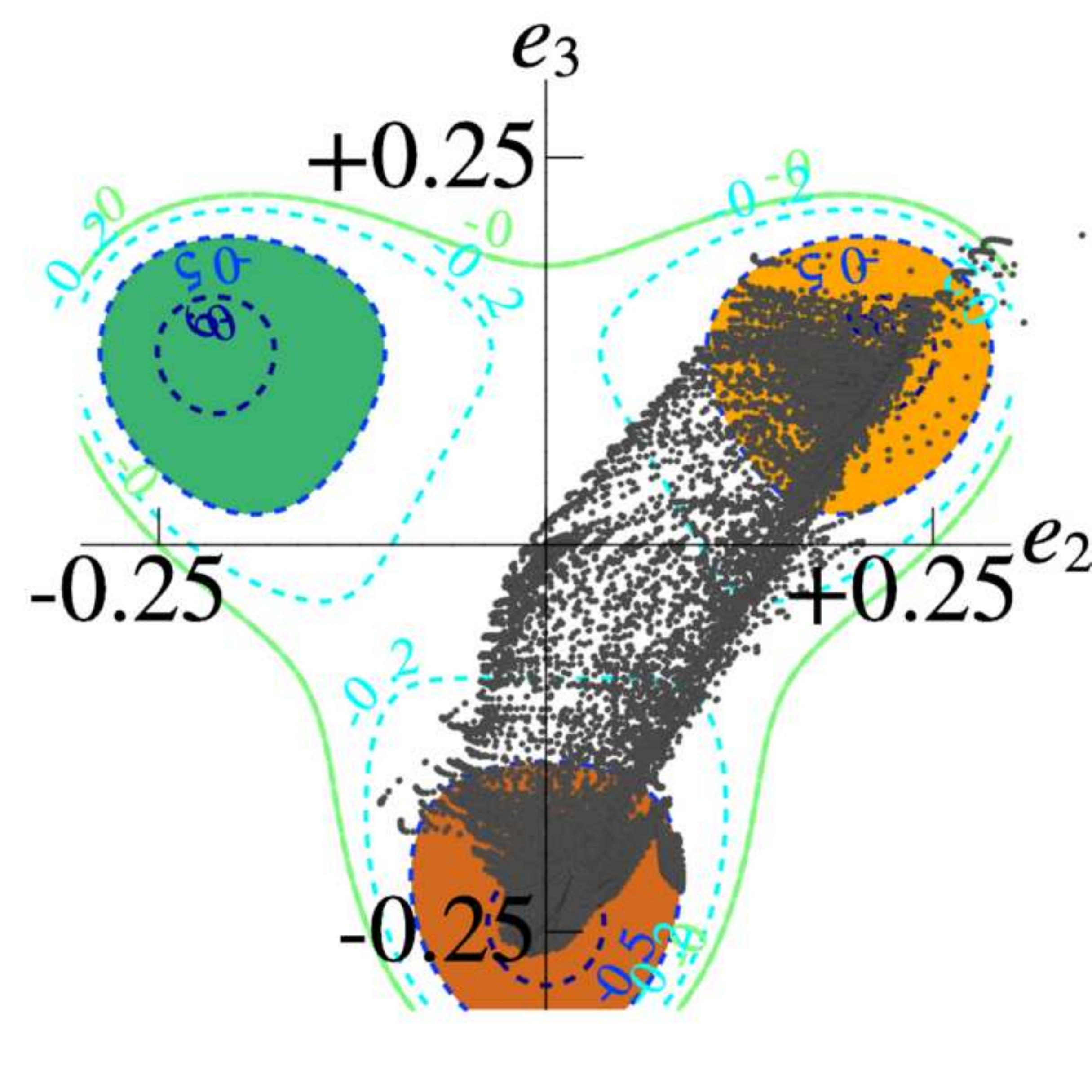} &
                \includegraphics[scale=0.12]{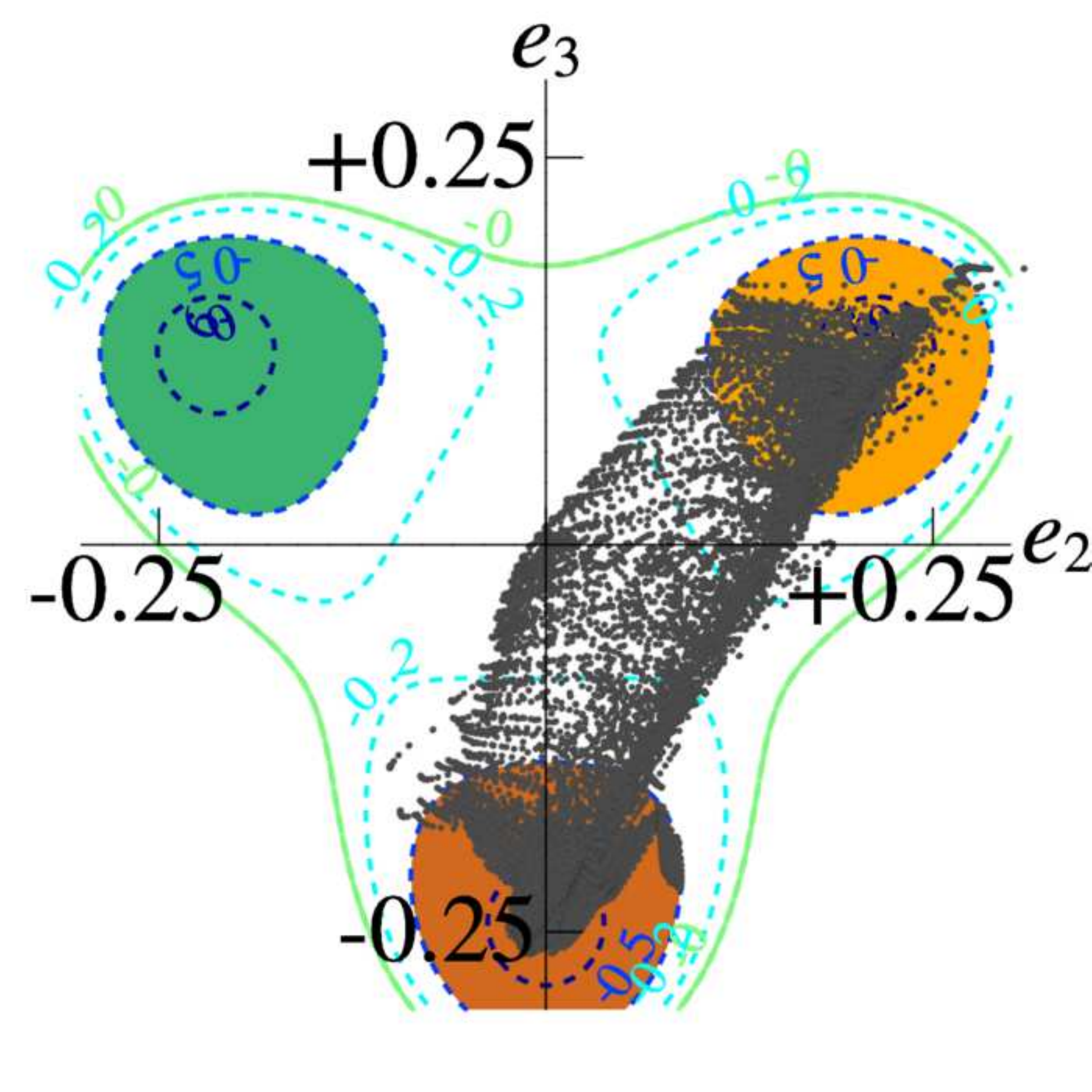} &
                \includegraphics[scale=0.12]{plot/plot_colorbar_phase.pdf}
            \end{tabular}
        \end{tabular}
    \end{center}
    \caption{Contours of $(e_2,e_3)$ for branches A - E computed from the solutions on the $128^3$ mesh at $33^3$ uniformly spaced points. Selected values of $l$ are indicated.   
    The three-well contour diagram, Fig. \ref{Fi:3well_diagram}, is also shown, underlying each plot as an indicator of the formation of tetragonal variants. }
    \label{Fi:e1-e2}
\end{figure}

\begin{figure}
    \begin{center}
        \begin{tabular}{rp{16.5cm}}
            \parbox[t]{0.5cm}{ }&
            \begin{tabular}{p{2.8cm}p{2.8cm}p{2.8cm}p{2.8cm}p{2.8cm}p{2.5cm}}
                \hspace{0.8cm}$l\!=\!0.0625$ & 
                \hspace{0.8cm}$l\!=\!0.0750$ & 
                \hspace{0.8cm}$l\!=\!0.1000$ &
                \hspace{0.8cm}$l\!=\!0.1500$ &
                \hspace{0.8cm}$l\!=\!0.2000$ &
                \vspace{0.5\baselineskip}
            \end{tabular}  \\
            \parbox[t]{0.5cm}{$E$} &
            \begin{tabular}{p{2.8cm}p{2.8cm}p{2.8cm}p{2.8cm}p{2.8cm}p{2.5cm}}
                &
                \includegraphics[scale=0.12]{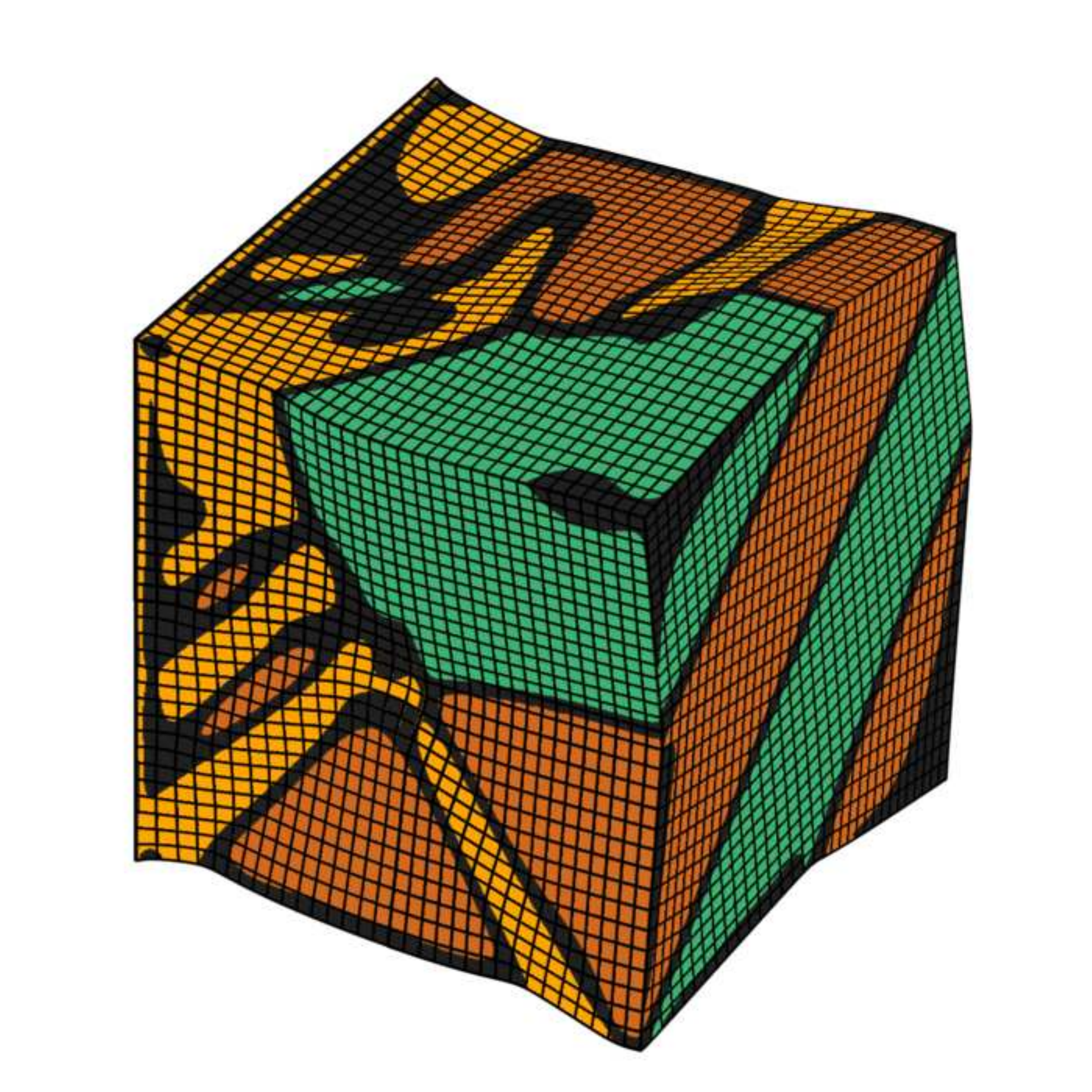} & 
                \includegraphics[scale=0.12]{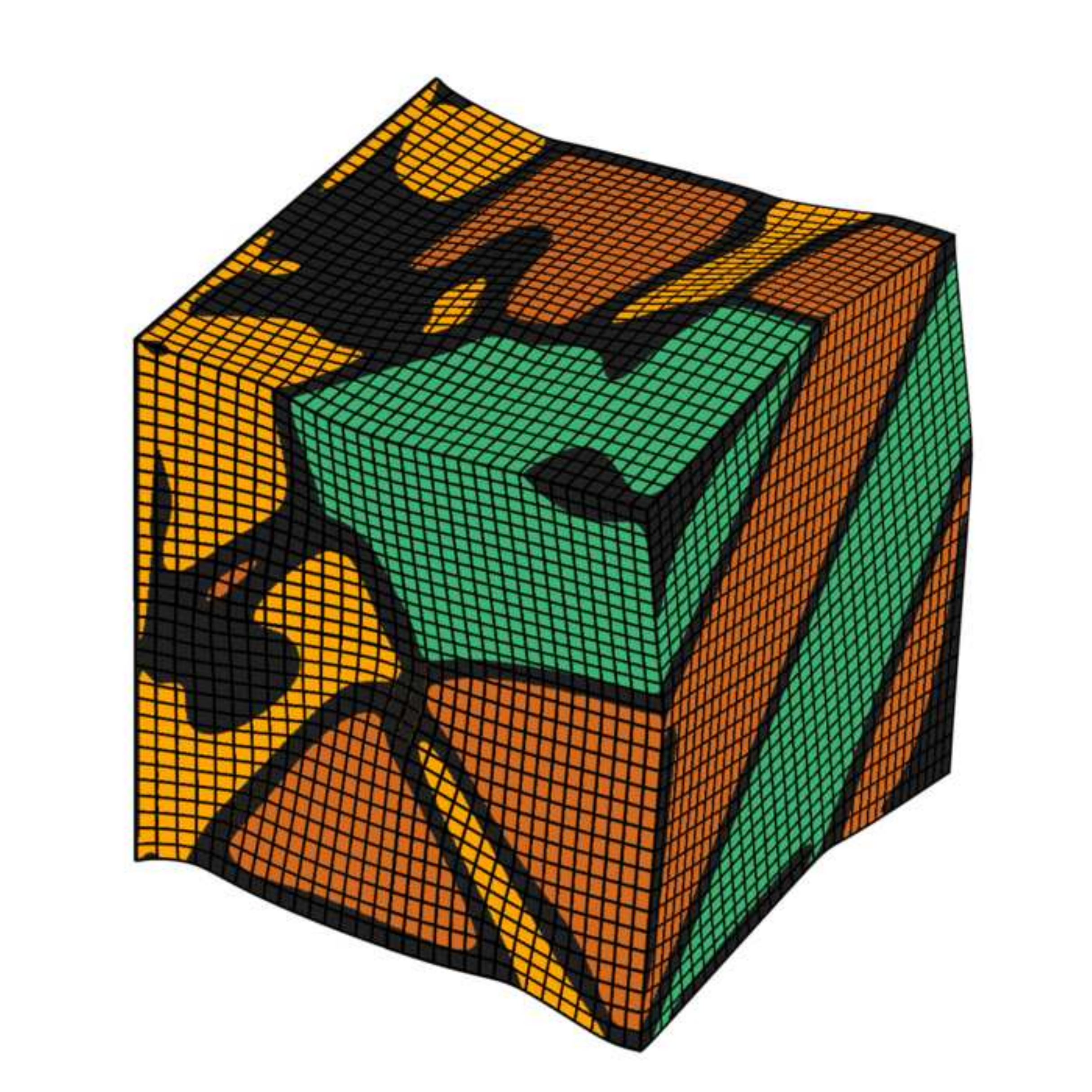} &
                \includegraphics[scale=0.12]{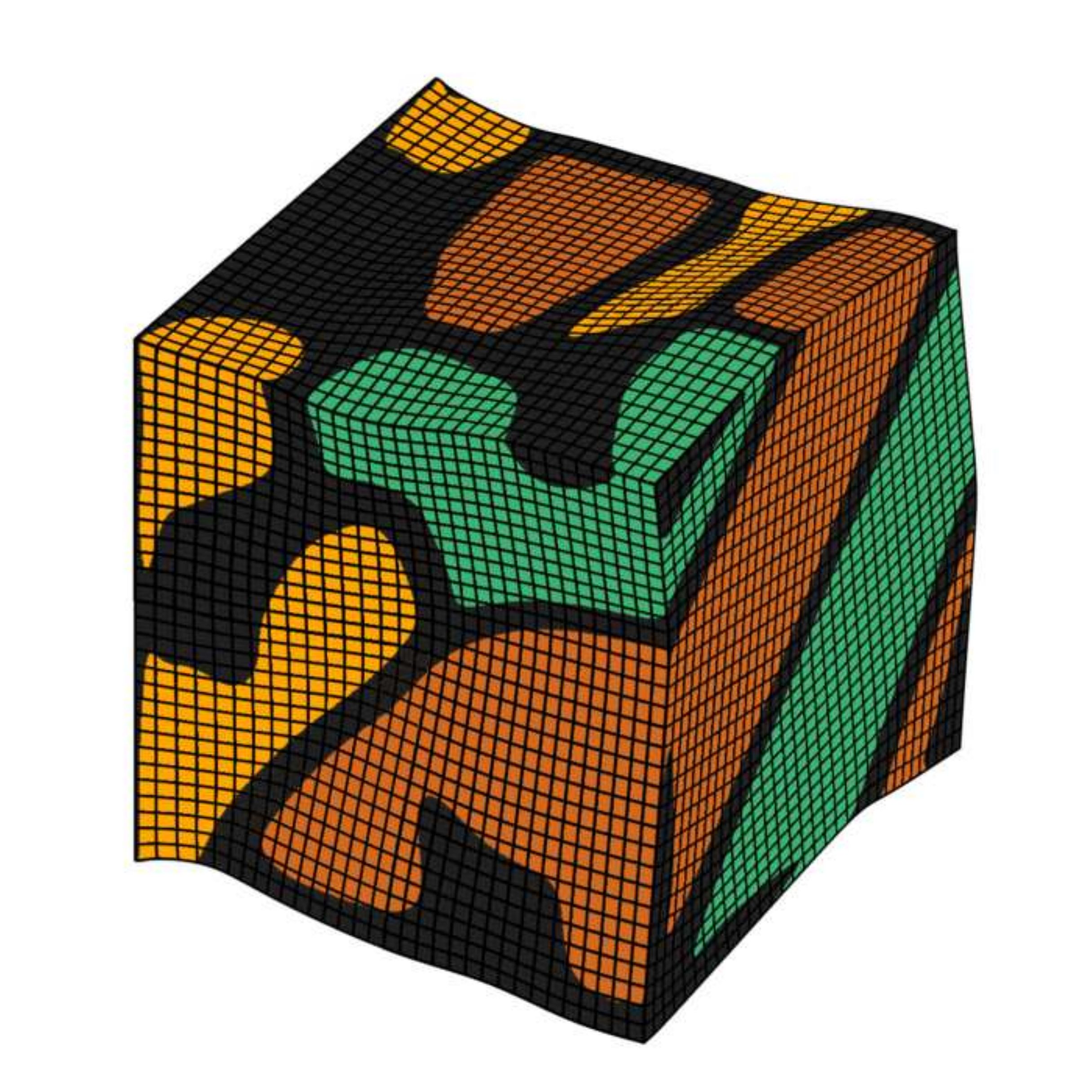} &
                \includegraphics[scale=0.12]{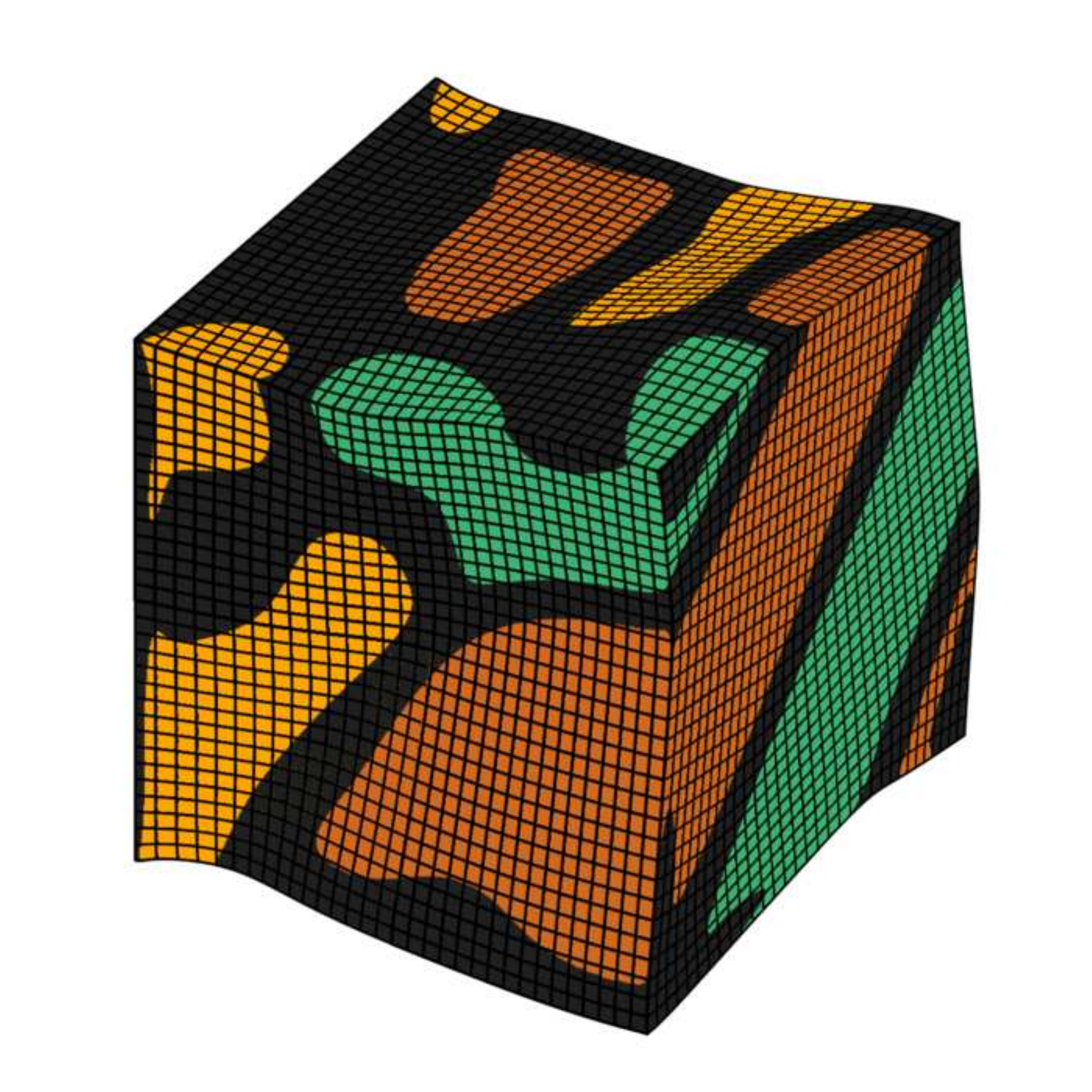} &
            \end{tabular}  \\
            \parbox[t]{0.5cm}{ $D$ } &
            \begin{tabular}{p{2.8cm}p{2.8cm}p{2.8cm}p{2.8cm}p{2.8cm}p{2.5cm}}
                \includegraphics[scale=0.12]{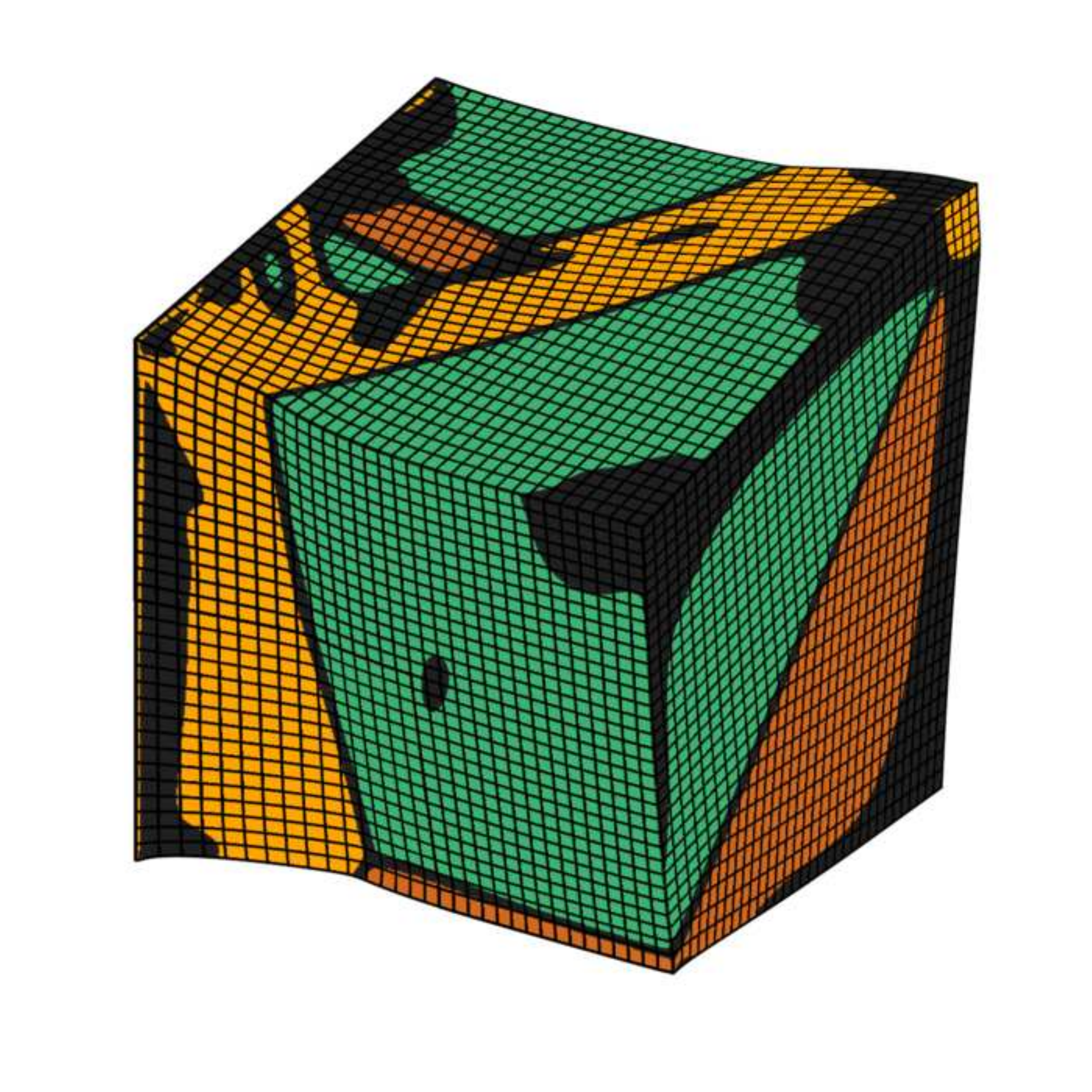} &
                \includegraphics[scale=0.12]{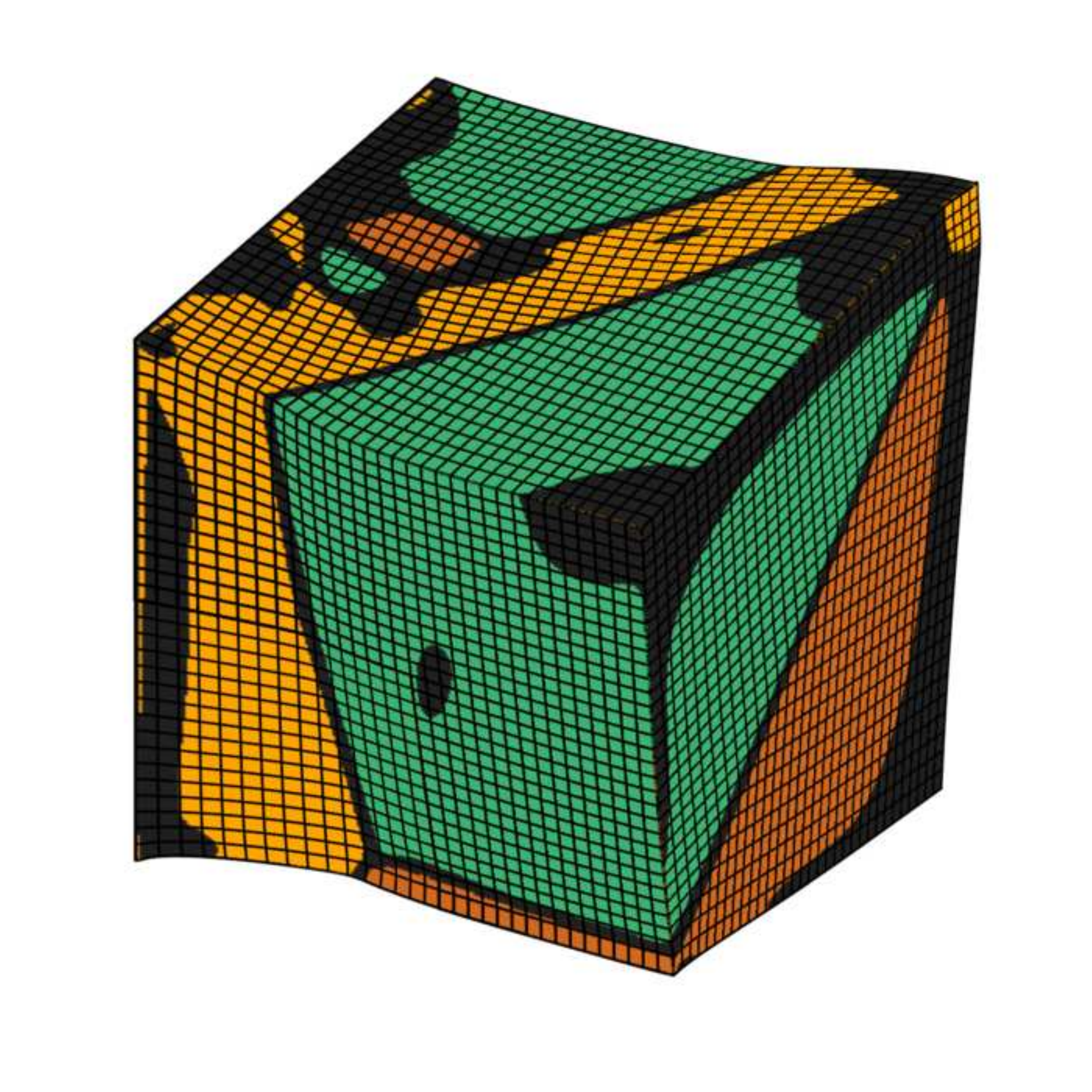} & 
                \includegraphics[scale=0.12]{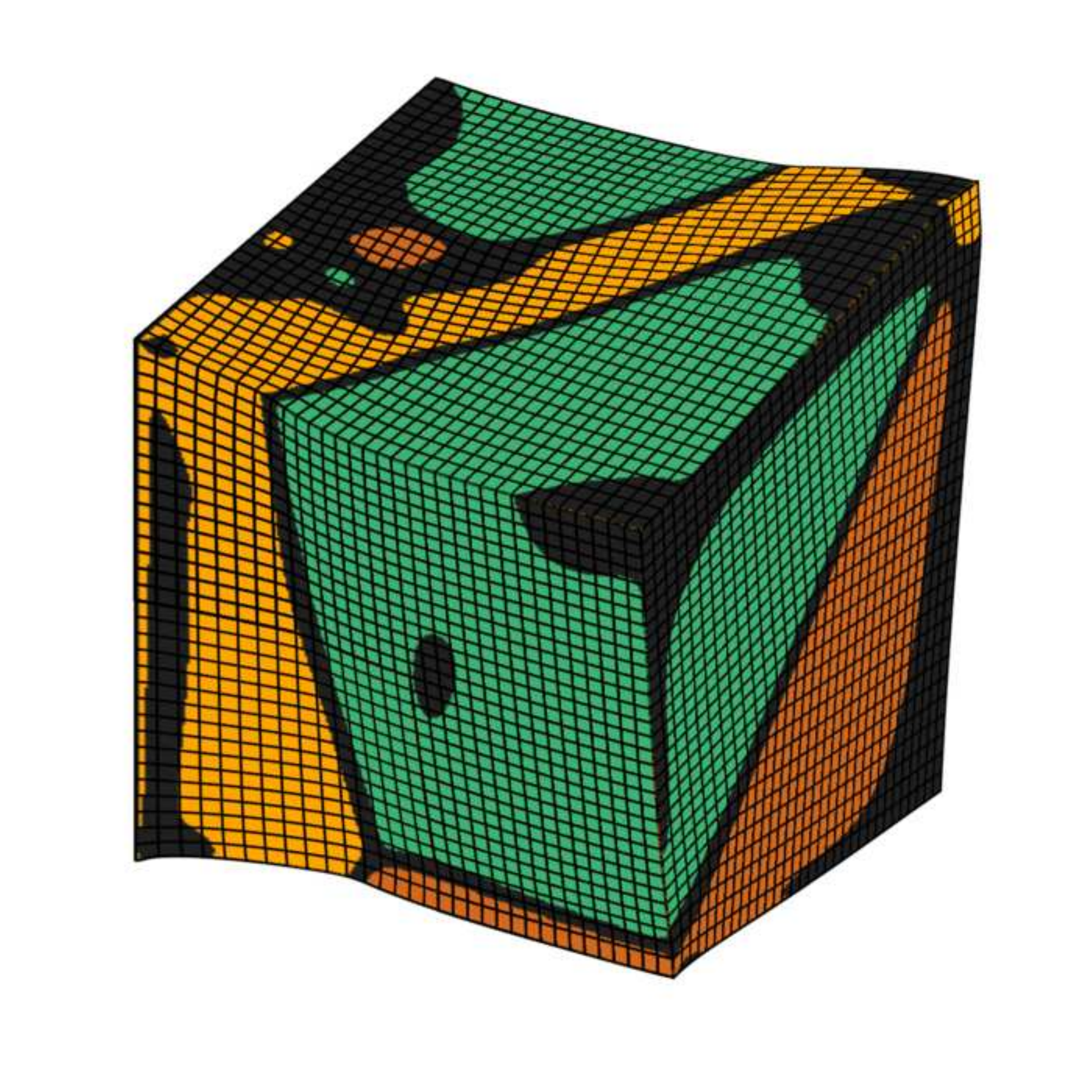} &
                \includegraphics[scale=0.12]{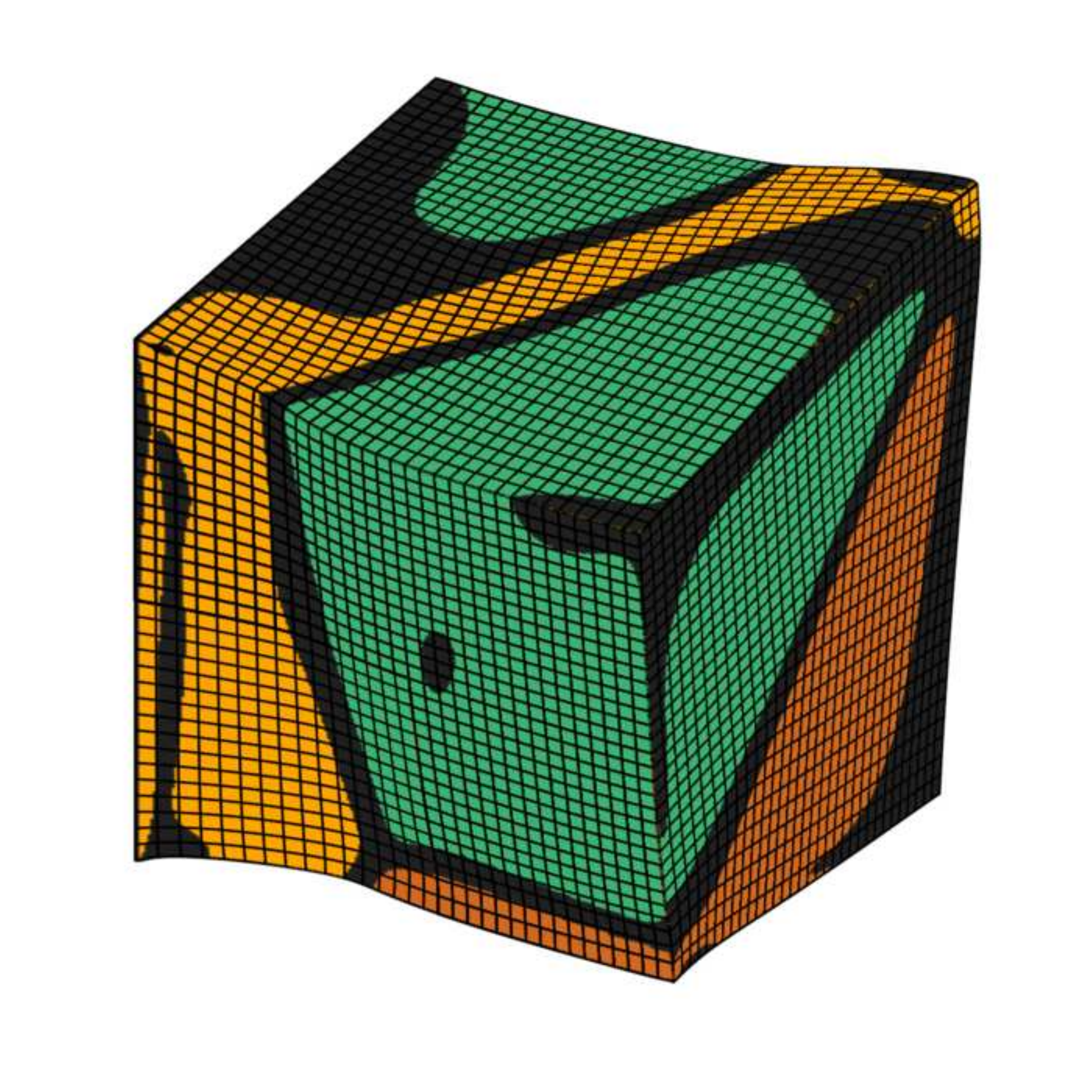} &
                \includegraphics[scale=0.12]{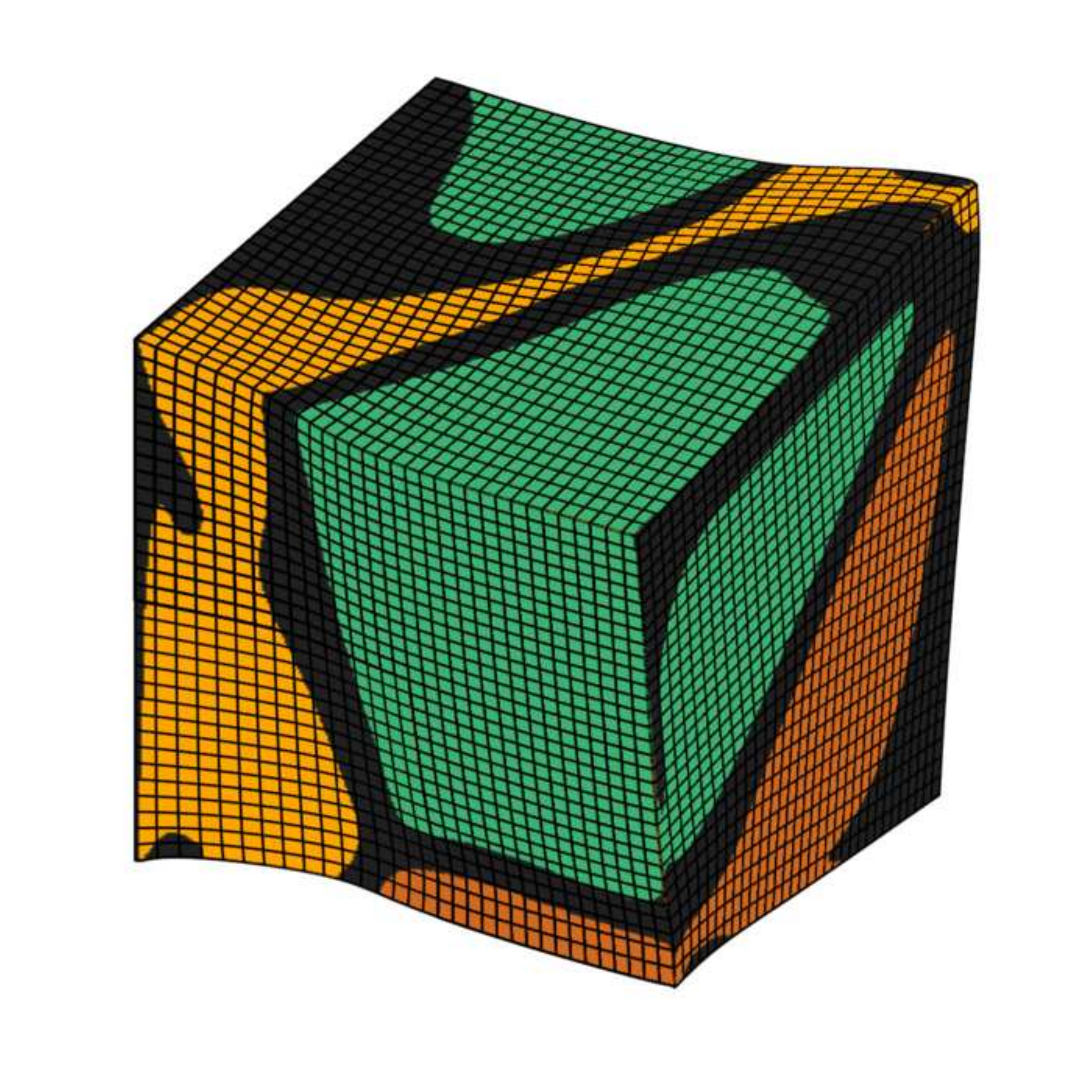} &
            \end{tabular}  \\
            \parbox[t]{0.5cm}{$C$} &
            \begin{tabular}{p{2.8cm}p{2.8cm}p{2.8cm}p{2.8cm}p{2.8cm}p{2.5cm}}
                \includegraphics[scale=0.12]{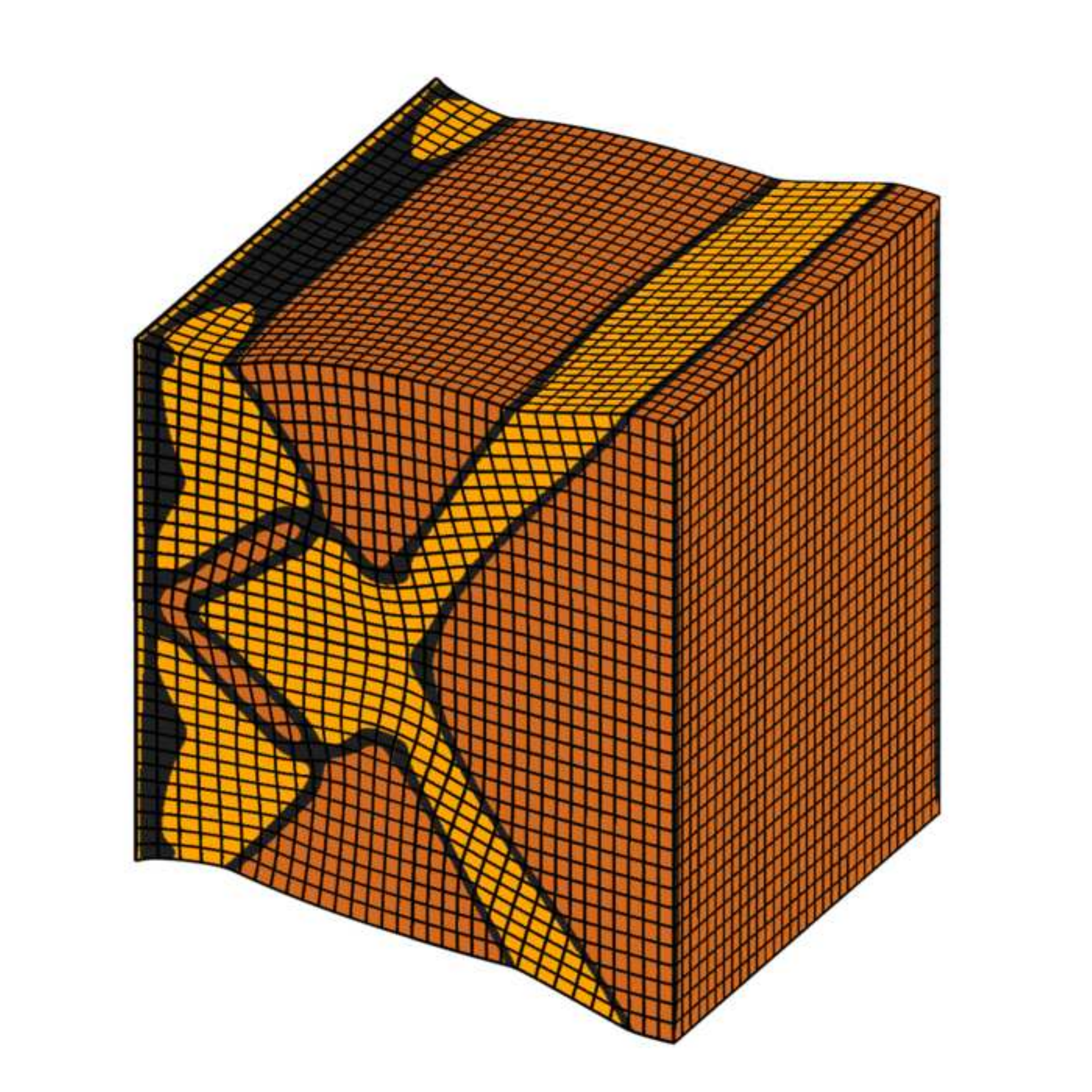} &
                \includegraphics[scale=0.12]{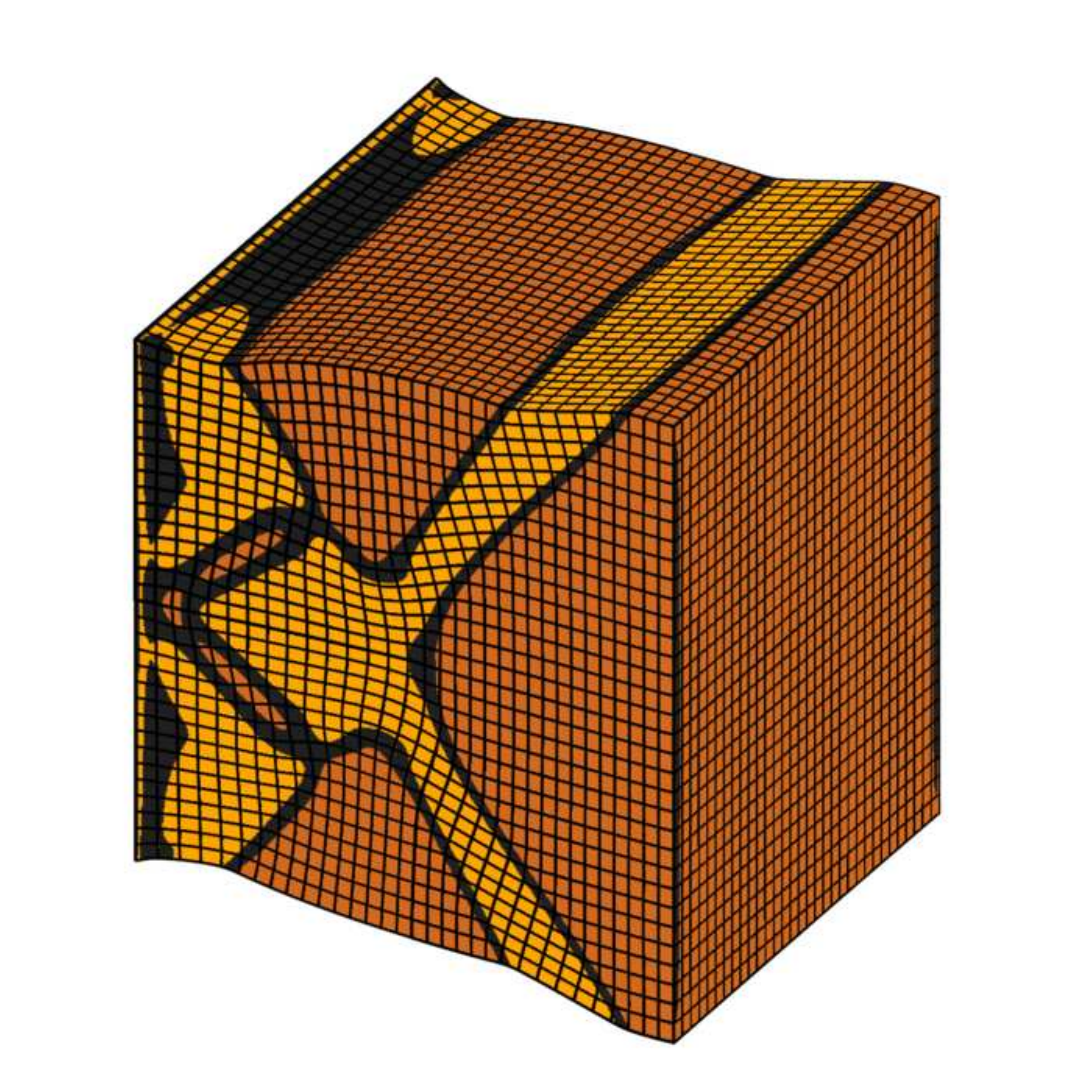} & 
                \includegraphics[scale=0.12]{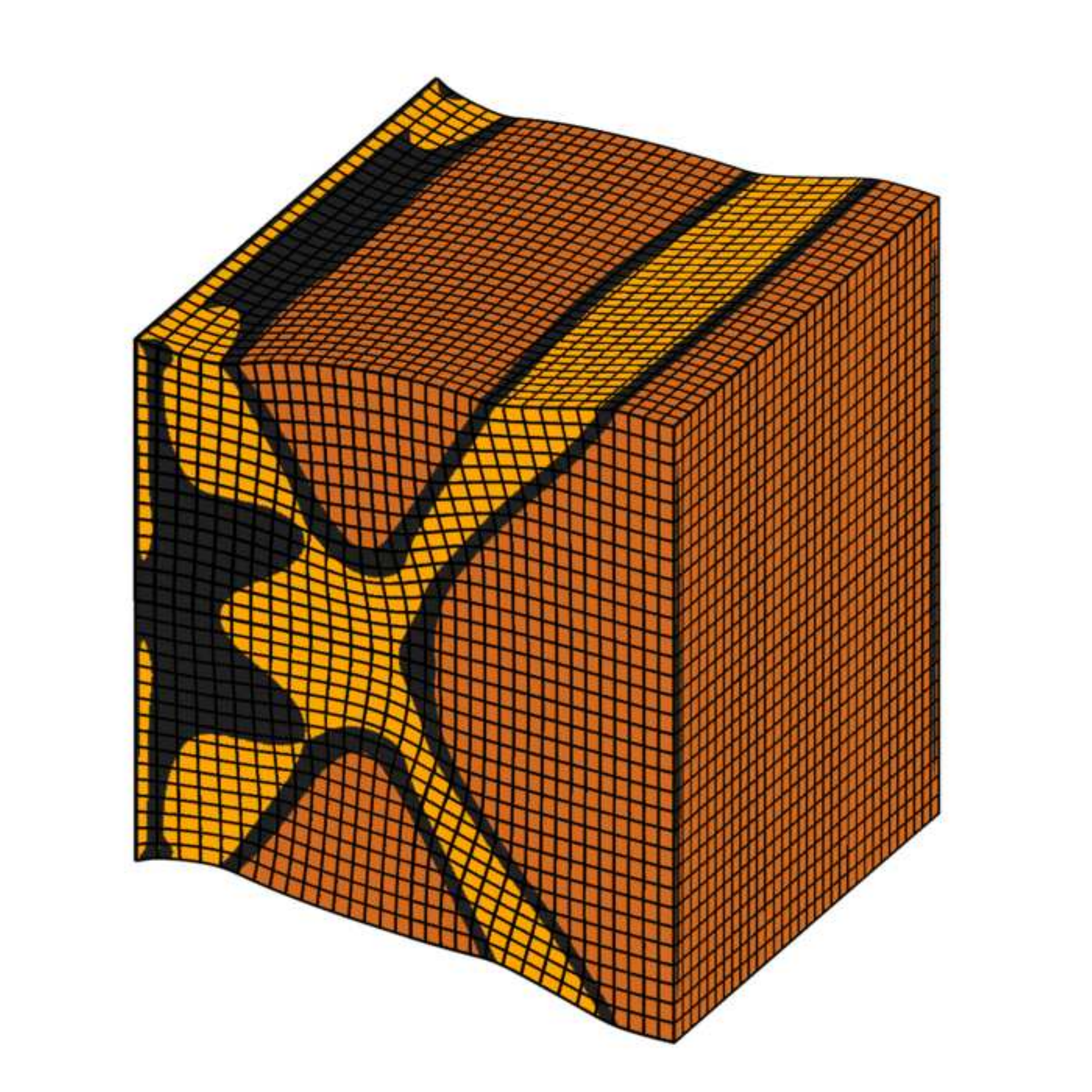} &
                \includegraphics[scale=0.12]{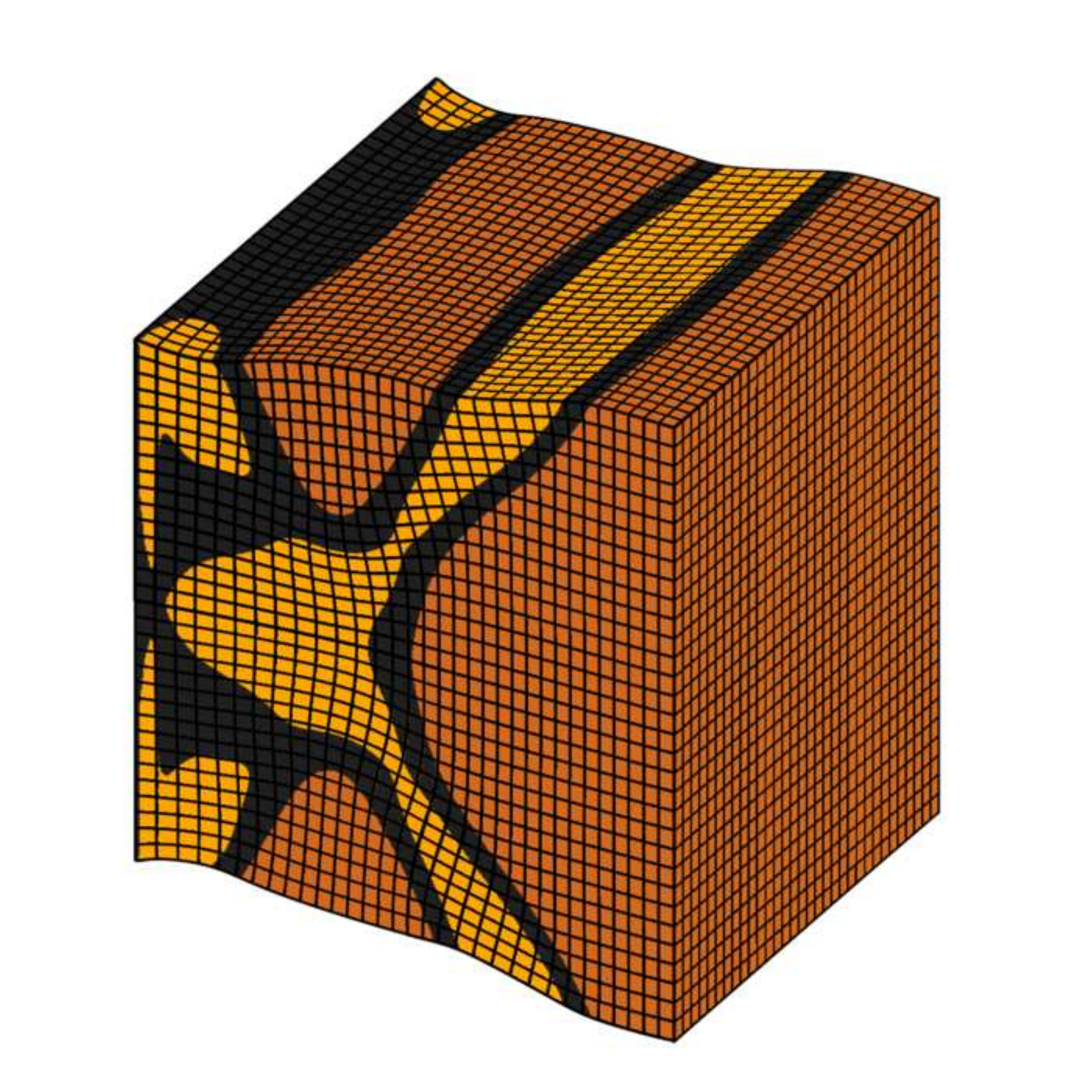} &
                \includegraphics[scale=0.12]{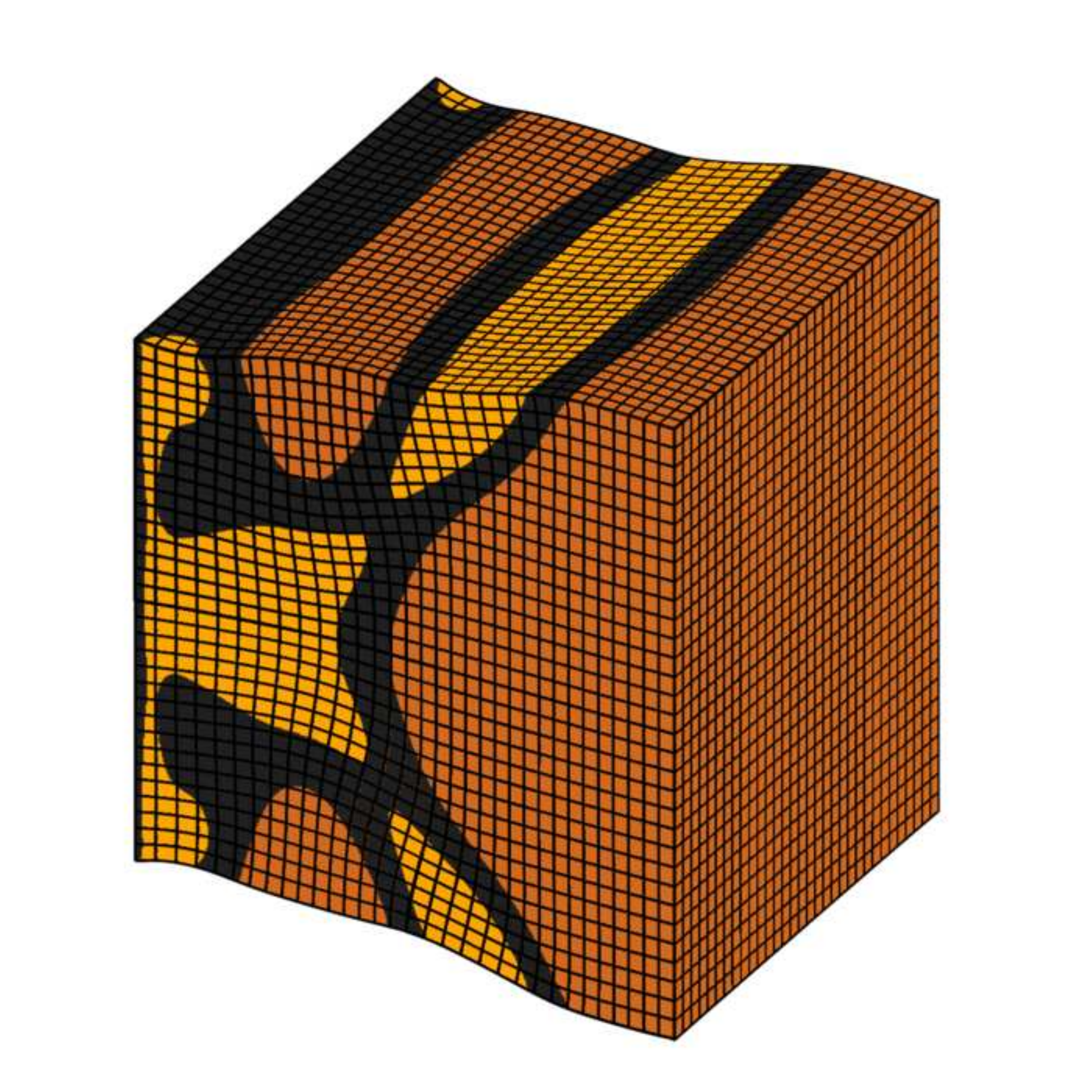} &
            \end{tabular}  \\
            \parbox[t]{0.5cm}{$B$} &
            \begin{tabular}{p{2.8cm}p{2.8cm}p{2.8cm}p{2.8cm}p{2.8cm}p{2.5cm}}
                \includegraphics[scale=0.12]{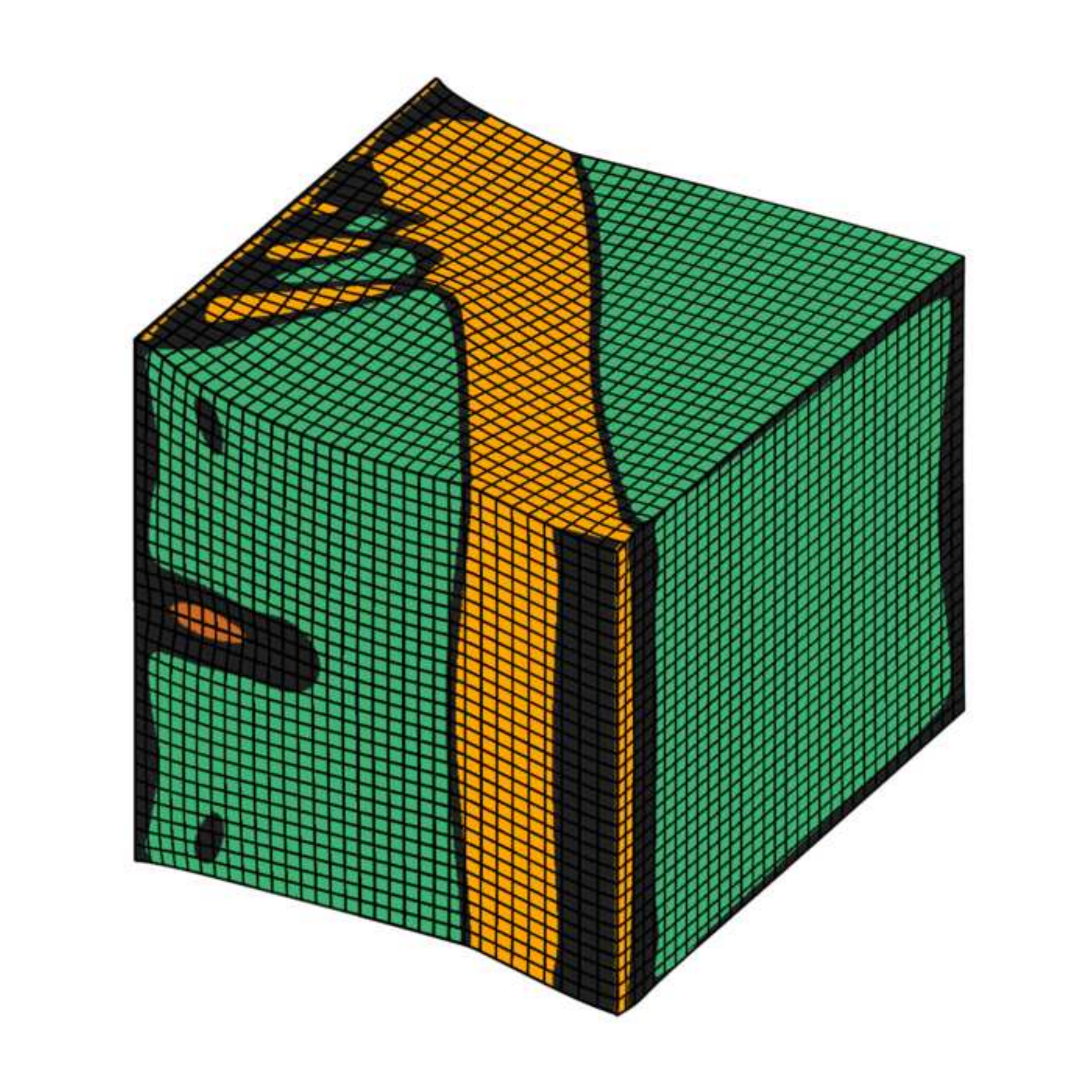} &
                \includegraphics[scale=0.12]{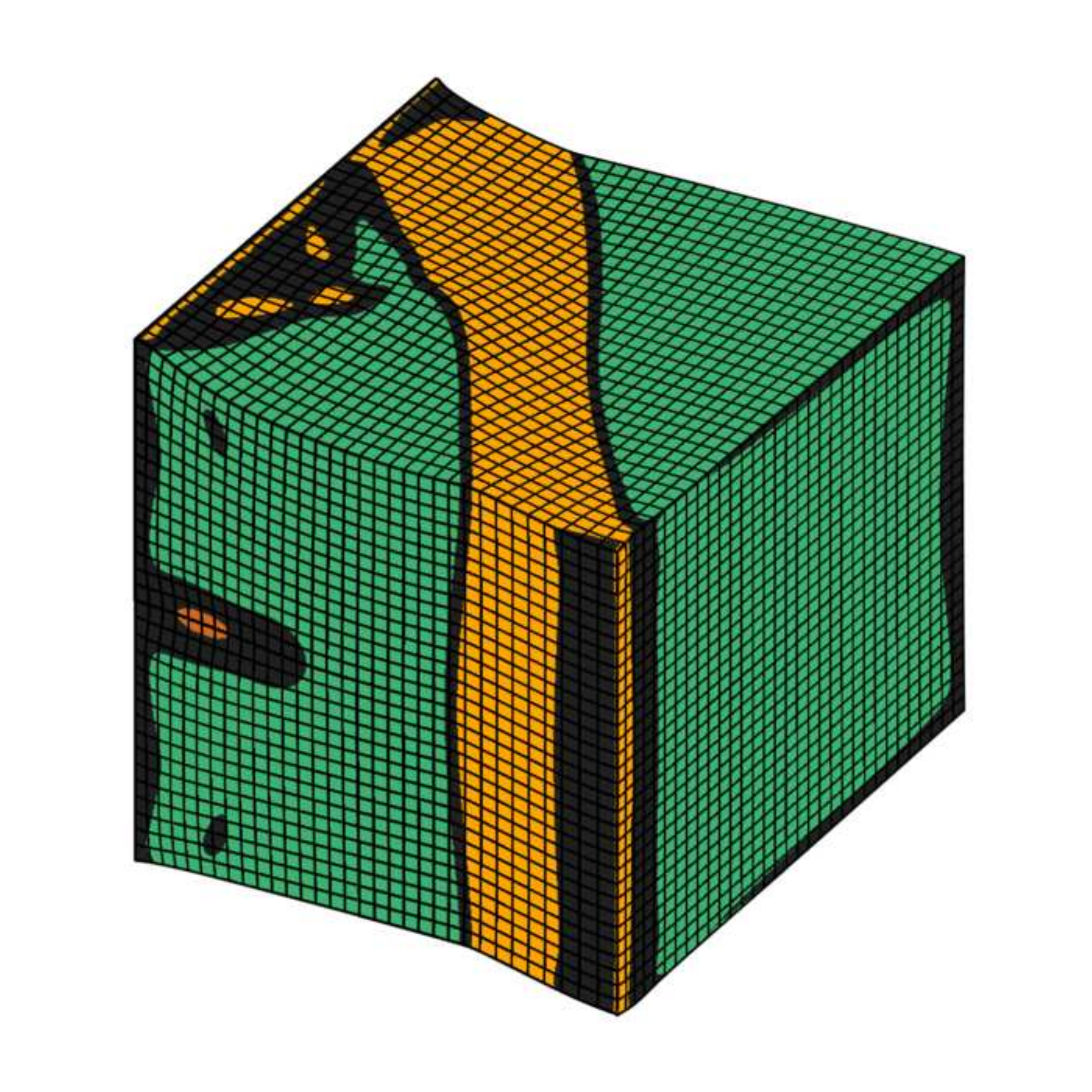} & 
                \includegraphics[scale=0.12]{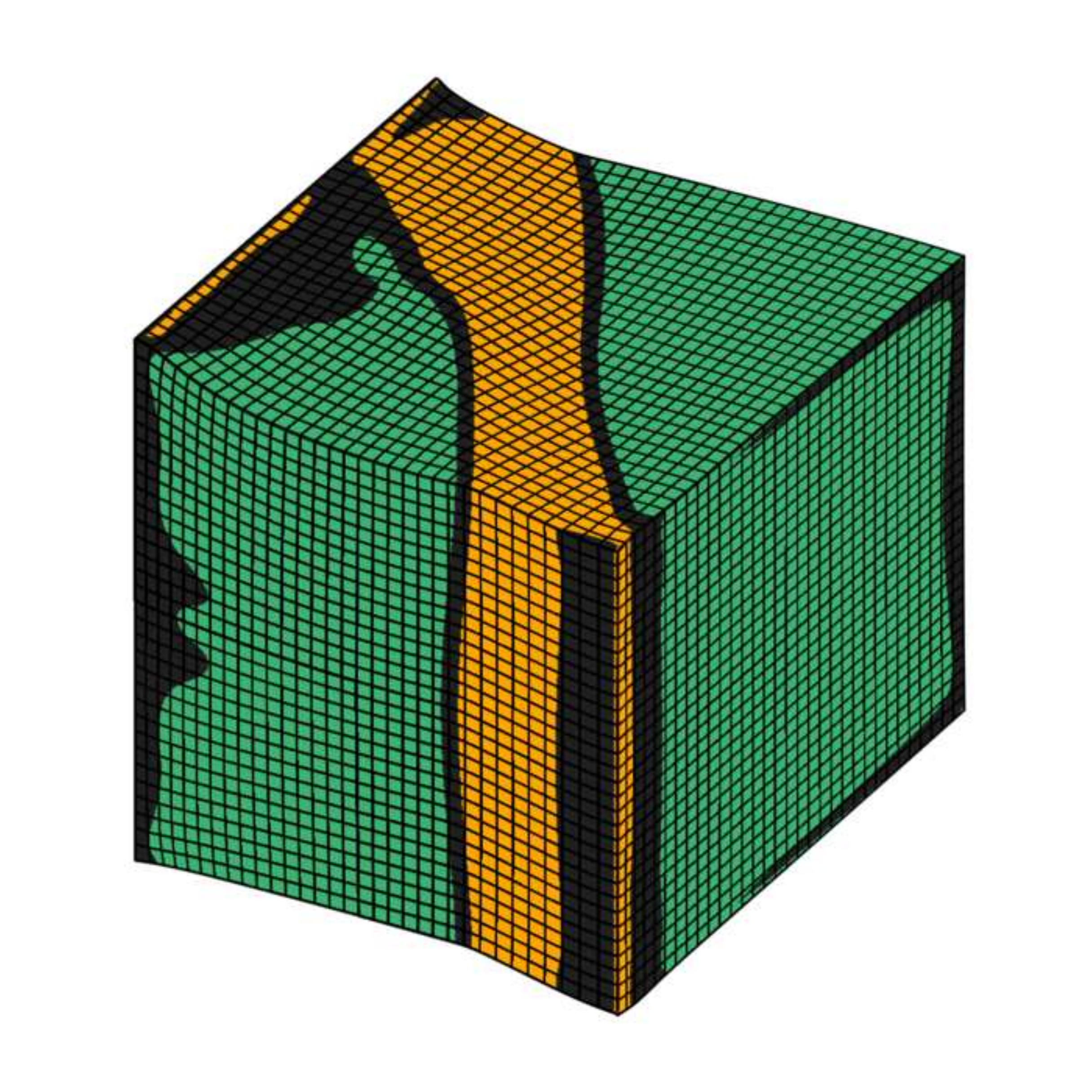} &
                \includegraphics[scale=0.12]{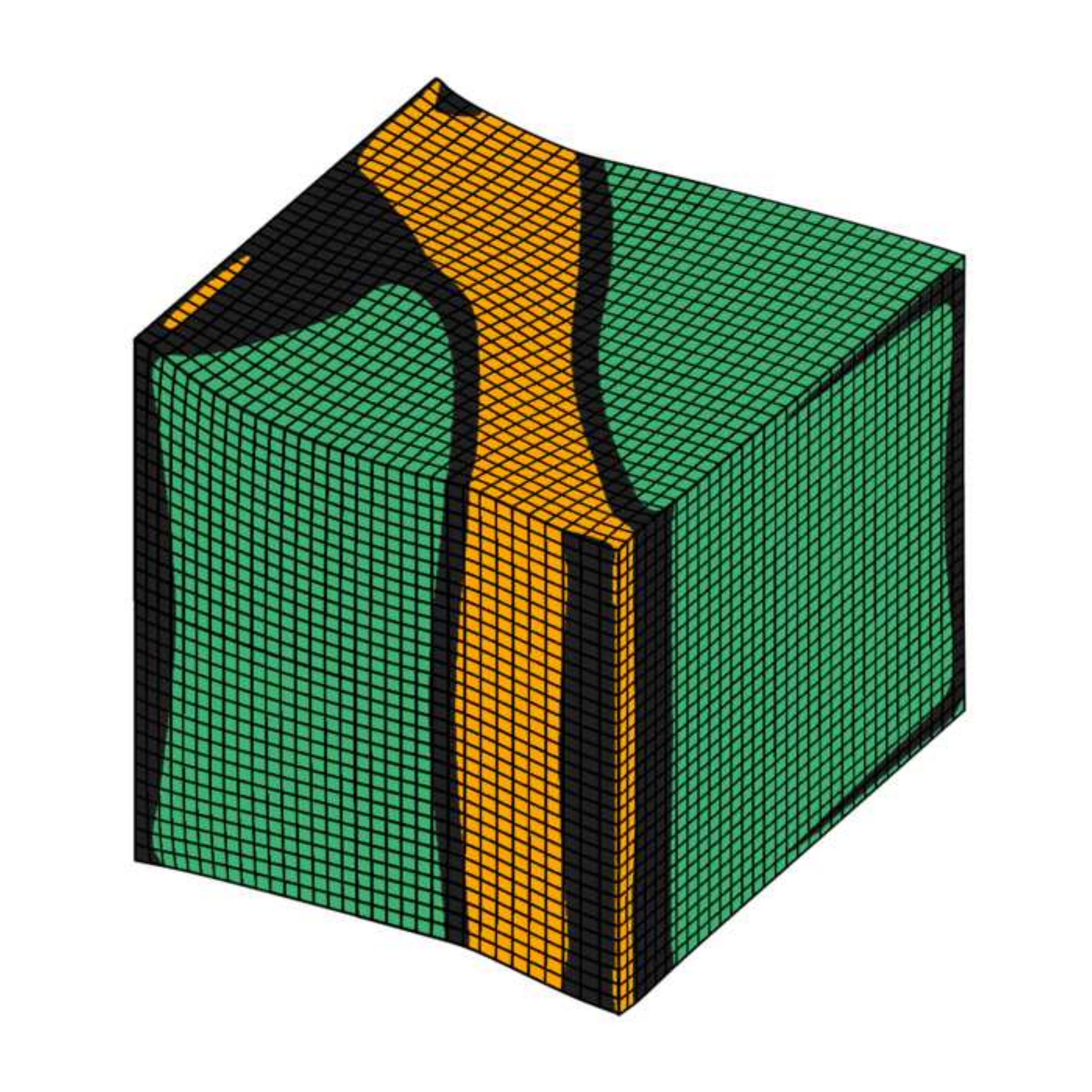} &
                \includegraphics[scale=0.12]{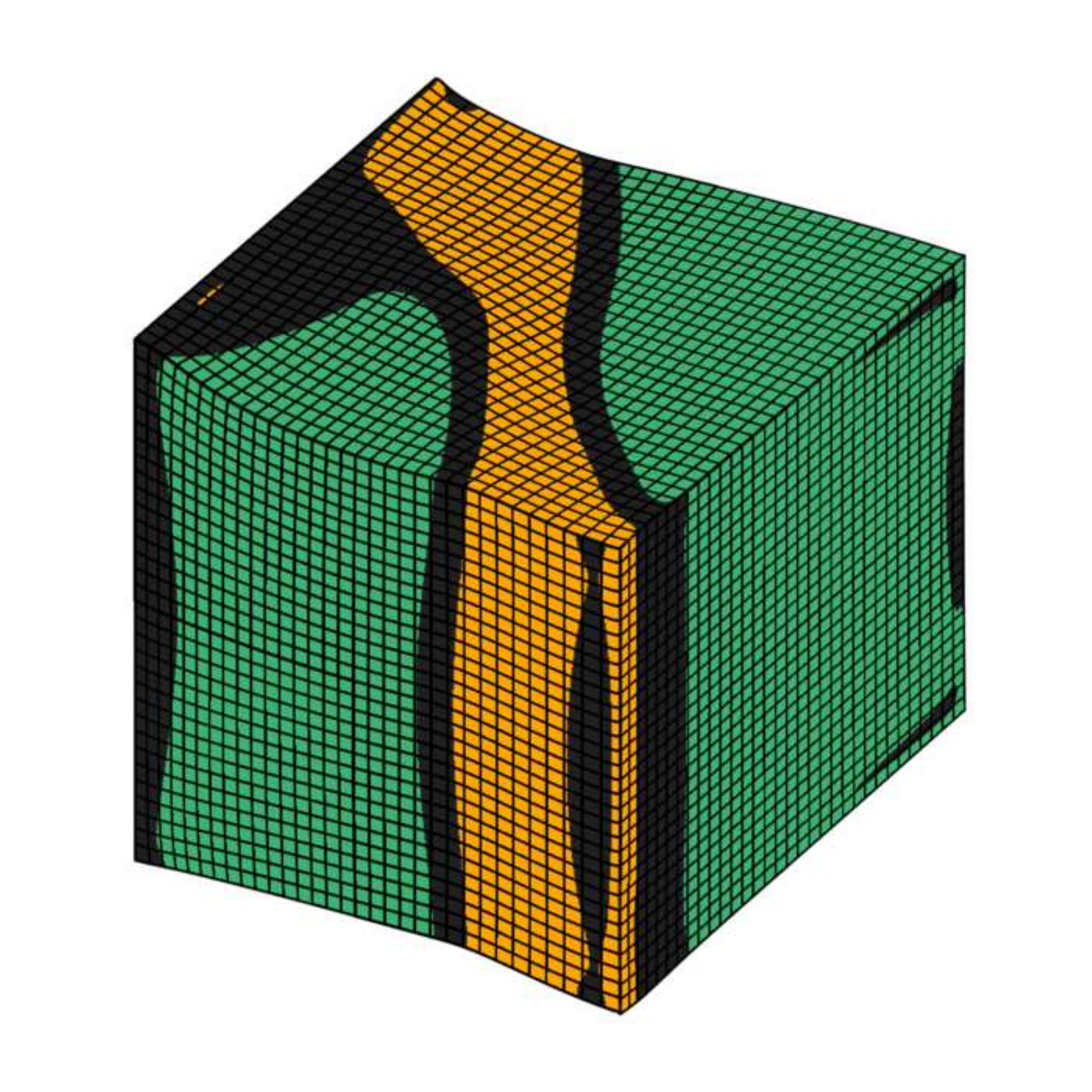} &
                \includegraphics[scale=0.12]{plot/plot_xyz.pdf}
            \end{tabular}  \\
            \parbox[t]{0.5cm}{ $A$ } &
            \begin{tabular}{p{2.8cm}p{2.8cm}p{2.8cm}p{2.8cm}p{2.8cm}p{2.5cm}}
                &
                \includegraphics[scale=0.12]{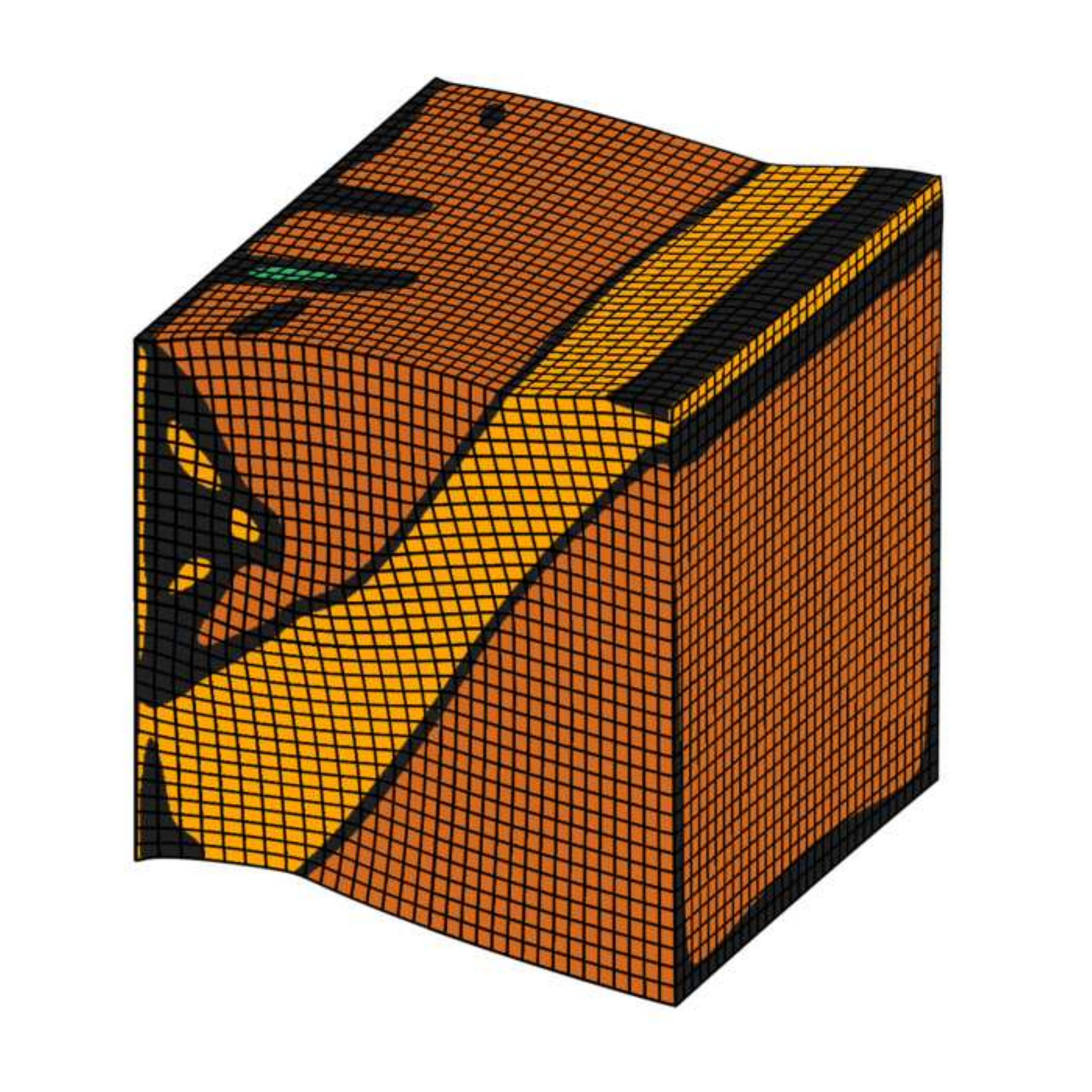} & 
                \includegraphics[scale=0.12]{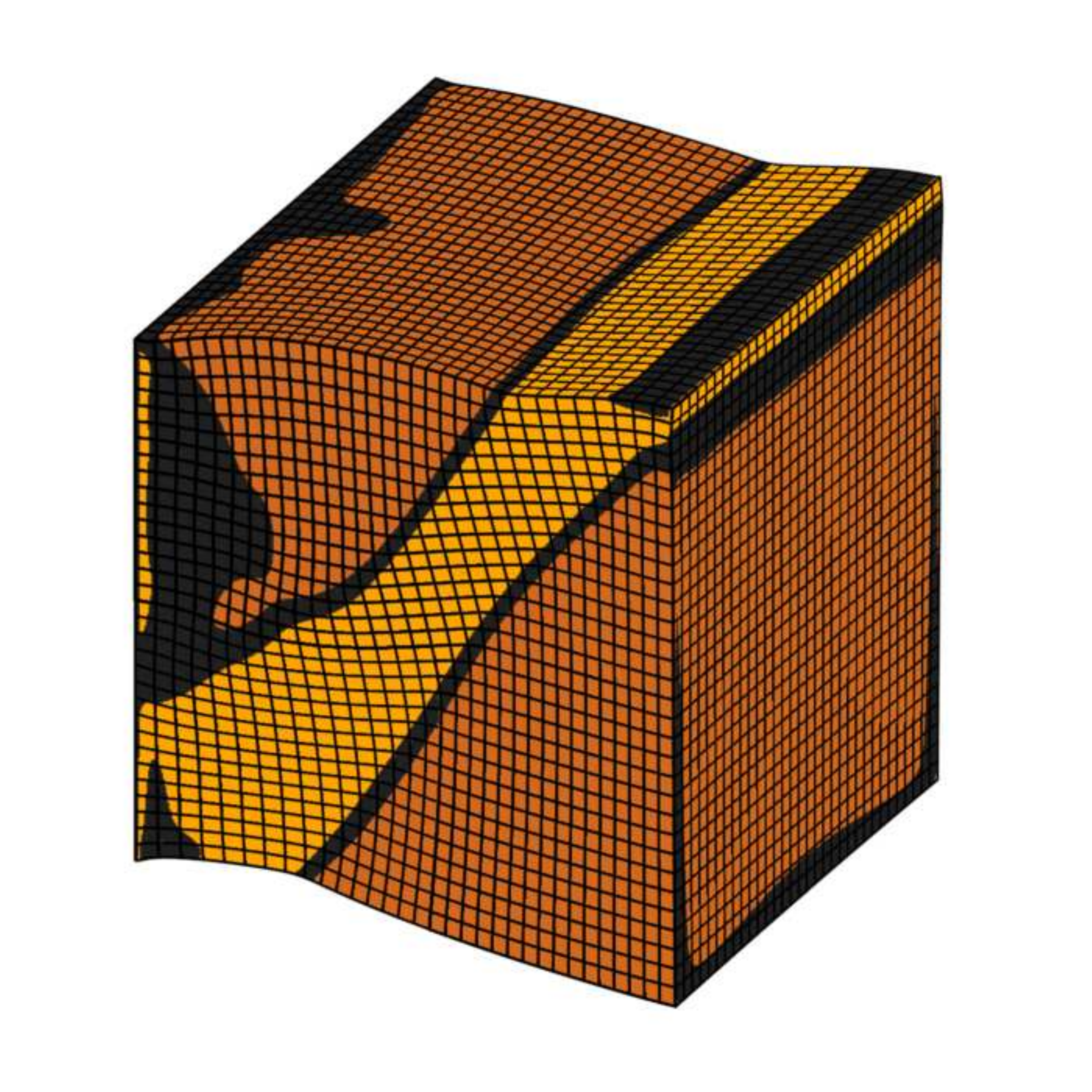} &
                \includegraphics[scale=0.12]{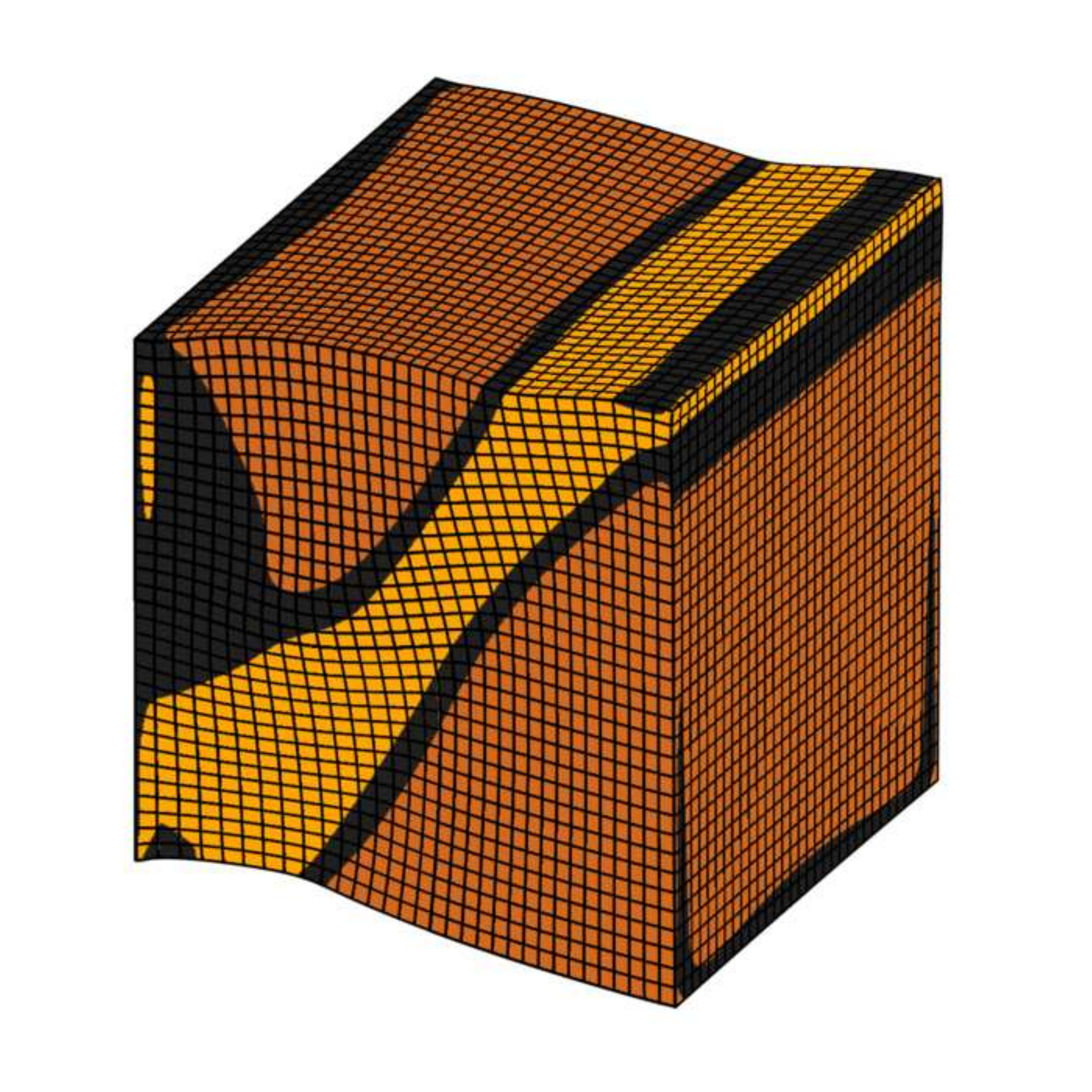} &
                \includegraphics[scale=0.12]{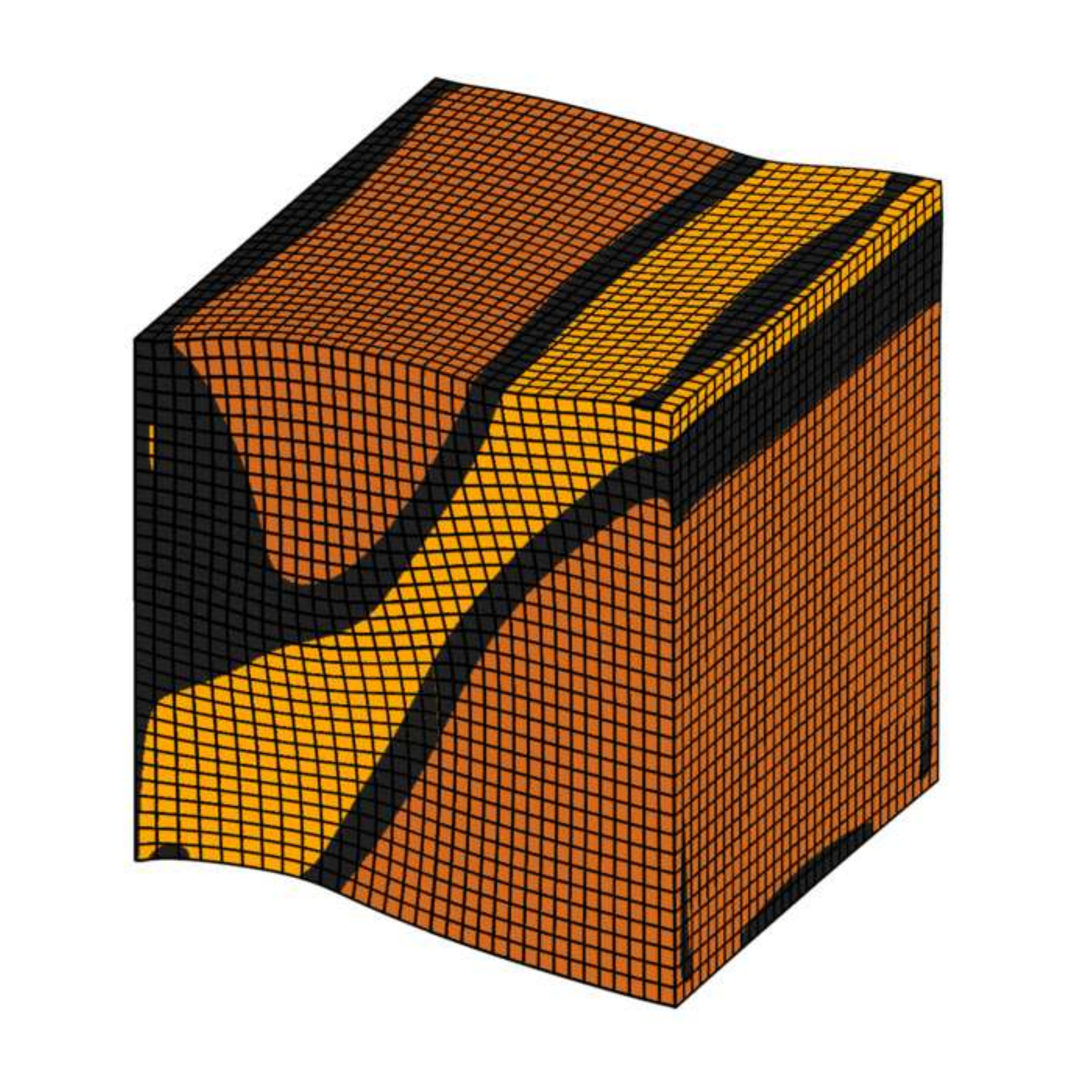} &
                \includegraphics[scale=0.12]{plot/plot_colorbar_phase.pdf}
            \end{tabular}
        \end{tabular}
    \end{center}
    \caption{Distribution of the three tetragonal variants for branches A - E on deformed configurations for selected values of $l$.  
    Solutions on the $128^3$ mesh have been overlaid with a $32^3$-plotting mesh.  
    The three tetragonal variants are indicated by different colors; see Fig. \ref{Fi:3well_diagram}.  }
    \label{Fi:phase}
\end{figure}

\begin{figure}
\begin{center}
        \begin{subfigure}[b]{2.0cm}
            \centering
            \includegraphics[scale=0.15]{plot/plot_xy.pdf}
        \end{subfigure}
        ~
        \begin{subfigure}[b]{3.8cm}
            \centering
            \includegraphics[scale=0.17]{plot/twin_1-2.pdf}
            \caption{}
            \label{Fi:twin2}
        \end{subfigure}
        ~
        \begin{subfigure}[b]{3.8cm}
            \centering
            \includegraphics[scale=0.19]{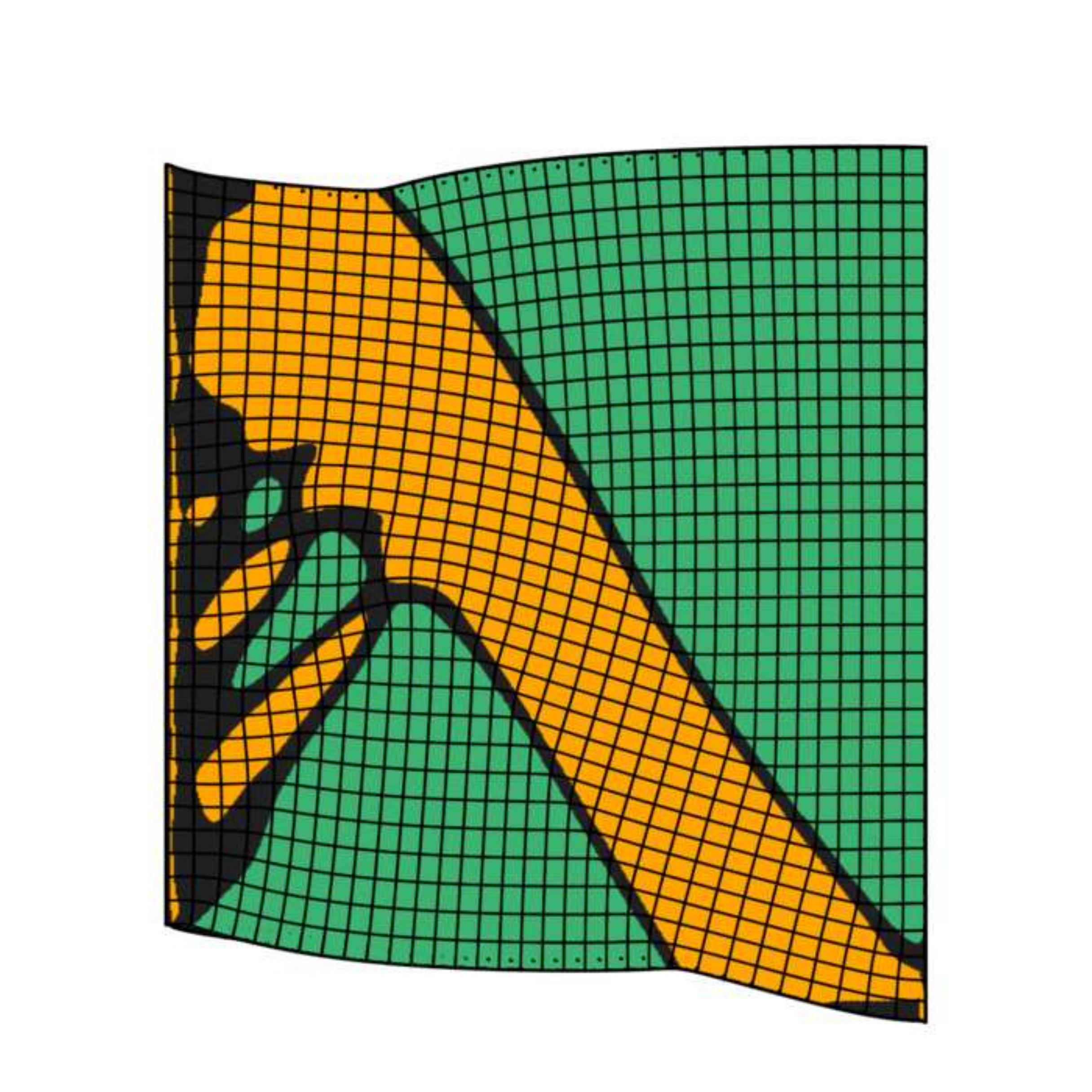}
            \caption{}
            \label{Fi:B_xy}
        \end{subfigure}
        ~
        \begin{subfigure}[b]{3.8cm}
            \centering
            \includegraphics[scale=0.19]{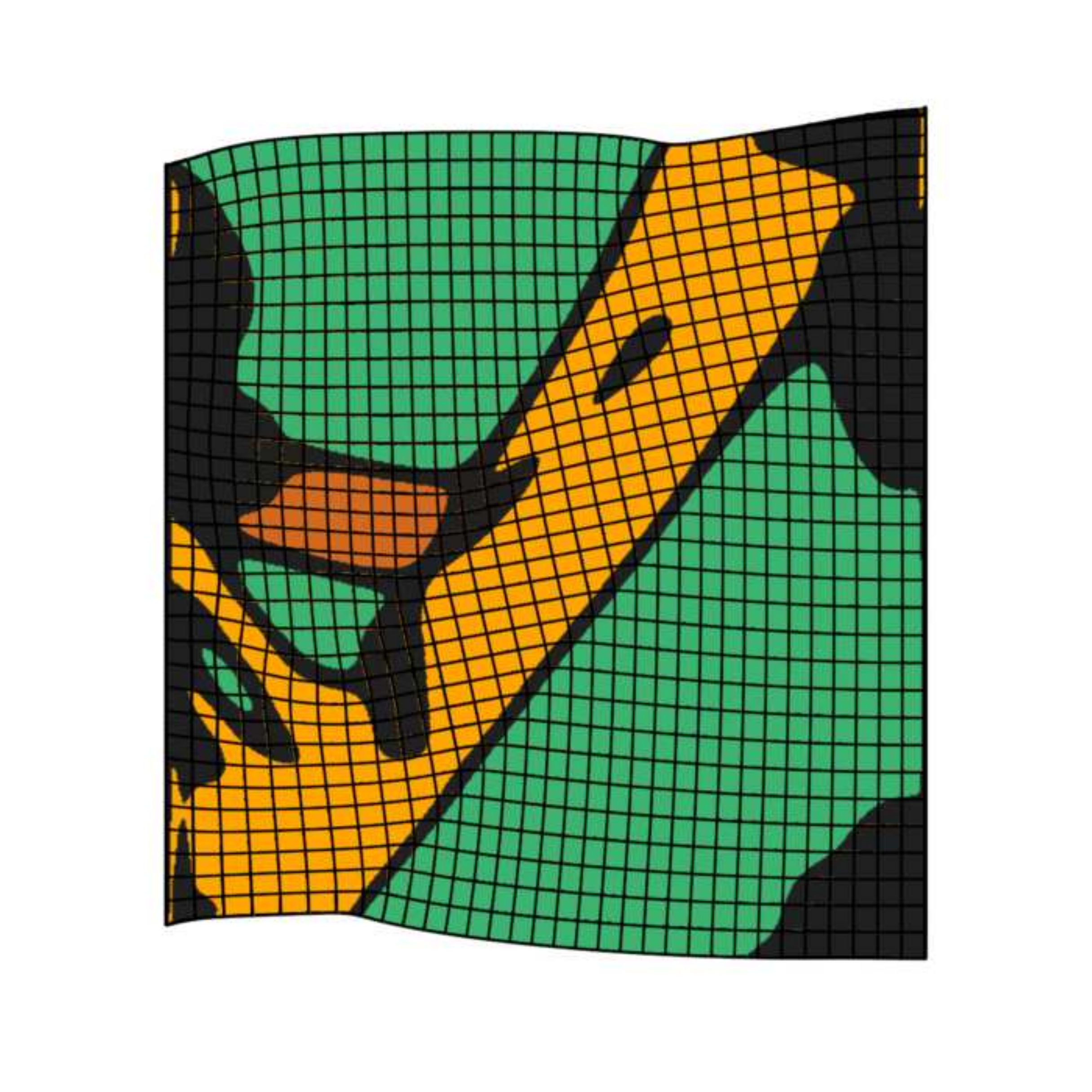}
            \caption{}
            \label{Fi:D_xy}
        \end{subfigure}
        ~
        \begin{subfigure}[b]{2.0cm}
            \centering
            \includegraphics[scale=0.15]{plot/plot_colorbar_phase.pdf}
        \end{subfigure}
\caption{(\subref{Fi:twin2}) One typical solution expected from the pure energy minimization argument for $l=0$, showing twinning between tetragonal variants 1 and 2 (reproduced from Fig. \ref{Fi:twin}).  Top-views of the microstructures with twinning for (\subref{Fi:B_xy}) branch B and (\subref{Fi:D_xy}) branch D obtained from the boundary value problem at $l=0.0625$.  }
\label{Fi:top}
\end{center}
\end{figure}

One can make several important observations for the solutions obtained in this section.  
As was also observed in one dimension, the interfaces between tetragonal variants, which are represented by the dark gray regions in Fig. \ref{Fi:phase}, become sharper as the length scale parameter $l$ decreases, the state of much of the material descending into the energy wells as indicated in Figs. \ref{Fi:e1-e2}.  These plots also highlight how a richer microstructure develops at lower values of $l$, with more material points being localized to those wells that are sparsely populated at higher $l$.
Fig. \ref{Fi:phase} also suggests that the interface thickness is proportional to $l$ and that, at a fixed value of $l$, the thickness is virtually the same over different branches -- an observation that can also be made for the one-dimensional problem in Sec. \ref{S:One-dimensional_example} from Fig. \ref{Fi:uX-X}.  
One can further infer from Figs. \ref{Fi:PI-l} and \ref{Fi:phase} that, regardless of stability, the equilibrium solutions achieving relatively low total free energy form via tetragonal variants with twin interfaces. This is as shown for the limiting case of $l=0$ in Fig. \ref{Fi:twin} according to the pure energy minimization argument.  
Fig. \ref{Fi:top}, for instance, shows the twin-structures observed in branches B and D at $l=0.0625$; Fig. \ref{Fi:twin} is repeated here to ease comparison.

\begin{table}
\begin{center}
\begin{tabular}[t]{|c|r|r|r|r|r|r|}
\hline
Branch & mesh & $l=0.0625$ & $l=0.0750$ & $l=0.1000$ & $l=0.1500$ & $l=0.2000$\\
\hline
  &  $64^3$ &           & -0.088755 & -0.065765 & -0.040403 & -0.027748 \\ 
E & $128^3$ &           & -0.020364 & -0.013910 & -0.007705 & -0.004932 \\
  & $256^3$ &           & -0.003886 & -0.002482 & -0.001261 & -0.009806 \\
\hline
  &  $64^3$ & -0.000301 & -0.000285 & -0.000284 & -0.000283 & -0.000206 \\ 
D & $128^3$ & -0.000041 & -0.000038 & -0.000038 & -0.000037 & -0.000027 \\
  & $256^3$ & -0.000005 & -0.000005 & -0.000005 & -0.000005 & -0.011602 \\
\hline
  &  $64^3$ & -0.046564 & -0.022921 & -0.000227 & -0.000264 & -0.000169 \\ 
C & $128^3$ & -0.006671 & -0.003216 & -0.000030 & -0.000034 & -0.000023 \\%
  & $256^3$ & -0.000861 & -0.000408 & -0.000004 & -0.000004 & -0.006915 \\
\hline
  &  $64^3$ & -0.007714 & -0.008696 & -0.008139 &  0.000352 &  0.000350 \\ 
B & $128^3$ & -0.001112 & -0.001184 & -0.001133 &  0.000046 &  0.000046 \\
  & $256^3$ & -0.000143 & -0.000152 & -0.000147 &  & \\
\hline
  &  $64^3$ &           & -0.003220 &  0.000313 &  0.000352 &  0.000350 \\ 
A & $128^3$ &           & -0.000522 &  0.000041 &  0.000046 &  0.000046 \\
  & $256^3$ &           & -0.000070 &  &  & \\
\hline
\end{tabular}
\end{center}
\caption{Smallest eigenvalues of the Hessians corresponding to the discretized counterpart of the second variation \eqref{E:D2_PI} for branches A - E at different refinement levels for selected values of $l$.  Positive values indicate stable/metastable solutions.  }
\label{Ta:eig}
\end{table}

\section{Conclusion and future works}\label{S:Conclusion}
We have considered martensitic phase transformations in three dimensions that are modeled by a free energy density function that is non-convex in strain space, and is regularized by Toupin's theory of gradient elasticity at finite strain.  
There exist three minima in the non-convex free energy density in strain space, corresponding to three, symmetrically equivalent, tetragonal, martensitic variants.  
The single maximum represents the cubic, austenite.
Our primary interest was to establish numerical procedures to obtain solution branches corresponding to the extrema/saddle points of the total free energy, and to assess their stability.
To this end, we have employed a simple branch-tracking technique to continuously follow a solution branch along the strain gradient length-scale parameter starting from the first solution computed using either a random initial guess or the homogeneous initial guess.  
The stability of each solution was then investigated in terms of the positive definiteness, or lack thereof, of the second variation of the total free energy. 
Each solution branch corresponds to a distinct twinned microstructure for the same boundary value problem. The phase interfaces between energetically stable tetragonal variants in a microstructure become sharper as the length scale parameter decreases.  The microstructures of certain branches themselves become richer with variants missing at higher values of $l$ emerging at lower $l$. To our knowledge this is the first work that comprehensively studies branching of solutions and observed twin structures in three-dimensional diffuse-interface problems based on a non-convex density function regularized by strain gradient terms.  

This work forms a foundation to study shape-memory alloys under loading, where different microstructures are experimentally observed for the same set of boundary conditions. A proper investigation of that class of problems also requires the incorporation of elastodynamics, in which case, the variations from initial conditions lead to different solution branches and therefore different microstructures for the same set of boundary conditions. A more direct comparison with experiments also needs a treatment of plasticity coupled with twinning as modelled here.
This work also provides a basis to study the homogenized response of a material exhibiting the microstructures corresponding to different solution branches. From such a study it may be possible to develop reduced order, effective constitutive models that also incorporate the evolution of martensitic microstructures.

\section*{Acknowledgments}
The numerical formulation and computations have been carried out as part of research supported by the U.S. Department of Energy, Office of Basic Energy Sciences, Division of Materials Sciences and Engineering under Award \#DE-SC0008637 that funds the PRedictive Integrated Structural Materials Science (PRISMS) Center at University of Michigan.  
This work used the Extreme Science and Engineering Discovery Environment (XSEDE), which is supported by National Science Foundation grant number ACI-1053575.  We used XSEDE resources \cite{xsede} through the Campus Champions program.  
The numerical computations in three dimensions presented here also made intensive use of resources of the National Energy Research Scientific Computing Center, a DOE Office of Science User Facility supported by the Office of Science of the U.S. Department of Energy under Contract No. DE-AC02-05CH11231.  
Finally, this research was supported in part through computational resources and services provided by Advanced Research Computing at the University of Michigan, Ann Arbor.

\bibliographystyle{plain}
\bibliography{reference}

\begin{thebibliography}{10}

\bibitem{petsc-user-ref}
Satish Balay, Shrirang Abhyankar, Mark~F. Adams, Jed Brown, Peter Brune, Kris
  Buschelman, Lisandro Dalcin, Victor Eijkhout, William~D. Gropp, Dinesh
  Kaushik, Matthew~G. Knepley, Dave~A. May, Lois~Curfman McInnes, Karl Rupp,
  Patrick Sanan, Barry~F. Smith, Stefano Zampini, Hong Zhang, and Hong Zhang.
\newblock {PETS}c users manual.
\newblock Technical Report ANL-95/11 - Revision 3.8, Argonne National
  Laboratory, 2017.

\bibitem{petsc-web-page}
Satish Balay, Shrirang Abhyankar, Mark~F. Adams, Jed Brown, Peter Brune, Kris
  Buschelman, Lisandro Dalcin, Victor Eijkhout, William~D. Gropp, Dinesh
  Kaushik, Matthew~G. Knepley, Dave~A. May, Lois~Curfman McInnes, Karl Rupp,
  Barry~F. Smith, Stefano Zampini, Hong Zhang, and Hong Zhang.
\newblock {PETS}c {W}eb page.
\newblock \url{http://www.mcs.anl.gov/petsc}, 2017.

\bibitem{petsc-efficient}
Satish Balay, William~D. Gropp, Lois~Curfman McInnes, and Barry~F. Smith.
\newblock Efficient management of parallelism in object oriented numerical
  software libraries.
\newblock In E.~Arge, A.~M. Bruaset, and H.~P. Langtangen, editors, {\em Modern
  Software Tools in Scientific Computing}, pages 163--202. Birkh{\"{a}}user
  Press, 1997.

\bibitem{Ball1987}
J.~M. Ball and R.~D. James.
\newblock Fine phase mixtures as minimizers of energy.
\newblock {\em Archive for Rational Mechanics and Analysis}, 100:13--52, 1987.

\bibitem{Carr1984}
Jack Carr, Morton~E. Gurtin, and Marshall Slemrod.
\newblock Structured phase transitions on a finite interval.
\newblock {\em Archive for Rational Mechanics and Analysis}, 86:317--351, 1984.

\bibitem{Chipot1988}
Michel Chipot and David Kinderlehrer.
\newblock Equilibrium configurations of crystals.
\newblock {\em Archive for Rational Mechanics and Analysis}, 103:237--277,
  1988.

\bibitem{Cottrell2009}
J.~Austin Cottrell, Thomas J.~R. Hughes, and Yuri Bazilevs.
\newblock {\em {I}sogeometric {A}nalysis}.
\newblock John Wiley \& Sons, Ltd, 2009.

\bibitem{Healey2007}
Timothy~J Healey and Ulrich Miller.
\newblock Two-phase equilibria in the anti-plane shear of an elastic solid with
  interfacial effects via global bifurcation.
\newblock {\em Proceedings of the Royal Society of London A: Mathematical,
  Physical and Engineering Sciences}, 463:1117--1134, 2007.

\bibitem{slepc2}
V.~Hernandez, J.~E. Roman, and V.~Vidal.
\newblock {SLEPc}: {S}calable {L}ibrary for {E}igenvalue {P}roblem
  {C}omputations.
\newblock {\em Lect. Notes Comput. Sci.}, 2565:377--391, 2003.

\bibitem{slepc1}
Vicente Hernandez, Jose~E. Roman, and Vicente Vidal.
\newblock {SLEPc}: A scalable and flexible toolkit for the solution of
  eigenvalue problems.
\newblock {\em {ACM} Trans. Math. Software}, 31(3):351--362, 2005.

\bibitem{Polizzi2009}
Eric Polizzi.
\newblock Density-matrix-based algorithm for solving eigenvalue problems.
\newblock {\em Phys. Rev. B}, 79:115112, Mar 2009.

\bibitem{slepc3}
J.~E. Roman, C.~Campos, E.~Romero, and A.~Tomas.
\newblock {SLEPc} users manual.
\newblock Technical Report DSIC-II/24/02 - Revision 3.9, D. Sistemes
  Inform\`atics i Computaci\'o, Universitat Polit\`ecnica de Val\`encia, 2018.

\bibitem{Rudraraju2014}
S.~Rudraraju, A.~Van der Ven, and K.~Garikipati.
\newblock Three-dimensional isogeometric solutions to general boundary value
  problems of {T}oupin's gradient elasticity theory at finite strains.
\newblock {\em Computer Methods in Applied Mechanics and Engineering}, 278:705
  -- 728, 2014.

\bibitem{IGAP4GradElast}
K.~{Sagiyama}.
\newblock {\tt IGAP4GradElast}: {An isogeometric analysis program for gradient
  elasticity written in \textsf{C} and aided by \textsf{Mathematica}. }.
\newblock \url{https://github.com/mechanoChem/IGAP4GradElast}, 2017.

\bibitem{Toupin1962}
R.A. Toupin.
\newblock Elastic materials with couple-stress.
\newblock {\em Archive for Rational Mechanics and Analysis}, 11:385--414, 1962.

\bibitem{Toupin1964}
R.A. Toupin.
\newblock {T}heories of elasticity with couple-stress.
\newblock {\em Archive for Rational Mechanics and Analysis}, 17:85--112, 1964.

\bibitem{xsede}
John Towns, Timothy Cockerill, Maytal Dahan, Ian Foster, Kelly Gaither, Andrew
  Grimshaw, Victor Hazlewood, Scott Lathrop, Dave Lifka, Gregory~D. Peterson,
  Ralph Roskies, J.~Ray Scott, and Nancy Wilkins-Diehr.
\newblock Xsede: Accelerating scientific discovery.
\newblock {\em Computing in Science \& Engineering}, 16(5):62--74, 2014.

\bibitem{Truskinovsky1995}
Lev Truskinovsky and Giovanni Zanzotto.
\newblock Finite-scale microstructures and metastability in one-dimensional
  elasticity.
\newblock {\em Meccanica}, 30:577--589, 1995.

\bibitem{Truskinovsky1996}
Lev Truskinovsky and Giovanni Zanzotto.
\newblock Ericksen's bar revisited : Energy wiggles.
\newblock {\em Journal of the Mechanics and Physics of Solids}, 44:1371 --
  1408, 1996.

\bibitem{Vainchtein1998}
A.~Vainchtein, T.~Healey, P.~Rosakis, and L.~Truskinovsky.
\newblock The role of the spinodal region in one-dimensional martensitic phase
  transitions.
\newblock {\em Physica D: Nonlinear Phenomena}, 115:29 -- 48, 1998.

\bibitem{Vainchtein1999}
Anna Vainchtein, Timothy~J. Healey, and Phoebus Rosakis.
\newblock Bifurcation and metastability in a new one-dimensional model for
  martensitic phase transitions.
\newblock {\em Computer Methods in Applied Mechanics and Engineering}, 170:407
  -- 421, 1999.

\end{thebibliography}


\end{document}